\documentclass [11pt]{article}
\usepackage {amsfonts}
\usepackage {mathrsfs}
\usepackage {amsmath}
\usepackage {latexsym}
\usepackage {amssymb}
\usepackage {amsthm}
\usepackage {amscd}
\date{}
\renewcommand{\baselinestretch}{1.2}
\title{\bf Classification of quasi-affine Generalized Dynkin Diagrams with Rank $ 4$ }
\author{ \small Zhengtang Tan $^{a}$,   Shouchuan Zhang $^{b}$   \\
\small $a$.School of Engineering and Design,  Hunan Normal University\\ Changsha  410081,    P.R. China \\
\small $b$. Department  of Mathematics,    Hunan University\\
Changsha  410082,   P.R. China\\
\small {\tt Emails:  z9491@sina.cn (SZ);    1843186255@qq.com (ZTT)} }
\date{}

\begin{document}
\newtheorem{Proposition}{Proposition}[section]
\newtheorem{Theorem}[Proposition]{Theorem}
\newtheorem{Definition}[Proposition]{Definition}
\newtheorem{Corollary}[Proposition]{Corollary}
\newtheorem{Lemma}[Proposition]{Lemma}
\newtheorem{Example}[Proposition]{Example}
\newtheorem{Remark}[Proposition]{Remark}

\maketitle


\begin {abstract}   All quasi-affine  connected Generalized Dynkin Diagram with rank $= 4$ are found. All quasi-affine  Nichols (Lie braided) algebras with rank $ 4$  are also found.
\end {abstract}

\section {Introduction and Preliminaries}\label {s0}

Nichols algebras play a fundamental role in the classification of finite-dimensional complex pointed Hopf algebras by means of the lifting method developed by Andruskiewitsch and Schneider \cite {AS02,  AS10,  AHS08}.
Heckenberger \cite {He06a,  He05} classified arithmetic root systems.   Heckenberger \cite {He06b} proved a {\rm  GDD} is  arithmetic  if and only if corresponding matrix is a finite Cartan matrix for {\rm  GDD}s of  Cartan types. W. Wu,    S. Zhang and   Y.-Z. Zhang \cite {WZZ15b} proved a Nichols Lie braided algebra is a finite dimensional  if and only if its {\rm  GDD},  which   fixed
parameter is of finite order,  is  arithmetic.

In fact every {\rm  GDD} is a weighted  graph. For every {\rm  GDD}, we can obtain two  braided tensor
categories $^{kG}_{kG} {\mathcal YD}$ and ${\mathcal T} (V).$ Braided diagrams in   braided tensor
categories can be used to represent highway traffic maps, aircraft route maps, circuit diagrams, and chip circuit diagrams. All these diagrams require the use of braids, as the intersections of two roads, two aircraft routes, and two chip circuits can all be represented by the braids of GDD

In order to classify finite dimensional Nichols algebras,  Heckenberger \cite {He06a,  He05} introduced {\rm  GDD}s and  classified arithmetic {\rm  GDD}s.
$\mathfrak B(V)$ is finite dimensional if and only if {\rm  GDD} of $V$ is arithmetic. Zhengtang Tan and   Shouchuan Zhang \cite {TZ24}
introduced quasi-affine {\rm  GDD}s  and  found all quasi-affine  connected Generalized Dynkin Diagrams with rank $> 5$.
The relationship between the connected components of Yang-Baxter  graphs and Nichols braided Lie algebras is given in \cite {WZT}.
We use notation in \cite {TZ24}.

We now recall some basic concepts of the graph theory (see \cite {Ha}).
Let $\Gamma _1$ be a non-empty set and $\Gamma _2 \subseteq \{  \{ u,  v\} \mid u,  v \in \Gamma _1,  \hbox { with } u \not= v \} \subseteq 2 ^{\Gamma_1}.$ Then $\Gamma = (\Gamma _1,  \Gamma_2)$ is called a graph;   $\Gamma_1$ is called the vertex set of  $\Gamma$;  $\Gamma_2$ is called the edge set of  $\Gamma$;  Element $\{u,  v\} \in \Gamma_2$ is called an edge. Let $F$ be an algebraically closed field of characteristic zero and $F^* := \{x \mid x\in F,  x \not= 0 \}$.  If   $\{x_1,    \cdots,    x_n\}$ is  a basis of  vector space  $V$ and
$C(x_i\otimes  x_j) = q_{ij} x_j\otimes x_i$ with $q_{ij} \in F^*$,
then $V$  is called a braided vector space of diagonal type,  $\{x_1,    \cdots,    x_n\}$
is called  canonical basis and $(q_{ij})_{n\times n}$  is called braided matrix. Let
$\widetilde{q}_{i j}:= q_{ij}q_{ji}$ for  $i,  j \in \{1,  2,  \cdots,  n \}$ with $i \not=j$. If let $\Gamma _1= \{1,  2,  \cdots,  n \}$ and $\Gamma _2 = \{  \{ u,  v\} \mid \widetilde{q}_{u v} \not= 1,    u,  v \in \Gamma _1,  \hbox { with } u \not= v \}$,  then $\Gamma = (\Gamma _1,  \Gamma_2)$ is a graph, called a Yang-Baxter graph. Set $q_{ii}$ over vertex $i$ and $\widetilde{q}_{i j}$ over edge $\{i,  j\}$ for $i,  j \in \Gamma _1$ with $i \not=j$.
Then $\Gamma = (\Gamma _1,  \Gamma_2)$ is called a Generalized Dynkin Diagram of braided vector space $V, $ written as {\rm  GDD} in short(see \cite [Def. 1.2.1] {He05}).
If $\Delta (\mathfrak{B}(V))$ is an arithmetic root system.   then we call its {\rm  GDD}  an arithmetic {\rm  GDD}  for convenience. All arithmetic {\rm  GDD}s are listed in Table $A1$,  $A2$,  $B$ and $C$ in \cite {He06a,  He05}

Let $(q_{ij})_{n \times n }$ be a braided  matrix. If  $q_{ij} q_{ji}  \left \{  \begin{array}{ll}
 \not= 1,  & \mbox {when }  \mid j-i \mid = 1\\
  =1,   & \mbox {when}  \mid j-i \mid \not= 1    \\
\end{array}\right. $ for any $1\le i\not= j \le n,$ then  $(q_{ij})_{n \times n }$ is called a  chain or labelled  chain.
If $(q_{ij})_{n \times n }$ is  a  chain and \begin {eqnarray} \label {ppe1}(q_{11} q_{1, 2} q_{2, 1} -1)(q_{11} +1)=0;
 (q_{n, n } q_{n, n-1} q_{n-1, n} -1)(q_{n,n} +1)=0;\end {eqnarray} i.e.
  \begin {eqnarray} \label {ppe2}
   q_{ii} +1= q_{i, i-1} q_{i-1, i}q_{i, i+1} q_{i+1, i}-1=0
   \end  {eqnarray}
    \begin {eqnarray} \label {ppe3}
   \mbox { or  }
  q_{ii}q_{i, i-1} q_{i-1, i}=q_{ii} q_{i, i+1} q_{i+1, i}=1,\end  {eqnarray}  $1<i < n$,
   then  the braided matrix's {\rm  GDD} is called a simple chain (see \cite [Def.1]{He06a}).  Conditions (\ref {ppe1}), (\ref {ppe2}) and (\ref {ppe3}) are called simple chain conditions.

Every   arithmetic {\rm  GDD} in Row 1-10 in Table C is  called a classical {\rm  GDD}.

A  {\rm  GDD}  is called a quasi-arithmetic  {\rm  GDD}  if its every sub-{\rm  GDD}   by deleting every vertex is arithmetic.
A  {\rm  GDD}  is called a quasi-affine  {\rm  GDD}  if it is quasi-arithmetic and is not arithmetic. In this case,  Nichols algebra $\mathfrak B(V)$
and Nichols Lie braded algebra $\mathfrak L(V)$ are said to be quasi-affine.
In other word,  if a {\rm  GDD} is quasi-affine  of a braided vector space $V$ which   fixed
parameter is of finite order,  then Nichols algebra and Nichols Lie braded algebra
of every proper sub{\rm  GDD} are finite dimensional with $\dim   \mathfrak B { (V)} = \infty$
and $\dim   \mathfrak L { (V)} = \infty$.

In this paper,   using Table $A1$,  $A2$,  $B$ and $C$ in \cite {He06a,  He05},  we find all quasi-affine  connected Generalized Dynkin Diagram with rank $=4$. We also  find all quasi-affine   Nichols  algebras and quasi-affine  Nichols Lie braided algebras with rank $ 4$.

\section {Properties about arithmetic {\rm  GDD} } \label {s2}


\begin {Lemma} \label {mainlemma} A {\rm  GDD} is a classical {\rm  GDD} if and only if it is one of classical Type 1-7,  Here
classical types are  listed as follows:\\

 {\ }\ \ \ \ \ \ \ \ \ \ \ \ \ \ \ \ \ \ \ \ \ \ \ \   \ \ \ \ \ \  $\begin{picture}(100,       15)

\put(-125,       -1){ {\rm   Type   1},   $2\le n$.}
\put(80,       1){\makebox(0,       0)[t]{$\bullet$}}

\put(48,       -1){\line(1,       0){33}}
\put(10,       1){\makebox(0,      0)[t]{$C_{n-1,  q,  i_1,  i_2,  \cdots,  i_j }$}}

\put(-18,      10){$$}
\put(0,       5){$$}
\put(22,      10){$$}
\put(50,       5){$q^{-2}$}

\put(68,       10){$q^2$}

 \ \ \ \ \ \ \ \ \ \ \ \ \ \ \ \ \ \ \
  \ \ \ \ \ \ \ \ \ \ \ \ \ \ \ \ \ \ \ { }$q \in F^{*}\setminus \{1,  -1\}$,  $0\le j\le n-1.$

 \put(220,       -1) {}
\end{picture}$\\

 {\ }\ \ \ \ \ \ \ \ \ \ \ \ \ \ \ \ \ \ \ \ \ \ \ \   \ \ \ \ \ \  $\begin{picture}(100,       15)

\put(-125,       -1){ {\rm   Type   2}. $2\le n$.}

\put(90,       1){\makebox(0,       0)[t]{$\bullet$}}

\put(58,       -1){\line(1,       0){33}}
\put(27,       1){\makebox(0,      0)[t]{$C_{n-1,  q^2,  i_1,  i_2,  \cdots,  i_j }$}}

\put(-8,      10){$$}
\put(0,       5){$$}
\put(22,      10){$$}
\put(70,       5){$q^{-2}$}

\put(88,       10){$q$}

  \ \ \ \ \ \ \ \ \ \ \ \ \ \ \ \ \ \ \
  \ \ \ \ \ \ \ \ \ \ \ \ \ \ \ \ \ \ \ $q \in F^{*}\setminus \{1,  -1\}$.  $0\le j\le n-1.$
\end{picture}$\\

 {\ }\ \ \ \ \ \ \ \ \ \ \ \ \ \ \ \ \ \ \ \ \ \ \ \   \ \ \ \ \ \  $\begin{picture}(100,       15)

\put(-125,       -1){ {\rm   Type   3},  $2\le n$.}

\put(80,       1){\makebox(0,       0)[t]{$\bullet$}}

\put(48,       -1){\line(1,       0){33}}
\put(27,       1){\makebox(0,      0)[t]{$C_{n-1,  q^{-2},  i_1,  i_2,  \cdots,  i_j }$}}

\put(-18,      10){$$}
\put(0,       5){$$}
\put(22,      10){$$}
\put(60,       5){$q^{2}$}

\put(78,       10){$-q^{-1}$}

  \ \ \ \ \ \ \ \ \ \ \ \ \ \ \ \ \ \ \
  \ \ \ \ \ \ \ \ \ \ \ \ \ \ \ \ \ \ \ $q \in F^{*}\setminus \{1,  -1\}$,  $0\le j\le n-1.$

\end{picture}$\\

 {\ }\ \ \ \ \ \ \ \ \ \ \ \ \ \ \ \ \ \ \ \ \ \ \ \   \ \ \ \ \ \  $\begin{picture}(100,       15)

\put(-125,       -1){{\rm   Type   4}. $2\le n$. }

\put(80,       1){\makebox(0,       0)[t]{$\bullet$}}

\put(48,       -1){\line(1,       0){33}}
\put(27,       1){\makebox(0,      0)[t]{$C_{n-1,  -q^{-1},  i_1,  i_2,  \cdots,  i_j }$}}

\put(-18,      10){$$}
\put(0,       5){$$}
\put(22,      10){$$}
\put(60,       5){$-q$}

\put(78,       10){$q$}

  \ \ \ \ \ \ \ \ \ \ \ \ \ \ \ \ \ \ \
  \ \ \ \ \ \ \ \ \ \ \ \ \ \ \ \ \ \ \ $q^3 =1$. $0\le j\le n-1.$

\end{picture}$\\ \\

\ \ \ \     \ \ \ \ \ \  $\begin{picture}(100,       15)

\put(-45,       -1){ {\rm   Type   5}. $3\le n$.}
\put(124,       1){\makebox(0,       0)[t]{$C_{n-2,  q^{},  i_1,  i_2,  \cdots,  i_j }$}}
\put(190,      -11){\makebox(0,      0)[t]{$\bullet$}}
\put(190,     15){\makebox(0,      0)[t]{$\bullet$}}
\put(162,     -1){\line(2,      1){27}}
\put(190,       -14){\line(-2,      1){27}}

\put(120,       10){$$}

\put(135,       5){$$}

\put(155,      10){$$}

\put(160,      -20){$q^{-1}$}
\put(165,       15){$q^{-1}$}

\put(193,       -12){$q$}
\put(193,       18){$q$}

 \ \ \ \ \ \ \ \ \ \ \ \ \ \ \ \ \ \ \ \ \ \ \ \ \ \ \ \ \ \ \ \ \ \ \ \ \ \
  \ \ \ \ \ \ \ \ \ \ \ \ \ \ \ \ \ \ \ $q\not= 1, $    $0\le j\le n-2.$

\put(215,         -1)  { }
\end{picture}$\\  \\

\ \ \ \      \ \ \ \ \ \  $\begin{picture}(100,       15)

\put(-45,       -1){{\rm   Type   6}. $3\le n$. }

\put(104,       1){\makebox(0,       0)[t]{$C_{n-2,  q^{},  i_1,  i_2,  \cdots,  i_j }$}}
\put(170,      -11){\makebox(0,      0)[t]{$\bullet$}}
\put(170,     15){\makebox(0,      0)[t]{$\bullet$}}

\put(142,     -1){\line(2,      1){27}}
\put(170,       -14){\line(-2,      1){27}}

\put(170,       -14){\line(0,      1){27}}

\put(100,       10){$$}
\put(120,       5){$$}

\put(127,      10){$$}

\put(140,      -20){$q^{-1}$}
\put(145,       15){$q^{-1}$}

\put(178,       -20){$-1$}

\put(178,       0){$q^{2}$}

\put(175,       10){$-1$}

 \ \ \ \ \ \ \ \ \ \ \ \ \ \ \ \ \ \ \ \ \ \ \ \ \ \ \ \ \ \ \ \ \ \ \ \ \ \
  \ \ \ \ \ \ \ \ \ \ \ \ \ \ \ \ \ \ \  $q^2 \not= 1$. $0\le j\le n-2.      $
\put(195,         -1)  {,   }
\end{picture}$\\

{\ }\!\!\!\!\!\!
{\rm   Type   7}. $1\le n$.
$C_{n,  q^{-1},  i_1,  i_2,  \cdots,  i_j }.$
$q \not= 1$,   $0\le j\le n$.

\end {Lemma}

The sub{\rm  GDD} of the right hand of  Type i in above Lemma is called head type i,  written as {\rm h}-type i with $i=1, 2, \cdots,  7.$ {\rm h}-type 7 has two forms. Right leaf vertex  is $-1$:    $\begin{picture}(100,       15)  \put(0,       -1){ }

\put(60,       1){\makebox(0,       0)[t]{$\bullet$}}

\put(28,       -1){\line(1,       0){33}}



\put(22,      10){}
\put(38,       5){$q^{-1}$}

\put(60,       10){$-1^{}$}

  \ \ \ \ \ \ \ \ \ \ \ \ \ \ \ \ \ \ \ {  }
\end{picture}$\\
with  $q \not= -1$, written $T6$;  right leaf vertex  is $q$:    $\begin{picture}(100,       15)  \put(-20,       -1){ }

\put(60,       1){\makebox(0,       0)[t]{$\bullet$}}

\put(28,       -1){\line(1,       0){33}}



\put(22,      10){}
\put(38,       5){$q^{-1}$}

\put(60,       10){$q^{}$}

  \ \ \ \ \ \ \ \ \ \ \ \ \ \ \ \ \ \ \ { written  T5. }
\end{picture}$\\ The left leaf vertex  of   above every {\rm  GDD}    is called the classical  tail  of the {\rm  GDD}.

A {\rm  GDD} is called a simple circle  if  it is a circle and its every sub-GDD by deleting every vertex
is a simple chain.
A {\rm  GDD} with rank 4 is called a semi-simple circle  if  it is a circle and its every sub-GDD by deleting every vertex
is a simple chain or a Type 5.

A vertex is called a leaf if it has only an edge connected it.

\begin {Definition} \label {2.1''}
{ {\rm (i)}} A connected arithmetic chain   is called a quasi(semi)-classical {\rm  GDD}  if the sub-{\rm GDD} by deleting a leaf  is a classical {\rm  GDD}( Type 7 ) which classical tail is other leaf. Further more,  the other leaf  is called the tail of the quasi(semi)-classical {\rm  GDD}.

{ {\rm (ii)}}
A connected  arithmetic non chain is called a quasi(semi)-classical {\rm  GDD} if there exist two vertexes
 such that following conditions hold:  sub-{\rm GDD} by deleting  a vertex is a connected classical (simple) {\rm  GDD};  sub-{\rm GDD} by deleting  an other vertex is also a connected classical (simple) {\rm  GDD};  the two classical tails  are the same  and  the two fixed
parameters  of the two classical {\rm  GDD}s  are the same. The two classical tails are called the tail of the quasi(semi)-classical  {\rm  GDD}.
\end {Definition}

A quasi-classical  {\rm  GDD} which is not  semi-classical is called a strict quasi-classical  {\rm  GDD}. A semi-classical {\rm  GDD} which is not classical is called a strict semi-classical {\rm  GDD}.

\begin {Definition} \label {6.1.-3}
{\rm (i)} If adding T5 on tail of a  quasi-classical   {\rm  GDD} is an arithmetic {\rm  GDD}   then the  quasi-classical  {\rm  GDD} is called continual on the tail via T5.  Otherwise  adding T5 on tail of a  quasi-classical   {\rm  GDD} is called a discontinuous {\rm  GDD} on the tail via T5.  For T6  we can similarly define these.

 {\rm (ii)}
 {\rm h}-Type i  adding on tail of  semi-classical {\rm  GDD}, which is not a Type 7, is called classical + semi-classical,  $1\le i \le 4.$
{\rm h}-Type 5  adding on tail of  semi-classical {\rm  GDD} is called classical + semi-classical  when  semi-classical {\rm  GDD} is continual via T5 ;   {\rm h}-Type 6  adding on tail of  semi-classical {\rm  GDD} is called classical + semi-classical when  semi-classical {\rm  GDD} is continual via T6.

 \end {Definition}

  {\rm h}-Type i ( $1\le i \le 6)$ adding on classical tail of  classical {\rm  GDD}, which is not a Type 7, is called classical + classical or bi-classical.

For a semi-classical {\rm  GDD}   $\alpha $  with rank 3, which is not Type 7,
 {\rm  GDD} $\beta $ by adding a vertex on the middle vertex of  $\alpha $  such that  tail of  $\alpha $ is  contained in {\rm h}-Type 5  if the tail of  $\alpha $  is in T5;
 {\rm  GDD} $\gamma$ by   adding a vertex on the middle vertex  of  $\alpha $  and tail of   $\alpha $  such that  tail  of   $\alpha $ is contained in {\rm h}- Type 6  if the tail of $\alpha $  is in T6. $\beta $ is called  Type 5 union semi-classical  $\alpha $.    $\gamma $ is called  Type 6 union semi-classical  $\alpha $.

\begin {Lemma}\label {1.1.5}
 If {\rm  GDD}  with $3 \le \hbox { rank } n \le 4$ is an arithmetic non chain,   then sub-{\rm GDD} by deleting  some vertex is connected simple chain.

\end {Lemma}

\begin {Lemma}\label {2.63}
{\rm (i)}  All  arithmetic {\rm GDD}s with rank $n=3$ are quasi-classical  except   {\rm GDD} $3$ of Row $16$,    {\rm GDD} $8$ of Row $17$,
 {\rm GDD} $3$ of Row $7$ when   $q \notin  R_{4}$,   {\rm GDD} $3$ of Row $9.$

{\rm (ii)}  All  arithmetic  {\rm GDD}s with rank $n=4$ are quasi-classical   except  {\rm GDD} $4$ of Row $9$ when   $q \notin  R_{4} \cup   R_{5}$,
 {\rm GDD} $6$ of Row $9$ when   $q \notin    R_{5}$;  {\rm GDD} $2$ of Row $14$,  $q_{}^4 \not= 1$;  {\rm GDD} $3$ of Row $14$,  $q_{}^4 \not= 1$;  {\rm GDD} $4$ of Row $14$,  $q_{}^4 \not= 1$.
\end {Lemma}

\begin {Lemma}\label {2.63}
{\rm (i)}  All  arithmetic {\rm GDD}s with rank $n=3$ are quasi-classical  except   {\rm GDD} $3$ of Row $16$,    {\rm GDD} $8$ of Row $17$,
 {\rm GDD} $3$ of Row $7$ when   $q \notin  R_{4}$,   {\rm GDD} $3$ of Row $9.$

{\rm (ii)}  All  arithmetic  {\rm GDD}s with rank $n=4$ are quasi-classical   except  {\rm GDD} $4$ of Row $9$ when   $q \notin  R_{4} \cup   R_{5}$,
 {\rm GDD} $6$ of Row $9$ when   $q \notin    R_{5}$;  {\rm GDD} $2$ of Row $14$,  $q_{}^4 \not= 1$;  {\rm GDD} $4$ of Row $14$,  $q_{}^4 \not= 1$.
\end {Lemma}

\begin {Proposition} \label {1.1.10}Assume rank $n =4.$

{\rm (i)}  Classical  + semi-classical {\rm  GDD}s are quasi-affine.

{\rm (ii)}  Type 5, 6  union  semi-classical {\rm  GDD}s are quasi-affine.

\end {Proposition}
{\bf Proof.}
We only need prove that  these {\rm  GDD}s are not arithmetic.  If  a non chain {\rm  GDD}  satisfied conditions is  an arithmetic {\rm  GDD},  then we obtain a contradiction by Lemma \ref {1.1.5} and \\ \\

$\begin{picture}(100,       15) \put(-58,        -1){ }

\put(60,       1){\makebox(0,       0)[t]{$\bullet$}}

\put(28,       -1){\line(1,       0){33}}
\put(27,       1){\makebox(0,      0)[t]{$\bullet$}}
\put(-14,       1){\makebox(0,      0)[t]{$\bullet$}}

\put(-14,      -1){\line(1,       0){50}}

\put(58,       -12){$q$}

\put(40,       -12){$q^{-1}$}

\put(22,      -12){$-1$}
\put(0,       -12){$-1{}$}

\put(-18,      -12){${-1}$}

\put(27,     38){\makebox(0,      0)[t]{$\bullet$}}

\put(27,       0){\line(0,      1){35}}

\put(30,       30){${-1}$}

\put(30,      15){$-1$}

\ \ \ \ \ \ \ \ \ \ \ \ \ \ \ \ \ \ \ \ \ \ {,    order of $q$ $> 2.$ If  a chain {\rm  GDD}  satisfied conditions is  an arithmetic {\rm  GDD},}
\put(80,         -1)  {    } \end{picture}$\\ \\
then we obtain a contradiction by \cite [Lemma 2.9]{TZ22}.  \hfill $\Box$

\begin {Proposition} \label {1.1.14}  Adding a vertex on tail of  strict quasi-classical  {\rm  GDD}  with rank 3 is not  quasi-affine except discontinuous {\rm  GDD}s.

\end{Proposition}
{\bf Proof.} It follows from \cite [Lemma 2.4]{TZ22}.  \hfill $\Box$

\begin {Definition} \label {1.1.11} A {\rm  GDD}  is called a quasi-arithmetic {\rm  GDD}  if sub-{\rm GDD} by deleting  every vertex of the {\rm  GDD}  is arithmetic.
\end {Definition}

Obviously,  adding an h-Type 7 on a tail of quasi-classical  {\rm  GDD}  is quasi-arithmetic. Adding a head of classical {\rm  GDD}s on a tail of semi-classical {\rm  GDD}s is quasi-arithmetic. The two {\rm  GDD}S are called near-classical  {\rm  GDD}s,  written nc in short.

Every  quasi-affine  circle  \\ \\ \\

 {\ }\ \ \ \ \ \ \ \ \ \ \ \    \ \ \ \ \ \  $\begin{picture}(100,       15) \put(-85,       -1){(a)}

\put(60,       1){\makebox(0,       0)[t]{$\bullet$}}

\put(28,       -1){\line(1,       0){33}}
\put(27,       1){\makebox(0,      0)[t]{$\bullet$}}

\put(-14,       1){\makebox(0,      0)[t]{$\bullet$}}

\put(-14,      -1){\line(1,       0){50}}

\put(26,     38){\makebox(0,      0)[t]{$\bullet$}}

\put(-18,     - 15){$q_{11}$}
\put(0,       -15){$\widetilde{q}_{12}$}
\put(22,     - 15){$q_{22}$}
\put(40,       -15){$\widetilde{q}_{32}$}

\put(58,      - 15){$q_{33}$}


\put(30,       40){$q_{44}$}

\put(-12,       20){$\widetilde{q}_{14}$}

\put(58,       10){$\widetilde{q}_{34}$}

\put(-12,      1){\line(1,      1){35}}

\put(60,      1){\line(-1,      1){35}}

\put(80,         -1)  { over chain } \end{picture}$ {\ }\ \ \ \ \ \ \ \ \ \ \ {\ }\ \ \ \ \ \ \ \ \ \ \ {\ }\ \ \ \ \ \ \ \ \ \ \
 $\begin{picture}(100,       15)  \put(-45,       -1){ (b)}

\put(60,       1){\makebox(0,       0)[t]{$\bullet$}}

\put(28,       -1){\line(1,       0){33}}
\put(27,       1){\makebox(0,      0)[t]{$\bullet$}}

\put(-14,       1){\makebox(0,      0)[t]{$\bullet$}}

\put(-14,      -1){\line(1,       0){50}}

\put(-18,      10){$q_{11}$}
\put(0,       5){$\widetilde{q}_{12}$}
\put(22,      10){$q_{22}$}

\put(38,       5){$\widetilde{q}_{32}$}

\put(60,       10){$q_{33}$}

  \ \ \ \ \ \ \ \ \ \ \ \ \ \ \ \ \ \ \ \  {}
\end{picture}$\\ \\
consists of  two quasi-arithmetic {\rm  GDD}s:

{\ }\ \ \ \ \ \ \ \   $\begin{picture}(100,       15)
\put(-48,       -1){(c) }

\put(60,       1){\makebox(0,       0)[t]{$\bullet$}}

\put(28,       -1){\line(1,       0){33}}
\put(27,       1){\makebox(0,      0)[t]{$\bullet$}}

\put(-14,       1){\makebox(0,      0)[t]{$\bullet$}}

\put(-14,      -1){\line(1,       0){50}}
\put(104,       1){\makebox(0,     0)[t]{$\bullet$}}

\put(62,      -1){\line(1,       0){40}}

\put(-18,      10){$q_{44}$}
\put(0,       5){$\widetilde{q}_{43}$}
\put(22,      10){$q_{11}$}
\put(40,       5){$\widetilde{q}_{21}$}
\put(58,       10){$q_{22}$}

\put(75,      5){$\widetilde{q}_{23}$}

\put(100,       10){$q_{33}$}
 \ \ \ \ \ \ \ \ \
   \ \ \ \ \ \ \ \ \ \ \ \ \ \ \ \ \ \ \ {   }
\end{picture}${\ }\ \ \ \ \ \ \ \ \ \ \ \ {\ }\ \ \ \ \ \ \ \ \ \ \ \ {\ }\ \
 $\begin{picture}(100,       15)
\put(-68,       -1){and \ \ \ (d) }

\put(60,       1){\makebox(0,       0)[t]{$\bullet$}}

\put(28,       -1){\line(1,       0){33}}
\put(27,       1){\makebox(0,      0)[t]{$\bullet$}}

\put(-14,       1){\makebox(0,      0)[t]{$\bullet$}}

\put(-14,      -1){\line(1,       0){50}}
\put(104,       1){\makebox(0,     0)[t]{$\bullet$}}

\put(62,      -1){\line(1,       0){40}}

\put(-18,      10){$q_{11}$}
\put(0,       5){$\widetilde{q}_{12}$}
\put(22,      10){$q_{22}$}
\put(40,       5){$\widetilde{q}_{23}$}
\put(58,       10){$q_{33}$}

\put(75,      5){$\widetilde{q}_{34}$}

\put(100,       10){$q_{44}$}
 \ \ \ \ \ \ \ \ \ \ \ \ \ \ \ \ \ \ \
   \ \ \ \ \ \ \ \ \ \ \ \ \ \ \ \ \  {,  }
\end{picture}$\\
written as (c)(d) in short.

A  circle is called  a near-classical circle if it consists of two near-classical  {\rm  GDD}s.


\begin {Lemma} \label {3.1.1} Rank  $n=3$. {\rm (I)}   Non classical   {\rm  GDD}  is not an arithmetic {\rm  GDD} if there exist two vertexes  $i$  and $j$ with $q_{ii} \not= -1$
and $\widetilde{q}_{ij} \not= 1 $ such that  $q_{ii} \widetilde{q}_{ij} \not=1$ and $i$ is not a middle vertex when {\rm  GDD}  is a chain  except
GDD $5$ of Row $16$;
GDD $9$ of Row $17$;
 {\rm  GDD}  $1$ of Row $18$;
 {\rm  GDD}  $2$ of Row $18$;
 {\rm  GDD}  $8$ of Row $17$;
 {\rm  GDD}  $4$ of Row $7$;
 {\rm  GDD}  $4$ of Row $16$;
 {\rm  GDD}  $3$ of Row $7$;
GDD $3$ of Row $16$;
 {\rm  GDD}  $3$ of Row $17$;

{\rm (II)} Non classical  Chain {\rm  GDD}  is not  an arithmetic {\rm  GDD}   if
 $\widetilde{q}_{21} \widetilde{q}_{23} \not=1$,
 $q_{22} =-1$ and there exists another vertex another vertex  $i\not= 2$ such that $q_{ii} =-1$ and $\widetilde{q}_{i2} \not=-1$ except
  {\rm  GDD}  $2$ of Row $7$;
 {\rm  GDD}  $2$ of Row $15$;
  {\rm  GDD}  $2$ of Row $16$;
  {\rm  GDD}  $1$ of Row $17$;
 {\rm  GDD}  $4$ of Row $17$.

 {\rm (III)} Non classical chain   {\rm  GDD}  is not an arithmetic  {\rm  GDD}   if $q_{22} =-1$,
 $\widetilde{q}_{21} \widetilde{q}_{23} \not=1$  and
  $\widetilde{q}_{23} {q}_{33} =1$,   $q_{11}=-1$
 except
    {\rm  GDD}  $2$ of Row $7$;
   {\rm  GDD}  $2$ of Row $16$;
 {\rm  GDD}  $1$ of Row $9$;    $q^2 \not=1$,  $r= -1.$

 {\rm (IV)}  {\rm  GDD}  is an arithmetic  {\rm  GDD}   if $\widetilde{q}_{21} q_{11}=1$  and $\widetilde{q}_{23} q_{33}=1$ with $q_{22}=-1$.

 {\rm (V)}   {\rm  GDD}  is not an arithmetic  {\rm  GDD}   if there exist different three $i,  j,  k$ such that    $q_{ii} \not= -1$,  $q_{ii} \widetilde{q}_{ij} \not=1$
and   $q_{ii} \widetilde{q}_{ik} \not=1$ with  $\widetilde{q}_{ik} \not=1 $ and  $\widetilde{q}_{ij} \not=1 $ except
   {\rm  GDD}  $2$ of Row $13$;
    {\rm  GDD}  $4$ of Row $15$;
      {\rm  GDD}  $6$ of Row $17$.

 {\rm (VI)}
Chain  {\rm  GDD}  is not an arithmetic  {\rm  GDD}   if $q_{11} =-1$,  $\widetilde{q}_{12} =q_{22}\not= -1$  or $q_{33} =-1$,  $\widetilde{q}_{32} =q_{22}\not= -1$ except
 {\rm  GDD}  $1$ of Row $15$;
    {\rm  GDD}  $6$ of Row $17$;
  {\rm  GDD}  $1$ of Row $13, $ $q \in R_3$;
 {\rm  GDD}  $2$ of Row $13$,  $q \in R_6$.

 {\rm (VII)} Non classical
 {\rm  GDD}  is not an arithmetic  {\rm  GDD}   if  there exist two different $i,  j$ such that    $q_{ii} =-1=\widetilde{ q}_{ij} $ with the fixed parameter $q\in R_3$ in Table A2 except
  {\rm  GDD}  $1$ of Row $17$;
   {\rm  GDD}  $2$ of Row $17$;
     {\rm  GDD}  $9$ of Row $17$;
       {\rm  GDD}  $8$ of Row $17$;
 {\rm  GDD}  $1$ of Row $9$,   $q \in  R_{3}$,   $r =-1; $
  {\rm  GDD}  $3$ of Row $9$, $rq =-1, $   $s =-1, $  $r \not=q, $ $q \in R_3$;
         {\rm  GDD}  $3$ of Row $16$;
          {\rm  GDD}  $4$ of Row $16$.

 {\rm (VIII)} Chain
 {\rm  GDD}  is not an arithmetic  {\rm  GDD}   if   $q_{22} =- \widetilde{ q}_{21}\notin R_4 $  with $q_{11} =-1$ or $q_{22} =- \widetilde{ q}_{23}\notin R_4  $ with $q_{33} =-1$  except
  {\rm  GDD}  $2$ of Row $13$,  $q \in  R_{3}$;   {\rm  GDD}  $4$ of Row $15$;
 {\rm  GDD}  $5$ of Row $17$;    {\rm  GDD}  $6$ of Row $17$;   {\rm  GDD}  $7$ of Row $17$.

 {\rm (IX)} Chain
 {\rm  GDD}  is not an arithmetic  {\rm  GDD}   if   $q_{22} =q_{33}= -1$,  $\widetilde{q}_{32}\not=-1$ and $q_{11}\widetilde{q}_{12}=1 $ or    $q_{22} =q_{11}= -1$,  $\widetilde{q}_{12}\not=-1$ and $q_{33}\widetilde{q}_{32}=1 $ except
   {\rm  GDD}  $2$ of Row $4$;    {\rm  GDD}  $2$ of Row $6$;   {\rm  GDD}  $2$ of Row $7$;   {\rm  GDD}  $2$ of Row $16$;   {\rm  GDD}  $1$ of Row $17$;
  {\rm  GDD}  $4$ of Row $7$,  $q\in R_6$.

{\rm (X)}    Non classical  chain   {\rm  GDD}  is not  an arithmetic  {\rm  GDD}   if
 $q_{11}=q_{22} =q_{33}= -1$
except
  {\rm  GDD}  $2$ of Row $15$;    {\rm  GDD}  $4$ of Row $17$;
 {\rm  GDD}  $2$ of Row $7$,  $q \in R_6$.

 {\rm (XI)} Chain  {\rm  GDD}  is not an arithmetic  {\rm  GDD}   if $q_{22} \not=-1$  and $q_{11} =-1$,  $q_{22}\widetilde{q}_{12} \not= 1 $ or
 $q_{22} \not=-1$  and $q_{33} =-1$ and
$q_{22}\widetilde{q}_{32} \not= 1 $   except
  {\rm  GDD}  $1$ of Row $13$;   {\rm  GDD}  $2$ of Row $13$,   $q \in  R_{3}$;   {\rm  GDD}  $4$ of Row $15$;   {\rm  GDD}  $2$ of Row $17$;   {\rm  GDD}  $5$ of Row $17$;   {\rm  GDD}  $6$ of Row $17$;    {\rm  GDD}  $7$ of Row $17$;   {\rm  GDD}  $1$ of Row $15$.

 {\rm (XII)} Circle   {\rm  GDD}  is not an arithmetic  {\rm  GDD}   if $q_{ii} \not= -1$ for $ 1\le i \le 3.$

 {\rm (XIII)}  {\rm  GDD}  is not an arithmetic  {\rm  GDD}   if there exist distinct two vertexes  $i$  and $j$ with $q_{ii} \not= -1$
and $\widetilde{q}_{ij} \not= 1 $ such that  $q_{ii} \widetilde{q}_{ij} \not=1$ and $i$ is not a middle vertex when  {\rm  GDD}  is a chain,  as well as there exists other place where  simple chain condition (\ref {ppe3}) is not satisfied   except
   {\rm  GDD}  $2$ of Row $18$.

{\rm (XIV)}  Non classical chain  {\rm  GDD}    is not an arithmetic  {\rm  GDD}   if  $q_{22} ^{-2}= \widetilde{q}_{12} $  or  $q_{22} ^{-2}= \widetilde{q}_{32} $   except
  {\rm  GDD}  $2$ of Row $18$;    {\rm  GDD}  $2$ of Row $13$;   {\rm  GDD}  $1$ of Row $13$;    {\rm  GDD}  $4$ of Row $15$;    {\rm  GDD}  $5$ of Row $17$;    {\rm  GDD}  $6$ of Row $17$;    {\rm  GDD}  $1$ of Row $15$.

 {\rm (XV)} Circle  {\rm  GDD}    is not an arithmetic  {\rm  GDD}   if  $q_{ii} \not= -1$  for $1\le i \le 3$

 {\rm (XVI)} Circle  {\rm  GDD}    is not an arithmetic  {\rm  GDD}   if  only one vertex is $-1$ except
  {\rm  GDD}  $3$ of Row $15$;    {\rm  GDD}  $3$ of Row $17$.

 {\rm (XVII)} If $srq =1$ and $r= q^{-2}$,  then $ s = q.$

\end {Lemma}

{\bf Proof.} We only prove (IV) since others are clear.  We show it by three cases. Case one. If $\widetilde{q}_{21}\widetilde{q}_{23} =1$,  then  {\rm  GDD}  is an arithmetic  {\rm  GDD}    since it is a simple chain.  Case two. If $\widetilde{q}_{21}=\widetilde{q}_{23}$,  then  {\rm  GDD}  is an arithmetic  {\rm  GDD}   by   {\rm  GDD}  $1$ of Row $10$ ( $, q \in F^{*}\setminus \{1,  -1\}$ and $q\notin R_3$ ) and   {\rm  GDD}  $1$ of Row $11$,  $q\in R_3.$
Case three.  If $\widetilde{q}_{21} \not= \widetilde{q}_{23},   \widetilde{q}_{23}^{-1} $.
 then  {\rm  GDD}  is an arithmetic  {\rm  GDD}   by   {\rm  GDD}  $1$ of Row $9$.              \hfill $\Box$

\begin {Lemma} \label {3.1.2} Rank $n=3$.
 {\rm (I)} Non classical
   {\rm  GDD}  is not  an arithmetic  {\rm  GDD}   if
 $\widetilde{q}_{23} q_{33} = 1$ and  $q_{22} =-1$ with $q_{33} \not= -1$
 except
  {\rm  GDD}  $2$ of Row $7$,  $q^6 \not=1; $
  {\rm  GDD}  $3$ of Row $7$;   {\rm  GDD}  $4$ of Row $7$,   $q^6 \not=1; $   {\rm  GDD}  $1$ of Row $9$,  $r =-1, $ $q^2 \not=1; $   {\rm  GDD}  $3$ of Row $15$;    {\rm  GDD}  $2$ of Row $16$;   {\rm  GDD}  $3$ of Row $16$;   {\rm  GDD}  $4$ of Row $16$;    {\rm  GDD}  $3$ of Row $17$;     {\rm  GDD}  $8$ of Row $17$.

 {\rm (II)} Non classical   {\rm  GDD}  is not  an arithmetic  {\rm  GDD}  if
  $\widetilde{q}_{23} q_{33} = \widetilde{q}_{21} q_{22} = 1$ and $q_{33} \not= -1$
except
 {\rm  GDD}  $1$ of Row $13$;   {\rm  GDD}  $3$ of Row $15$;   {\rm  GDD}  $5$ of Row $16$;   {\rm  GDD}  $1$ of Row $18$.

 {\rm (III)}   Non classical     {\rm  GDD}  is not  an arithmetic  {\rm  GDD}   if there exist $i,  j$ such that $q_{ii} \not= -1$,  $q_{jj} = -1$,
 $\widetilde{q}_{ij} q_{ii} =  1$
 except
   {\rm  GDD}  $1$ of Row $7$;     {\rm  GDD}  $3$ of Row $7$;       {\rm  GDD}  $1$ of Row $15$;   {\rm  GDD}  $3$ of Row $15$;    {\rm  GDD}  $1$ of Row $16$;   {\rm  GDD}  $3$ of Row $16$;   {\rm  GDD}  $2$ of Row $17$;   {\rm  GDD}  $3$ of Row $17$;   {\rm  GDD}  $5$ of Row $17$;   {\rm  GDD}  $7$ of Row $17$;     {\rm  GDD}  $8$ of Row $17$.

 {\rm (IV)} Non classical   {\rm  GDD}  is not  an arithmetic  {\rm  GDD}   if $q_{22} =q_{33}= -1$  and $\widetilde{q}_{32} \not=-1$  except
    {\rm  GDD}  $2$ of Row $7$;    {\rm  GDD}  $3$ of Row $7$;   {\rm  GDD}  $3$ of Row $9$;       {\rm  GDD}  $2$ of Row $15$;    {\rm  GDD}  $2$ of Row $16$;    {\rm  GDD}  $3$ of Row $16$;   {\rm  GDD}  $1$ of Row $17$; {\rm  GDD}  $4$ of Row $17$.

 {\rm (V)} Non classical    {\rm  GDD}  is not  an arithmetic  {\rm  GDD}   if
 $\widetilde{q}_{23} =q_{22} =q_{33} =- 1$
 except
  {\rm  GDD}  $2$ of Row $7$, $q\in R_6; $
  {\rm  GDD}  $1$ of Row $17$,
 {\rm  GDD}  $9$ of Row $17$;    {\rm  GDD}  $8$ of Row $17$;   {\rm  GDD}  $4$ of Row $7$,  $q\in R_6$;
   {\rm  GDD}  $1$ of Row $9$,    $r= -1.$ $q^2 \not= 1$;
 {\rm  GDD}  $3$ of Row $9$,  $r,  q \not=-1, $ $q,  q^{-1} \not=r$,       $s= -1$;
 {\rm  GDD}  $3$ of Row $7$,  $q\in R_6$.

\end {Lemma}


\begin {Lemma} \label {3.1.3} Rank $n=4$.
Non classical  {\rm  GDD}  is not an arithmetic  {\rm  GDD}   if there exist two different vertexes  $i$  and $j$ with $q_{ii} \not= -1$
and $\widetilde{q}_{ij} \not= 1 $ such that  $q_{ii} \widetilde{q}_{ij} \not=1$ and $i$ is not a middle vertex (i.e. $i \not= 2, 3$) when  {\rm  GDD}  is a chain  except
  {\rm  GDD}  $1$ of Row $17$  in Table B;   {\rm  GDD}  $2$ of Row $17$  in Table B;   {\rm  GDD}  $5$ of Row $18$ in Table B;   {\rm  GDD}  $10$ of Row $20$ in Table B;   {\rm  GDD}  $4$ of Row $22$ in Table B;   {\rm  GDD}  $3$ of Row $17$ in Table B.

\end {Lemma}

\begin {Lemma} \label {3.1.4} Assume rank
$n=4, $ and $q^4,  q^6 \not =1, $ where $q$ is a fix  parameter in Table B.
Non classical
 {\rm  GDD}  is not an arithmetic  {\rm  GDD}   if there exist two different vertexes  $i$ and  $j$  with $q_{ii} = -1$,  $\widetilde{q}_{ij} = -1$ and ${q}_{jj} = -1$ except
 {\rm  GDD}  $2$ of Row $14$ in Table B;   {\rm  GDD}  $3$ of Row $14$ in Table B.

\end {Lemma}

\begin {Lemma} \label {3.1.7} Rank
$n=4.$
 {\rm  GDD}   is not an arithmetic  {\rm  GDD}   if  {\rm  GDD}  is\\

    {\ }\!\!\!\!\!\!\!\!\!\!\!\!\!\!\!\!\!\!\!\!\!\!\!\!\!\!\!\!\!  {\ }\!\!\!\!\!\!\!\!\!\!\!\!\!\!\!\!\!\!\!\!\!\!\!\!\!\!\!\!\!  {\ }\!\!\!\!\!\!\!\!\!\!\!\!\!\!\!\!\!\!\!\!\!   $\begin{picture}(100,       15)\put(-68,       -1){ }

\put(111,      1){\makebox(0,      0)[t]{$\bullet$}}
\put(144,       1){\makebox(0,       0)[t]{$\bullet$}}
\put(170,      -11){\makebox(0,      0)[t]{$\bullet$}}
\put(170,     15){\makebox(0,      0)[t]{$\bullet$}}
\put(113,      -1){\line(1,      0){33}}
\put(142,     -1){\line(2,      1){27}}
\put(170,       -14){\line(-2,      1){27}}

\put(170,       -14){\line(0,      1){27}}

\put(100,       10){$$}
\put(120,       5){$$}

\put(140,      10){$$}

\put(140,      -20){$$}
\put(145,       15){$$}

\put(178,       -20){$$}

\put(178,       0){$$}

\put(175,       10){$$}

\put(225,         -1)  {, with   $\widetilde{q}_{12}q_{22} \not=1, $ and $q_{22} \not= -1$,
except {\rm  GDD}  $4$ of Row $22$ in Table B. }

\end{picture}$

\end {Lemma}

\begin {Lemma} \label {3.1.8} Rank
$n=4.$ Non classical
chain  {\rm  GDD}  is not an arithmetic  {\rm  GDD}   if  there exists a vertex  $i$  with $q_{ii} = -1$,  $\widetilde{q}_{i,  i+1} \widetilde{q}_{i,  i-1}\not= 1$ and there exists other vertex which does not satisfy  simple chain condition (\ref {ppe2}) or condition (\ref {ppe3}). except
   {\rm  GDD}  $6$ of Row $9$  in Table B;   {\rm  GDD}  $2$ of Row $17$;   {\rm  GDD}  $5$ of Row $17$;   {\rm  GDD}  $6$ of Row $22$;   {\rm  GDD}  $7$ of Row $22$;   {\rm  GDD}  $4$ of Row $14$.

\end {Lemma}

\begin {Lemma} \label {3.1.9} Rank
$n=4. $
 Non classical
 {\rm  GDD}
is not an arithmetic  {\rm  GDD}   if  {\rm  GDD}  is

 {\ }\!\!\!\!\!\!\!\!\!\!\!\!\!\!\!\!\!\!\!\!\!\!\!\!\!\!\!\!\!  {\ }\!\!\!\!\!\!\!\!\!\!\!\!\!\!\!\!\!\!\!\!\!\!\!\!\!\!\!\!\!  {\ }\!\!\!\!\!\!\!\!\!\!\!\!\!\!\!\!\!\!\!\!\! $\begin{picture}(100,       15)\put(-68,       -1){ }
\put(111,      1){\makebox(0,      0)[t]{$\bullet$}}
\put(144,       1){\makebox(0,       0)[t]{$\bullet$}}
\put(170,      -11){\makebox(0,      0)[t]{$\bullet$}}
\put(170,     15){\makebox(0,      0)[t]{$\bullet$}}
\put(113,      -1){\line(1,      0){33}}
\put(142,     -1){\line(2,      1){27}}
\put(170,       -14){\line(-2,      1){27}}

\put(100,       10){$$}

\put(115,       5){$$}

\put(135,      10){$$}

\put(140,      -20){$$}
\put(145,       15){$$}

\put(173,       -12){$$}
\put(173,       18){$$}

\end{picture}$
  {\ }\ \ \ \ \ \ \ \ \ \ \ \
     $\begin{picture}(100,       15)\put(68,       -1){ or}

\put(111,      1){\makebox(0,      0)[t]{$\bullet$}}
\put(144,       1){\makebox(0,       0)[t]{$\bullet$}}
\put(170,      -11){\makebox(0,      0)[t]{$\bullet$}}
\put(170,     15){\makebox(0,      0)[t]{$\bullet$}}
\put(113,      -1){\line(1,      0){33}}
\put(142,     -1){\line(2,      1){27}}
\put(170,       -14){\line(-2,      1){27}}

\put(170,       -14){\line(0,      1){27}}

\put(100,       10){$$}
\put(120,       5){$$}

\put(127,      10){$$}

\put(140,      -20){$$}
\put(145,       15){$$}

\put(178,       -20){$$}

\put(178,       0){$$}

\put(175,       10){$$}

\put(245,         -1)  {   with $q_{ii} \not=-1$ for $i =1,  3,  4.$  }

\end{picture}$

\end {Lemma}

\begin {Lemma} \label {3.1.10}
$n=4, $ Non classical
 {\rm  GDD}  is not an arithmetic  {\rm  GDD}   if $q_{ii} = -1$ for $1\le i \le 4$,  except
  {\rm  GDD}  $5$ of Row $22$;   {\rm  GDD}  $4$ of Row $9$,    $q \in  R_{4}$;   {\rm  GDD}  $2$ of Row $9$,    $q \in  R_{4}$;   {\rm  GDD}  $3$ of Row $20$.

\end {Lemma}

\begin {Lemma} \label {3.1.12}
$n=4, $
Non classical chain  {\rm  GDD}  is not an arithmetic  {\rm  GDD}   if $q_{44} = q_{22}=-1$,  $\widetilde{q}_{12} \widetilde{q}_{23} \not=1$ and $q_{11} \widetilde{q}_{12} =1$,  or $q_{11} = q_{33}=-1$,  $\widetilde{q}_{43} \widetilde{q}_{23} \not=1$ and $q_{44} \widetilde{q}_{43} =1$,
 except
 {\rm  GDD}  $3$ of Row $22$  in Table B;    {\rm  GDD}  $5$ of Row $9$,   $q \in  R_{4}$;   {\rm  GDD}  $6$ of Row $9$,   $q \in  R_{4}$;   {\rm  GDD}  $6$ of Row $22$;   {\rm  GDD}  $7$ of Row $22$;   {\rm  GDD}  $2$ of Row $22$.

\end {Lemma}

\begin {Lemma} \label {3.1.13}

 {\rm  GDD}
is not an arithmetic  {\rm  GDD}   if  {\rm  GDD}  is\\

{\ }\!\!\!\!\!\!\!\!\!\!\!\!\!\!\!\!\!\!\!\!\!\!\!\!\!\!\!\!\!\!\!\!\!\!\!\!\!\!\!\!\!\!\!\!\!\!\!\!
    $\begin{picture}(100,       15)\put(-68,       -1){ }

\put(111,      1){\makebox(0,      0)[t]{$\bullet$}}
\put(144,       1){\makebox(0,       0)[t]{$\bullet$}}
\put(170,      -11){\makebox(0,      0)[t]{$\bullet$}}
\put(170,     15){\makebox(0,      0)[t]{$\bullet$}}
\put(113,      -1){\line(1,      0){33}}
\put(142,     -1){\line(2,      1){27}}
\put(170,       -14){\line(-2,      1){27}}

\put(170,       -14){\line(0,      1){27}}

\put(100,       10){$$}
\put(120,       5){$$}

\put(127,      10){$-1$}

\put(140,      -20){$$}
\put(145,       15){$$}

\put(178,       -20){$-1$}

\put(178,       0){$$}

\put(175,       10){$q_{33}$}

\put(195,         -1)  {,   }
\end{picture}$
 {\ }\ \ \ \ \ \ \ \ \ \ \ \   {\ }\ \ \ \ \ \ \
 {\ }\ \ \ \ \ \ \ \ \ \ \ \   $\begin{picture}(100,       15)\put(28,       -1){  $q_{33} \not=-1$ \ \ \  or}

\put(111,      1){\makebox(0,      0)[t]{$\bullet$}}
\put(144,       1){\makebox(0,       0)[t]{$\bullet$}}
\put(170,      -11){\makebox(0,      0)[t]{$\bullet$}}
\put(170,     15){\makebox(0,      0)[t]{$\bullet$}}
\put(113,      -1){\line(1,      0){33}}
\put(142,     -1){\line(2,      1){27}}
\put(170,       -14){\line(-2,      1){27}}

\put(170,       -14){\line(0,      1){27}}

\put(100,       10){$$}
\put(120,       5){$$}

\put(127,      10){$-1$}

\put(140,      -20){$$}
\put(145,       15){$$}

\put(178,       -20){$q_{44}$}

\put(178,       0){$$}

\put(175,       10){$-1$}

\put(195,         -1)  {$q_{44} \not=-1$ except }
\end{picture}$\\ \\
 {\rm  GDD}  $3$ of Row $9$ in Table B;   {\rm  GDD}  $3$ of Row $17$ in Table B;   {\rm  GDD}  $3$ of Row $18$;   {\rm  GDD}  $8$ of Row $20$;   {\rm  GDD}  $8$ of Row $22$.

\end {Lemma}

\begin {Lemma} \label {3.1.16}

 {\rm  GDD}
is not an arithmetic  {\rm  GDD}   if  {\rm  GDD}  is\\

   \ \ \ \ \ \  {\ }\ \ \ \ \ \ \ \ \ \ \ \   $\begin{picture}(100,       15)\put(-68,       -1){ }
\put(111,      1){\makebox(0,      0)[t]{$\bullet$}}
\put(144,       1){\makebox(0,       0)[t]{$\bullet$}}
\put(170,      -11){\makebox(0,      0)[t]{$\bullet$}}
\put(170,     15){\makebox(0,      0)[t]{$\bullet$}}
\put(113,      -1){\line(1,      0){33}}
\put(142,     -1){\line(2,      1){27}}
\put(170,       -14){\line(-2,      1){27}}

\put(100,       10){$r$}

\put(115,       5){$r^{-1}$}

\put(135,      10){$-1$}

\put(140,      -20){$q$}
\put(145,       15){$$}

\put(173,       -12){$-1$}
\put(173,       18){$$}

\put(195,         -1) { $r \not= q,  q^{-1}$. }
\end{picture}$\\

\end {Lemma}

\begin {Lemma} \label {3.1.18}
Non classical
 {\rm  GDD}
is not an arithmetic  {\rm  GDD}   if  {\rm  GDD}  is
\\

 {\ }\!\!\!\!\!\!\!\!\!\!\!\!\!\!\!\!\!\!\!\!\!\!\!\!\!\!\!\!\!\!\!\!\!\!\!\!\!\!\!\!\!\!\!\!\!\!\!\!
 $\begin{picture}(100,       15)\put(-68,       -1){ }

\put(111,      1){\makebox(0,      0)[t]{$\bullet$}}
\put(144,       1){\makebox(0,       0)[t]{$\bullet$}}
\put(170,      -11){\makebox(0,      0)[t]{$\bullet$}}
\put(170,     15){\makebox(0,      0)[t]{$\bullet$}}
\put(113,      -1){\line(1,      0){33}}
\put(142,     -1){\line(2,      1){27}}
\put(170,       -14){\line(-2,      1){27}}

\put(170,       -14){\line(0,      1){27}}

\put(100,       10){$r$}
\put(114,       5){$r^{-1}$}

\put(130,      10){$-1$}

\put(140,      -20){$q$}
\put(155,       15){$$}

\put(178,       -20){$-1$}

\put(178,       0){$$}

\put(175,       10){$$}

\put(195,         -1)  { $r \not= q,  q^{-1}$ except  {\rm  GDD}  $4$ of Row $9$;
  {\rm  GDD}  $2$ of Row $14$; }
\end{picture}$\\ \\
  {\rm  GDD}  $3$ of Row $14$;
  {\rm  GDD}  $4$ of Row $17$;   {\rm  GDD}  $8$ of Row $22$.

\end {Lemma}

\begin {Lemma} \label {3.1.19} $n=4, $
Non classical chain  {\rm  GDD}
is not an arithmetic  {\rm  GDD}   if  {\rm  GDD}  is of  $q_{22}=-1$ and $q_{11} \not=-1$ with $\widetilde{q}_{12}= \widetilde{q}_{23}$ or $q_{33}=-1$ and $q_{44} \not=-1$ with $\widetilde{q}_{43}= \widetilde{q}_{23}$ except
 {\rm  GDD}  $3$ of Row $22$;
 {\rm  GDD}  $2$ of Row $22$;    {\rm  GDD}  $5$ of Row $9$,    $q \in  R_{5}$;   {\rm  GDD}  $6$ of Row $9$,   $q \in  R_{4}$;   {\rm  GDD}  $1$ of Row $13$,  $q \in  R_{3}$;    {\rm  GDD}  $1$ of Row $14$,   $q\in R_4; $    {\rm  GDD}  $5$ of Row $14$,   $q\in R_4; $      {\rm  GDD}  $6$ of Row $17$.

\end {Lemma}

\begin {Lemma} \label {3.1.21} Rank $n=4, $
Non classical chain  {\rm  GDD}
is not an arithmetic  {\rm  GDD}   if there exist $i$  such that  $\widetilde{q}_{i+1,  i} \widetilde{q}_{i,  i-1}\not=1$ and $q_{ii} =-1$ with  $q_{i-1,  i-1} =-1$,  or  $q_{i+1,  i+1} =-1$,
 except
  {\rm  GDD}  $5$ of Row $9$,   $q \in  R_{4}$;   {\rm  GDD}  $3$ of Row $20$;   {\rm  GDD}  $4$ of Row $20$;   {\rm  GDD}  $2$ of Row $22$;   {\rm  GDD}  $7$ of Row $22$;   {\rm  GDD}  $2$ of Row $9$;
   {\rm  GDD}  $2$ of Row $14$,
  $q \in  R_{3}$;
   {\rm  GDD}  $6$ of Row $21$,
  $q \in  R_{4}$;    {\rm  GDD}  $5$ of Row $9$,   $q \in  R_{6}$;
  {\rm  GDD}  $6$ of Row $9$,   $q \in  R_{6}$;   {\rm  GDD}  $2$ of Row $18$.

\end {Lemma}

\begin {Lemma} \label {3.1.22} Rank  $n=4, $ Non classical
 {\rm  GDD}
is not an arithmetic  {\rm  GDD}   if  {\rm  GDD}
\\

 {\ }\!\!\!\!\!\!\!\!\!\!\!\!\!\!\!\!\!\!\!\!\!\!\!\!\!\!\!\!\!\!\!\!\!\!\!\!\!\!\!\!\!\!\!\!\!\!\!\! {\ }\!\!\!\!\!\!\!\!\!\!\!\!\!\!\!\!\!\!\!\!\!\!\!\!\!\!
 $\begin{picture}(100,       15)\put(-68,       -1){ }

\put(111,      1){\makebox(0,      0)[t]{$\bullet$}}
\put(144,       1){\makebox(0,       0)[t]{$\bullet$}}
\put(170,      -11){\makebox(0,      0)[t]{$\bullet$}}
\put(170,     15){\makebox(0,      0)[t]{$\bullet$}}
\put(113,      -1){\line(1,      0){33}}
\put(142,     -1){\line(2,      1){27}}
\put(170,       -14){\line(-2,      1){27}}

\put(170,       -14){\line(0,      1){27}}

\put(100,       10){$$}
\put(120,       5){$$}

\put(127,      10){$-1$}

\put(140,      -20){$$}
\put(155,       15){$$}

\put(178,       -20){$-1$}

\put(178,       0){$$}

\put(175,       10){$-1$}

\put(195,         -1)  {,  with $\widetilde{q}_{12}\widetilde{q}_{23}\not=1$ or $\widetilde{q}_{12}\widetilde{q}_{24}\not=1$,  \ \ \  or}
\end{picture}$\\ \\

  {\ }\!\!\!\!\!\!\!\!\!\!\!\!\!\!\!\!\!\!\!\!\!\!\!\!\!\!\!\!\!\!\!\!\!\!\!\!\!\!\!\!\!\!\!\!\!\!\!\! {\ }\!\!\!\!\!\!\!\!\!\!\!\!\!\!\!\!\!\!\!\!\!\!\!\!\!\!\!
 $\begin{picture}(100,       15)\put(-68,       -1){ }
\put(111,      1){\makebox(0,      0)[t]{$\bullet$}}
\put(144,       1){\makebox(0,       0)[t]{$\bullet$}}
\put(170,      -11){\makebox(0,      0)[t]{$\bullet$}}
\put(170,     15){\makebox(0,      0)[t]{$\bullet$}}
\put(113,      -1){\line(1,      0){33}}
\put(142,     -1){\line(2,      1){27}}
\put(170,       -14){\line(-2,      1){27}}

\put(100,       10){$$}

\put(115,       5){$$}

\put(135,      10){$-1$}

\put(140,      -20){$$}
\put(145,       15){$$}

\put(173,       -12){$-1$}
\put(173,       18){$-1$}

\put(195,         -1)  { with $\widetilde{q}_{12}\widetilde{q}_{23}\not=1$ or $\widetilde{q}_{12}\widetilde{q}_{24}\not=1$ except  {\rm  GDD}  $4$ of Row $17$;
 {\rm  GDD}  $4$ of Row $9$;  }
\end{picture}$\\ \\
 {\rm  GDD}  $3$ of Row $14$;   {\rm  GDD}  $2$ of Row $14$;
 {\rm  GDD}  $4$ of Row $18$;
 {\rm  GDD}  $5$ of Row $22$.
\end {Lemma}

\begin {Lemma} \label {3.1.23} Rank
$n=4. $
Non classical
 {\rm  GDD}  is not an arithmetic  {\rm  GDD}   if  there exist three different  vertexex  $i$,  $j$ and $k$  such that  $q_{ii} \not= -1$,  $\widetilde{q}_{i,  j} {q}_{i,  i}\not= 1$ and
$\widetilde{q}_{i,  k} {q}_{i,  i}\not= 1$ with $\widetilde{q}_{i,  j} \not= 1$ and
$\widetilde{q}_{i,  k} \not= 1$ except
 {\rm  GDD}  $4$ of Row $21$;   {\rm  GDD}  $2$ of Row $20$;   {\rm  GDD}  $2$ of Row $21$.

\end {Lemma}

\begin {Lemma} \label {3.1.25} Rank
$n=4. $
Non classical
 {\rm  GDD}  is not an arithmetic  {\rm  GDD}   if   $q_{11} = -1$,  $\widetilde{q}_{12} = {q}_{22}\not= -1$ or  $q_{44} = -1$,  $\widetilde{q}_{43} = {q}_{33}$,
 except
 {\rm  GDD}  $1$ of Row $20$;   {\rm  GDD}  $1$ of Row $21$;   {\rm  GDD}  $3$ of Row $21$.

\end {Lemma}

\begin {Lemma} \label {3.1.26} Rank
$n=4,  $ Non classical
chain
 {\rm  GDD}  is not an arithmetic  {\rm  GDD}   if   $q_{11} = -1$,  $\widetilde{q}_{12} = -{q}_{22}$ or  $q_{44} = -1$,  $\widetilde{q}_{43} = -{q}_{33}$ with  $q ^4\not=1$,
 except
 {\rm  GDD}  $2$ of Row $20$;   {\rm  GDD}  $2$ of Row $21$;
 {\rm  GDD}  $4$ of Row $21$.

\end {Lemma}
\begin {Lemma} \label {3.1.27} Rank
$n=4. $ Non classical
 {\rm  GDD}  is not an arithmetic  {\rm  GDD}   if\\
  {\ }\!\!\!\!\!\!\!\!\!\!\!\!\!\!\!\!\!\!\!\!\!\!\!\!\!\!\!\!\!\!\!\!\!\!\!\!\!\!\!\!\!\!\!\!\!\!\!\! {\ }\!\!\!\!\!\!\!\!\!\!\!\!\!\!\!\!\!\!\!\!\!\!\!\!\!\!\!\!\!\!\!\!\!\!\!\!\!\!\!\!\!\!\!\!\!\!\!\! {\ }\!\!\!\!\!\!\!\!\!\!\!\!\!\!\!\!\!\!\!\!\!\!\!\!\!\!\!\!\!\!\!\!\!\!\!\!\!\!\!\!\!\!\!\!\!\!\!\!
  $\begin{picture}(100,       15)\put(-68,       -1){ }
\put(111,      1){\makebox(0,      0)[t]{$\bullet$}}
\put(144,       1){\makebox(0,       0)[t]{$\bullet$}}
\put(170,      -11){\makebox(0,      0)[t]{$\bullet$}}
\put(170,     15){\makebox(0,      0)[t]{$\bullet$}}
\put(113,      -1){\line(1,      0){33}}
\put(142,     -1){\line(2,      1){27}}
\put(170,       -14){\line(-2,      1){27}}

\put(100,       10){$$}

\put(115,       5){$$}

\put(130,      10){$$}

\put(140,      -20){$$}
\put(145,       15){$$}

\put(173,       -12){$$}
\put(173,       18){$$}

\put(195,         -1)  { }
\end{picture}$ {\ }\ \ \ \ \ \ \ \ \ \ \ \  {\ }\ \ \ \ \ \ \ \ \ \ \ \
 {\ }\ \ \ \ \ \ \ \ \ \ \ \   $\begin{picture}(100,       15)\put(-68,       -1){ or}

\put(111,      1){\makebox(0,      0)[t]{$\bullet$}}
\put(144,       1){\makebox(0,       0)[t]{$\bullet$}}
\put(170,      -11){\makebox(0,      0)[t]{$\bullet$}}
\put(170,     15){\makebox(0,      0)[t]{$\bullet$}}
\put(113,      -1){\line(1,      0){33}}
\put(142,     -1){\line(2,      1){27}}
\put(170,       -14){\line(-2,      1){27}}

\put(170,       -14){\line(0,      1){27}}

\put(100,       10){$$}
\put(120,       5){$$}

\put(127,      10){$$}

\put(140,      -20){$$}
\put(155,       15){$$}

\put(178,       -20){$$}

\put(178,       0){$$}

\put(175,       10){$$}

\put(195,         -1)  {  }
\end{picture}$\\
\\
 with
  $q_{ii} = -1$ for $1\le i \le 4$ and  $\widetilde{q}_{12} = -\widetilde{q}_{23}$ or  $\widetilde{q}_{24} = -\widetilde{q}_{12}$,
 except
 {\rm  GDD}  $5$ of Row $22$.

\end {Lemma}

\begin {Lemma} \label {3.1.28} Rank
$n=4. $ Non classical
 {\rm  GDD}  is not an arithmetic  {\rm  GDD}   if  {\rm  GDD}  contains a subGDD    {\ }\ \ \ \ \ \ \ \ \ \ \ \
    $\begin{picture}(100,       15)  \put(-68,       -1){ }

\put(60,       1){\makebox(0,       0)[t]{$\bullet$}}

\put(28,       -1){\line(1,       0){33}}
\put(27,       1){\makebox(0,      0)[t]{$\bullet$}}

\put(-14,       1){\makebox(0,      0)[t]{$\bullet$}}

\put(-14,      -1){\line(1,       0){50}}

\put(-18,      10){$-1$}
\put(0,       5){$q$}
\put(22,      10){$-1$}
\put(40,       5){$q$}

\put(58,       10){$-1$}

  \ \ \ \ \ \ \ \ \ \ \ \ \ \ \ \ \ \ \ except  {\rm  GDD}  $3$ of Row $20$;   {\rm  GDD}  $2$ of Row $9$,   $q \in  R_{4}$;
\end{picture}$
\\  {\rm  GDD}  $6$ of Row $21$;
 {\rm  GDD}  $5$ of Row $9$,   $q \in  R_{4}$;  {\rm  GDD}  $3$ of Row $9$,   $q \in  R_{4}$;   {\rm  GDD}  $4$ of Row $9$,   $q \in  R_{4}$.

\end {Lemma}

\begin {Lemma} \label {3.1.30} Rank
$n=4. $ Non classical   chain  {\rm  GDD}  is not an arithmetic  {\rm  GDD}   if there exists $i$ such that  $q_{ii} = -1$ and $q_{jj} \not= -1$ for $j \not= i$ with $\widetilde{q}_{i,  i-1} = \widetilde{q}_{i,  i+1}$,  except
  {\rm  GDD}  $5$ of Row $9$,   $q \in  R_{5}$;    {\rm  GDD}  $1$ of Row $14$;  $q \in  R_{4}$;   {\rm  GDD}  $5$ of Row $14$,  $q \in  R_{4}$;   {\rm  GDD}  $6$ of Row $17$.

 \end {Lemma}

  \begin {Lemma} \label {3.1.33} Rank
$n=4. $ Non classical   chain  {\rm  GDD}  is not an arithmetic  {\rm  GDD}   if
$q_{22} =-1, $ $\widetilde{q}_{12}\widetilde{q}_{23} \not=1$ and $\widetilde{q}_{12}\not= \widetilde{q}_{23},  $ or $q_{33} =-1, $ $\widetilde{q}_{43}\widetilde{q}_{23} \not=1$ and $\widetilde{q}_{34}\not= \widetilde{q}_{23},  $ except
 {\rm  GDD}  $5$ of Row $9$;   {\rm  GDD}  $1$ of Row $14$,  $q \notin  R_{2}\cup R_{4}$;   {\rm  GDD}  $5$ of Row $14$,  $q \notin  R_{2}\cup R_{4}$;     {\rm  GDD}  $2$ of Row $9$;   {\rm  GDD}  $6$ of Row $9$,   $q^4,  q^3 \not=1; $    {\rm  GDD}  $2$ of Row $17$;   {\rm  GDD}  $5$ of Row $17$;   {\rm  GDD}  $6$ of Row $22$;   {\rm  GDD}  $7$ of Row $22$;   {\rm  GDD}  $4$ of Row $14$.

 \end {Lemma}

  \begin {Lemma} \label {3.1.34} Rank
$n=4. $

  {\rm  GDD}    {\ }\!\!\!\!\!\!\!\!\!\!\!\!\!\!\!\!\!\!\!\!\!\!\!\!\!\!\!\!\!\!\!\!\!\!\!\!\!\!\!\!\!\!\!\!\!\!\!\! {\ }\!\!\!\!\!\!\!\!\!\!\!\!\!\!\!\!\!\!\!\!\!\!\!\!\!\!\!{\ }\ \ \ \ \ \ \ \ \ \ \ \   $\begin{picture}(100,       15)\put(-68,       -1){ }

\put(111,      1){\makebox(0,      0)[t]{$\bullet$}}
\put(144,       1){\makebox(0,       0)[t]{$\bullet$}}
\put(170,      -11){\makebox(0,      0)[t]{$\bullet$}}
\put(170,     15){\makebox(0,      0)[t]{$\bullet$}}
\put(113,      -1){\line(1,      0){33}}
\put(142,     -1){\line(2,      1){27}}
\put(170,       -14){\line(-2,      1){27}}

\put(170,       -14){\line(0,      1){27}}

\put(100,       10){$$}
\put(120,       5){$$}

\put(127,      10){$$}

\put(140,      -20){$$}
\put(155,       15){$$}

\put(178,       -20){$$}

\put(178,       0){$$}

\put(175,       10){$$}

\put(195,         -1)  { is not an arithmetic  {\rm  GDD}   if
  $q_{33}\not= -1$ and $q_{44}\not= -1$.}
\end{picture}$

\end {Lemma}

 \begin {Lemma} \label {3.1.35} Rank
$n=4. $ Non classical
  {\rm  GDD}  {\ }\!\!\!\!\!\!\!\!\!\!\!\!\!\!\!\!\!\!\!\!\!\!\!\!\!\!\!\!\!\!\!\!\!\!\!\!\!\!\!\!\!\!\!\!\!\!\!\! {\ }\!\!\!\!\!\!\!\!\!\!\!\!\!\!\!\!\!\!  $\begin{picture}(100,       15)\put(-68,       -1){ }
\put(111,      1){\makebox(0,      0)[t]{$\bullet$}}
\put(144,       1){\makebox(0,       0)[t]{$\bullet$}}
\put(170,      -11){\makebox(0,      0)[t]{$\bullet$}}
\put(170,     15){\makebox(0,      0)[t]{$\bullet$}}
\put(113,      -1){\line(1,      0){33}}
\put(142,     -1){\line(2,      1){27}}
\put(170,       -14){\line(-2,      1){27}}

\put(100,       10){$$}

\put(115,       5){$$}

\put(135,      10){$$}

\put(140,      -20){$$}
\put(145,       15){$$}

\put(173,       -12){$$}
\put(173,       18){$$}

\put(195,         -1)  { is not an arithmetic  {\rm  GDD}  }
\end{picture}$\\ \\
 if there are two among $q_{11},   q_{33}$ and $q_{44}$ which are not equal to $-1, $ except
  {\rm  GDD}  $9$ of Row $20$;   {\rm  GDD}  $10$ of Row $20$.

 \end {Lemma}

 \begin {Lemma} \label {3.1.36} Rank
$n=4. $  Non classical
 {\rm  GDD}  is not an arithmetic  {\rm  GDD}
  if there exist two different $i$ and $j$ such that
$\widetilde{q}_{ij} = -1$ and $q^4 \not=1$,  where $q$ is a fixed  parameter in Table B except
 {\rm  GDD}  $2$ of Row $17$;   {\rm  GDD}  $3$ of Row $17$;   {\rm  GDD}  $2$ of Row $14$,   $q^4 \not=1; $  {\rm  GDD}  $3$ of Row $14$,  $q^4 \not=1; $ {\rm  GDD}  $4$ of Row $14$,   $q \in  R_{3}$;   {\rm  GDD}  $4$ of Row $9$,    $q \in  R_{6}$;   {\rm  GDD}  $5$ of Row $9$,   $q \in  R_{6}$;   {\rm  GDD}  $6$ of Row $9$,   $q \in  R_{6}$.

 \end {Lemma}

\begin {Lemma} \label {6.2.2'} Rank
 $n=4$. Non classical  chain  {\rm  GDD}  is not  an arithmetic  {\rm  GDD}   if
 $\widetilde{q}_{i,  i+1} \widetilde{q}_{i,  i-1} =-1$,  and
$q_{ii} =-1$,  except
 {\rm  GDD}  $6$ of Row $9$,   $q \in  R_{4}$;   {\rm  GDD}  $1$ of Row $14$;   {\rm  GDD}  $5$ of Row $14$;   {\rm  GDD}  $2$ of Row $22$;   {\rm  GDD}  $3$ of Row $22$;   {\rm  GDD}  $6$ of Row $9$,   $q \in  R_{4}$.

\end  {Lemma}

\begin  {Lemma} \label {6.2.2''}  Rank $n=4$.   Non classical   {\rm  GDD}  is not  an arithmetic  {\rm  GDD}   if  {\rm  GDD}  contains sub- {\rm  GDD}
  {\ }\ \ \ \ \ \ \ \ \ \ \ \    $\begin{picture}(100,       15)  \put(-68,       -1){ }

\put(60,       1){\makebox(0,       0)[t]{$\bullet$}}

\put(28,       -1){\line(1,       0){33}}
\put(27,       1){\makebox(0,      0)[t]{$\bullet$}}

\put(-14,       1){\makebox(0,      0)[t]{$\bullet$}}

\put(-14,      -1){\line(1,       0){50}}

\put(-18,      10){$q^{}$}
\put(0,       5){$q^{-1}$}
\put(22,      10){$-1{}$}
\put(40,       5){$q^{-1}$}

\put(58,       10){$q^{}$}

\ \ \ \ \ \ \ \ \ \ \ \
  \ \ \ \ \ \
except  {\rm  GDD}  $6$ of Row $17$;   {\rm  GDD}  $5$ of Row $9$,   $q \in  R_{5}$;

\end{picture}$\\  {\rm  GDD}  $6$ of Row $9$,   $q \in  R_{4}$;
 {\rm  GDD}  $3$ of Row $12$,   $q \in  R_{3}$;   {\rm  GDD}  $1$ of Row $13$,   $q \in  R_{3}$;   {\rm  GDD}  $1$ of Row $14$,   $q \in  R_{4}$;   {\rm  GDD}  $5$ of Row $14$,   $q \in  R_{4}$;   {\rm  GDD}  $3$ of Row $22$;   {\rm  GDD}  $7$ of Row $21$.

\end  {Lemma}

\begin  {Lemma} \label {3.2.5'''} Rank
 $n=4$.    Non classical  {\rm  GDD}  is not  an arithmetic  {\rm  GDD}   if  the fixed parameter $q\notin R_4$  and  {\rm  GDD}  contains sub {\rm  GDD}
    {\ }\ \ \ \ \ \ \ \ \ \ \ \  $\begin{picture}(100,       15)  \put(-85,       -1){}

\put(27,       1){\makebox(0,      0)[t]{$\bullet$}}

\put(-14,       1){\makebox(0,      0)[t]{$\bullet$}}

\put(-14,      -1){\line(1,       0){45}}

\put(-18,      10){$-1$}
\put(0,       5){$-1$}
\put(22,      10){$-1$}
\put(40,       5){$$}

\put(58,       10){$$}

  \ \ \ \ \ \ \ \ \ \ \ \ \ \ \ \ \ \ \ {except  {\rm  GDD}  $2$ of Row $14$; }
\end{picture}$
\\   {\rm  GDD}  $3$ of Row $14$;   {\rm  GDD}  $4$ of Row $14$;
 {\rm  GDD}  $5$ of Row $9$,   $q \in  R_{6}$;   {\rm  GDD}  $6$ of Row $9$,   $q \in  R_{6}$;   {\rm  GDD}  $4$ of Row $9$,   $q \in  R_{6}$.

\end  {Lemma}

\begin {Lemma} \label {3.2.32}

 $n=4$. Non classical     chain  {\rm  GDD}  is not  an arithmetic  {\rm  GDD}   if
 there exists $i$ such that $q_{ii} = -1$ and
 $\widetilde{q}_{i,  i+1} = \widetilde{q}_{i,  i-1} \not=-1$,  except
 {\rm  GDD}  $1$ of Row $14$,   $q \in  R_{4}$;   {\rm  GDD}  $5$ of Row $14$,   $q \in  R_{4}$;   {\rm  GDD}  $6$ of Row $17$;   {\rm  GDD}  $2$ of Row $18$;   {\rm  GDD}  $3$ of Row $20$;   {\rm  GDD}  $4$ of Row $20$;   {\rm  GDD}  $6$ of Row $21$;   {\rm  GDD}  $2$ of Row $22$;   {\rm  GDD}  $3$ of Row $22$;   {\rm  GDD}  $5$ of Row $9$,   $q \in  R_{5}$;   {\rm  GDD}  $6$ of Row $9$,   $q \in  R_{4}$.

\end  {Lemma}

\begin {Lemma} \label {3.2.33} Rank
 $n=4$. Non classical  chain  {\rm  GDD}  is not  an arithmetic  {\rm  GDD}   if there exists only one place where simple chain condition   is not satisfied and
$q_{22} \not=-1$,  $q_{33} \not=-1$ with
 $\widetilde{q}_{23}q_{22} \not= 1$ or $ q_{33} \widetilde{q}_{34} \not=1$,
 or $ q_{22} \widetilde{q}_{21} \not=1$,  or $ q_{33} \widetilde{q}_{32} \not=1$,
except
  {\rm  GDD}  $1$ of Row $9$;   {\rm  GDD}  $1$ of Row $21$;   {\rm  GDD}  $1$ of Row $18$;   {\rm  GDD}  $6$ of Row $20$.

\end  {Lemma}

\begin {Lemma} \label {3.2.35} Rank
 $n=4$. Non classical   {\rm  GDD}  which is is not  an arithmetic  {\rm  GDD}   if
 there exists only one place where simple chain conditions  are not satisfied,  as well as,
 $\widetilde{q}_{12} \widetilde{q}_{23} \not=1$ and
$q_{22} =-1$ with $q_{11},  \widetilde{q}_{12},  \widetilde{q}_{23} \not= -1$,  or $\widetilde{q}_{32} \widetilde{q}_{34} \not=1$ and
$q_{33} =-1$ with $q_{44},  \widetilde{q}_{32},  \widetilde{q}_{43} \not= -1$,  except
  {\rm  GDD}  $5$ of Row $9$;   {\rm  GDD}  $5$ of Row $14$;   {\rm  GDD}  $6$ of Row $17$;   {\rm  GDD}  $2$ of Row $22$;   {\rm  GDD}  $3$ of Row $22$;   {\rm  GDD}  $1$ of Row $14$.

\end  {Lemma}

\begin {Lemma} \label {3.2.36} Rank
 $n=4$.   {\rm  GDD}  is not  an arithmetic  {\rm  GDD}   if this  {\rm  GDD}  is not classical and every vertex is not $-1$,
 except
 {\rm  GDD}  $1$ of Row $17$.

\end {Lemma}

\begin {Lemma} \label {3.2.37} Rank
 $n=4$. Chain  {\rm  GDD}  is not  an arithmetic  {\rm  GDD}   if $q_{11} = -1$,  $q_{22} \widetilde{q}_{12} \not= 1$ or $q_{44} = -1$,  $q_{33} \widetilde{q}_{34} \not= 1$,  as well as only one vertex is $-1$.

\end {Lemma}

\begin {Lemma} \label {3.2.38} Rank
  $n=4$. Non classical chain  {\rm  GDD}  is not  an arithmetic  {\rm  GDD}   if there exists only one place where simple chain condition   is not satisfied with
 $\widetilde{q}_{21}q_{22} \not= 1$,  $q_{22}\not=-1$ or $ q_{33} \widetilde{q}_{34} \not=1$,   $q_{33}\not=-1$,  or $ q_{22} \widetilde{q}_{23} \not=1$,   $q_{22}\not=-1$
 or $ q_{33} \widetilde{q}_{32} \not=1$, \emph{} $q_{33}\not=-1$.
except
 {\rm  GDD}  $1$ of Row $9$;   {\rm  GDD}  $1$ of Row $21$;   {\rm  GDD}  $1$ of Row $18$;   {\rm  GDD}  $6$ of Row $20$;   {\rm  GDD}  $1$ of Row $20$;  {\rm  GDD}  $5$ of Row $20$;   {\rm  GDD}  $3$ of Row $21$;   {\rm  GDD}  $5$ of Row $21$;   {\rm  GDD}  $1$ of Row $22$.

\end {Lemma}

\section {  Main result} \label {s3}

In this section main theorem is given. That is,   all quasi-affine  connected generalized Dynkin diagram with rank 4 are found.

\begin {Lemma} \label {3.1}
All Type 5,6 union  semi-classical and  classical + semi-classical, which is not a bi-classical  {\rm GDD}s, with rank $n=4$ are listed.\\ \\

 {\ }\ \ \ \ \ \ \ \ \ \ \
 $
$

\end {Lemma}


{\bf Proof.}
 (2.1.1) is quasi-affine  by Lemma \ref {3.2.38}. (2.1.2) is quasi-affine  by Lemma \ref {3.1.3}. (2.1.1) is quasi-affine  by Lemma \ref {3.1.3}.
(3.1.1) is quasi-affine  by Lemma \ref {3.1.3}. (3.1.2) is quasi-affine  by Lemma \ref {3.2.38}.
(3.1.3) is quasi-affine  by Lemma \ref {3.2.38}.
(3.1.4) is quasi-affine  by Lemma \ref {3.1.3}.
(4.1.1) is quasi-affine  step by step.
(5.1.1) is quasi-affine  by Lemma \ref {3.1.3}.
(5.1.2) is quasi-affine  by Lemma \ref {3.2.38}.
(5.1.3) is quasi-affine  by Lemma \ref {3.2.38}.
(5.1.4) is quasi-affine  by Lemma \ref {3.1.3}.
(5.2.1) is quasi-affine  by Lemma \ref {3.1.1}.
(5.1.1) is quasi-affine  by Lemma \ref {3.1.3}.
(5.2.2) is quasi-affine  by Lemma \ref {3.2.5'''}.
(5.2.3) is quasi-affine  by Lemma \ref {3.2.5'''}.
(3.3.2) is quasi-affine  by Lemma \ref {3.1.3}.
(6.1.1) is quasi-affine  by Lemma \ref {3.1.3}.
(3.1.2) is quasi-affine  by Lemma \ref {3.2.38}.
(6.1.3) is quasi-affine  by Lemma \ref {3.2.38}.
(6.1.4) is quasi-affine  by Lemma \ref {3.1.3}.
(6.2.1) is quasi-affine  by Lemma \ref {3.1.33}.
(6.2.2) is quasi-affine  by Lemma \ref {3.1.16}.
(6.2.3) is quasi-affine  by Lemma \ref {3.1.16}.
(6.3.1) is quasi-affine  by Lemma \ref {3.1.13}.
(8.3.1) is quasi-affine  step by step.
(8.3.2) is quasi-affine  step by step.
(10.1.1) is quasi-affine   by Lemma \ref {6.2.2''}.
(10.1.2) is quasi-affine   by Lemma \ref {6.2.2''}.
(10.2.1) is quasi-affine   by Lemma \ref {3.1.33}.
(10.2.2) is quasi-affine   by Lemma \ref {3.1.33}.
(10.2.3) is quasi-affine   by Lemma \ref {3.1.9}.
(10.2.4) is quasi-affine   by Lemma \ref {3.1.1} (IX)  and   Lemma \ref {3.1.16}.
(11.1.1) is quasi-affine   by Lemma \ref {6.2.2''}.
(11.1.2) is quasi-affine   by Lemma \ref {6.2.2''}.
(12.1.1) is quasi-affine   by Lemma \ref {3.1.3}.
(12.1.2) is quasi-affine   by Lemma \ref {3.1.3}. (14.1.1) is quasi-affine   by Lemma \ref {3.1.3}.
(14.1.2) is quasi-affine   by Lemma \ref {3.1.3}.

\begin  {Lemma} \label {3.3}
Near-classical  circles, which are quasi-affine,  with rank $n=4$ are listed.\\ \\ \\

{\ }\ \ \ \ \ \ \ \ \ \ \ \ $
$   \\

\end  {Lemma}

{\bf Proof.}
 We write the proof according  to the following method.

{\rm  (1)   } If a quasi-affine {\rm GDD} is over two {\rm GDD}s, then we only write the quasi-affine {\rm GDD} over the latter one.

{\rm  (2)   }
We do not consider  two-type 6 and four-Type 6.
(6.2.3)
is quasi-affine by    {\rm  GDD}  $1$ of Row $13$.

 (6.2.5) is quasi-affine by  Lemma \ref {3.1.1}  {\rm  (IX)}.
  (8.3.1) is quasi-affine by  Lemma \ref {3.1.1}  {\rm  (I)}.
   (8.3.2) is quasi-affine by  Lemma \ref {3.1.1}  {\rm  (X)}.
 (8.3.3) is quasi-affine by   {\rm  GDD}  $1$ of Row $15$.
  (8.3.6) is quasi-affine by  Lemma \ref {3.1.1}  {\rm  (X)}.
 (15.2.2) is quasi-affine by  Lemma \ref {3.1.2}  {\rm  (III)} (h). (15.2.4) is quasi-affine by  Lemma \ref {3.1.1}  {\rm  (X)} (e).
(15.4.2) is quasi-affine by  Lemma \ref {3.1.1}  {\rm  (X)} (e).
 (17.4.1) is quasi-affine by   {\rm  GDD}  $1$ of Row $15$.
 (17.4.2) is quasi-affine by   {\rm  GDD}  $6$ of Row $17$.
(17.4.4) is quasi-affine by   {\rm  GDD}  $6$ of Row $17$.
(17.5.1) is quasi-affine by   {\rm  GDD}  $6$ of Row $17$.
(17.6.1) is quasi-affine by    Lemma \ref {3.1.1}  {\rm  (X)}.
(17.7.1) is quasi-affine by   {\rm  GDD}  $5$ of Row $17$.

(nc)(nc) which is not quasi-affine is listed: (nc)(nc) over GDD $1$ of Row $4$ is not  quasi-affine by Lemma \ref {3.1.1}(I), or Lemma \ref {3.1.1}(V). (nc)(nc) over GDD $2$ of Row $6$ is not  quasi-affine by Lemma \ref {3.1.1}(I), or Lemma \ref {3.1.2}(II).
 (nc)(nc) over GDD $2$ of Row $7$ is not  quasi-affine by Lemma \ref {3.1.1}(I), or Lemma \ref {3.1.2}(II),  or Lemma \ref {3.1.1}(V).
(nc)(nc) over GDD $1$ of Row $8$ is not  quasi-affine Lemma \ref {3.1.1}(I).
(nc)(nc) over GDD $2$ of Row $8$ is not  quasi-affine Lemma \ref {3.1.2}(III).
(nc)(nc) over GDD $1$ of Row $9$ is not  quasi-affine by Lemma \ref {3.1.1}(I).
(nc)(nc) over GDD $1$ of Row $10$ is not  quasi-affine by Lemma \ref {3.1.1}(I).
(nc)(nc) over GDD $2$ of Row $10$ is not  quasi-affine by Lemma \ref {3.1.1}(I).
(nc)(nc) over GDD $1$ of Row $11$ is not  quasi-affine by Lemma \ref {3.1.1}(I),  Lemma \ref {3.1.1}(V).
(nc)(nc) over GDD $2$ of Row $16$ is not  quasi-affine by Lemma \ref {3.1.1}(I), Lemma \ref {3.1.1}(IX),  Lemma \ref {3.1.1}(V),  Lemma \ref {3.1.1}(III).
(nc)(nc) over GDD $5$ of Row $17$ is not  quasi-affine by Lemma \ref {3.1.1}(X).
\hfill$\Box$

\begin {Theorem} \label {Main}
        All connected  quasi-affine {\rm GDD}s with rank $= 5$ are  listed.
{\rm (i)}  Bi-classical {\rm GDD}s.
{\rm (ii)} Simple circle and semi-simple circles.
{\rm (iii)}  Discontinuous {\rm GDD}s.
{\rm (iv)} Classical + semi-classical {\rm GDD}s. {\rm (v)} Type 5, 6 union semi-classical {\rm GDD}s.
{\rm (vi)} Other cases
\end {Theorem}

{\ }\\

\ \ \ \  $
$\\  \\
$(q,r,s) \in A\cap B\cap C$, where $A:= \{ (q,r,s) \mid  $  $ s =q^{2}$ or $ s =q^{-2}$, $q\in R_6$ or
 $ s =-q^{-2}$, $q \in  R_{12} \cup  R_6  $ or $ s =-1,$ $q^{} \in  R_4 \cup R_6 \cup R_{12}  \cup R_{8}  $
  or $ q^{-2} =-1,$ $s^{} \in  R_2 \cup R_3 \cup R_4 \cup R_6  $    $\}$,   $B:= \{ (q,r,s) \mid  $ $ r^{} =-q^{}$,$q\in R_6$,  or $ r^{} ={-1}$, $q \in  R_4\cup R_6$  $\}$, $C:= \{ (q,r,s) \mid  $  $ qrs =1,  q\not= r,  q\not= s,  s\not= r.$ $\}$. \\ \\ \\ \\ \\ \\

  {\ }\ \ \ \ \ \ {\ }\!\!\!\!\!\!\!\!\!\!\!\!\!\!\!\!\!\!\!\!\!\!\!\!\!\!\!\!\!\!\!\!\!\!\!\!
\!\!\!\!\!\!\!\!\!\!\!\!\!\!\!\!\!\!\!\!\!\!\!\!\!\!\!\!\!\!\!\!\!\!\!\!\!\!\!\!\!\!\!\!\!\!\!\! $\begin{picture}(100,       15)\put(98,       -1){ 10.3.11 }

\put(170,     10){\makebox(0,      0)[t]{$\bullet$}}

\put(170,     70){\makebox(0,      0)[t]{$\bullet$}}

\put(230,     10){\makebox(0,      0)[t]{$\bullet$}}
\put(230,     70){\makebox(0,      0)[t]{$\bullet$}}

\put(170,       10){\line(0,      1){60}}

\put(170,       10){\line(1,      1){60}}

\put(230,       10){\line(0,      1){60}}


\put(170,       10){\line(1,       0){60}}
\put(170,       70){\line(1,       0){60}}

\put(150,     10){$-1$}
\put(150,     30){$q^{-1}$}

\put(150,     70){$q$}

\put(180,     30){$q^{3}$}
\put(190,     80){$q^{-2}$}
\put(190,     -10){$q^{3}$}
\put(220,     30){$$}

\put(250,     10){$-1$}
\put(250,     30){$q^{-6}$}

\put(250,     70){$-1$}

\ \ \ \ \ \ \ \ \ \ \ \ \ \ \ \ \ \ \ \ \ \ \
  \ \ \ \ \ \ \ \ \ \ \ \ \ \  \ \ \ \ \ \ \ \ \ \ \
 \ \ \ \ \ \ \ \ \ \ \ \ \ \ \ \ \ \ \ \ \ \ \ \ \ \ {  $q\in R_8.$.}
\put(80,         -1)  {    } \end{picture}$\\ \\ \\ \\ \\

{\ }\!\!\!\!\!\!\!\!\!\!\!\!\!\!\!\!\!\!\!\!\!{\ }\!\!\!\!\!\!\!\!\!\!\!\!\!\!\!\!\!\!\!\!\!{\ }\!\!\!\!\!\!\!\!\!\!\!\!\!\!\!\!\!\!\!\!\!{\ }\!\!\!\
  $\begin{picture}(100,       15)\put(88,       -1){ 10.3.12}

\put(170,     10){\makebox(0,      0)[t]{$\bullet$}}

\put(170,     70){\makebox(0,      0)[t]{$\bullet$}}

\put(230,     10){\makebox(0,      0)[t]{$\bullet$}}
\put(230,     70){\makebox(0,      0)[t]{$\bullet$}}

\put(170,       10){\line(0,      1){60}}

\put(170,       10){\line(1,      1){60}}

\put(230,       10){\line(0,      1){60}}


\put(170,       10){\line(1,       0){60}}
\put(170,       70){\line(1,       0){60}}

\put(150,     10){$-1$}
\put(150,     30){$s^{}$}

\put(150,     70){$-1$}

\put(180,     30){$q$}
\put(190,     80){$r^{}$}
\put(190,     -10){$q$}
\put(220,     30){$$}

\put(250,     10){$-1$}
\put(250,     30){$q^{-2}$}

\put(250,     70){$-1$}
 \put(280,     -1){ }
\put(80,         -1)  {    } \end{picture}$\\
$(q,r,s) \in A\cap B\cap C$, where $A:= \{ (q,r,s) \mid  $  $ s =-q^{}$, $ r =q^{2}$, $q\in R_6$  or
 $ s =-q^{}$, $ r = -q^{2}$, $q\in R_6$ or
  $ s =-q^{}$, $ r =-1$, $q\in R_6$   $\}$,   $B:= \{ (q,r,s) \mid  $
   $ s =-1^{}$, $ r =q^{2}$, $q\in R_6 \cup  R_4,$
  or  $ s =-1^{}$, $ r =-q^{2}$, $q\in R_6 \cup  R_4,$
  or  $ s =-1^{}$,  $r\in R_6 \cup  R_4 \cup   R_2 \cup   R_3,$ $q\in R_4,$
  or  $ s =-1^{}$, $ r =-1^{2}$, $q\in R_6 \cup  R_4,$
     $\}$.
\\ \\ \\ \\ \\

  {\ }\ \ \ \ \ \ \ \  \ {\ }\!\!\!\!\!\!\!\!\!\!\!\!\!\!\!\!\!\!\!\!\!\!\!\!\!\!\!\!\!\!\!\!\!\!\!\!
\!\!\!\!\!\!\!\!\!\!\!\!\!\!\!\!\!\!\!\!\!\!\!\!\!\!\!\!\!\!\!\!\!\!\!\!\!\!\!\!\!\!\!\!\!\!\!\! $
$

\subsection*
{  Two -Type 6.   That is,   quasi-affine  {\rm  GDD}s consist of two Type 6. }\label {ss2.3}

{\ }\\ \\ \\ \\

\ \ \ \ \ \ {\ }\!\!\!\!\!\!\!\!\!\!\!\!\!\!\!\!\!\!\!\!\!\!\!\!\!\!\!\!\!
  $%
$

\subsection*
{ Four-Type 6. That is,   quasi-affine  {\rm  GDD}s consist of four Type 6. }\label {ss2.4}

{\ }\\ \\ \\ \\

 {\ }\!\!\!\!\!\!\!\!\!\!\!\!\!\!\!\!\!\!\!\!\!\!\!\!\!\!\!\!\!  {\ }\!\!\!\!\!\!\!\!\!\!\!\!\!\!\!\!\!\!\!\!\!\!\!\!\!\!\!\!\!  {\ }\!\!\!\!\!\!\!\!\!\!\!\!\!\!\!\!\!\!\!\!\!\!\!\!\!\!\!\!\!  $%
$\\




{\bf Proof.}
 We write the proof according  to the following method.

{\rm  (1)   } If a {\rm GDD} is quasi-affine over two {\rm GDD}s,  we write the quasi-affine {\rm GDD} over  the latter {\rm GDD}.

{\rm  (2)   } If a non circle  {\rm GDD} is  quasi-affine with circle,  we write it in ones over {\rm GDD} with circle.

{\rm  (3)   }
We do not consider   near-classical circles and near-classical  {\rm  GDD}s.

{\rm  (4)   }
We do not consider  two- type 6 and four-Type 6.

\subsection* {Quasi-affine over
 {\rm  GDD}  $1$ of Row $4$}

   {\rm (ii) } All quasi-arithmetic {\rm  GDD}s by  adding a vertex on  Vertex $2$  are listed.
 According  to Lemma \ref {3.1.2} (I) we have to consider following cases.\\ \\

  {\ }\ \ \ \ \ \ \ \ \ \ \ \ $\begin{picture}(100,       15) \put(-68,        -1){ (a)}

\put(60,       1){\makebox(0,       0)[t]{$\bullet$}}
\put(58,       -12){$q^2$}

\put(40,       -12){$q^{-2}$}
\put(28,       -1){\line(1,       0){33}}
\put(27,       1){\makebox(0,      0)[t]{$\bullet$}}

\put(22,      -12){$q^2$}
\put(0,       -12){$q^{-2}$}

\put(-14,       1){\makebox(0,      0)[t]{$\bullet$}}

\put(-14,      -1){\line(1,       0){50}}

\put(-18,      -12){$q$}

\put(27,     38){\makebox(0,      0)[t]{$\bullet$}}

\put(27,       0){\line(0,      1){35}}

\put(30,       15){$q^{-2}$}

\put(30,       30){$-1$}

\ \ \ \ \ \ \ \ \ \ \ \ \ \ \ \ \ \ \ \ \ \ \ \ \ \ \ \ \
  \ \ \ \ \ \ \ \ \ \ \ \ \ \ \ \ \ \ \ {$, q \in F^{*}\setminus \{1,  -1\}$,  by  {\rm  GDD}  $1$ of Row $2$.}
\put(80,         -1)  {    } \end{picture}$\\ \\
 It is quasi-affine  by Lemma \ref {3.1.3}.\\ \\

  {\ }\ \ \ \ \ \ \ \ \ \ \ \ $\begin{picture}(100,       15) \put(-68,        -1){(b) }

\put(60,       1){\makebox(0,       0)[t]{$\bullet$}}
\put(58,       -12){$q$}

\put(40,       -12){$q^{-1}$}
\put(28,       -1){\line(1,       0){33}}
\put(27,       1){\makebox(0,      0)[t]{$\bullet$}}

\put(22,      -12){$q$}
\put(0,       -12){$q^{-2}$}

\put(-14,       1){\makebox(0,      0)[t]{$\bullet$}}

\put(-14,      -1){\line(1,       0){50}}

\put(-18,      -12){$q^2$}

\put(27,     38){\makebox(0,      0)[t]{$\bullet$}}

\put(27,       0){\line(0,      1){35}}

\put(30,       15){$q^{-1}$}

\put(30,       30){$-1$}

\ \ \ \ \ \ \ \ \ \ \ \ \ \ \ \ \ \ \ \ \ \ \ \ \ \ \ \ \   \ \ \ \ \ \ \ \ \ \ \ \ \ \ \ \ \ \ \ {$, q \in F^{*}\setminus \{1,  -1\}$,  by   {\rm  GDD}  $1$ of Row $3$.}
\put(80,         -1)  {    } \end{picture}$\\ \\
 It is quasi-affineby Lemma \ref {3.1.3} when   $q^2,   q^3 \not=1$. It is  {\rm GDD} $10$ of Row $20$ when   $q \in  R_{3}$.
\subsection* {Quasi-affine over
 {\rm  GDD}  $2$ of Row $5$}
  {\rm (ii) } All quasi-arithmetic {\rm  GDD}s by  adding a vertex on  Vertex $2$  are listed.
 According  to  Type   2 and  {\rm  GDD}  $4$ of Row $7$ we have to consider following cases.

 (a) i.e.   (5.2.1)
  by   {\rm  GDD}  $2$ of Row $5$ or  Type   2.
The sub-{\rm GDD} by deleting  Vertex  4  is an arithmetic GDD  when $q\in R_3 \cup R_{4} \cup R_{6}$ by Lemma \ref {3.1.1} {\rm (X)}.  It is quasi-affine when $q\in R_3 \cup R_{4} \cup R_{6}$ by Lemma \ref {3.1.3}.

 (b)   i.e.  (5.2.2)    by  Type   2.
The sub-{\rm GDD} by deleting Vertex  4  is an arithmetic GDD when $q^2\in R_3\cup R_4\cup R_6\cup R_8 $    by Lemma \ref {3.1.1} (II).  It is quasi-affine $q\in R_3 \cup R_{4} \cup R_{6} \cup R_{8}$ by Lemma \ref {3.1.3} .
\subsection* {Quasi-affine over
 {\rm  GDD}  $3$ of Row $5$}
 {\rm (ii) } All quasi-arithmetic {\rm  GDD}s by  adding  a vertex on  Vertex $3$  are listed.
 According  to  Type   3 and  {\rm  GDD}  $4$ of Row $7$ we have to consider following cases.

 (a) i.e.  (5.3.1) by   Type   3,    $q^2\not=1. $
Sub-{\rm GDD} by deleting Vertex  4  is an arithmetic  {\rm  GDD}   when $q \in  R_4 \cup R_6 \cup R_8  $ by Lemma \ref {3.1.1} {\rm (IX)}.
 It is quasi-affine  by Lemma \ref {3.1.3}.
\subsection* {Quasi-affine over
 {\rm  GDD}  $1$ of Row $6$}
   {\rm (ii) } All quasi-arithmetic {\rm  GDD}s by  adding  a vertex on  Vertex $2$ are listed.
 According  to   Type   1 we have to consider following cases.

 (a) i.e.  (6.1.1)
 by  {\rm  GDD}  $1$ of Row $3$ in Table A2, $q^2\not=1$.
The sub-{\rm GDD} by deleting  4  is an arithmetic  {\rm  GDD}.    It  is  an arithmetic  {\rm  GDD}   when $q\in R_3$ by Lemma \ref {3.1.3}.
   It  is quasi-affine  when $q\notin R_3$ by Lemma \ref {3.1.3}.

 (b) i.e.  (6.1.2)  by   {\rm  GDD}  $1$ of Row $6$ in Table A2.
The sub-{\rm GDD} by deleting Vertex  4  is an arithmetic GDD.
   It  is quasi-affine  by Lemma \ref {3.1.3} with $q\notin R_3$.  It  is an arithmetic GDD  by Lemma \ref {3.1.3}  with $q\in R_3$.

\subsection* {Quasi-affine over
 {\rm  GDD}  $2$ of Row $6$ } \label {sub3.10}
   {\rm (ii) } All quasi-arithmetic {\rm  GDD}s by  adding  a vertex on  Vertex $2$  are listed.
 According  to Lemma \ref {3.1.2} (IV) we have to consider following cases.

 (a)   i.e.  (6.2.1) by   {\rm  GDD}  $2$ of Row $4$, $q^2 \not=1$.
  The sub-{\rm GDD} by deleting Vertex 4  is an arithmetic GDD  by Lemma \ref {3.1.1} {\rm (IV)}.  It  is an arithmetic   when  $q\in R_3$ by  {\rm  GDD} $9$ of Row $20$.  It  is quasi-affine   when  $q\notin R_3 \cup R_2$ by Lemma \ref {3.1.16}.

  (b)  i.e.  (6.2.2)  by    {\rm  GDD}  $2$ of Row $5$, $q \in F^{*}\setminus \{1,  -1\}$.
 The sub-{\rm GDD} by deleting  Vertex 4  is  an arithmetic  {\rm  GDD}   when  $q\in R_7$ by Lemma \ref {3.1.1}(I).
  The sub-{\rm GDD} by deleting Vertex  4  is  an arithmetic  {\rm  GDD}   when  $q\in R_6$ by  Type   2.
  It  is quasi-affine   when  $q\in R_6\cup R_7$ by Lemma \ref {3.1.3}.\\ \\

  {\ }\ \ \ \ \ \ \ \ \ \ \ \ $\begin{picture}(100,       15) \put(-68,        -1){(c) }

\put(60,       1){\makebox(0,       0)[t]{$\bullet$}}
\put(58,       -12){$-1$}

\put(40,       -12){$q$}
\put(28,       -1){\line(1,       0){33}}
\put(27,       1){\makebox(0,      0)[t]{$\bullet$}}

\put(22,      -12){$-1$}
\put(0,       -12){$q^{-2}$}

\put(-14,       1){\makebox(0,      0)[t]{$\bullet$}}

\put(-14,      -1){\line(1,       0){50}}

\put(-18,      -12){$q^2$}

\put(27,     38){\makebox(0,      0)[t]{$\bullet$}}

\put(27,       0){\line(0,      1){35}}

\put(30,       30){$q^{2}$}

\put(30,       15){$q^{-2}$}

\ \ \ \ \ \ \ \ \ \ \ \ \ \ \ \ \   \ \ \ \ \ \   \ \ \ \ \ \ \ \ \ \ \ \ \ \ \ \ \ \ \ {$, q \in F^{*}\setminus \{1,  -1\}$,  by   {\rm  GDD}  $2$ of Row $6$.}
\put(80,         -1)  {    } \end{picture}$\\ \\
 The sub-{\rm GDD} by deleting  Vertex 4 is an arithmetic  {\rm  GDD}   by Lemma \ref {3.1.1} {\rm (IV)}.
   It  is quasi-affine  by Lemma \ref {3.1.16}.

\subsection* {Quasi-affine over
 {\rm  GDD}  $3$ of Row $6$ } \label {sub3.11}

  {\rm (i) } All quasi-arithmetic {\rm  GDD}s by  adding  a vertex on  Vertex $1$  are listed.

  (a)  i.e.  (6.3.1),   $q^2,  q^3 \not=1, $  by    {\rm  GDD}  $1$ of Row $7$. It is quasi-affine  by Lemma \ref {3.1.3}.

   (b)  i.e.  (6.3.2) $ ,q \in R_3$,  by    {\rm  GDD}  $1$ of Row $15$. It is
quasi-affine  by Lemma \ref {3.1.3}.

  (c)  i.e.  (6.3.3)$ ,q \in R_3$,  by   {\rm  GDD}  $1$ of Row $16$. It is quasi-affine  by Lemma \ref {3.1.3}.

  (d)  i.e.  (6.3.4)
  {$, q \in R_3$,  by   {\rm  GDD}  $2$ of Row $17$. It is quasi-affine  by Lemma \ref {3.1.3}.

  (e)  i.e.  (6.3.5)
  $ ,q \in R_3$,  by  {\rm  GDD}  $5$ of Row $17$. It is quasi-affine  by Lemma \ref {3.1.3}.

  (f)  i.e.  (6.3.6)
  $ ,q \in R_3$,  by  {\rm  GDD}  $7$ of Row $17$. It is quasi-affine  by Lemma \ref {3.1.3}.

  {\rm (ii) } All quasi-arithmetic {\rm  GDD}s by  adding  a vertex on  Vertex $2$ of {\rm  GDD}  $3$ of Row $6$ are listed.
 According  to Lemma \ref {3.1.2} (I) we have to consider following cases.\\

{\ }\ \ \ \ \ \ \ \ \ \ \ \ $\begin{picture}(100,       15) \put(-68,        -1){ $(a)$}

\put(60,       1){\makebox(0,       0)[t]{$\bullet$}}

\put(28,       -1){\line(1,       0){33}}
\put(27,       1){\makebox(0,      0)[t]{$\bullet$}}
\put(-14,       1){\makebox(0,      0)[t]{$\bullet$}}

\put(-14,      -1){\line(1,       0){50}}

\put(58,       -12){$q$}

\put(40,       -12){$q^{-1}$}

\put(22,      -12){$-1$}
\put(0,       -12){${q}$}

\put(-18,      -12){${-1}$}

\put(59,       0){\line(-1,      1){17}}

\put(28,      -1){\line(1,       1){17}}

\put(18,     12){$q^{2}$}

\put(36,     22){$-1$}
\put(59,     12){$q^{-1}$}

\put(43,     18){\makebox(0,      0)[t]{$\bullet$}}

 \ \ \ \ \ \ \ \ \ \ \ \ \ \ \ \ \ \ \ \ \ \ \ \  \ \ \ \ \ \ \ \
 \ \ \ \ \ \ \ \ \ \ \ \ \ \ \ \ \ \ \ { ,  by   {\rm  GDD}  $2$ of Row $4$.}
\put(80,         -1)  {    } \end{picture}$\\ \\
The sub-{\rm GDD} by deleting  Vertex  4  is not arithmetic  {\rm  GDD}   when $q^3, q^4 \not=1$ by Lemma \ref {3.1.2} {\rm (IV)}.
 It  is an arithmetic  {\rm  GDD}    when $q^3 =1$ by  {\rm  GDD}  $8$ of Row $20$ or by Lemma \ref {3.1.13}.

(b) i.e.  (6.3.7)   by   {\rm  GDD}  $2$ of Row $6$. The sub-{\rm GDD} by deleting  Vertex  4  is
 an arithmetic  {\rm  GDD}    when $q \in  R_3\cup R_5\cup R_6 $  by Lemma \ref {3.1.1} {\rm (IX)},   It  is quasi-affine  when $q \in  R_3\cup R_5\cup R_6 $  by Lemma \ref {3.1.13}.

 (c)  i.e.  (6.3.8)
  $,q^3\notin R_2$,  by   {\rm  GDD}  $2$ of Row $7$. The sub-{\rm GDD} by deleting  Vertex  4. It
is an arithmetic  {\rm  GDD}    when $q \in R_7$  by Lemma \ref {3.1.1} {\rm (IX)}.   It  is quasi-affine  by Lemma \ref {3.1.13} when $q \in R_7$ .

 (d)  i.e.  (6.3.9) ,$q\notin R_3$,  by   {\rm  GDD}  $4$ of Row $7$. The sub-{\rm GDD} by deleting  Vertex  4   is
 an arithmetic  {\rm  GDD}   when   $q\in R_8$ by  Type   3.  The sub-{\rm GDD} by deleting  Vertex  4 is not an arithmetic  {\rm  GDD}   when   $q\notin R_8$  by Lemma \ref {3.1.1}{\rm (I)}.   It  is quasi-affine  when $q\in R_8$ by Lemma \ref {3.1.3}.

 (e)  i.e.  (6.3.10) ,by   {\rm  GDD}  $1$ of Row $8$.
 The sub-{\rm GDD} by deleting  Vertex  4 is an arithmetic  {\rm  GDD}    when $q \in  R_3 \cup R_5 \cup R_7 \cup R_6 \cup R_4$   by Lemma \ref {3.1.1} {\rm (IX)} and   by Lemma \ref {3.1.1} {\rm (X)}.
  It  is quasi-affine  when  $q \in  R_4 \cup R_5 \cup R_7 \cup R_6 $  by Lemma \ref {3.1.13};
 It  is an arithmetic  {\rm  GDD}    when  $q \in  R_3$  by  {\rm  GDD}  $3$ of Row $18$.

  (f)  i.e.  (6.3.11)    by   {\rm  GDD}  $1$ of Row $9$.
 The sub-{\rm GDD} by deleting  Vertex  4  is an arithmetic  {\rm  GDD}    when $ q^{2} =r^{}$,    $q\notin R_3$ or $ q^{4} =r^{}$,  $q\notin R_3\cup   R_5\cup R_4$ or $ q^{6} =r^{}$,   $q\notin R_5\cup   R_7\cup  R_6\cup R_4$ or $ -q^{2} =r^{}$,   $q\in R_6$,  or $r=-1, $ $q \in  R_4 \cup R_6 \cup R_8 \cup R_{12}\cup R_{3} $,  or $q\in R_4$ by Lemma \ref {3.1.1} {\rm (IV)},   by Lemma \ref {3.1.1} {\rm (IX)} and   by Lemma \ref {3.1.1} {\rm (X)}.  It  is quasi-affine   when $r\not=  q^{2}$ by Lemma \ref {3.1.13}.  It  is an arithmetic  {\rm  GDD}   when $r=  q^{2}$   by Lemma \ref {3.1.13}.

 (g)  i.e.  (6.3.12)
 $, q \notin R_3$,  by   {\rm  GDD}  $2$ of Row $10$.
 The sub-{\rm GDD} by deleting  Vertex  4  is an arithmetic  {\rm  GDD}    when $q \in  R_5 \cup R_9 \cup R_{8} \cup R_{13}$  by Lemma \ref {3.1.1} {\rm (IX)}.
  It  is an arithmetic  {\rm  GDD}   when $q \in  R_5$   by Lemma \ref {3.1.13}.  It  is quasi-affine  when $q\notin R_5$
 by Lemma \ref {3.1.13}.

  (h)  i.e.  (6.3.13)$, q \notin R_3$,  by   {\rm  GDD}  $1$ of Row $10$.
 The sub-{\rm GDD} by deleting  Vertex  4 is  an arithmetic  {\rm  GDD}    when $q\in R_5\cup R_4$ by Lemma \ref {3.1.2}{\rm (I)} and by Lemma \ref {3.1.1} {\rm (IV)}.
  It  is quasi-affine  by Lemma \ref {3.1.13}.

 (i)  i.e.  (6.3.14)$, q \notin R_3$,  by   {\rm  GDD}  $2$ of Row $10$. The sub-{\rm GDD} by deleting  Vertex  4.
It is an  arithmetic  {\rm  GDD}    when  $q \in  R_4 \cup R_6 \cup R_8 $  by Lemma \ref {3.1.1} {\rm (IX)}.   It  is an arithmetic  {\rm  GDD}   when $q\in R_4$     by Lemma \ref {3.1.13}.
  It  is quasi-affine   when $q\notin R_4.$    by Lemma \ref {3.1.13}.

 (j)  i.e.  (6.3.15)$, q \in R_3$,  by   {\rm  GDD}  $1$ of Row $11$.   The sub-{\rm GDD}
by deleting Vertex  4  is an arithmetic  {\rm  GDD}   when   $q\in R_3$.   It  is quasi-affine  by Lemma \ref {3.1.13}.\\

{\ }\ \ \ \ \ \ \ \ \ \ \ \  $\begin{picture}(100,       15) \put(-68,        -1){(k) }

\put(60,       1){\makebox(0,       0)[t]{$\bullet$}}

\put(28,       -1){\line(1,       0){33}}
\put(27,       1){\makebox(0,      0)[t]{$\bullet$}}
\put(-14,       1){\makebox(0,      0)[t]{$\bullet$}}

\put(-14,      -1){\line(1,       0){50}}

\put(58,       -12){$-q^{-1}$}

\put(40,       -12){$-q$}

\put(22,      -12){$-1$}
\put(0,       -12){$-q^{-1}{}$}

\put(-18,      -12){$q^{-1}$}

\put(59,       0){\line(-1,      1){17}}

\put(28,      -1){\line(1,       1){17}}

\put(18,     12){$q^{-2}$}

\put(36,     22){$-1$}
\put(59,     12){$-q^{}$}

\put(43,     18){\makebox(0,      0)[t]{$\bullet$}}

 \ \ \ \ \ \ \ \ \ \ \ \ \ \ \ \ \ \ \ \ \ \ \ \  \ \ \ \ \ \ \ \
  \ \ \ \ \ \ \ \ \ \ \ \ \ \ \ \ \ \ \ {$, q \in R_3$,  by    {\rm  GDD}  $3$ of Row $14$.}
\put(80,         -1)  {    } \end{picture}$\\ \\
 The sub-{\rm GDD} by deleting  Vertex  4  is not arithmetic  {\rm  GDD}    by Lemma \ref {3.1.1}{\rm (I)}.
 \\

 (l)  i.e.  (6.3.16)$, q \in R_3$,  by    {\rm  GDD}  $2$ of Row $16$. The sub-{\rm GDD} by deleting  Vertex  4  is
an arithmetic  {\rm  GDD}    when $q \in  R_3 $   by Lemma \ref {3.1.1} {\rm (X)}.
  It  is quasi-affine  by Lemma \ref {3.1.13}.
\\

{\ }\ \ \ \ \ \ \ \ \ \ \ \  \ \  $\begin{picture}(100,       15) \put(-68,        -1){(m) }

\put(60,       1){\makebox(0,       0)[t]{$\bullet$}}

\put(28,       -1){\line(1,       0){33}}
\put(27,       1){\makebox(0,      0)[t]{$\bullet$}}
\put(-14,       1){\makebox(0,      0)[t]{$\bullet$}}

\put(-14,      -1){\line(1,       0){50}}

\put(58,       -12){$-q$}

\put(35,       -12){$-q^{-1}$}

\put(22,      -12){$-1$}
\put(0,       -12){${-1}$}

\put(-18,      -12){${q}$}

\put(59,       0){\line(-1,      1){17}}

\put(28,      -1){\line(1,       1){17}}

\put(18,     12){$q^{2}$}

\put(36,     22){$-1$}
\put(59,     12){$-q^{-1}$}

\put(43,     18){\makebox(0,      0)[t]{$\bullet$}}

 \ \ \ \ \ \ \ \ \ \ \ \ \ \ \ \ \ \ \ \ \ \ \ \  \ \ \ \ \ \ \ \
  \ \ \ \ \ \ \ \ \ \ \ \ \ \ \ \ \ \ \ {$, q \in R_3$,  by   {\rm  GDD}  $4$ of Row $16$.}
\put(80,         -1)  {    } \end{picture}$\\ \\
 The sub-{\rm GDD} by deleting  Vertex  4  is not an arithmetic  {\rm  GDD}   by Lemma \ref {3.1.1}{\rm (I)}.

\subsection* {Quasi-affine over
 {\rm  GDD}  $1$ of Row $7$ }
  {\rm (i) }  All quasi-arithmetic {\rm  GDD}s by  adding  a vertex on  Vertex $1$  are listed.
 According  to Lemma \ref {3.1.2} (I) we have to consider following cases.

 (a)  i.e.  (7.1.1),
$q^2\notin R_3$,   $q^2\notin R_2$,  by   {\rm  GDD}  $2$ of Row $6$. It
is quasi-affine  by Lemma \ref {3.1.8} and by
Lemma \ref {3.1.21}.

 {\rm (ii) }  All quasi-arithmetic {\rm  GDD}s by  adding  a vertex on  Vertex $2$  are listed.

 (a)  i.e.  (7.1.2) by  {\rm  GDD}  $1$ of Row $7$ in Table A2. The sub-{\rm GDD}
 by deleting  Vertex  4  is an arithmetic  {\rm  GDD}
 by  Type   7.   It  is  quasi-affine  by Lemma \ref {3.1.3}.

\subsection* {Quasi-affine over
 {\rm  GDD}  $2$ of Row $7$ }
  {\rm (i) }  All quasi-arithmetic {\rm  GDD}s by  adding  a vertex on  Vertex $1$  are listed.
 According  to Lemma \ref {3.1.2} (IV) we have to consider following cases.

 (a)  i.e.  (7.2.1)$, q \in F^{*}\setminus \{1,  -1\}$. $ q^2 \notin R_3\cup R_2$,  by  {\rm  GDD}  $2$ of Row $7$.
It is  quasi-affine  by Lemma \ref {3.1.8} and  by Lemma \ref {3.1.21}.

 (b)  i.e.  (7.2.2)$, q \in F^{*}\setminus \{1,  -1\}$. $ q \notin R_3$,  by  {\rm  GDD}  $2$ of Row $7$.
It is  quasi-affine  by Lemma \ref {3.1.21}.

  {\rm (ii) }  All quasi-arithmetic {\rm  GDD}s by  adding a vertex on  Vertex $2$ of {\rm  GDD}  $2$ of Row $7$ are listed.
According  to Lemma \ref {3.1.2} (IV) we have to consider following cases.

 (a)  i.e.  (7.2.1)
$, q \in F^{*}\setminus \{1,  -1\}$, $q^3 \not=1.$ by   {\rm  GDD}  $2$ of Row $4$.

 The sub-{\rm GDD} by deleting  Vertex  4  is an arithmetic  {\rm  GDD}
 by Lemma \ref {3.1.1} {\rm (IV)}.   It  is
quasi-affine  when $q\notin R_4$ by Lemma \ref {3.1.16}.
 It  is quasi-affine  when $q\in R_4$ by
Lemma \ref {3.1.35}.

 (b)  i.e.  (7.2.3),  $q^2 \notin  R_2 \cup R_3 $,  $q^2,  q^4 \not=1, $   by   {\rm  GDD}  $2$ of Row $5$.
 The sub-{\rm GDD} by deleting  Vertex  4 is an arithmetic  {\rm  GDD}
when $q\in R_8$ by  Type   2.  The sub-{\rm GDD} by deleting  Vertex  4  is
an arithmetic  {\rm  GDD}   when $q\in R_{18}$
by Lemma \ref {3.1.1}{\rm (I)}.
  It  is quasi-affine  when   $q\in R_8\cup  R_{18}$
 by Lemma \ref {3.1.3}.

(c)  i.e.  (7.2.4)
$, q \in F^{*}\setminus \{1,  -1\}$,  by    {\rm  GDD}  $2$ of Row $6$. The sub-{\rm GDD} by
deleting  Vertex  4 is an arithmetic  {\rm  GDD}
 by Lemma \ref {3.1.1} {\rm (IV)}.
  It  is quasi-affine   by Lemma \ref {3.1.16}.

(d)  i.e.  (7.2.5), $q^2, q^3 \not=1,$ by   {\rm  GDD}  $2$ of Row $7$. The sub-{\rm GDD} by deleting  Vertex  4
 is an arithmetic  {\rm  GDD}
 by Lemma \ref {3.1.1} {\rm (IV)}.
  It  is quasi-affine   by Lemma \ref {3.1.35}.

  {\rm (iii) }  All quasi-arithmetic {\rm  GDD}s by  adding a vertex on  Vertex $3$ are listed.
 According  to Lemma \ref {3.1.2} (III) we have to consider following cases.

(a) i.e.  (7.2.6), $q^9, q^6 \not=1,$  by    {\rm  GDD}  $1$ of Row $7$. It is quasi-affine
by Lemma   \ref {2.63}.

  {\rm (iv) } All quasi-affine circles are listed.

  (nc) (a)   is not quasi-affine since\\

    $\begin{picture}(100,       15) \put(-68,        -1){}

\put(60,       1){\makebox(0,       0)[t]{$\bullet$}}

\put(28,       -1){\line(1,       0){33}}
\put(27,       1){\makebox(0,      0)[t]{$\bullet$}}

\put(-14,       1){\makebox(0,      0)[t]{$\bullet$}}

\put(-14,      -1){\line(1,       0){50}}

\put(-18,      10){$-1$}
\put(0,       5){$q^{-1} $}
\put(22,      10){$q$}
\put(40,       5){$q^{-9} $}

\put(58,       10){$q^{3} $}

\put(80,         -1)  { with  $q =q^{9} $ is not arithmetic  {\rm  GDD}
by Lemma \ref {3.1.1} (I).   } \end{picture}$

(nc)  (a)  is not quasi-affine since\\

  $\begin{picture}(100,       15) \put(-68,        -1){ }

\put(60,       1){\makebox(0,       0)[t]{$\bullet$}}

\put(28,       -1){\line(1,       0){33}}
\put(27,       1){\makebox(0,      0)[t]{$\bullet$}}

\put(-14,       1){\makebox(0,      0)[t]{$\bullet$}}

\put(-14,      -1){\line(1,       0){50}}

\put(-18,      10){$-1$}
\put(0,       5){$q^{-2} $}
\put(22,      10){$q$}
\put(40,       5){$q^{-18} $}

\put(62,       10){$q^{6} $}

\put(80,         -1)  {   with $q =q^{18} $ is not an arithmetic  {\rm  GDD}
by Lemma \ref {3.1.1} (I). } \end{picture}$

 (a) (a)  i.e.  (7.2.7)   is quasi-affine  since  the sub-{\rm GDD} by deleting   Vertex 2

  with $q^{3} =q^{9} $ is an  arithmetic  {\rm  GDD}
 by Lemma \ref {3.1.1} {\rm (IV)}.

 (a) (nc)  i.e.  (7.2.8)
  is quasi-affine  since  the sub-{\rm GDD} by deleting  Vertex 2
  with $-1 =q^{3} $ is an  arithmetic  {\rm  GDD}.

 (a) (nc)  i.e.  (7.2.9)   is quasi-affine  since the sub-{\rm GDD} by deleting   Vertex 2
 is an arithmetic  {\rm  GDD}.

 (a) (nc)  i.e.  (7.2.10) is quasi-affine  when $ q^{12} \not=1$   since  the sub-{\rm GDD} by deleting    Vertex 2
 is an arithmetic  {\rm  GDD}   by  Type   1.

(nc) (a)  is not quasi-affine since\\

 $\begin{picture}(100,       15) \put(-68,        -1){ }

\put(60,       1){\makebox(0,       0)[t]{$\bullet$}}

\put(28,       -1){\line(1,       0){33}}
\put(27,       1){\makebox(0,      0)[t]{$\bullet$}}

\put(-14,       1){\makebox(0,      0)[t]{$\bullet$}}

\put(-14,      -1){\line(1,       0){50}}

\put(-18,      10){$-1$}
\put(0,       5){$q^{-1}$}
\put(22,      10){$-1$}
\put(40,       5){$q^{-9}$}

\put(58,       10){$q^{3}$}

\put(80,         -1)  {  with $-1 =q^{9} $ is not arithmetic  {\rm  GDD}
by Lemma \ref {3.1.1} (I).  } \end{picture}$

\subsection* {Quasi-affine over
 {\rm  GDD}  $3$ of Row $7$ }
   {\rm (i) }  All quasi-arithmetic {\rm  GDD}s by  adding a vertex on  Vertex $1$ of {\rm  GDD}  $3$ of Row $7$ are listed.
 According  to Lemma \ref {3.1.2} (III) we have to consider following cases.\\

  {\ }\ \ \ \ \ \ \ \ \ \ \ \ $\begin{picture}(100,       15) \put(-68,        -1){ (a)}

\put(60,       1){\makebox(0,       0)[t]{$\bullet$}}
\put(58,       -12){$-1$}

\put(40,       -12){$q^{-1}$}
\put(28,       -1){\line(1,       0){33}}
\put(27,       1){\makebox(0,      0)[t]{$\bullet$}}

\put(22,      -12){$q$}
\put(0,       -12){$q^{-2}$}

\put(-14,       1){\makebox(0,      0)[t]{$\bullet$}}

\put(-14,      -1){\line(1,       0){50}}

\put(-18,      -12){$q^2$}

\put(59,       0){\line(-1,      1){17}}

\put(28,      -1){\line(1,       1){17}}

\put(43,     18){\makebox(0,      0)[t]{$\bullet$}}

\put(18,     12){$q^{-2}$}

\put(36,     22){$-1$}
\put(59,     12){$q^3$}

 \ \ \ \ \ \ \ \ \ \ \ \  \ \ \ \ \ \ \ \
  \ \ \ \ \ \ \ \ \ \ \ \ \ \ \ \ \ \ \ {$, q \in F^{*}\setminus \{1,  -1\}$ by    {\rm  GDD}  $1$ of Row $6$.}
\put(80,         -1)  {    } \end{picture}$\\ \\
 The sub-{\rm GDD} by deleting  Vertex  4  is not an arithmetic  {\rm  GDD}    by Lemma \ref {3.1.1} {\rm (V)}.

 (b)     i.e.  (7.3.2)$, q \in F^{*}\setminus \{1,  -1\}$,  by   {\rm  GDD}  $3$ of Row $6$. The sub-{\rm GDD} by deleting
Vertex  2 is an arithmetic  {\rm  GDD}   by Lemma \ref {3.1.2} {\rm (III)} with $q\in R_6$.   The sub-{\rm GDD} by deleting  Vertex  2 is not an arithmetic  {\rm  GDD}   by Lemma \ref {3.1.2} {\rm (III)} with $q\notin R_6$.  The sub-{\rm GDD} by deleting  Vertex  1 is  an arithmetic  {\rm  GDD}   by Lemma \ref {3.1.1} {\rm (X)} with $q\in R_6$.   It  is quasi-affine   with $q\in R_6$. \\

  {\ }\ \ \ \ \ \ \ \ \ \ \ \ $\begin{picture}(100,       15) \put(-68,        -1){ (c)}

\put(60,       1){\makebox(0,       0)[t]{$\bullet$}}
\put(58,       -12){$-1$}

\put(40,       -12){$q^{-1}$}
\put(28,       -1){\line(1,       0){33}}
\put(27,       1){\makebox(0,      0)[t]{$\bullet$}}

\put(22,      -12){$q$}
\put(0,       -12){$q^{-3}$}

\put(-14,       1){\makebox(0,      0)[t]{$\bullet$}}

\put(-14,      -1){\line(1,       0){50}}

\put(-18,      -12){$q^3$}

\put(59,       0){\line(-1,      1){17}}

\put(28,      -1){\line(1,       1){17}}

\put(43,     18){\makebox(0,      0)[t]{$\bullet$}}

\put(18,     12){$q^{-2}$}

\put(36,     22){$-1$}
\put(59,     12){$q^3$}

  \ \ \ \ \ \ \ \ \ \ \ \ \ \ \ \ \ \  \ \ \ \ \ \ \ \
  \ \ \ \ \ \ \ \ \ \ \ \ \ \ \ \ \ \ \ {$, q \in F^{*}\setminus \{1,  -1\}$ $, q \notin R_3$,  by   {\rm  GDD}  $1$ of Row $7$.}
\put(80,         -1)  {    } \end{picture}$\\
\\
 The sub-{\rm GDD} by deleting  Vertex  4   is not an arithmetic  {\rm  GDD}    by Lemma \ref {3.1.1} {\rm (V)}.

(d)     i.e.  (7.3.3)$, q \in F^{*}\setminus \{1,  -1\}$ $, q \notin R_3$,  by    {\rm  GDD}  $3$ of Row $7$. The sub-{\rm GDD}
by deleting  Vertex  2   is  an  arithmetic  {\rm  GDD}    by Lemma \ref {3.1.1} {\rm (V)} when  $q\in R_6$.  The sub-{\rm GDD} by deleting  Vertex  1  is  an  arithmetic  {\rm  GDD}    by  Type   7.   It  is quasi-affine  when  $q\in R_6$ by Lemma \ref {3.1.3} .

(e)     i.e.  (7.3.4)$, q \in F^{*}\setminus \{1,  -1\}$,  by  {\rm  GDD}  $3$ of Row $8$. The sub-{\rm GDD} by
deleting  Vertex  4 is  arithmetic  {\rm  GDD}   by Lemma \ref {3.1.2} (III).   It  is quasi-affine  by Lemma \ref {3.1.3}.\\

  {\ }\ \ \ \ \ \ \ \ \ \ \ \ \ $\begin{picture}(100,       15) \put(-68,        -1){ (f)}

\put(60,       1){\makebox(0,       0)[t]{$\bullet$}}
\put(58,       -12){$-1$}

\put(40,       -12){$-q$}
\put(28,       -1){\line(1,       0){33}}
\put(27,       1){\makebox(0,      0)[t]{$\bullet$}}

\put(18,      -12){$-q^{-1}$}
\put(-4,       -12){$-q$}

\put(-14,       1){\makebox(0,      0)[t]{$\bullet$}}

\put(-14,      -1){\line(1,       0){50}}

\put(-18,      -12){$q$}

\put(59,       0){\line(-1,      1){17}}

\put(28,      -1){\line(1,       1){17}}

\put(43,     18){\makebox(0,      0)[t]{$\bullet$}}

\put(18,     12){$q^{2}$}

\put(36,     22){$-1$}
\put(59,     12){$-q^{-3}$}

 \ \ \ \ \ \ \ \ \ \ \ \ \ \ \ \ \ \ \ \ \ \ \ \  \ \ \ \ \ \ \ \
  \ \ \ \ \ \ \ \ \ \ \ \ \ \ \ \ \ \ \ {$, q \in R_3$,  by   {\rm  GDD}  $1$ of Row $14$.}
\put(80,         -1)  {    } \end{picture}$\\ \\
 The sub-{\rm GDD} by deleting  Vertex  4 is not an arithmetic  {\rm  GDD}   by Lemma \ref {3.1.1}{\rm (I)}.

(g)     i.e.  (7.3.6)
  by   {\rm  GDD}  $1$ of Row $4$. The sub-{\rm GDD} by deleting  Vertex  4  is arithmetic
{\rm  GDD} when $q\in R_6$ by Lemma \ref {3.1.2}{\rm (II)}.   It  is quasi-affine  when $q\in R_6$ by Lemma \ref {3.1.3}.\\

{\ }\ \ \ \ \ \ \ \ \ \ \ \ $\begin{picture}(100,       15) \put(-68,        -1){ $(h)$}

\put(60,       1){\makebox(0,       0)[t]{$\bullet$}}

\put(28,       -1){\line(1,       0){33}}
\put(27,       1){\makebox(0,      0)[t]{$\bullet$}}
\put(-14,       1){\makebox(0,      0)[t]{$\bullet$}}

\put(-14,      -1){\line(1,       0){50}}

\put(58,       -12){$-1$}

\put(40,       -12){$q^{-2}$}

\put(22,      -12){$q^2$}
\put(0,       -12){${q^{-2}}$}

\put(-18,      -12){${q}$}

\put(59,       0){\line(-1,      1){17}}

\put(28,      -1){\line(1,       1){17}}

\put(43,     18){\makebox(0,      0)[t]{$\bullet$}}

\put(18,     12){$q^{-4}$}

\put(36,     22){$-1$}
\put(59,     12){$q^6$}

 \ \ \ \ \ \ \ \ \ \ \ \ \ \ \ \ \ \ \ \ \ \ \ \  \ \ \ \ \ \ \ \
  \ \ \ \ \ \ \ \ \ \ \ \ \ \ \ \ \ \ \ { ,  by   {\rm  GDD}  $1$ of Row $5$.}
\put(80,         -1)  {    } \end{picture}$\\ \\
 The sub-{\rm GDD} by deleting  Vertex  4  is not an  arithmetic  {\rm  GDD}    by Lemma \ref {3.1.1}{\rm (I)}.

(i)     i.e.  (7.3.7)  by   {\rm  GDD}  $2$ of Row $7$.  The sub-{\rm GDD} by deleting  Vertex  4
is  an  arithmetic  {\rm  GDD}    $q \in   R_{10}  $ by Lemma \ref {3.1.1} (X).  It  is quasi-affine  by Lemma \ref {3.1.3}.
\\ \\ \\ \\ \\

{\ }\!\!\!\!\!\!\!\!\!\!\!\!\!\!\!\!\!\!\!\!\!{\ }\!\!\!\!\!\!\!\!\!\!\!\!\!\!\!\!\!\!\!\!\!{\ }\!\!\!\!\!\!\!\!\!\!\!\!\!\!\!\!\!\!\!\!\!{\ }\!\!\!\
\!\!\! \!\!\! \!\!\! \!\!\!    \!\!       $\begin{picture}(100,       15)\put(98,       -1){ (j)}

\put(170,     10){\makebox(0,      0)[t]{$\bullet$}}

\put(170,     70){\makebox(0,      0)[t]{$\bullet$}}

\put(230,     10){\makebox(0,      0)[t]{$\bullet$}}
\put(230,     70){\makebox(0,      0)[t]{$\bullet$}}

\put(170,       10){\line(0,      1){60}}

\put(170,       10){\line(1,      1){60}}

\put(230,       10){\line(0,      1){60}}


\put(170,       10){\line(1,       0){60}}
\put(170,       70){\line(1,       0){60}}

\put(150,     10){$-1$}
\put(150,     30){$q^{3}$}

\put(150,     70){$-1$}

\put(180,     30){$q^{-1}$}
\put(190,     80){$q^{-2}$}
\put(190,     -10){$q^{3}$}
\put(220,     30){$$}

\put(250,     10){$-1$}
\put(250,     30){$q^{-2}$}

\put(250,     70){$q$}

\ \ \ \ \ \ \ \ \ \ \ \
 \ \ \ \ \ \ \ \ \ \ \ \ \ \ \ \ \ \ \ \ \ \ \ \  \ \ \ \ \ \ \ \ \ \ \  \ \ \ \ \ \ \ \
 \ \ \ \ \ \ \ \ \ \ \ \ \ \ \ \ \ \ \ \ \ \ \ \ { ,  by     {\rm  GDD}  $3$ of Row $7$.}
\put(80,         -1)  {    } \end{picture}$\\ \\
 The sub-{\rm GDD} by deleting  Vertex  2    is  an  arithmetic  {\rm  GDD}    when $q \in  R_9 $ by Lemma \ref {3.1.1} (X).
 The sub-{\rm GDD} by deleting  Vertex  1    is not  an arithmetic  {\rm  GDD}     by Lemma \ref {3.1.1} (V).

(k)     i.e.  (7.3.8)   by   {\rm  GDD}  $4$ of Row $7$. The sub-{\rm GDD} by deleting  Vertex  4
  is an arithmetic  {\rm  GDD}    when $q\in R_{11}$ by  Type   3.
 The sub-{\rm GDD} by deleting  Vertex  4    is not an  arithmetic  {\rm  GDD}    when $q\notin R_{11}$  by Lemma \ref {3.1.1}{\rm (I)}.
 It  is   quasi-affine  by Lemma \ref {3.1.3}.

   (l)     i.e.  (7.3.5)$, q \in R_3$,  by   {\rm  GDD}  $5$ of Row $17$. The sub-{\rm GDD} by deleting  Vertex  4
  is an arithmetic  {\rm  GDD}   by Lemma \ref {3.1.1} {\rm (V)}.  It  is quasi-affine  by Lemma \ref {3.1.3}.

  {\rm (ii) }   All quasi-arithmetic {\rm  GDD}s by  adding a vertex on  Vertex $2$ of  are listed.
According  to Lemma \ref {3.1.2} (I) we have to consider following cases.

(a)     i.e.  (7.3.9)
  by    {\rm  GDD}  $2$ of Row $6$. The sub-{\rm GDD} by deleting  Vertex  4    is an arithmetic
{\rm  GDD}    when $q\in R_7$ by Lemma \ref {3.1.1} {\rm (X)}.  It  is quasi-affine     when $q\in R_7$ by Lemma \ref {3.1.3}.

(b)     i.e.  (7.3.10)   by   {\rm  GDD}  $2$ of Row $7$.  The sub-{\rm GDD} by deleting  Vertex  4
is  an  arithmetic  {\rm  GDD}    $q \in   R_{10}  $ by Lemma \ref {3.1.1} (X).  It  is quasi-affine  by Lemma \ref {3.1.3}.

(c)     i.e.  (7.3.11)   by   {\rm  GDD}  $4$ of Row $7$. The sub-{\rm GDD} by deleting  Vertex  4
   is an arithmetic  {\rm  GDD}    when $q\in R_{11}$ by  Type   3.
 The sub-{\rm GDD} by deleting  Vertex  4    is not an  arithmetic  {\rm  GDD}    when $q\notin R_{11}$  by Lemma \ref {3.1.1}{\rm (I)}.
 It  is   quasi-affine  by Lemma \ref {3.1.3}.

(d)     i.e.  (7.3.12)   by    {\rm  GDD}  $1$ of Row $8$.
 The sub-{\rm GDD} by deleting  Vertex  4   is an arithmetic  {\rm  GDD}   $q \in  R_4 \cup R_7 \cup R_{10}  \cup R_6. $
by Lemma \ref {3.1.2} {\rm (IV)} and by Lemma \ref {3.1.1} {\rm (IV)}.
It  is quasi-affine  by Lemma \ref {3.1.3}.

(e)     i.e.  (7.3.13)  by   {\rm  GDD}  $1$ of Row $9$. The sub-{\rm GDD} by deleting  Vertex  4
  is an arithmetic  {\rm  GDD}   when $ r =q^{3}$,   $q\notin R_6$ or $ r =q^{6}$,   $q\notin R_{12}$ or
  $ r =-q^{-3}$,    $q\in R_9$   or  $  q^{3} ={-1}$,  or    $ r ={-1}$
by Lemma \ref {3.1.1} {\rm (IV)} (IX) (X). It  is  quasi-affine  by Lemma \ref {3.1.3}.

(f)     i.e.  (7.3.14)$, q \notin R_3$,  by  {\rm  GDD}  $2$ of Row $10$.
 The sub-{\rm GDD} by deleting  Vertex  4   is an arithmetic  {\rm  GDD}   when $q \in  R_7 \cup R_{13} \cup R_{19} \cup R_{12}
 $
by Lemma \ref {3.1.1} {\rm (IV)} {\rm (IX)}.
It is quasi-affine  by Lemma \ref {3.1.3}.

(g)     i.e.  (7.3.15)$, q \notin R_3$,  by   {\rm  GDD}  $1$ of Row $10$.
 The sub-{\rm GDD} by deleting  Vertex  4    is arithmetic  {\rm  GDD}    when $q \in  R_5 \cup R_{8}
\cup R_6 $
by Lemma \ref {3.1.2} {\rm (IV)} and by Lemma \ref {3.1.1} {\rm (IV)}.
 It  is  quasi-affine   by Lemma \ref {3.1.3}.

(h)     i.e.  (7.3.16)$, q \notin R_3$,  by   {\rm  GDD}  $2$ of Row $10$.
 The sub-{\rm GDD} by deleting  Vertex  4    is an arithmetic  {\rm  GDD}    when $q \in  R_4 \cup R_5 \cup R_{8} \cup R_{11}
\cup R_6 $
by Lemma \ref {3.1.2} {\rm (IV)} and  by Lemma \ref {3.1.1} {\rm (IV)}.
 It  is  quasi-affine   by Lemma \ref {3.1.3}.

(i)     i.e.  (7.3.17)$, q \in R_3$,  by   {\rm  GDD}  $3$ of Row $14$. The sub-{\rm GDD} by
 deleting  Vertex  4  is  an arithmetic  {\rm  GDD}    by Lemma \ref {3.1.1}{\rm (I)}. It  is  quasi-affine  by Lemma \ref {3.1.3}.

(j)     i.e.  (7.3.18)$, q \in R_3$,  by   {\rm  GDD}  $2$ of Row $16$. The sub-{\rm GDD} by deleting
 Vertex  4 is  an arithmetic  {\rm  GDD}    by Lemma \ref {3.1.2}V). It  is  quasi-affine  by Lemma \ref {3.1.3}.

  {\rm (iii) }   All quasi-arithmetic {\rm  GDD}s by  adding a vertex on  Vertex $3$  are listed.
 According  to Lemma \ref {3.1.2} (IV) we have to consider following cases.

(a)     i.e.  (7.3.19)$, q^{3} \not= 1$,  by    {\rm  GDD}  $2$ of Row $4$. The sub-{\rm GDD} by deleting Vertex
  4 is  arithmetic {\rm  GDD}     when  $q \in  R_5 $ by Lemma \ref {3.1.1} (I) or  Type   2.  The sub-{\rm GDD} by deleting  Vertex  4    is  arithmetic  {\rm  GDD}     when  $q \in  R_4 $ by  Lemma \ref {3.1.1}{\rm (I)}.
 It  is quasi-affine  by Lemma \ref {3.1.3}.\\

  {\ }\ \ \ \ \ \ \ \ \ \ \ \ $\begin{picture}(100,       15) \put(-68,        -1){ $(b)$}

\put(60,       1){\makebox(0,       0)[t]{$\bullet$}}
\put(58,       -12){$-1$}

\put(40,       -12){$q^6$}
\put(28,       -1){\line(1,       0){33}}
\put(27,       1){\makebox(0,      0)[t]{$\bullet$}}

\put(22,      -12){$-1$}
\put(0,       -12){$q^{-6}$}

\put(-14,       1){\makebox(0,      0)[t]{$\bullet$}}

\put(-14,      -1){\line(1,       0){50}}

\put(-18,      -12){$q^{3}$}

\put(59,       0){\line(-1,      1){17}}

\put(28,      -1){\line(1,       1){17}}

\put(43,     18){\makebox(0,      0)[t]{$\bullet$}}

\put(18,     12){$q^{-4}$}

\put(36,     22){$q^2$}
\put(59,     12){$q^{-2}$}

 \ \ \ \ \ \ \ \ \ \ \ \ \ \ \ \ \ \ \ \ \ \ \ \  \ \ \ \ \ \ \ \
  \ \ \ \ \ \ \ \ \ \ \ \ \ \ \ \ \ \ \ {$, q^{3} \in F^{*}\setminus \{1,  -1\}$,  by    {\rm  GDD}  $2$ of Row $5$.}
\put(80,         -1)  {    } \end{picture}$\\ \\
 The sub-{\rm GDD} by deleting  Vertex  4    is not  an arithmetic  {\rm  GDD}   by Lemma \ref {3.1.1}{\rm (I)}.

(c)     i.e.  (7.3.20)$, q^{3} \in F^{*}\setminus \{1,  -1\}$,  by   {\rm  GDD}  $2$ of Row $6$. The sub-{\rm GDD} by deleting
Vertex  4  is an arithmetic  {\rm  GDD}   when $ q^{8} =1$  by Lemma \ref {3.1.1}{\rm (I)}.
 The sub-{\rm GDD} by deleting  Vertex  4   is an arithmetic  {\rm  GDD}   when  $q\in R_{18}$ by Lemma \ref {3.1.1}{\rm (I)}.
 It  is quasi-affine  by Lemma \ref {3.1.3}.

(d)     i.e.  (7.3.23)$, q^{3} \in F^{*}\setminus \{1,  -1\}$,  by    {\rm  GDD}  $3$ of Row $6$. The sub-{\rm GDD}
by deleting  Vertex  2  is an arithmetic  {\rm  GDD}   when $q\in R_5$  by Lemma \ref {3.1.1}{\rm (I)} or  Type   2.
 The sub-{\rm GDD} by deleting  Vertex  2  is an arithmetic  {\rm  GDD}   when $q\in R_{12}$  by Lemma \ref {3.1.1}{\rm (I)}.
 The sub-{\rm GDD} by deleting  Vertex  3  is an arithmetic  {\rm  GDD}    by Lemma \ref {3.1.1} {\rm (IV)}.
 It  is quasi-affine.

(e)     i.e.  (7.3.21)$, q^{3} \in F^{*}\setminus \{1,  -1\}$. $ q \notin R_3$,  by   {\rm  GDD}  $2$ of Row $7$.
 The sub-{\rm GDD} by deleting  Vertex  4   is an arithmetic  {\rm  GDD}    when $q^{11} =1$ or $q^{12} =-1$   by Type   2 and by Lemma \ref {3.1.1}{\rm (I)}.
 It  is quasi-affine  by Lemma \ref {3.1.3}.

(f)     i.e.  (7.3.24)$, q^{3} \in F^{*}\setminus \{1,  -1\}$. $ q^{} \notin R_3$,  by     {\rm  GDD}  $3$ of Row $7$.
 The sub-{\rm GDD} by deleting  Vertex  2    is  an arithmetic  {\rm  GDD}   by Lemma \ref {3.1.1}{\rm (I)}.
 The sub-{\rm GDD} by deleting  Vertex  1    is the same as  The sub-{\rm GDD} by deleting  Vertex  1.  It  is quasi-affine.

(g)     i.e.  (7.3.22)$, q^{3} \in F^{*}\setminus \{1,  -1\}$,  by  {\rm  GDD}  $2$ of Row $8$. The sub-{\rm GDD}
 by deleting Vertex  4  is an  arithmetic  {\rm  GDD}   when $q \in R_5$ by  Type   2.
 It  is quasi-affine  when $q \in R_5$ by Lemma \ref {3.1.3}.

  {\rm (iv) } All quasi-affine  {\rm  GDD}s which are complete diagrams are listed.
\\ \\ \\ \\ \\

{\ }\!\!\!\!\!\!\!\!\!\!\!\!\!\!\!\!\!\!\!\!\!{\ }\!\!\!\!\!\!\!\!\!\!\!\!\!\!\!\!\!\!\!\!\!{\ }\!\!\!\!\!\!\!\!\!\!\!\!\!\!\!\!\!\!\!\!\!{\ }\!\!\!\
\ \ \ \ \ $\begin{picture}(100,       15)\put(70,       -1){(b) in  Case {\rm (i) } }

\put(170,     10){\makebox(0,      0)[t]{$\bullet$}}

\put(170,     70){\makebox(0,      0)[t]{$\bullet$}}

\put(230,     10){\makebox(0,      0)[t]{$\bullet$}}
\put(230,     70){\makebox(0,      0)[t]{$\bullet$}}

\put(170,       10){\line(0,      1){60}}

\put(170,       10){\line(1,      1){60}}

\put(230,       10){\line(0,      1){60}}

\put(230,       10){\line(-1,      1){60}}

\put(170,       10){\line(1,       0){60}}
\put(170,       70){\line(1,       0){60}}

\put(150,     10){$q$}
\put(150,     30){$q^{-1}$}

\put(150,     70){$-1$}

\put(180,     30){$q^{-1}$}
\put(190,     80){$q^{2}$}
\put(190,     -10){$q^{-2}$}
\put(220,     30){$q^{3}$}

\put(250,     10){$-1$}
\put(250,     30){$q^{3}$}

\put(250,     70){$-1$}

\put(290,         -1)  {   with  $q\notin R_3.$} \end{picture}$\\ \\
 The sub-{\rm GDD} by deleting   Vertex 2 is not an arithmetic  {\rm  GDD}   by Lemma \ref {3.1.1}{\rm (I)}.  \\ \\ \\ \\ \\

{\ }\!\!\!\!\!\!\!\!\!\!\!\!\!\!\!\!\!\!\!\!\!{\ }\!\!\!\!\!\!\!\!\!\!\!\!\!\!\!\!\!\!\!\!\!{\ }\!\!\!\!\!\!\!\!\!\!\!
 $\begin{picture}(100,       15)\put(68,       -1){(d) in  Case {\rm (i) } }

\put(170,     10){\makebox(0,      0)[t]{$\bullet$}}

\put(170,     70){\makebox(0,      0)[t]{$\bullet$}}

\put(230,     10){\makebox(0,      0)[t]{$\bullet$}}
\put(230,     70){\makebox(0,      0)[t]{$\bullet$}}

\put(170,       10){\line(0,      1){60}}

\put(170,       10){\line(1,      1){60}}

\put(230,       10){\line(0,      1){60}}

\put(230,       10){\line(-1,      1){60}}

\put(170,       10){\line(1,       0){60}}
\put(170,       70){\line(1,       0){60}}

\put(150,     10){$q$}
\put(150,     30){$q^{-2}$}

\put(150,     70){$-1$}

\put(180,     30){$q^{-1}$}
\put(190,     80){$q^{3}$}
\put(190,     -10){$q^{-2}$}
\put(220,     30){$$}

\put(250,     10){$-1$}
\put(250,     30){$q^{3}$}

\put(250,     70){$-1$}

\put(80,         -1)  {    } \end{picture}$\\ \\
 The sub-{\rm GDD} by deleting   Vertex 2 is not an arithmetic  {\rm  GDD}   by Lemma \ref {3.1.1}{\rm (I)}.
\\ \\ \\ \\ \\

{\ }\!\!\!\!\!\!\!\!\!\!\!\!\!\!\!\!\!\!\!\!\!{\ }\!\!\!\!\!\!\!\!\!\!\!\!\!\!\!\!\!\!\!\!\!{\ }\!\!\!\!\!\!\!\!\!\!\!\!\!\!\!
   $\begin{picture}(100,       15)\put(68,       -1){
(c) in  Case {\rm (ii) }}

\put(170,     10){\makebox(0,      0)[t]{$\bullet$}}

\put(170,     70){\makebox(0,      0)[t]{$\bullet$}}

\put(230,     10){\makebox(0,      0)[t]{$\bullet$}}
\put(230,     70){\makebox(0,      0)[t]{$\bullet$}}

\put(170,       10){\line(0,      1){60}}

\put(170,       10){\line(1,      1){60}}

\put(230,       10){\line(0,      1){60}}

\put(230,       10){\line(-1,      1){60}}

\put(170,       10){\line(1,       0){60}}
\put(170,       70){\line(1,       0){60}}

\put(150,     10){$-1$}
\put(150,     30){$q^{2}$}

\put(150,     70){$-1$}

\put(180,     30){$q^{-1}$}
\put(190,     80){$q^{-1}$}
\put(190,     -10){$q^{3}$}
\put(220,     30){$q^{3}$}

\put(250,     10){$-1$}
\put(250,     30){$q^{-2}$}

\put(250,     70){$q$}

\put(80,         -1)  {    } \end{picture}$\\ \\
with  $q\notin R_3.$
 the sub-{\rm GDD} by deleting  Vertex 2 is not an  arithmetic  {\rm  GDD}   by Lemma \ref {3.1.1}{\rm (I)}.
\\ \\ \\ \\ \\

{\ }\!\!\!\!\!\!\!\!\!\!\!\!\!\!\!\!\!\!\!\!\!{\ }\!\!\!\!\!\!\!\!\!\!\!\!\!\!\!\!\!\!\!\!\!{\ }\!\!\!\!\!\!\!  $\begin{picture}(100,       15)\put(68,       -1){(d) in  Case {\rm (iii) } }

\put(170,     10){\makebox(0,      0)[t]{$\bullet$}}

\put(170,     70){\makebox(0,      0)[t]{$\bullet$}}

\put(230,     10){\makebox(0,      0)[t]{$\bullet$}}
\put(230,     70){\makebox(0,      0)[t]{$\bullet$}}

\put(170,       10){\line(0,      1){60}}

\put(170,       10){\line(1,      1){60}}

\put(230,       10){\line(0,      1){60}}

\put(230,       10){\line(-1,      1){60}}

\put(170,       10){\line(1,       0){60}}
\put(170,       70){\line(1,       0){60}}

\put(150,     10){$-1$}
\put(150,     30){$q^{-3}$}

\put(150,     70){$q^{3}$}

\put(180,     30){$q^{3}$}
\put(190,     80){$q^{-3}$}
\put(190,     -10){$q^{-2}$}
\put(220,     30){$$}

\put(250,     10){$q$}
\put(250,     30){$q^{-1}$}

\put(250,     70){$-1$}

\put(280,         -1)  {    the sub-{\rm GDD} by deleting  Vertex 2 is not } \end{picture}$
\\ \\  arithmetic  {\rm  GDD} by Lemma \ref {3.1.1}{\rm (I)}.

(f)  in  Case {\rm (iii) } i.e.  (7.3.1)
 If  the sub-{\rm GDD} by deleting  Vertex 2 is  arithmetic
 {\rm  GDD},   then and $ q^{-2} =q^{3}$ and $\widetilde{q}_{14} =q^{-2}$
by Lemma \ref {3.1.1}{\rm (I)}. The sub-{\rm GDD} by deleting  Vertex 3  is  arithmetic  {\rm  GDD}   by Lemma \ref {3.1.1}{\rm (I)},  It  is quasi-affine  when $ \widetilde{q}_{14} =q^{-2}$ and $q^{}\in R_5.$

\subsection* {Quasi-affine over
 {\rm  GDD}  $4$ of Row $7$ }
   {\rm (i) } All quasi-arithmetic {\rm  GDD}s by  adding  a vertex on  Vertex $1$  are listed.  According  to Lemma \ref {3.1.2} (III) we have to consider following cases.

 (a)  i.e.  (7.4.3), $q^{3}\notin R_3\cup R_2$,   $q^{}\notin R_3\cup R_2$,  by    {\rm  GDD}  $1$ of Row $7$. It
is quasi-affine  by Lemma \ref {3.1.3}.

  {\rm (ii) }  All quasi-arithmetic {\rm  GDD}s by  adding a vertex on  Vertex $3$  are listed. According  to  Type   3 we have to consider following cases.

(a) i.e. (7.4.2),
 $q^{}\notin R_3\cup R_2$,  by  Type   3.
 The sub-{\rm GDD} by deleting  Vertex  4   is an arithmetic  {\rm  GDD}   when $q \in  R_5 \cup R_7 \cup R_9 \cup R_6\cup R_4 $ by Lemma \ref {3.1.1} {\rm (IV)} {\rm (IX)}.
 It  is quasi-affine  by Lemma \ref {3.1.3}.

 (b)  i.e.  (7.4.1),   $q^{}\notin R_3\cup R_2$,  by   Type   3. The sub-{\rm GDD} by deleting  Vertex  4
  is an arithmetic  {\rm  GDD}    by Lemma \ref {3.1.1} {\rm (IV)}.
 It  is quasi-affine  by Lemma \ref {3.1.3}.

\subsection* {Quasi-affine over
 {\rm  GDD}  $1$ of Row $8$ }
   {\rm (i) } All quasi-arithmetic {\rm  GDD}s by  adding a vertex on  Vertex $1$  are listed..
According  to Lemma \ref {3.1.2} (III) we have to consider following cases.\\

  {\ }\ \ \ \ \ \ \ \ \ \ \ \ $\begin{picture}(100,       15) \put(-68,        -1){ (a)}

\put(60,       1){\makebox(0,       0)[t]{$\bullet$}}
\put(58,       10){$-1$}

\put(40,       5){$q^{-1}$}
\put(28,       -1){\line(1,       0){33}}
\put(27,       1){\makebox(0,      0)[t]{$\bullet$}}

\put(22,      10){$q$}
\put(0,       5){$q^{-2}$}

\put(-14,       1){\makebox(0,      0)[t]{$\bullet$}}

\put(-14,      -1){\line(1,       0){50}}

\put(-18,      10){$q^2$}

\put(104,       1){\makebox(0,     0)[t]{$\bullet$}}

\put(62,      -1){\line(1,       0){40}}

\put(75,      5){$q$}

\put(100,       10){$q^{-1}$}\ \ \ \ \ \ \ \ \ \ \ \ \ \ \ \ \ \ \ \ \ \ \ \ \ \ \ \ \ \ \ \ \ \   \ \ \ \ \ \ \ \ \ \ \ \ \ \ \ \ \ \ \ { $, q \in F^{*}\setminus \{1,  -1\}$,  by   {\rm  GDD}  $1$ of Row $6$.}
\put(80,         -1)  {    } \end{picture}$\\
It is an arithmetic  {\rm  GDD}   by  Type   1.

   (b)  i.e.  (8.1.1),   $q^2,  q^3 \not=1, $ by  {\rm  GDD}  $1$ of Row $7$. It is  an
arithmetic
  {\rm  GDD}  when $q\in R_4$ by  {\rm  GDD}  1 of Row 22;  Otherwise it is quasi-affine
by Lemma \ref {3.2.38}.

  {\rm (ii) }  All quasi-arithmetic {\rm  GDD}s by  adding  a vertex on  Vertex $2$  are listed.
According  to Lemma \ref {3.1.2} (I) we have to consider following cases.

 (a)  i.e.  (8.1.2)  by   {\rm  GDD}  $2$ of Row $4$. The sub-{\rm GDD} by deleting  Vertex  4    is
arithmetic  {\rm  GDD}   when $q\in R_3\cup R_4$ by Lemma \ref {3.1.1} {\rm (IX)}.  It  is quasi-affine  by Lemma \ref {3.1.35}
 when $q\in R_4$.  It  is  {\rm  GDD}  $9$ of Row $20$
 when $q\in R_3$.

 (b)  i.e.  (8.1.3)
   by   {\rm  GDD}  $2$ of Row $6$. The sub-{\rm GDD} by deleting  Vertex  4    is  an
arithmetic {\rm  GDD}    when $q \in  R_3 \cup R_4 \cup R_5  $ by Lemma \ref {3.1.1} (IX). It  is quasi-affine  when $q\notin R_3$
by Lemma \ref {3.1.16} or by Lemma \ref {3.1.35}.    It  is an arithmetic  {\rm  GDD}   when $q\in R_3$
by Lemma \ref {3.1.35}.

 (c)  i.e.  (9.1.4)
  by   {\rm  GDD}  $2$ of Row $7$. The sub-{\rm GDD} by deleting  Vertex 4    is  an
arithmetic {\rm  GDD}     when $q \in  R_4 \cup R_6 \cup R_5 $ by Lemma \ref {3.1.2} (IV). It  is quasi-affine    by Lemma \ref {3.1.35}.

 (d)  i.e.  (8.1.5)
,  $q\notin R_2 \cup R_3$,  by   {\rm  GDD}  $4$ of Row $7$.
 The sub-{\rm GDD} by deleting  Vertex  4    is  an arithmetic  {\rm  GDD}     when $q\in R_5 \cup R_6$  by  Type   3 and by Lemma \ref {3.1.1}{\rm (I)}.  The sub-{\rm GDD} by deleting  Vertex  4    is not  an arithmetic  {\rm  GDD}     when $q\notin R_5 \cup R_6$  by Lemma \ref {3.1.1}{\rm (I)}.  It  is   quasi-affine  by Lemma \ref {3.1.3}.\\
  \\

 {\ }\ \ \ \ \ \ \ \ \ \ $\begin{picture}(100,       15) \put(-58,        -1){(e) }

\put(60,       1){\makebox(0,       0)[t]{$\bullet$}}

\put(28,       -1){\line(1,       0){33}}
\put(27,       1){\makebox(0,      0)[t]{$\bullet$}}
\put(-14,       1){\makebox(0,      0)[t]{$\bullet$}}

\put(-14,      -1){\line(1,       0){50}}

\put(58,       -12){$q$}

\put(40,       -12){$q^{-1}$}

\put(22,      -12){$-1$}
\put(0,       -12){$q{}$}

\put(-18,      -12){$q^{-1}$}

\put(27,     38){\makebox(0,      0)[t]{$\bullet$}}

\put(27,       0){\line(0,      1){35}}

\put(30,       30){$q^{-1}$}

\put(30,      15){$q$}

\ \ \ \ \ \ \ \ \ \ \ \ \ \ \ \ \ \ \ \ \ \ {,  by   {\rm  GDD}  $1$ of Row $8$. The sub-{\rm GDD} by deleting  Vertex  4   is an }
\put(80,         -1)  {    } \end{picture}$\\ \\ arithmetic
{\rm  GDD}   by Lemma \ref {3.1.1} {\rm (IV)}.
 It  is an arithmetic  {\rm  GDD}    by  {\rm  GDD}  6 Row 12 or  Type   5.

 (f)  i.e.  (8.1.6),
  $q \notin  R_2 \cup R_4  $,  by  {\rm  GDD} $3$ of Row $5$.
 The sub-{\rm GDD} by deleting  Vertex  4  is   {\rm  GDD} $4$ of Row $7$  when $q \in  R_5 $;
  is not an arithmetic {\rm  GDD}   when $q \notin  R_5 $  by Lemma \ref {3.1.35}.
It is quasi-affine  by Lemma \ref {3.1.3}.

  {\rm (iii) }  All quasi-arithmetic {\rm  GDD}s by  adding a vertex on  Vertex $3$  are listed.
It is the same as   {\rm (i)}.

  {\rm (iv) }
 All quasi-affine circles are listed. \\ \\  \\

 {\ }\ \ \ \ \ \ \ \ \ \ \ \    \ \ \ \ \ \  $\begin{picture}(100,       15) \put(-85,       -1){ (a)  (nc)}

\put(60,       1){\makebox(0,       0)[t]{$\bullet$}}

\put(28,       -1){\line(1,       0){33}}
\put(27,       1){\makebox(0,      0)[t]{$\bullet$}}

\put(-14,       1){\makebox(0,      0)[t]{$\bullet$}}

\put(-14,      -1){\line(1,       0){50}}

\put(26,     38){\makebox(0,      0)[t]{$\bullet$}}

\put(-18,     - 15){$q$}
\put(0,       -15){$q^{-1}$}
\put(22,     - 15){$-1$}
\put(40,       -15){$q^{}$}

\put(58,      - 15){$q^{-1}$}


\put(30,       40){$q^{2}$}

\put(-12,       20){$q^{-2}$}

\put(58,       10){$q^{}$}

\put(-12,      1){\line(1,      1){35}}

\put(60,      1){\line(-1,      1){35}}

\put(80,         -1)  { is quasi-affine   since the sub-{\rm GDD} by deleting  Vertex 2
   } \end{picture}$\\ \\

  $\begin{picture}(100,       15) \put(-125,       -1){  }

\put(60,       1){\makebox(0,       0)[t]{$\bullet$}}

\put(28,       -1){\line(1,       0){33}}
\put(27,       1){\makebox(0,      0)[t]{$\bullet$}}

\put(-14,       1){\makebox(0,      0)[t]{$\bullet$}}

\put(-14,      -1){\line(1,       0){50}}

\put(-18,      10){$q$}
\put(0,       5){$q^{-2} $}
\put(22,      10){$q^{2} $}
\put(40,       5){$q^{} $}

\put(58,       10){$q^{-1} $}

\put(80,         -1)  {   with  $q^{-1} =q^{2} $ is an  arithmetic  {\rm  GDD}   by  Type   2. } \end{picture}$

\subsection* {Quasi-affine over
 {\rm  GDD}  $2$ of Row $8$  }
  {\rm (i) }  All quasi-arithmetic {\rm  GDD}s by  adding  a vertex on  Vertex $1$  are listed.
According  to Lemma \ref {3.1.2} (IV) we have to consider following cases.

   (a)  i.e.  (8.2.1)$,  q^2, q^3 \not=1$  by   {\rm  GDD}  $2$ of Row $7$. It is an
arithmetic
{\rm  GDD}    when $q\in R_4$ by  {\rm  GDD}  2 of Row 22;   It  is quasi-affine  when $q\notin R_4$ by Lemma \ref {3.1.21}.

  {\rm (ii) } All quasi-arithmetic {\rm  GDD}s by  adding a vertex on  Vertex $2$  are listed.
According  to Lemma \ref {3.1.2} (IV) we have to consider following cases.

 (a)  i.e.  (8.2.2)
 $, q \in F^{*}\setminus \{1,  -1\}$,  by   {\rm  GDD}  $2$ of Row $4$. The sub-{\rm GDD} by deleting
Vertex  4    is  an arithmetic  {\rm  GDD}    when  $q\in R_3\cup R_4$  by Lemma \ref {3.1.1} {\rm (IX)}. It  is quasi-affine   when    $q\in  R_4$ by Lemma \ref {3.1.22}. It  is an arithmetic  {\rm  GDD}    when    $q\in  R_3$ by Lemma \ref {3.1.22}.

  (b)  i.e.  (8.2.3)$, q \in F^{*}\setminus \{1,  -1\}$,  by   {\rm  GDD}  $2$ of Row $6$. The sub-{\rm GDD} by deleting
 Vertex  4  is
an arithmetic  {\rm  GDD}    when $q \in  R_4 \cup R_5 \cup R_3  $ by Lemma \ref {3.1.1} {\rm (IX)}.  It  is an arithmetic  {\rm  GDD}    when $q \in  R_3  $ by Lemma \ref {3.1.22}.
 It  is  quasi-affine   when $q \in  R_4 \cup R_5 $ by Lemma \ref {3.1.22}.

 (c)  i.e.  (8.2.4)
 $, q \in F^{*}\setminus \{1,  -1\}$. $ q \notin R_3$,  by  {\rm  GDD}  $2$ of Row $7$.
 The sub-{\rm GDD} by deleting  Vertex  4    is an arithmetic  {\rm  GDD}    when $q \in  R_4 \cup R_5 \cup R_6 $ by Lemma \ref {3.1.1} {\rm (IX)}.  It  is quasi-affine     when $q \in  R_5 \cup R_6$ by Lemma \ref {3.1.22}. It  is {\rm  GDD}  $4$ of Row $18$    when $q \in  R_4. $

 (d)  i.e.  (8.2.5)$, q \in F^{*}\setminus \{1,  -1\}$,  by    {\rm  GDD}  $2$ of Row $8$.
 The sub-{\rm GDD} by deleting  Vertex  4    is an arithmetic  {\rm  GDD}    when  $q\in R_3$ by  {\rm  GDD}  2 of Row 15 or by Lemma \ref {3.1.1} {\rm (X)}.  It  is quasi-affine  by Lemma \ref {3.1.22}  when  $q\notin R_2\cup R_3$.

  {\rm (iii) } All quasi-arithmetic {\rm  GDD}s by  adding a vertex  on  Vertex $3$  are listed.

It is the same as  Case {\rm (i) }when we change q in  Case {\rm (i) }into $q^{-1}$.

  {\rm (iv) } All quasi-affine circles are listed.

(nc) (a)  i.e.  (8.2.6)
  is quasi-affine  since the sub-{\rm GDD} by deleting  Vertex 2
  with  $q =q^{-3} $ is an  arithmetic  {\rm  GDD}.

 (nc)  (a)  i.e.  (8.2.7) is quasi-affine  since the sub-{\rm GDD} by deleting  Vertex 2 with $q^{2} =q^{-3} $ is  arithmetic  {\rm  GDD}.

  (a)  (a)  i.e.  (8.2.8)  is quasi-affine  since the sub-{\rm GDD} by deleting  Vertex 2 with  $q ^3=q^{-3} $ is an  arithmetic  {\rm  GDD}.

 (a) (nc)  i.e.  (8.2.9)
 is quasi-affine  since the sub-{\rm GDD} by deleting  Vertex 2
 with  $q ^3={-1} $ is  arithmetic  {\rm  GDD}   by Lemma \ref {3.1.1} {\rm (X)}.

\subsection* {Quasi-affine over
 {\rm  GDD}  $3$ of Row $8$ }
  {\rm (i) }  All quasi-arithmetic {\rm  GDD}s by  adding a vertex  on  Vertex $1$  are listed.
According  to Lemma \ref {3.1.2} (I) we have to consider following cases.

 (a)  i.e.  (8.3.2)
  by   {\rm  GDD}  $2$ of Row $7$. It is quasi-affine   by Lemma \ref {3.1.21}.

 (b)  i.e.  (8.3.3)
, $q \notin  R_2 \cup R_3 $,  by   {\rm  GDD}  $4$ of Row $7$. It is quasi-affine by Lemma \ref {3.1.3}.

  {\rm (ii) }  All quasi-arithmetic {\rm  GDD}s by  adding a vertex  on  Vertex $2$  are listed..

  (a)  i.e.  (8.3.4), $q^2, q^3 \not=1,$   by   {\rm  GDD}  $1$ of Row $7$. It is quasi-affine  by Lemma \ref {3.1.3}.

   {\rm (iv) }  All quasi-affine circles are listed.

 (nc)(a)  i.e.  (8.3.6)
     is quasi-affine  since the sub-{\rm GDD} by deleting  Vertex 2
is an  arithmetic  {\rm  GDD}   when  $ q^{3} =q^{-1}$ or  $ q^{3} ={-1}$ by Lemma \ref {3.1.1} (X).

 (nc)(a)  i.e.  (8.3.7)
is quasi-affine since the sub-{\rm GDD} by deleting  Vertex 2
  is  an arithmetic  {\rm  GDD}   when $q^{3}= q^{-2}$. \\ \\

  {\ }\ \ \ \ \ \ \ \ \ \ \ \ $\begin{picture}(100,       15) \put(-68,        -1){  (a) (nc)}

\put(60,       1){\makebox(0,       0)[t]{$\bullet$}}

\put(28,       -1){\line(1,       0){33}}
\put(27,       1){\makebox(0,      0)[t]{$\bullet$}}

\put(-14,       1){\makebox(0,      0)[t]{$\bullet$}}

\put(-14,      -1){\line(1,       0){50}}

\put(26,     38){\makebox(0,      0)[t]{$\bullet$}}

\put(-18,     - 15){$-1$}
\put(0,       -15){$q^{-3}$}
\put(22,     - 15){$q^{3}$}
\put(40,       -15){$q^{-3}$}

\put(58,      - 15){$-1$}


\put(30,       40){$-1$}

\put(-12,       20){$q^{}$}

\put(58,       10){$q^{6}$}

\put(-12,      1){\line(1,      1){35}}

\put(60,      1){\line(-1,      1){35}}

\put(80,         -1)  {  is quasi-affine since the sub-{\rm GDD} by deleting  Vertex 2
 } \end{picture}$\\ \\

  $\begin{picture}(100,       15) \put(-68,        -1){ }

\put(60,       1){\makebox(0,       0)[t]{$\bullet$}}

\put(28,       -1){\line(1,       0){33}}
\put(27,       1){\makebox(0,      0)[t]{$\bullet$}}

\put(-14,       1){\makebox(0,      0)[t]{$\bullet$}}

\put(-14,      -1){\line(1,       0){50}}

\put(-18,      10){$-1$}
\put(0,       5){$q$}
\put(22,      10){$-1$}
\put(40,       5){$q^{6}$}

\put(58,       10){$-1$}

\put(80,         -1)  { with $-1 =q^{-6} $ is an arithmetic  {\rm  GDD}   by Lemma \ref {3.1.1} {\rm (X)}.   } \end{picture}$

 (c) (nc) i.e.  (8.3.5)
    is quasi-affine since the sub-{\rm GDD} by deleting  Vertex 2
  with $q^{-3} =-q^{-1} $ is  arithmetic  {\rm  GDD}    by  Type   1.

\subsection* {Quasi-affine over
 {\rm  GDD}  $1$ of Row $9$  }
   {\rm (i) } All quasi-arithmetic {\rm  GDD}s by  adding a vertex  on  Vertex $1$  are listed.
 According  to Lemma \ref {3.1.2} (III) we have to consider following cases.

 (a)  i.e.  (9.1.1)
 $, q \in F^{*}\setminus \{1,  -1\}$ $, q \notin R_3$,  by   {\rm  GDD}  $1$ of Row $7$.
is
quasi-affine  by Lemma \ref {3.1.8}.

 (ii)  All quasi-arithmetic {\rm  GDD}s by  adding a vertex  on  Vertex $2$  are listed.
 According  to Lemma \ref {3.1.2} (I) we have to consider following cases.

 (a)  i.e.  (9.1.2)
  by   {\rm  GDD}  $2$ of Row $4$.
 The sub-{\rm GDD} by deleting  Vertex  4    is arithmetic  {\rm  GDD}    when $r= -q, $ and $q\in R_3$ or  $r= q^3, $ and $q^4 \not=1$ or
 $r= q^2$ and $r\notin R_3$ or $r=-1$ and $q\in R_3\cup R_4\cup R_6$   by Lemma \ref {3.1.1} {\rm (IX)}.  It  is quasi-affine  by Lemma \ref {3.1.16}.

 (b)  i.e.  (9.1.3)
  by    {\rm  GDD}  $3$ of Row $5$. The sub-{\rm GDD} by deleting  Vertex  4    is arithmetic  {\rm  GDD}    when $ r^{-2} = q^{-3}$ and $q^3 \not=1$  or $ r^{-2} = q^{-2}$   by Lemma \ref {3.1.1}{\rm (I)}  and    Type   3.
 It  is  quasi-affine  by Lemma \ref {3.1.3} when   $ r^{-2} = q^{-3}$ and $q^3 \not=1$  or $ r^{-2} = q^{-2}$.

 (c)  i.e.  (9.1.4)
  by    {\rm  GDD}  $2$ of Row $6$. The sub-{\rm GDD} by deleting  Vertex  4    is arithmetic  {\rm  GDD}    when $ r^{2} = q$,  with $q^2,  q^3 \not=1; $
 or $ r^{2} = q^3$,  $q^3,  q^5 \not=1; $ or $ r^{2} = -q$, $q\in R_3; $
 or $ r^{2} = -1$,  $q\in R_3\cup R_6.$ It  is
  quasi-affine  by Lemma \ref {3.1.16}.

 (d)  i.e.  (9.1.5),
$q^2,  q^3 \not=1, $ by   {\rm  GDD}  $2$ of Row $7$. The sub-{\rm GDD} by deleting  Vertex  4
is arithmetic  {\rm  GDD}    when $r^3 =q$ or $r^3 =q^2$,  $q^5 \not=1$
or   $q^5 =-1$,  $q\in R_4$  by Lemma \ref {3.1.1} {\rm (IX)} {\rm (X)}.
 It  is   quasi-affine   when $r^3 =q$ or $r^3 =q^2$,  $q^5 \not=1$ or   $q^5 =-1$,  $q\in R_4$  by Lemma \ref {3.1.16}.

  (e)  i.e.  (9.1.6), $ q^3 \not=1$,  by  {\rm  GDD}  $4$ of Row $7$.
 The sub-{\rm GDD} by deleting  Vertex  4    is   arithmetic  {\rm  GDD}   when  $ r^{3} =q^{2}$,  $q\notin R_5$ by Lemma \ref {3.1.1}{\rm (I)}  or  Type   3.  It  is quasi-affine  by Lemma \ref {3.1.3}\\
\\

 {\ }\ \ \ \ \ \ \ \ \ \ $\begin{picture}(100,       15) \put(-58,        -1){ (f) }

\put(60,       1){\makebox(0,       0)[t]{$\bullet$}}

\put(28,       -1){\line(1,       0){33}}
\put(27,       1){\makebox(0,      0)[t]{$\bullet$}}
\put(-14,       1){\makebox(0,      0)[t]{$\bullet$}}

\put(-14,      -1){\line(1,       0){50}}

\put(58,       -12){$q$}

\put(40,       -12){$q^{-1}$}

\put(22,      -12){$-1$}
\put(0,       -12){$q{}$}

\put(-18,      -12){$q^{-1}$}

\put(27,     38){\makebox(0,      0)[t]{$\bullet$}}

\put(27,       0){\line(0,      1){35}}

\put(30,       30){$r$}

\put(30,       15){$r^{-1}$}

\ \ \ \ \ \ \ \ \ \ \ \ \ \ \ \ \ \ \ \ \ \ \ \ \ \ \ \   {,  by    {\rm  GDD}  $1$ of Row $8$. The sub-{\rm GDD} by deleting  Vertex  4    is   }
\put(80,         -1)  {    } \end{picture}$\\ \\
arithmetic {\rm  GDD}    by Lemma \ref {3.1.1} {\rm (IV)}. It  is quasi-affine  by Lemma \ref {3.1.9} and  by Lemma \ref {3.1.4}.

 (g)  i.e.  (9.1.7)
   by    {\rm  GDD}  $1$ of Row $9$. The sub-{\rm GDD} by deleting  Vertex  4    is
 arithmetic {\rm  GDD}    by Lemma \ref {3.1.1} {\rm (IV)}. It  is quasi-affine  by Lemma \ref {3.1.9} and  by Lemma \ref {3.1.4}.

\subsection* {Quasi-affine over
 {\rm  GDD}  $2$ of Row $9$ }

 {\rm  GDD}  $2$ of Row $9$ is the same  as  {\rm  GDD}  $1$ of Row $9$.

\subsection* {Quasi-affine over
 {\rm  GDD}  $3$ of Row $9$   }

 Let $C= \{ (q, s, r) \mid    qrs =1,  q\not= r,  q\not= s,  s\not= r     \}$.

  {\rm (i) } All quasi-arithmetic {\rm  GDD}s by  adding a vertex  on  Vertex $1$  are listed.
 According  to Lemma \ref {3.1.2} (IV) we have to consider following cases.

  (a)  i.e.  (9.3.1)$, q \in F^{*}\setminus \{1,  -1\}$,  by    {\rm  GDD}  $2$ of Row $4$. The sub-{\rm GDD}
by deleting  Vertex  4    is an arithmetic  {\rm  GDD}    when $q =s^2, $ $ s^3 \not=1, $ or $q =s^3, $ with  $s^2,  s^3,    s^4 \not=1, $ or $q =-s^{}$,   $s\in R_3$  or $s ={-1}$,  $q^2,  r^2  \not=1$ by Lemma \ref {3.1.1} {\rm (IX)}.
It is  GDD $2$ of Row $14$ when  $s = -1, $ $r = -q^{-1}. $ It  is quasi-affine  by Lemma \ref {3.1.18}, $(q, s, r)\in (A-B) \cap C, $   where $A = \{ (q, s, r) \mid  $
$q =s^2, $ $ s^3 \not=1, $ or $q =s^3; $ $s^2,  s^3,    s^4 \not=1 $
 or $q =-s^{}$,   $s\in R_3$  or $s ={-1}$,  $q^2,  r^2  \not=1 \}$, $B= \{   (q, s, r) \mid $  $s =-1^{}$, $r =-s^{-1}$   $ \}.$

  (b)  i.e.  (9.3.2)
$, q \in F^{*}\setminus \{1,  -1\}$,  by    {\rm  GDD}  $2$ of Row $5$. The sub-{\rm GDD} by deleting
Vertex  4    is an arithmetic  {\rm  GDD}    when $s=-q $ by  Type   2.  It  is quasi-affine  by Lemma \ref {3.1.3}.

   (c)  i.e.  (9.3.3), $q^2 \not=1,$   by    {\rm  GDD}  $2$ of Row $6$. The sub-{\rm GDD} by deleting Vertex  4    is an
arithmetic  {\rm  GDD}     when $q^2=s$, or $s^2 = q^2$,  or $s^3 = q^2$,  or  $q^2 = -s$,  $s\in R_3$ or $s=-1$ or
$q^2 = -1$,    $s \in  R_2 \cup R_3 \cup R_4 \cup R_6 $
  by Lemma \ref {3.1.1} (IX)  (IV) and (X).
 It  is  quasi-affine  by Lemma \ref {3.1.22}. $(q, s, r) \in (A-B)\cap C,$
 where $A = \{ (q, s, r) \mid $ $q^2=s$, or $s^2 = q^2$,  or $s^3 = q^2$,  or  $q^2 = -s$,  $s\in R_3$ or $s=-1$
                or
$q^2 = -1$,    $s \in  R_2 \cup R_3 \cup R_4 \cup R_6 $;  $q \in F^{*}\setminus \{1,  -1\}$;
$q,  r,  s \in F^{*}\setminus \{1\},  qrs =1,  q\not= r,  q\not= s,  s\not= r.$
            $ \}$, $B= \{   (q, s, r) \mid $  $s =q^{2}$, $r =q^{-3}$   $ \}.$

 (d)  i.e.  (9.3.4)
  $, q \in F^{*}\setminus \{1,  -1\}$,  by     {\rm  GDD}  $3$ of Row $6$.
 The sub-{\rm GDD} by deleting  Vertex  2    is an  arithmetic  {\rm  GDD}     when $q = s^2$,
$1 \not= s^3$,  or $q = s^3$ with $s^2,  s^3,  s^4 \not=1, $
    or  $s = -1$  or $s = -q$,  $s\in R_3$
 by Lemma \ref {3.1.1} (IX).  The sub-{\rm GDD} by deleting  Vertex  1    is  an arithmetic  {\rm  GDD}     when $q = r^2$,
$1 \not= r^3$,  or $q = r^3$ with $r^2,  r^3,  r^4 \not=1, $
    or  $r = -1$  or $r = -q$,  $r\in R_3$
 by Lemma \ref {3.1.1} (IX).
 It  is  quasi-affine.

 (e)  i.e.  (9.3.5)$, q \in F^{*}\setminus \{1,  -1\}$. $ q \notin R_3$,  by   {\rm  GDD}  $2$ of Row $7$.
 The sub-{\rm GDD} by deleting  Vertex  4    is  an arithmetic  {\rm  GDD}     when $q^3=s$  or $s^2 = q^3$,  or $q^3 = s^3$ or $-s^{} = q^{3}$,  $s\in R_3$  or $q^3 = -1$,   $s \in  R_2 \cup R_3 \cup R_4 \cup R_6 $;  $q\notin R_2\cup R_3$
 by Lemma \ref {3.1.1} (IX) and by Lemma \ref {3.1.2} (V).
 It  is quasi-affine  by Lemma \ref {3.1.22}.

 (f)  i.e.  (9.3.6), $q\notin R_2\cup R_3$,  by      {\rm  GDD}  $3$ of Row $7$.
 The sub-{\rm GDD} by deleting  Vertex  2    is an arithmetic  {\rm  GDD}     when  $ s^{3} =q^{}$
or  $ s^{6} =q^{}$ or  $ s^{9} =q^{}$  or  $- s^{3} =q^{}$,  $s\in R_9$
or $ s^{3} =-1^{}$ by Lemma \ref {3.1.1} (IX)(IV).  The sub-{\rm GDD} by deleting  Vertex  1    is an arithmetic  {\rm  GDD}     when $ r^{3} =q^{2}$ by  Type   2.
 It  is  quasi-affine.

 (g)  i.e.  (9.3.7)
  $, q \in F^{*}\setminus \{1,  -1\}$,  by   {\rm  GDD}  $2$ of Row $8$.
 The sub-{\rm GDD} by deleting  Vertex  4   is an arithmetic  {\rm  GDD}   when  $ s ={-1}$,   $q \in  R_3 \cup R_4 \cup R_6  $ or $s =- q^{-1}$ with  $s\in R_3\cup R_6$ by Lemma \ref {3.1.1} {\rm (X)}. It  is quasi-affine  by Lemma \ref {3.1.22}.

 (h)  i.e.  (9.3.8)$, q,  r,  s \in F^{*}\setminus \{1\},  qrs =1,  q\not= r,  q\not= s,  s\not= r$
 by    {\rm  GDD}  $3$ of Row $9$. $ q',  r',  s' \in F^{*}\setminus \{1\},  q'r's' =1,  q'\not= r',  q'\not= s',  s'\not= r'.$ $ q^{} =q'$ The sub-{\rm GDD} by deleting  Vertex  2    is  an arithmetic  {\rm  GDD}     when $ s' =s, $ $s\in R_3$ or $ s' =s^{-1}, $ or $ s' =-s^{}$ with  $s\in R_3\cup  R_6$,   or $ s' ={-1},  s \in  R_2 \cup R_3 \cup R_4 \cup R_6 $  or $ s ={-1}, $  $s' \in  R_2 \cup R_3 \cup R_4  \cup R_6 $ by Lemma \ref {3.1.1} {\rm (X)}.
 The sub-{\rm GDD} by deleting  Vertex  1    is  an arithmetic  {\rm  GDD}     when $ r' =r$,  $r\in R_3$ or
$ r' =r^{-1}$ or $ r' =-r$ with  $r\in R_3\cup  R_6$,  or $ r' =-1$,   $r\in R_2 \cup R_3 \cup R_4 \cup R_6 $ or $ r =-1, $ $r'\in R_2 \cup R_3 \cup R_4 \cup R_6 $ by Lemma \ref {3.1.1} {\rm (X)}.
It  is  quasi-affine. \\

  {\ }\ \ \ \ \ \ \ \ \ \ \ \ $\begin{picture}(100,       15) \put(-68,        -1){ (i) }

\put(60,       1){\makebox(0,       0)[t]{$\bullet$}}
\put(58,       -12){$-1$}

\put(37,       -12){$-q^{-1}$}
\put(28,       -1){\line(1,       0){33}}
\put(27,       1){\makebox(0,      0)[t]{$\bullet$}}

\put(22,      -12){$-1$}
\put(0,       -12){$-q$}

\put(-14,       1){\makebox(0,      0)[t]{$\bullet$}}

\put(-14,      -1){\line(1,       0){50}}

\put(-18,      -12){$q$}

\put(59,       0){\line(-1,      1){17}}

\put(28,      -1){\line(1,       1){17}}

\put(43,     18){\makebox(0,      0)[t]{$\bullet$}}

\put(10,     12){$-s^{-1}$}

\put(36,     22){$-1$}
\put(60,     12){$-r^{-1}$}

 \ \ \ \ \ \ \ \ \ \ \ \ \ \ \ \ \ \ \ \ \ \ \ \  \ \ \ \ \ \ \ \
  \ \ \ \ \ \ \ \ \ \ \ \ \ \ \ \ \ \ \ {$, q \in R_3$,  by    {\rm  GDD}  $2$ of Row $14$.}
\put(80,         -1)  {    } \end{picture}$\\ \\
 The sub-{\rm GDD} by deleting  Vertex  4    is not an arithmetic  {\rm  GDD}  by Lemma \ref {3.1.1} (I).

 (j)  i.e.  (9.3.9)$, q \in R_3$,  by  {\rm  GDD}  $2$ of Row $15$. The sub-{\rm GDD} by deleting  Vertex 4    is
 an arithmetic  {\rm  GDD}     when $q=s$  or $s = q^{-1}$,  or $s = -q$ ($q\in R_3$) or  $s = -1$($q\in R_3$) by Lemma \ref {3.1.1} (X).
  It  is quasi-affine  by Lemma \ref {3.1.22}    when $q=-s$,  $q\in R_3$ or $s = q^{-1}$, $q\in R_3$.

 (k)  i.e.  (9.3.10),  $q \in R_3$,  by    {\rm  GDD}  $2$ of Row $16$.
 The sub-{\rm GDD} by deleting  Vertex  4    is  an arithmetic  {\rm  GDD}   when $ s^{} =-q^{}$ or $ -q^{} =s^{2}$ or $-q^{} =s^{3}$,   $ s^3 \not=1, $ or $ s^{} =-1$  by Lemma \ref {3.1.1} (III).  It  is quasi-affine  by Lemma \ref {3.1.22}.

  (l)  i.e.  (9.3.11) $, t \in R_3$,  by    {\rm  GDD}  $1$ of Row $17$.
   The sub-{\rm GDD} by deleting  Vertex  4    is  an arithmetic  {\rm  GDD}    when $t= s$,  or $t= s^{-1}$, or $t= -s$,  or  $ s ={-1}$   by Lemma \ref {3.1.1} (IX).  It  is quasi-affine   by Lemma \ref {3.1.22}.

  (m)  i.e.  (9.3.12) $, q \in R_3$,  by    {\rm  GDD}  $1$ of Row $17$. The sub-{\rm GDD} by deleting  Vertex  4
is an arithmetic  {\rm  GDD}   when   $s \in  R_2 \cup  R_3 \cup R_4 \cup R_6 $  by Lemma \ref {3.1.1} (X).
 It  is quasi-affine  by Lemma \ref {3.1.22}.

 (n)  i.e.  (9.3.13) $, q \in R_3$,  by    {\rm  GDD}  $4$ of Row $17$.
 The sub-{\rm GDD} by deleting  Vertex  4    is  an arithmetic {\rm  GDD}     when $ q^{} =-s^{-1}$ or $s^{} =-1$ by Lemma \ref {3.1.1} (X).
  It  is quasi-affine  by Lemma \ref {3.1.22}.

  (o)  i.e.  (9.3.14) $, q \in R_3$,  by     {\rm  GDD}  $4$ of Row $17$.
 The sub-{\rm GDD} by deleting  Vertex  4    is  an arithmetic  {\rm  GDD}    when $ q =-s^{}$ or $ q^{} =-s^{-1}$ or $s^{} =-1^{}$   by Lemma \ref {3.1.1} (X).
  It  is quasi-affine  by Lemma \ref {3.1.22}.

  (p)  i.e.  (9.3.15) $, q \in R_3$,  by    {\rm  GDD}  $9$ of Row $17$.  The sub-{\rm GDD} by deleting  Vertex  4
 is an arithmetic  {\rm  GDD}   when $ s^{} =-q^{}$ by Lemma \ref {3.1.1} (I).
 It  is quasi-affine  by Lemma \ref {3.1.3}.

  {\rm (ii) } All quasi-arithmetic {\rm  GDD}s by  adding  a vertex on  Vertex $2$  are listed. It is the same as  {\rm (i) }.
  {\rm (iii) }  All quasi-arithmetic {\rm  GDD}s by  adding  a vertex on  Vertex $3$  are listed. It is the same  as  {\rm (i) }

  {\rm (iv) } All quasi-affine  {\rm  GDD}s which are complete diagrams are listed.
\\ \\ \\ \\ \\

{\ }\!\!\!\!\!\!\!\!\!\!\!\!\!\!\!\!\!\!\!\!\!{\ }\!\!\!\!\!\!\!\!\!\!\!\!\!\!\!\!\!\!\!\!\!{\ }\!\!\!\!\!\!\!\!\!\!\!\!
  $\begin{picture}(100,       15)\put(68,       -1){(a) }

\put(170,     10){\makebox(0,      0)[t]{$\bullet$}}

\put(170,     70){\makebox(0,      0)[t]{$\bullet$}}

\put(230,     10){\makebox(0,      0)[t]{$\bullet$}}
\put(230,     70){\makebox(0,      0)[t]{$\bullet$}}

\put(170,       10){\line(0,      1){60}}

\put(170,       10){\line(1,      1){60}}

\put(230,       10){\line(0,      1){60}}

\put(230,       10){\line(-1,      1){60}}

\put(170,       10){\line(1,       0){60}}
\put(170,       70){\line(1,       0){60}}

\put(150,     10){$-1$}
\put(150,     30){$q^{-1}$}

\put(150,     70){$q$}

\put(180,     30){$q$}
\put(190,     80){$q^{-1}$}
\put(190,     -10){$s$}
\put(220,     30){$$}

\put(250,     10){$-1$}
\put(250,     30){$r$}

\put(250,     70){$-1$}

\put(270,         -1)  {  If the sub-{\rm GDD} by deleting  Vertex 2  } \end{picture}$\\ \\
is an arithmetic  {\rm  GDD}   then considering  Lemma \ref {3.1.1}{\rm (I)} we have
$a_{1}.$ $\widetilde{q}_{34} = q^{-1} $,   $ s^{} =q^{2}$
or
$a_{2}.$ $\widetilde{q}_{34} =-1.$ $ s^{} = -q^{}$,  $q\in R_3$,  by Lemma \ref {3.1.1}{\rm (I)}.
or
$a_{3}.$ $\widetilde{q}_{34} =q^{-2}, $ $ s^{} = q^{3}$ by Lemma \ref {3.1.1}{\rm (I)}.
$a_{4}.$ $\widetilde{q}_{34} =-q.$ $ s^{} = -1^{}$,  $q\in R_3$.
$a_{5}.$ $\widetilde{q}_{34} = q^{-1},$ $ s^{} = q^{2}$, $ q^{2} \not= 1^{}$.
For $a_{3}, $ if the sub-{\rm GDD} by deleting  Vertex 1 is an arithmetic  {\rm  GDD} ,   then  $ r^{} =q^{3}$ by Lemma \ref {3.1.1}{\rm (I)}. It is not quasi-affine.
For $a_{1}, $ if the sub-{\rm GDD} by deleting  Vertex 1  is an arithmetic  {\rm  GDD} ,   then     $ r^{} =q^{2}$. It is not quasi-affine.
For
$a_{2}.$ if the sub-{\rm GDD} by deleting  Vertex 1  is an arithmetic  {\rm  GDD} ,   then     $ r^{} =-q^{}$. It is not quasi-affine.
For $a_{4}$,  $r ={-1}$,  It is not quasi-affine. For $a_{5}$,   It is not quasi-affine since  $ s^{} = q^{2}$.
\\ \\ \\ \\ \\

 \ \ \ \ \ \  $\begin{picture}(100,     15)\put(-45,      -1){ (b) }

\put(170,    10){\makebox(0,     0)[t]{$\bullet$}}

\put(170,    70){\makebox(0,     0)[t]{$\bullet$}}

\put(230,    10){\makebox(0,     0)[t]{$\bullet$}}
\put(230,    70){\makebox(0,     0)[t]{$\bullet$}}

\put(170,      10){\line(0,     1){60}}

\put(170,      10){\line(1,     1){60}}

\put(230,      10){\line(0,     1){60}}

\put(230,      10){\line(-1,     1){60}}

\put(170,      10){\line(1,      0){60}}
\put(170,      70){\line(1,      0){60}}

\put(150,    10){$-1$}
\put(150,    30){$q^{-1}$}

\put(150,    70){$q$}

\put(180,    30){$q^{3}$}
\put(190,    80){$q^{-2}$}
\put(190,    -10){$s^{3}$}
\put(220,    30){$$}

\put(250,    10){$-1$}
\put(250,    30){$r^{3}$}

\put(250,    70){$-1$}

\put(80,        -1)  {    } \end{picture}$\\ \\
with  $q^2, q^3 \not=1.$
If The sub-{\rm GDD} by deleting Vertex 1 is an arithmetic {\rm GDD},  then $\widetilde{q}_{34} =q^{-1}, $   $ q^{3} =r^{3}$ by Lemma \ref {3.1.1}{\rm (I)}. It is a contradiction.

(c)   i.e.  (9.3.16)
$a_{2}.$ If the sub-{\rm GDD} by deleting  Vertex 2 is  {\rm  GDD}  $3$ of Row $9$, then
$\widetilde{q}_{34} =s^{-1}(s')^{-1}$ and  $ s^{} \not=(s')^{-1}.$
  $a_{22}.$ If the sub-{\rm GDD} by deleting  Vertex 1  is by  Type   6, then
$ r' =r,  $  $ s' =s,$ $\widetilde{q}_{34} = r^{-2}= s^{-2}, $ $r =-s.$ It is not quasi-affine.
   $a_{23}. $ If the sub-{\rm GDD} by deleting  Vertex 1  is  {\rm  GDD}  $3$ of Row $9$, then
$s^{-1}(s')^{-1}r^{}(r')^{} =1$,  i.e. $ (r')^{2} =s^{2}$,  $ (s')^{2} =r^{2}$, $ r' =\pm s$,  $ s' =\pm r$.
 $a_{231}.$
 $ r' = s,  $ $ s'=r$,  $\widetilde{q}_{34}= q. $  It  is quasi-affine.
  $a_{232}.$ $ r' = s,  $ $ s'=-r$. It is a contradiction since $r's' q^{} =s (-r)q^{}=-1$.
   $a_{233}.$ $ r' =- s,  $ $ s'=r$. It is not quasi-affine.
   $a_{234}.$ $ r' = -s,  $ $ s'=-r$.  It  is quasi-affine  when $\widetilde{q}_{34} = -q.$
$a_{3}.$  The sub-{\rm GDD} by deleting  Vertex 2 is  Type   6.
 $a_{31}.$ $ s' = s,  $ then
$\widetilde{q}_{34} =s^{-2}$. Consequently  $ r =r'.$ If the sub-{\rm GDD} by deleting  Vertex 1  is an arithmetic  {\rm  GDD} ,   then  $ s^{-2} =r^{-2}$,   $ s =-r.$ It is  quasi-affine.
   $a_{32}.$ $\widetilde{q}_{34} =s^{}$. $s' =s^{-2}. $
     $a_{322}.$   The sub-{\rm GDD} by deleting  Vertex 1 is  {\rm  GDD}  $3$ of Row $9$,  $ r' =q^{}$.  It is not quasi-affine.
      $a_{323}.$ If  the sub-{\rm GDD} by deleting  Vertex 1  is  Type   6,  then $ r' =r^{}$,  $ s =r^{-2}$. It is a contradiction by Lemma \ref {3.1.1}{\rm (XVII)}.
  $a_{33}.$ $\widetilde{q}_{34} =s'$, $s =(s')^{-2}. $
  $a_{332}.$  The sub-{\rm GDD} by deleting  Vertex 1  is  {\rm  GDD}  $3$ of Row $9$,  then     $ s'= (r')^{-1}r^{-1}$,  i.e.
  $ r= (s')^{-1} (r')^{-1} = q$,  which is a contradiction.
  $a_{333}.$ If  the sub-{\rm GDD} by deleting  Vertex 1 is  Type   6.
   $a_{3331}.$ $ r =r'$,  $ s' =r^{-2}$,  $ rs =r's',$ $ rs =r^{-2}r. $  It is a contradiction by Lemma \ref {3.1.1}{\rm (XVII)}.
    $a_{3332}.$ If $ r =s'$,  $ r' =r^{-2}$, then it is a contradiction by Lemma \ref {3.1.1}{\rm (XVII)}.
  $a_{3333}.$ $ r' =s'$,  $ r =(r')^{-2}$,$ s=(r')^{-2}=r$. It is not quasi-affine.

 It  is quasi-affine  when    $\widetilde{q}_{34} =r^{-1}s^{-1}$, $s =r',$ $s' =r;$ or
  $\widetilde{q}_{34} =-r^{-1}s^{-1}$. $r' =-s, $ $s' =-r;$  or  $s' =s, $ $r' =r,$ $\widetilde{q}_{34} =s^{-2}$,
 $s=-r.$

\subsection* {Quasi-affine over
 {\rm  GDD}  $1$ of Row $10$ }

  {\rm (i) }  All quasi-arithmetic {\rm  GDD}s by  adding a vertex  on  Vertex $1$  are listed.
 According  to Lemma \ref {3.1.2} {\rm  (III)}  we have to consider following cases.

  {\rm (ii) }  All quasi-arithmetic {\rm  GDD}s by  adding a vertex  on  Vertex $2$  are listed.
 According  to Lemma \ref {3.1.2}  {\rm  (I)}  we have to consider following cases.

 The sub-{\rm GDD} by deleting  Vertex  1   is an arithmetic  {\rm  GDD}   in this case by symmetry.

(a)   i.e. (10.1.1) $,q \in F^{*}\setminus \{1, -1\}$ $,q \notin R_3$.  GDD $1$ of Row $7$.
It is quasi-affine  $q\notin   R_2\cup   R_3$ by Lemma \ref {3.1.8}.

 (b)  i.e.  (10.1.3)
 $q^2\notin R_2\cup R_3$,  by   {\rm  GDD}  $2$ of Row $6$. It  is quasi-affine  by Lemma \ref {6.2.2''}.

(c)  i.e.  (10.1.4), $q^3\notin R_2\cup R_3$,  by  {\rm  GDD}  $2$ of Row $7$. It  is quasi-affine  by Lemma \ref {6.2.2''}.

(d)  i.e.  (10.1.5)  $q^3\notin R_2\cup R_3$,  by    {\rm  GDD}  $4$ of Row $7$.
 It  is quasi-affine    $q\notin R_2\cup R_3$ by Lemma \ref {3.1.3}. \\ \\

 {\ }\ \ \ \ \ \
 $\begin{picture}(100,       15) \put(-45,        -1){(e) }

\put(60,       1){\makebox(0,       0)[t]{$\bullet$}}

\put(28,       -1){\line(1,       0){33}}
\put(27,       1){\makebox(0,      0)[t]{$\bullet$}}
\put(-14,       1){\makebox(0,      0)[t]{$\bullet$}}

\put(-14,      -1){\line(1,       0){50}}

\put(58,       -12){$q$}

\put(40,       -12){$q^{-1}$}

\put(22,      -12){$-1$}
\put(0,       -12){$q{}$}

\put(-18,      -12){$q^{-1}$}

\put(27,     38){\makebox(0,      0)[t]{$\bullet$}}

\put(27,       0){\line(0,      1){35}}

\put(30,       30){$q$}

\put(30,       20){$q^{-1}$}

\ \ \ \ \ \ \ \ \ \ \ \ \ \ \ \ \ \ \ \  {, $q\notin R_2\cup R_3$,  by   {\rm  GDD}  $1$ of Row $8$. It  is an
arithmetic  {\rm  GDD}   by  Type   5.}
\put(80,         -1)  {    } \end{picture}$\\

 (f)  i.e.  (10.1.6), $q\notin R_2\cup R_3$,  by    {\rm  GDD}  $1$ of Row $9$. It  is quasi-affine  by Lemma \ref {6.2.2''}.

 (g)  i.e.  (10.1.2)
$, q \notin R_2 \cup R_3$,  by    {\rm  GDD}  $1$ of Row $10$. It  is quasi-affine  by Lemma \ref {6.2.2''}.

\subsection* {Quasi-affine over
 {\rm  GDD}  $2$ of Row $10$  }
  {\rm (i) }  All    quasi-arithmetic {\rm  GDD}s by  adding a vertex  on  Vertex $1$  are listed.
 According  to Lemma \ref {3.1.2} (III) we have to consider following cases.

 (a)  i.e.  (10.2.1)$, q \in F^{*}\setminus \{1,  -1\}$ $, q \notin R_3$,  by    {\rm  GDD}  $1$ of Row $7$
is quasi-affine
 by Lemma \ref {3.1.8}.

  {\rm (ii) }  All quasi-arithmetic {\rm  GDD}s by  adding a vertex  on  Vertex $2$  are listed.
 According  to Lemma \ref {3.1.2} (I) we have to consider following cases.

 (a)  i.e.  (10.2.2) by    {\rm  GDD}  $2$ of Row $4$. The sub-{\rm GDD} by deleting  Vertex  4
is an arithmetic  {\rm  GDD}    when  $q \in  R_4 \cup R_5  $
by Lemma \ref {3.1.1} {\rm (IX)}.   It  is quasi-affine  by Lemma \ref {3.1.16}.

(b)  i.e.  (10.2.3)
 $q^2\notin R_2\cup R_3$,  by   {\rm  GDD}  $3$ of Row $5$.
 The sub-{\rm GDD} by deleting  Vertex  4    is  an arithmetic  {\rm  GDD}   when $q\in R_6$ by  Type   3;  when $q\in R_7$  by Lemma \ref {3.1.1}{\rm (I)}.   It  is quasi-affine  when $q\in   R_7$  by Lemma \ref {3.1.3}.

(c)  i.e.  (10.2.6)  by    {\rm  GDD}  $2$ of Row $6$. The sub-{\rm GDD} by deleting  Vertex  4 is  an arithmetic
{\rm  GDD}   when   $q \in  R_5 \cup R_6 \cup R_7$  by Lemma \ref {3.1.1} {\rm (IX)}.  It  is quasi-affine   by Lemma \ref {3.1.16}.

(d)  i.e.  (10.2.7)
$q^3\notin R_2\cup R_3$,  by   {\rm  GDD}  $2$ of Row $7$. The sub-{\rm GDD} by deleting  Vertex  4
is an arithmetic  {\rm  GDD}    when  $q \in  R_7 \cup R_8  $ by Lemma \ref {3.1.1} {\rm (IX)}.  It  is quasi-affine  by Lemma \ref {3.1.16}.

(e)  i.e.  (10.2.8), $q^3\notin R_2\cup R_3$,  by   {\rm  GDD}  $4$ of Row $7$. The sub-{\rm GDD} by deleting
 Vertex  4 is an arithmetic  {\rm  GDD}    when  $q \in  R_8 $ by  Type   3.
  The sub-{\rm GDD} by deleting  Vertex  4    is an arithmetic  {\rm  GDD}    when  $ q^{-6} =q^{3}$ by Lemma \ref {3.1.1}{\rm (I)}.
 It  is quasi-affine  by Lemma \ref {3.1.3}.

(f)  i.e.  (10.2.9)  by   {\rm  GDD}  $1$ of Row $8$. The sub-{\rm GDD} by
deleting  Vertex  4    is  an arithmetic  {\rm  GDD}   by Lemma \ref {3.1.1} {\rm (IV)}. It  is quasi-affine  by Lemma \ref {3.1.9}.

(g)  i.e.  (10.2.5)  by    {\rm  GDD}  $1$ of Row $9$. The sub-{\rm GDD} by
deleting  Vertex  4    is  an arithmetic  {\rm  GDD}   by Lemma \ref {3.1.1} {\rm (IV)}. It  is quasi-affine  by Lemma \ref {3.1.9}.

(h)  i.e.  (10.2.4)  $, q \notin R_3 \cup R_2,  $ $q^2 \notin R_3 \cup R_2$,  by   {\rm  GDD}  $2$ of Row $10$.
 The sub-{\rm GDD} by deleting  Vertex  4    is  an arithmetic  {\rm  GDD}    by Lemma \ref {3.1.1} {\rm (IV)}.
 It  is quasi-affine  by Lemma \ref {3.1.35}.

  {\rm (iii) }  All quasi-arithmetic {\rm  GDD}s by  adding a vertex  on  Vertex $3$  are listed.
 According  to Lemma \ref {3.1.2} (III) we have to consider following cases.

(a)    i.e.  (10.2.11),
$q^{6} \not=1$,  by   {\rm  GDD}  $1$ of Row $7$.
It is quasi-affine by by  \cite [Lemma 2.11]{TZ22}.

\subsection* {Quasi-affine over
 {\rm  GDD}  $3$ of Row $10$   }
  {\rm (i) }  All   quasi-arithmetic {\rm  GDD}s by  adding a vertex  on  Vertex $2$  are listed.
According  to Lemma \ref {3.1.2} (IV) we have to consider following cases.

(a)  i.e.  (10.3.1)
  $, q \in F^{*}\setminus \{1,  -1\}$,  by    {\rm  GDD}  $2$ of Row $4$. The sub-{\rm GDD}
by deleting  Vertex  4    is  an arithmetic  {\rm  GDD}   when $q \in  R_5 \cup R_7 \cup R_6 \cup R_4 $  by Lemma \ref {3.1.1} {\rm (IX)}.  It  is quasi-affine  when $q \in  R_7 \cup R_6 \cup R_4 $ by Lemma \ref {3.1.22}. It  is    {\rm  GDD}  $4$ of Row $9$.  when $q \in  R_5.  $

   (b)  i.e.  (10.3.2) $, q^2 \in F^{*}\setminus \{1,  -1\}$,  by   {\rm  GDD}  $2$ of
  Row $5$. The sub-{\rm GDD} by deleting
Vertex  4    is  an arithmetic  {\rm  GDD}   when $q \in  R_6 $  by  Type   2;   The sub-{\rm GDD} by deleting  Vertex  4    is not an arithmetic {\rm  GDD}     when $q \notin  R_6 $  by Lemma \ref {3.1.1}{\rm (I)}.
 It  is quasi-affine  by Lemma \ref {3.1.3}.

 (c)  i.e.  (10.3.3) $, q \in F^{*}\setminus \{1,  -1\}$,  by  {\rm  GDD}  $2$ of Row $6$. The sub-{\rm GDD} by
deleting  Vertex  4    is an arithmetic  {\rm  GDD}    when $q \in  R_4 \cup R_6 \cup R_8 $ by Lemma \ref {3.1.1} (IX).  It  is quasi-affinen $q \in   R_6 \cup R_8 $  by Lemma \ref {3.1.22}.  It  is  {\rm  GDD}  $4$ of Row $9$. $q \in  R_4.$

 (d)  i.e.  (10.3.4)  $, q \in F^{*}\setminus \{1,  -1\}$,  by    {\rm  GDD}  $3$ of Row $6$.
 The sub-{\rm GDD} by deleting  Vertex  2    is an arithmetic  {\rm  GDD}    when $q \in  R_5 \cup R_9 \cup R_{13}\cup R_8 $ by Lemma \ref {3.1.1} (IX).
 The sub-{\rm GDD} by deleting  Vertex  1    is an arithmetic  {\rm  GDD}    when $q \in  R_4 \cup R_5 \cup R_{3} $ by Lemma \ref {3.1.1} (X). It  is  quasi-affine  when $q \in  R_5 $.

(e)  i.e.  (10.3.5)$, q \in F^{*}\setminus \{1,  -1\}$. $ q \notin R_3$,  by   {\rm  GDD}  $2$ of Row $7$.
 The sub-{\rm GDD} by deleting  Vertex  4    is an arithmetic  {\rm  GDD}    when $q \in  R_5 \cup R_7 \cup R_9 \cup R_6 \cup R_4 $ by Lemma \ref {3.1.1} {\rm (IX)} and   Lemma \ref {3.1.2} {\rm (V)}.
 It  is quasi-affine  by Lemma \ref {3.1.22}.\\  \\ \\

  {\ }\\

{\ }\!\!\!\!\!\!\!\!\!\!\!\!\!\!\!\!\!\!\!\!\!{\ }\!\!\!\!\!\!\!\!\!\!\!\!\!\!\!\!\!\!\!\!\!{\ }\!\!\!\!\!\!\!\!\!\!\!\!\!\!\! \ \ \ \  $\begin{picture}(100,       15)\put(68,       -1){(f) }

\put(170,     10){\makebox(0,      0)[t]{$\bullet$}}

\put(170,     70){\makebox(0,      0)[t]{$\bullet$}}

\put(230,     10){\makebox(0,      0)[t]{$\bullet$}}
\put(230,     70){\makebox(0,      0)[t]{$\bullet$}}

\put(170,       10){\line(0,      1){60}}

\put(170,       10){\line(1,      1){60}}

\put(230,       10){\line(0,      1){60}}


\put(170,       10){\line(1,       0){60}}
\put(170,       70){\line(1,       0){60}}

\put(150,     10){$-1$}
\put(150,     30){$q^{-1}$}

\put(150,     70){$q$}

\put(180,     30){$q^{3}$}
\put(190,     80){$q^{-2}$}
\put(190,     -10){$q^{-6}$}
\put(220,     30){$$}

\put(250,     10){$-1$}
\put(250,     30){$q^{3}$}

\put(250,     70){$-1$}

\put(280,         -1)  { $ q \notin R_3$,  $ q ^2\notin R_3\cup   R_2$,  by    {\rm  GDD}  $3$ of Row $7$.   } \end{picture}$\\ \\
 The sub-{\rm GDD} by deleting  Vertex  1    is not an arithmetic  {\rm  GDD}   by Lemma \ref {3.1.1}{\rm (I)}. \\

  {\ }\ \ \ \ \ \ \ \ \ \ \ \ $\begin{picture}(100,       15) \put(-68,        -1){ (g)}

\put(60,       1){\makebox(0,       0)[t]{$\bullet$}}
\put(58,       -12){$-1$}

\put(40,       -12){$q$}
\put(28,       -1){\line(1,       0){33}}
\put(27,       1){\makebox(0,      0)[t]{$\bullet$}}

\put(22,      -12){$-1$}
\put(0,       -12){$q^{-1}$}

\put(-14,       1){\makebox(0,      0)[t]{$\bullet$}}

\put(-14,      -1){\line(1,       0){50}}

\put(-18,      -12){$-1$}

 \put(59,       0){\line(-1,      1){17}}

\put(28,      -1){\line(1,       1){17}}

\put(43,     18){\makebox(0,      0)[t]{$\bullet$}}

\put(18,     12){$q^{-2}$}

\put(36,     22){$-1$}
\put(59,     12){$q$}

 \ \ \ \ \ \ \ \ \ \ \ \ \ \ \ \ \ \ \ \ \ \ \ \  \ \ \ \ \ \ \ \
 \ \ \ \ \ \ \ \ \ \ \ \ \ \ \ \ \ \ \ {$, q \in F^{*}\setminus \{1,  -1\}$,  by   {\rm  GDD}  $2$ of Row $8$.}
\put(80,         -1)  {    } \end{picture}$\\ \\
 It is     {\rm  GDD}  $5$ of Row $22$.

(h)  i.e.  (10.3.10)
$, q,  r,  s \in F^{*}\setminus \{1\},  qrs =1,  q\not= r,  q\not= s,  s\not= r$,
  by   {\rm  GDD}  $3$ of Row $9$.  The sub-{\rm GDD} by deleting  Vertex  2    is an arithmetic  {\rm  GDD}   when  $ s =q^{2}$ or $ s =q^{-2}$,  $q\in R_3 \cup R_6$ or
 $ s =-q^{-2}$,  $s \in  R_3 \cup R_6  $ or $ s =-1, $ $q^{-2} \in  R_2 \cup R_3 \cup R_4 \cup R_6  $
  or $ q^{-2} =-1, $ $s^{} \in  R_2 \cup R_3 \cup R_4 \cup R_6  $ by Lemma \ref {3.1.1} {\rm (X)}.
 The sub-{\rm GDD} by deleting  Vertex  1    is an arithmetic  {\rm  GDD}
 when $ r^{} =-q^{}$, $q\in R_6$,   or $ r^{} ={-1}$,  $q \in  R_4\cup R_6$
   by Lemma \ref {3.1.1} {\rm (X)}. It  is  quasi-affine.
\\

  {\ }\ \ \ \ \ \ \ \ \ \ \ \ $\begin{picture}(100,       15) \put(-68,        -1){ (i)}

\put(60,       1){\makebox(0,       0)[t]{$\bullet$}}
\put(58,       -12){$-1$}

\put(37,       -12){$-q^{-1}$}
\put(28,       -1){\line(1,       0){33}}
\put(27,       1){\makebox(0,      0)[t]{$\bullet$}}

\put(22,      -12){$-1$}
\put(0,       -12){$-q$}

\put(-14,       1){\makebox(0,      0)[t]{$\bullet$}}

\put(-14,      -1){\line(1,       0){50}}

\put(-18,      -12){$q$}

 \put(59,       0){\line(-1,      1){17}}

\put(28,      -1){\line(1,       1){17}}

\put(43,     18){\makebox(0,      0)[t]{$\bullet$}}

\put(18,     12){$q^{2}$}

\put(36,     22){$-1$}
\put(59,     12){$-q^{-1}$}

 \ \ \ \ \ \ \ \ \ \ \ \ \ \ \ \ \ \ \ \ \ \ \ \  \ \ \ \ \ \ \ \
\ \ \ \ \ \ \ \ \ \ \ \ \ \ \ \ \ \ \ {$, q \in R_3$,  by    {\rm  GDD}  $2$ of Row $14$.}
\put(80,         -1)  {    } \end{picture}$\\ \\
 The sub-{\rm GDD} by deleting  Vertex  4    is not  an arithmetic  {\rm  GDD}   by Lemma \ref {3.1.1}{\rm (I)}.

(j)  i.e.  (10.3.6)
$,q \in R_3$,  by  GDD $4$ of Row $17$.
The sub-{\rm GDD} by deleting  Vertex  4  is an arithmetic GDD   by Lemma \ref {3.1.1} {\rm (IX)}.
It is quasi-affine  by Lemma \ref {3.1.22}.

  {\rm (ii) }  All quasi-arithmetic {\rm  GDD}s by  adding a vertex  on  Vertex $1$  are listed.
 According  to Lemma \ref {3.1.2} (IV) we have to consider following cases.

(a)  i.e.  (10.3.7), $q^2, q^3 \not=1,$  by   {\rm  GDD}  $2$ of Row $7$. It  is quasi-affine  by Lemma \ref {3.1.22}.

(b)  i.e.  (10.3.11),  $q^6,  q^9 \not=1, $  by    {\rm  GDD}  $3$ of Row $7$. The sub-{\rm GDD}
by deleting  Vertex  1    is an arithmetic  {\rm  GDD}    when $q\in R_5 \cup R_8  $ by Lemma \ref {3.1.1} {\rm (IX)}.
 The sub-{\rm GDD} by deleting  Vertex  2    is an arithmetic  {\rm  GDD}    when $q\in R_8 \cup R_4 $ by  Type   2. It  is
quasi-affine when $q\in R_8$.

(c)  i.e.  (10.3.12)$, q,  r,  s \in F^{*}\setminus \{1\},  qrs =1,  q\not= r,  q\not= s,  s\not= r$,
 by   {\rm  GDD}  $3$ of Row $9$.  The sub-{\rm GDD} by deleting  Vertex  2    is an arithmetic  {\rm  GDD}    when $ s =-q^{}$,  $q\in R_6$ or  $ -1 =s$,  $q \in  R_4 \cup R_6 $
by Lemma \ref {3.1.1} {\rm (X)}.
 The sub-{\rm GDD} by deleting  Vertex  1    is an arithmetic  {\rm  GDD}    when   $ r =q^{2}$ or  $ r =-q^{-2}$
$q\in R_6\cup R_{12}$ or $ -1 =q^{-2}$,  $r \in  R_2  \cup R_3 \cup R_4 \cup R_6 $  or $ -1 =r^{}$,  $ q^{-2}\in  R_2 \cup R_3 \cup R_4  \cup R_6 $. It  is quasi-affine when $(q,r,s) \in A\cap B\cap C$, where $A:= \{ (q,r,s) \mid  $  $ s =q^{2}$ or $ s =q^{-2}$, $q\in R_6$ or
 $ s =-q^{-2}$, $q \in  R_{12} \cup  R_6  $ or $ s =-1,$ $q^{} \in  R_4 \cup R_6 \cup R_{12}  \cup R_{8}  $
  or $ q^{-2} =-1,$ $s^{} \in  R_2 \cup R_3 \cup R_4 \cup R_6  $    $\}$,   $B:= \{ (q,r,s) \mid  $ $ r^{} =-q^{}$,$q\in R_6$,  or $ r^{} ={-1}$, $q \in  R_4\cup R_6$  $\}$, $C:= \{ (q,r,s) \mid  $  $ qrs =1,  q\not= r,  q\not= s,  s\not= r.$ $\}$..
\\ \\ \\

  {\ }\\

{\ }\!\!\!\!\!\!\!\!\!\!\!\!\!\!\!\!\!\!\!\!\!{\ }\!\!\!\!\!\!\!\!\!\!\!\!\!\!\!\!\!\!\!\!\!{\ }\!\!\!\!\!\!\!\!\!\!\!\!\!\!\!\!\!\!\!\!\!{\ }\!\!\!\
   $\begin{picture}(100,       15)\put(88,       -1){ (d)}

\put(170,     10){\makebox(0,      0)[t]{$\bullet$}}

\put(170,     70){\makebox(0,      0)[t]{$\bullet$}}

\put(230,     10){\makebox(0,      0)[t]{$\bullet$}}
\put(230,     70){\makebox(0,      0)[t]{$\bullet$}}

\put(170,       10){\line(0,      1){60}}

\put(170,       10){\line(1,      1){60}}

\put(230,       10){\line(0,      1){60}}


\put(170,       10){\line(1,       0){60}}
\put(170,       70){\line(1,       0){60}}

\put(150,     10){$-1$}
\put(150,     30){$q^{}$}

\put(150,     70){$-1$}

\put(180,     30){$q$}
\put(190,     80){$q^{-2}$}
\put(190,     -10){$q$}
\put(220,     30){$$}

\put(250,     10){$-1$}
\put(250,     30){$q^{-2}$}

\put(250,     70){$-1$}

\put(280,         -1)  { $, q^2,  q^3 \not= 1$,  by   {\rm  GDD}  $3$ of Row $10$.   } \end{picture}$\\ \\
 The sub-{\rm GDD} by deleting  Vertex  1    is not an arithmetic  {\rm  GDD}   by Lemma \ref {3.1.1} {\rm (X)}.

(e)  i.e.  (10.3.8)$, q \in R_3$, by   GDD $4$ of Row $17$.
  It is quasi-affine  by Lemma \ref {3.1.22}.

  {\rm (iii) }  All quasi-arithmetic {\rm  GDD}s by  adding a vertex  on  Vertex $3$  are listed. It is the same   as  Case {\rm (i) }.

  {\rm (iv) } All quasi-affine  {\rm  GDD}s which are complete diagrams are listed.\\ \\ \\ \\ \\ \\ \\

  \ \ \ \ \ \  $\begin{picture}(100,     15)\put(-45,      -1){(b) in  Case {\rm (ii) } }

\put(170,    10){\makebox(0,     0)[t]{$\bullet$}}

\put(170,    70){\makebox(0,     0)[t]{$\bullet$}}

\put(230,    10){\makebox(0,     0)[t]{$\bullet$}}
\put(230,    70){\makebox(0,     0)[t]{$\bullet$}}

\put(170,      10){\line(0,     1){60}}

\put(170,      10){\line(1,     1){60}}

\put(230,      10){\line(0,     1){60}}

\put(230,      10){\line(-1,     1){60}}

\put(170,      10){\line(1,      0){60}}
\put(170,      70){\line(1,      0){60}}

\put(150,    10){$-1$}
\put(150,    30){$q^{-1}$}

\put(150,    70){$q$}

\put(180,    30){$q^{3}$}
\put(190,    80){$q^{-2}$}
\put(190,    -10){$q^{3}$}
\put(220,    30){$$}

\put(250,    10){$-1$}
\put(250,    30){$q^{-6}$}

\put(250,    70){$-1$}

\put(80,        -1)  {    } \end{picture}$\\ \\   If the sub-{\rm GDD} by deleting  Vertex  1
 is an arithmetic {\rm GDD},  then  $a_1.$ $\widetilde{q}_{34} =q^{-1}$ and  $ q^{-6} =q^{3}$ by GDD $3$ of Row $7$, or $a_2.$   $\widetilde{q}_{34} =q^{-2}$, $ q^{-2} =q^{-1}$ and $ q^{-6} =q^{3}$.  It is clear that $a_2$ is a contradiction.
For  $a_1,$ if the sub-{\rm GDD} by deleting  Vertex  2  is an arithmetic {\rm GDD},  then $ q^{3} =q^{2}$ by  Type   6. It is a contradiction.

(f) in  Case {\rm (i) }  i.e.  (10.3.9)  If the sub-{\rm GDD} by deleting  Vertex 2 is  an arithmetic  {\rm  GDD}, then $a_1.$ $\widetilde{q}_{34} =q^{-1}$ and $ q^{-6} =q^{2}$  by  Type   6
or $a_2.$ $ \widetilde{q}_{34} =q^{-2}$ and $ q^{-6} =q^{3}$ by Lemma \ref {3.1.1}{\rm (I)}.
For $a_1, $ The sub-{\rm GDD} by deleting  Vertex 1  is an arithmetic  {\rm  GDD}   by Lemma \ref {3.1.1}{\rm (I)},  For $a_2, $ the sub-{\rm GDD} by deleting  Vertex 1  is not an arithmetic  {\rm  GDD}   by Lemma \ref {3.1.1}{\rm (I)}.  It  is quasi-affine  when $\widetilde{q}_{34} =q^{-1}$ and $ q^{8} =1$.

 (h) in  Case {\rm (i) } i.e.  (10.3.13)
$a_1$,   Assume  the sub-{\rm GDD} by deleting Vertex 2 is  Type   6.
$a_{11}$.  $ \widetilde{q}_{34} =q^{-2}$ and $ s =q^{4}$.
$a_{111}$. Assume   the sub-{\rm GDD} by deleting Vertex 1 is  Type   6.  Then $ r^{} =q^{-2}$ and $ q^{} =r^{-2}$.
Consequently $ q^{} =q^{4}$. It is not quasi-affine.
$a_{112}$. If the sub-{\rm GDD} by deleting  Vertex 1  is    {\rm  GDD}  $3$ of Row $9$, then $ r^{} =q^{}$. It is not quasi-affine.
$a_{12}$. $\widetilde{q}_{34} =s,$ $ q^{-2} =s^{-2}$, $ q^{} =-s^{}$. The sub-{\rm GDD} by deleting  Vertex 1  is an arithmetic
by   {\rm  GDD}  $3$ of Row $9.$
$a_{2}$. The sub-{\rm GDD} by deleting  Vertex 2 is    {\rm  GDD}  $3$ of Row $9$.
 $ \widetilde{q}_{34}  =q^{2}s^{-1}$.
$a_{21}$. If the sub-{\rm GDD} by deleting  Vertex 1  is   {\rm  GDD}  $3$ of Row $9$,  then $ q^{2} s^{-1} =r^{-1} q^{-1}$,  i, e,
$ q^{3}  =s^{}r^{-1}$.
$a_{22}$.  If the sub-{\rm GDD} by deleting  Vertex 1  is   Type   6,
then $ q^{2} s^{-1} =r^{} $ and
$ q^{}  =r^{-2}$.

 It  is quasi-affine  when $ \widetilde{q}_{34}  =s^{}$, $ q^{} =-s;$ or  $ \widetilde{q}_{34}  =
q^{2}s^{-1}$, $ q^{3} =sr^{-1};$ or $ \widetilde{q}_{34}  =q^{2}s^{-1}$, $ q^{2} s^{-1}=r,$ $ q^{} =r^{-2}$ with    $ q,  r,  s \in F^{*}\setminus \{1\},  qrs =1,  q\not= r,  q\not= s,  s\not= r.$

(b) in  Case {\rm (ii) } i.e. (10.3.11).  If the sub-{\rm GDD} by deleting  Vertex 2 is an
arithmetic  {\rm  GDD}, then  $a_1.$ $\widetilde{q}_{34} =q^{-1}$ and  $ q^{-6} =q^{3}$ by  {\rm  GDD}  $3$ of Row $7$,  or $a_2.$   $\widetilde{q}_{34} =q^{-2}$,  $ q^{-2} =q^{-1}$ and $ q^{-6} =q^{3}$. It is not quasi-affine.
For  $a_1, $ if the sub-{\rm GDD} by deleting  Vertex 1  is an arithmetic  {\rm  GDD} ,   then $ q^{3} =q^{2}$ by  Type   6. It is not quasi-affine.

 (c) in  Case {\rm (ii) } i.e.  (10.3.14)
$a_{1}.$ If the sub-{\rm GDD} by deleting  Vertex 1  is  {\rm  GDD}  $3$ of Row $9$,  then $\widetilde{q}_{34} =r.$
$a_{11}.$ The sub-{\rm GDD} by deleting  Vertex 2 is  Type   6. then  $ q^{-2} =r^{-2}$,  i.e. $ q^{} =-r^{}$.
$a_{2}.$ Assume the sub-{\rm GDD} by deleting  Vertex 1   is  Type   6. $a_{21}.$ $\widetilde{q}_{34} =s^{}$ and $ q^{} =s^{-2}$.
$a_{211}.$ If the sub-{\rm GDD} by deleting  Vertex 2  is  Type   6,  then
 $ s^{} =q^{-2}$ and $ r^{} =s^{-2}$. It is not quasi-affine.
$a_{212}. $ If the sub-{\rm GDD} by deleting  Vertex 2 is  {\rm  GDD}  $3$ of Row $9$,  then
 $ rs^{} q^{-2}=1$,  i.e.  $ sr^{} =q^{2}$.  $ sr =s^{-4}$. $ r^{} =s^{-5}$.
$a_{22}.$ $\widetilde{q}_{34} =q^{}$,   $ q^{-2} =s^{}$.
$a_{221}.$ If the sub-{\rm GDD} by deleting  Vertex 2  is  Type   6,  then it is not quasi-affine.
$a_{222}.$ If the sub-{\rm GDD} by deleting  Vertex 2  is  {\rm  GDD}  $3$ of Row $9$,   then $ rq^{-1} =1.$ It is not quasi-affine.
 It  is quasi-affine  when  $\widetilde{q}_{34} =r^{}$ and $ q^{} =-r$,  or  $\widetilde{q}_{34} =s^{}$ and $ q^{} =s^{-2}$. $rs q^{-2} =1.$
\\ \\ \\ \\ \\ \\ \\ \\

{\ }\!\!\!\!\!\!\!\!\!\!\!\!\!\!\!\!\!\!\!\!\!{\ }\!\!\!\!\!\!\!\!\!\!\!\!\!\!\!\!\!\!\!\!\!{\ }\!\!\!\!\!\!\!\!\!\!\!\!\!\!\!
  $\begin{picture}(100,       15)\put(68,       -1){
(d) in  Case {\rm (ii) } }

\put(170,     10){\makebox(0,      0)[t]{$\bullet$}}

\put(170,     70){\makebox(0,      0)[t]{$\bullet$}}

\put(230,     10){\makebox(0,      0)[t]{$\bullet$}}
\put(230,     70){\makebox(0,      0)[t]{$\bullet$}}

\put(170,       10){\line(0,      1){60}}

\put(170,       10){\line(1,      1){60}}

\put(230,       10){\line(0,      1){60}}

\put(230,       10){\line(-1,      1){60}}

\put(170,       10){\line(1,       0){60}}
\put(170,       70){\line(1,       0){60}}

\put(150,     10){$-1$}
\put(150,     30){$q^{}$}

\put(150,     70){$-1$}

\put(180,     30){$q$}
\put(190,     80){$q^{-2}$}
\put(190,     -10){$q$}
\put(220,     30){$$}

\put(250,     10){$-1$}
\put(250,     30){$q^{-2}$}

\put(250,     70){$-1$}

\put(80,         -1)  {    } \end{picture}$\\ \\
If the sub-{\rm GDD} by deleting  Vertex 1  is  Type   6,  then
 $\widetilde{q}_{34} =q^{-2}$.
 The sub-{\rm GDD} by deleting  Vertex 2 is not an arithmetic  {\rm  GDD}    checked  step by step.

\subsection* {Quasi-affine over
  {\rm  GDD}  $1$ of Row $11$ }
  {\rm (i) }  All quasi-affine  {\rm  GDD}s which are complete diagrams are listed.
 According  to Lemma \ref {3.1.2} (III) we have to consider following cases.

 (a)  i.e.  (11.1.1), by {\rm  GDD} $2$ of Row $7$,  It  is quasi-affine  by Lemma \ref {6.2.2''} with
      $q\in R_9$.

 (b)  i.e.  (11.1.2), by  GDD $4$ of Row $7$.  It  is quasi-affine  by Lemma \ref {3.2.36} with
      $q\in R_9$.

 (c)  i.e.  (11.1.3), by   {\rm  GDD}  $1$ of Row $9$.
 It  is quasi-affine  by Lemma \ref {6.2.2''} with
$q\in R_3$,  $ r \not=q^{-1}, q$.

\subsection* {Quasi-affine over
  {\rm  GDD}  $2$ of Row $11$ }
   {\rm (i) } All quasi-arithmetic {\rm  GDD}s by  adding a vertex  on  Vertex $1$  are listed.
  According  to Lemma \ref {3.1.2} (IV) we have to consider following cases. \\ \\ \\ \\ \\

{\ }\!\!\!\!\!\!\!\!\!\!\!\!\!\!\!\!\!\!\!\!\!{\ }\!\!\!\!\!\!\!\!\!\!\!\!\!\!\!\!\!\!\!\!\!{\ }\!\!\!\!\!\!\!\!\!\!\!\!\!\!\!\!\!\!\!\
{\ }\!\!\!\!\!\!\!\!\!\!\!\!\!\!\!\!\!\!\!\!\!\!\!\!\!
 $\begin{picture}(100,       15)\put(120,       -1){(a) }

\put(170,     10){\makebox(0,      0)[t]{$\bullet$}}

\put(170,     70){\makebox(0,      0)[t]{$\bullet$}}

\put(230,     10){\makebox(0,      0)[t]{$\bullet$}}
\put(230,     70){\makebox(0,      0)[t]{$\bullet$}}

\put(170,       10){\line(0,      1){60}}

\put(170,       10){\line(1,      1){60}}

\put(230,       10){\line(0,      1){60}}


\put(170,       10){\line(1,       0){60}}
\put(170,       70){\line(1,       0){60}}

\put(150,     10){$-1$}
\put(150,     30){$s$}

\put(150,     70){$-1$}

\put(180,     30){$q$}
\put(190,     80){$r$}
\put(190,     -10){$q$}
\put(220,     30){$$}

\put(250,     10){$-1$}
\put(250,     30){$q$}

\put(250,     70){$-1$}

    \put(280,     -1){$, q,  r,  s \in F^{*}\setminus \{1\},  qrs =1,  q\not= r,  q\not= s,  s\not= r $,}
\put(80,         -1)  {    } \end{picture}$\\ \\
 by   {\rm  GDD}  $3$ of Row $9$. The sub-{\rm GDD} by deleting  Vertex  1 is an arithmetic  {\rm  GDD}   when $-q=r$ or $r=-1$ by Lemma \ref {3.1.1} {\rm (X)}.
 The sub-{\rm GDD} by deleting  Vertex  2 is an arithmetic  {\rm  GDD}   when $-q=s$ or $s=-1$ by Lemma \ref {3.1.1} {\rm (X)}.
It  is a contradiction. \\ \\ \\ \\ \\

{\ }\!\!\!\!\!\!\!\!\!\!\!\!\!\!\!\!\!\!\!\!\!\!\!\!\!\!\!\!\!\!\!\!\!\!\!\!\!\!\!\!\!\!
  \ \ \ \ \ \  $\begin{picture}(100,     15)\put(45,      -1){ (b)}

\put(170,    10){\makebox(0,     0)[t]{$\bullet$}}

\put(170,    70){\makebox(0,     0)[t]{$\bullet$}}

\put(230,    10){\makebox(0,     0)[t]{$\bullet$}}
\put(230,    70){\makebox(0,     0)[t]{$\bullet$}}

\put(170,      10){\line(0,     1){60}}

\put(170,      10){\line(1,     1){60}}

\put(230,      10){\line(0,     1){60}}


\put(170,      10){\line(1,      0){60}}
\put(170,      70){\line(1,      0){60}}

\put(150,    10){$-1$}
\put(150,    30){$q^{-1}$}

\put(150,    70){$q$}

\put(180,    30){$q^{3}$}
\put(190,    80){$q^{-2}$}
\put(190,    -10){$q^{3}$}
\put(220,    30){$$}

\put(250,    10){$-1$}
\put(250,    30){$q^{3}$}

\put(250,    70){$-1$}

 \ \ \ \ \ \  \ \ \ \ \ \ \ \ \ \ \ \ \ \  \ \ \ \ \ \ \ \  \ \ \  \ \ \ \ \ \ \ \  \ \ \  \ \ \ \ \ \ \ \
 \ \ \ \ \ \ \ \ \ \ \ \ \ \ \ \ \ \ \ \ \ \ \ \ { $,q \in F^{*}\setminus \{1, -1\}$. $ q \notin R_3$,    GDD $3$ of Row $7$.}
\put(80,        -1)  {    } \end{picture}$\\ \\
 The sub-{\rm GDD} by deleting  Vertex  1 is not an arithmetic  {\rm  GDD}   when $q \in  R_{9}$ by Lemma \ref {3.1.1} {\rm (I)}.

(c) i.e.  (11.2.1), $q\in R_3,$ by {\rm  GDD}  $2$ of Row $15.$ It is affine by
   Lemma \ref {3.1.22}.

(d) i.e.  (11.2.2)
,   $q\in R_3,$ by {\rm  GDD}  $2$ of Row $16.$ It is affine by
   Lemma \ref {3.1.22}.

(e) i.e.  (11.2.3),
  $q\in R_3,$ by {\rm  GDD}  $1$ of Row $17.$ It is affine by
   Lemma \ref {3.1.22}.

(f) i.e.  (11.2.4)
, $q\in R_3,$  by {\rm  GDD}  $4$ of Row $17.$ It is affine by
   Lemma \ref {3.1.22}.

\subsection* {Quasi-affine over
  {\rm  GDD}  $1$ of Row $13$ }
  {\rm (i) }  All quasi-arithmetic {\rm  GDD}s by  adding a vertex  on  Vertex $1$  are listed.
 According  to Lemma \ref {3.1.2} (II) we have to consider following cases.

 (a)  i.e.  (13.1.1)$, q \in R_3 \cup R_6$,  by  {\rm  GDD}  $1$ of Row $13$. It  is
quasi-affine when  $q\in R_3$ by Lemma \ref {3.1.25}.
 It  is quasi-affine  when  $q\in R_6$ by Lemma \ref {3.2.38}.

   {\rm (ii) }  All    quasi-arithmetic {\rm  GDD}s by  adding a vertex  on  Vertex $2$  are listed.
 According  to Lemma \ref {3.1.2} (II) we have to consider following cases.

  (a)  i.e.  (13.1.2) $, q \in F^{*}\setminus \{1\}$,  by   {\rm  GDD}  $1$ of Row $1$.  The sub-{\rm GDD} by deleting
Vertex  4  is  {\rm  GDD}  $1$ of Row $13$. It  is quasi-affine  by Lemma \ref {3.1.3}.\\ \\

 {\ }\ \ \ \ \ \ \ \ \ \ \ \ $\begin{picture}(100,       15) \put(-68,        -1){(b) }

\put(60,       1){\makebox(0,       0)[t]{$\bullet$}}
\put(58,       -12){$q^2$}

\put(40,       -12){$q^{-2}$}
\put(28,       -1){\line(1,       0){33}}
\put(27,       1){\makebox(0,      0)[t]{$\bullet$}}

\put(22,      -12){$q^2$}
\put(0,       -12){$q^{-2}$}

\put(-14,       1){\makebox(0,      0)[t]{$\bullet$}}

\put(-14,      -1){\line(1,       0){50}}

\put(-18,      -12){$q$}

\put(27,     38){\makebox(0,      0)[t]{$\bullet$}}

\put(27,       0){\line(0,      1){35}}

\put(30,       30){$-1$}

\put(30,       20){$q^{-4}$}

 \ \ \  \ \ \ \ \ \ \ \ \ \ \ \ \ \ \ \ \ \ \ {     $q^2\in R_3\cup R_6$, $q\in R_3\cup R_6$,   $ q^2 \in F^{*}\setminus \{1,  -1\}$,  by   {\rm  GDD}  $1$ of Row $2$.}
\put(80,         -1)  {    } \end{picture}$\\ \\
 The sub-{\rm GDD} by deleting  Vertex  4    is not an arithmetic  {\rm  GDD}   by Lemma \ref {3.1.1}{\rm (I)}.\\
\\

  {\ }\ \ \ \ \ \ \ \ \ \ \ \ $\begin{picture}(100,       15) \put(-68,        -1){ (c)}

\put(60,       1){\makebox(0,       0)[t]{$\bullet$}}
\put(58,       -12){$q$}

\put(40,       -12){$q^{-1}$}
\put(28,       -1){\line(1,       0){33}}
\put(27,       1){\makebox(0,      0)[t]{$\bullet$}}

\put(22,      -12){$q$}
\put(0,       -12){$q^{-2}$}

\put(-14,       1){\makebox(0,      0)[t]{$\bullet$}}

\put(-14,      -1){\line(1,       0){50}}

\put(-18,      -12){$q^2$}

\put(27,     38){\makebox(0,      0)[t]{$\bullet$}}

\put(27,       0){\line(0,      1){35}}

\put(30,       30){$-1$}

\put(30,       20){$q^{-2}$}

  \ \ \ \ \ \ \ \ \ \ \ \ \ \ \ \ \ \ \ {$, q \in F^{*}\setminus \{1,  -1\}$,  by   {\rm  GDD}  $1$ of Row $3$.}
\put(80,         -1)  {    } \end{picture}$\\
\\
 The sub-{\rm GDD} by deleting  Vertex  4    is not an arithmetic  {\rm  GDD}   by Lemma \ref {3.1.1} {\rm (V)}.\\
\\

{\ }\ \ \ \ \ \ \ \ \ \ \ \ $\begin{picture}(100,       15) \put(-68,        -1){(d) }

\put(60,       1){\makebox(0,       0)[t]{$\bullet$}}

\put(28,       -1){\line(1,       0){33}}
\put(27,       1){\makebox(0,      0)[t]{$\bullet$}}

\put(-14,       1){\makebox(0,      0)[t]{$\bullet$}}

\put(-14,      -1){\line(1,       0){50}}

\put(27,     38){\makebox(0,      0)[t]{$\bullet$}}

\put(27,       0){\line(0,      1){35}}

\put(58,       -12){$-q^{}$}

\put(40,       -12){$-q^{-1}$}
\put(12,      -12){$-q^{}$}
\put(-10,       -12){$-q^{-1}$}

\put(-18,      -12){$q^{-1}$}

\put(30,       30){$-1$}

\put(30,       20){$q^{-2}$}

  \ \ \ \ \ \ \ \ \ \ \ \ \ \ \ \ \ \ \ \ \ \ \ \ \  {$, q \in R_3$,  by    {\rm  GDD}  $1$ of Row $12$.}
\put(80,         -1)  {    } \end{picture}$\\
\\
 The sub-{\rm GDD} by deleting  Vertex  4    is not an arithmetic  {\rm  GDD}   by Lemma \ref {3.1.1}{\rm (I)}.

(e) i.e.  (13.1.3)
  $,q \in F^{*}\setminus \{1, -1\}$. by GDD $1$ of Row $4$.
 The sub-{\rm GDD} by deleting  Vertex  4  is  an arithmetic GDD  when $q\in R_3$ by Lemma \ref {3.1.2} {\rm (III)}(h).  The sub-{\rm GDD} by deleting  Vertex  4  is not an arithmetic GDD  when $q\in R_6$ by Lemma \ref {3.1.2} {\rm (III)}.  It is quasi-affine  $q\in R_3$ by Lemma \ref {3.1.3}.

(f) i.e.  (13.1.4)
  $,q \in F^{*}\setminus \{1, -1\}$, by GDD $1$ of Row $4$.
 The sub-{\rm GDD} by deleting  Vertex  4    is an arithmetic  {\rm  GDD}   when $q\in R_3$ by Lemma \ref {3.1.2} {\rm (III)}.  The sub-{\rm GDD} by deleting  Vertex  4    is not an arithmetic  {\rm  GDD}   when $q\in R_6$ by Lemma \ref {3.1.2} {\rm (III)}.  It  is quasi-affine $q\in R_3$ by Lemma \ref {3.1.3}.\\ \\

  {\ }\ \ \ \ \ \ \ \ \ \ \ \ $\begin{picture}(100,       15) \put(-68,        -1){(g) }

\put(60,       1){\makebox(0,       0)[t]{$\bullet$}}
\put(58,       -12){$-1$}

\put(40,       -12){$q^{-1}$}
\put(28,       -1){\line(1,       0){33}}
\put(27,       1){\makebox(0,      0)[t]{$\bullet$}}

\put(22,      -12){$q$}
\put(0,       -12){$q^{-2}$}

\put(-14,       1){\makebox(0,      0)[t]{$\bullet$}}

\put(-14,      -1){\line(1,       0){50}}

\put(-18,      -12){$q^2$}

\put(27,     38){\makebox(0,      0)[t]{$\bullet$}}

\put(27,       0){\line(0,      1){35}}

\put(30,       30){$-q^{}$}

\put(30,      15){$-q^{-1}$}

\ \ \ \ \ \ \ \ \ \ \ \ \ \ \ \ \ \ \ \ \ \ \ \ \ \ \ \ \ \ \ \ \ \ \ \ \ \   \ \ \ \ \ \ \ \ \ \ \ \ \ \ \ \ \ \ \ \ {$, q \in R_3$,  by    {\rm  GDD}  $1$ of Row $6$.}
\put(80,         -1)  {    } \end{picture}$\\
\\
 The sub-{\rm GDD} by deleting  Vertex  4    is not an arithmetic  {\rm  GDD}    by Lemma \ref {3.1.1} {\rm (V)}.

  (h)  i.e.  (13.1.5)
 $, q \in R_3 \cup R_6$,  by   {\rm  GDD}  $1$ of Row $13$.  The sub-{\rm GDD} by deleting  Vertex  4    is an arithmetic  {\rm  GDD}   by Lemma \ref {3.1.1} {\rm (V)} with  $q\in R_6$.
 The sub-{\rm GDD} by deleting  Vertex  4    is not an arithmetic  {\rm  GDD}   by Lemma \ref {3.1.1} {\rm (V)} with  $q\in R_3$.
  It  is quasi-affine   by Lemma \ref {3.1.3}.

  {\rm (iv) }  All quasi-affine circles are listed.

  (a)
  (nc)  i.e.  (13.1.6)
is quasi-affine  since  the sub-{\rm GDD} by deleting  Vertex  4

\subsection* {Quasi-affine over
  {\rm  GDD}  $2$ of Row $13$  }

  {\rm (i) }  All quasi-arithmetic {\rm  GDD}s by  adding a vertex  on  Vertex $1$  are listed.
  According  to Lemma \ref {3.1.2} (II) we have to consider following cases.

   {\rm (ii) }  All quasi-arithmetic {\rm  GDD}s by  adding a vertex  on  Vertex $2$  are listed.

Checking  step by step we have to consider following cases.

  (a)  i.e.  (13.2.1)
  by   {\rm  GDD}  $2$ of Row $13$.  The sub-{\rm GDD} by deleting  Vertex  4   is an arithmetic  {\rm  GDD}
when $q\in R_3$ by Lemma \ref {3.1.1} {\rm (V)}.
 It  is quasi-affine $q\in R_3$ by Lemma \ref {3.1.3}.

  {\rm (iii) }  All quasi-arithmetic {\rm  GDD}s by  adding a vertex  on  Vertex $3$  are listed.
 According  to  Type   3 and  {\rm  GDD}  $4$ of Row $7$ we have to consider following cases.

  (a)  i.e.  (13.2.2),
   $q^2,  q^3 \not=1, $ by   {\rm  GDD}  $4$ of Row $7$  in Table A2. It is quasi-affine by   Lemma   \ref {2.63}.

  {\rm (iv) } All quasi-affine circles are listed.

  (nc)
  (a)  i.e.  (13.2.3)
  is  quasi-affine  since the sub-{\rm GDD} by deleting  Vertex 2.
with    $q^{3}=-1$ is an arithmetic  {\rm  GDD}   by Lemma \ref {3.1.2} (I) or by Lemma \ref {3.1.1} {\rm (IV)}.

\subsection* {Quasi-affine over
  {\rm  GDD}  $1$ of Row $14$   }
   {\rm (i) }  All quasi-arithmetic {\rm  GDD}s by  adding a vertex  on  Vertex $2$  are listed.
 According  to Lemma \ref {3.1.2} (I) we have to consider following cases.

  (a)  i.e.  (14.1.1)   by   {\rm  GDD}  $2$ of Row $6$.
 It  is quasi-affine   by Lemma \ref {3.1.3},  where $q^2 =- \xi,  \xi \in R_3$.

 (b)  i.e.  (14.1.4)
, $q\notin R_3$,  $q^3\in R_6$,  by   {\rm  GDD}  $2$ of Row $7$.
 It  is quasi-affine   by Lemma \ref {3.1.3}.

 (c)  i.e.  (14.1.5), $q\notin R_3\cup R_2$,  $q^3\in R_6$,  by   {\rm  GDD}  $4$ of Row $7$.
 It  is quasi-affine  when $q\in R_{18}$ by Lemma \ref {3.1.3}.

 (d)  i.e.  (14.1.3), by
 {Lemma \ref {3.1.1}(IV).
It is quasi-affine by Lemma \ref {3.1.3} when  $q\in R_6$.

 (e)  i.e.  (14.1.6)
  by    {\rm  GDD}  $1$ of Row $9$.
 It  is quasi-affine   by Lemma \ref {3.1.3} when  $q\in R_6$.

 (f)  i.e.  (14.1.2), $q\notin R_3\cup R_2$,  by {\rm  GDD}  $2$ of Row $10$.
 It  is quasi-affine   by Lemma \ref {3.1.3} when  $q\in R_6$.

  {\rm (ii) }  All quasi-arithmetic {\rm  GDD}s by  adding a vertex  on  Vertex $2$  are listed.
 According  to     Lemma  \ref {3.1.1}{\rm (I) }  we have to consider following cases.

  {\rm (iii) }  All quasi-arithmetic {\rm  GDD}s by  adding a vertex  on  Vertex $3$  are listed.
 According  to   Type   4  we have to consider following cases.

 (a)  i.e.  (14.1.7), $q\in R_3$ ,  by   {\rm  GDD}  $2$ of Row $13$  in Table A2.
 It  is quasi-affine by \cite [Lemma 2.4]{TZ22}.

  {\rm (iv) } All quasi-affine circles are listed.

(nc)  (a)   is not quasi-affine since\\

  $\begin{picture}(100,       15) \put(-68,        -1){}

\put(60,       1){\makebox(0,       0)[t]{$\bullet$}}

\put(28,       -1){\line(1,       0){33}}
\put(27,       1){\makebox(0,      0)[t]{$\bullet$}}

\put(-14,       1){\makebox(0,      0)[t]{$\bullet$}}

\put(-14,      -1){\line(1,       0){50}}

\put(-18,      10){$-1$}
\put(0,       5){$q$}
\put(22,      10){$-1$}
\put(40,       5){$q^{2}$}

\put(58,       10){$-q^{-1}$}

\put(80,         -1)  { is not  an arithmetic  {\rm  GDD}   by Lemma \ref {3.1.1}{\rm (I)} with $q\in R_6.$   } \end{picture}$

(a)(a)   is not quasi-affine since\\

  $\begin{picture}(100,       15) \put(-68,        -1){}

\put(60,       1){\makebox(0,       0)[t]{$\bullet$}}

\put(28,       -1){\line(1,       0){33}}
\put(27,       1){\makebox(0,      0)[t]{$\bullet$}}

\put(-14,       1){\makebox(0,      0)[t]{$\bullet$}}

\put(-14,      -1){\line(1,       0){50}}

\put(-18,      10){$-1$}
\put(0,       5){$q$}
\put(22,      10){$-1$}
\put(40,       5){$q^{2}$}

\put(58,       10){$-q^{-1}$}

\put(80,         -1)  {  is not an arithmetic  {\rm  GDD}   by Lemma \ref {3.1.1} (I) with $q\in R_6.$  } \end{picture}$

(b)  (a)  is not quasi-affine since\\

  $\begin{picture}(100,       15) \put(-68,        -1){}

\put(60,       1){\makebox(0,       0)[t]{$\bullet$}}

\put(28,       -1){\line(1,       0){33}}
\put(27,       1){\makebox(0,      0)[t]{$\bullet$}}

\put(-14,       1){\makebox(0,      0)[t]{$\bullet$}}

\put(-14,      -1){\line(1,       0){50}}

\put(-18,      10){$-1$}
\put(0,       5){$q$}
\put(22,      10){$-1$}
\put(40,       5){$q^{6}$}

\put(58,       10){$-q^{-3}$}

\put(80,         -1)  { with    $q\in R_{18}$ is not an arithmetic  {\rm  GDD}   by Lemma \ref {3.1.1} (I).   } \end{picture}$

 (d)  (a)  i.e.  (14.1.8)   is quasi-affine since the sub-{\rm GDD} by deleting  Vertex 2
  with $r=-1$ is  an arithmetic  {\rm  GDD}   by Lemma \ref {3.1.1}{\rm (I)} with $q\in R_6$.

\subsection* {Quasi-affine over
  {\rm  GDD}  $2$ of Row $14$   }
  {\rm (i) }  All quasi-arithmetic {\rm  GDD}s by  adding a vertex  on  Vertex $1$  are listed.
According  to Lemma \ref {3.1.2} (IV) we have to consider following cases.

 (a)  i.e.  (14.2.1), $q \in R_3 $. by {\rm  GDD}  $2$ of Row $7$. It  is quasi-affine  by Lemma \ref {3.1.3}.

 (b)  i.e.  (14.2.2),   $q\in R_6 $.  by Lemma \ref {3.1.1} {\rm (X)}.
It  is quasi-affine   by  {\rm  GDD}  $2$ of Row $7$.

   {\rm (ii) }  All quasi-arithmetic {\rm  GDD}s by  adding a vertex  on  Vertex $2$  are listed.
According  to Lemma \ref {3.1.2} (IV) we have to consider following cases.\\ \\

  {\ }\ \ \ \ \ \ \ \ \ \ \ \ $\begin{picture}(100,       15) \put(-68,        -1){$(a)$ }

\put(60,       1){\makebox(0,       0)[t]{$\bullet$}}
\put(58,       -12){$-1$}

\put(40,       -12){$q$}
\put(28,       -1){\line(1,       0){33}}
\put(27,       1){\makebox(0,      0)[t]{$\bullet$}}

\put(22,      -12){$-1$}
\put(0,       -12){$q^{-1}$}

\put(-14,       1){\makebox(0,      0)[t]{$\bullet$}}

\put(-14,      -1){\line(1,       0){50}}

\put(-18,      -12){$q$}

\put(27,     38){\makebox(0,      0)[t]{$\bullet$}}

\put(27,       0){\line(0,      1){35}}

\put(30,       35){$-q^{-1}$}

\put(30,       20){$q^{-1}$}

\ \ \ \ \ \ \ \ \ \ \ \ \   \   \ \ \ \ \ \ \ \  {$q\in R_6 $,  by   {\rm  GDD}  $2$ of Row $4$. The sub-{\rm GDD} by }
\put(80,         -1)  {    } \end{picture}$\\ \\
deleting  Vertex  4    is not an arithmetic  {\rm  GDD}   by Lemma \ref {3.1.1}{\rm (I)}.

\subsection* {Quasi-affine over
  {\rm  GDD}  $3$ of Row $14$ }
  {\rm (i) }  All quasi-arithmetic {\rm  GDD}s by  adding a vertex  on  Vertex $1$  are listed.
 According  to Lemma \ref {3.1.2} (III) we have to consider following cases.

 (a)  i.e.  (14.3.1),
 $ q \in R_6$,  by    {\rm  GDD}  $1$ of Row $7$. It  is quasi-affine  by Lemma \ref {3.1.3}.

 {\rm (ii) }  All quasi-arithmetic {\rm  GDD}s by  adding a vertex  on  Vertex $2$  are listed.\\

{\ }\!\!\!\!\!\!\!\!\!\!\!\!\!\!\!\!\!\!\!\!\!{\ }\!\!\!\!\!\!\!\!\!\!\!\!\!\!\!\!\!\!\!\!\!{\ }\!\!\!\!\!\!\!\!\!\!\!\!
  $\begin{picture}(55,      15)\put(68,       -1){ (a)}
\put(111,      1){\makebox(0,      0)[t]{$\bullet$}}
\put(144,       1){\makebox(0,       0)[t]{$\bullet$}}
\put(170,      -11){\makebox(0,      0)[t]{$\bullet$}}
\put(170,     15){\makebox(0,      0)[t]{$\bullet$}}
\put(113,      -1){\line(1,      0){33}}
\put(142,     -1){\line(2,      1){27}}
\put(170,       -14){\line(-2,      1){27}}

\put(100,       10){$-1$}

\put(115,       5){$-q$}

\put(135,      10){$-1$}

\put(140,      -20){$-q^{-1}$}
\put(145,       15){$-q^{}$}

\put(173,       -12){$q^{-1}$}
\put(173,       18){$-q^{-1}$}

\put(200,         -1)  {  $q\in R_3$ ,  by   {\rm  GDD}  $2$ of Row $14$ in Table A2 or  Type   4.  } \end{picture}$\\ \\
 The sub-{\rm GDD} by deleting  Vertex  4    is not an arithmetic  {\rm  GDD}     by Lemma \ref {3.1.1}{\rm (II)}.

\subsection* {Quasi-affine over
  {\rm  GDD}  $1$ of Row $15$  }
  {\rm (i) }  All quasi-arithmetic {\rm  GDD}s by  adding a vertex  on  Vertex $1$  are listed.
 According  to Lemma \ref {3.1.2} (I) we have to consider following cases.

 (a)  i.e.  (15.1.1),  $q^2 \not=1$,  by    {\rm  GDD}  $2$ of Row $6$.
 It  is quasi-affine  by Lemma \ref {3.1.8}.

 (b)  i.e.  (15.1.2),  $q ^ 3\in R_3$,  by  {\rm  GDD}  $2$ of Row $7$.
 It  is quasi-affine  when    $q ^ 3\in R_3$ by Lemma \ref {3.1.21}.

 (c)  i.e.  (15.1.3)  by   {\rm  GDD}  $4$ of Row $7$.
 It  is quasi-affine  when    $q ^ 3\in R_3$ by by Lemma \ref {3.1.3}.

 (d)  i.e.  (15.1.4)  by    {\rm  GDD}  $1$ of Row $9$.
 It  is quasi-affine  by Lemma \ref {3.1.8}.

 (e)  i.e.  (15.1.5)$, q \notin R_3$,  by   {\rm  GDD}  $2$ of Row $10$.
 It  is quasi-affine  by Lemma \ref {3.1.8} with $q\in R_6$.


   {\rm (ii) }  All quasi-arithmetic {\rm  GDD}s by  adding a vertex  on  Vertex $2$  are listed.
 According  to Lemma \ref {3.1.2} (III) we have to consider following cases.

 (a)  i.e.  (15.1.6)$, q \in F^{*}\setminus \{1,  -1\}$,  by   {\rm  GDD}  $3$ of Row $8$.
 The sub-{\rm GDD} by deleting  Vertex  4    is an arithmetic  {\rm  GDD}    by  {\rm  GDD}  1 of Row 15.  It  is quasi-affine  by Lemma \ref {3.1.3}.

 (b)  i.e.  (15.1.7)
 by   {\rm  GDD}  $1$ of Row $4$. The sub-{\rm GDD} by
deleting  Vertex  4 is   an arithmetic  {\rm  GDD}    by Lemma \ref {3.1.2} (III).  It  is quasi-affine  by Lemma \ref {3.1.3}.

  {\rm (iv) } All quasi-affine circles are listed.

 {\ }\ \ \ \ \ \ \ \ \ \ \ \ $\begin{picture}(100,       15) \put(-68,        -1){(b)  (nc)  }

\put(60,       1){\makebox(0,       0)[t]{$\bullet$}}

\put(28,       -1){\line(1,       0){33}}
\put(27,       1){\makebox(0,      0)[t]{$\bullet$}}

\put(-14,       1){\makebox(0,      0)[t]{$\bullet$}}

\put(-14,      -1){\line(1,       0){50}}

\put(-18,      10){$-1$}
\put(0,       5){$q$}
\put(22,      10){$-1$}
\put(40,       5){$q^{6}$}

\put(58,       10){$-1$}

\put(80,         -1)  {  is not  an arithmetic  {\rm  GDD}   by Lemma \ref {3.1.1} (X).
  } \end{picture}$

 (d)  (nc)  i.e.  (15.1.8) is quasi-affine since the sub-{\rm GDD} by deleting  Vertex 2
 with  $r =-1 $. It  is an arithmetic  {\rm  GDD}    by Lemma \ref {3.1.1} (X).

 (d)  (nc)  i.e.  (15.1.9)
 is quasi-affine  since the sub-{\rm GDD} by deleting  Vertex 2
 with  $r =q^{-2} $ is an arithmetic  {\rm  GDD}    by  Type   7.

 (e)   (nc)  i.e.  (15.1.10)  is quasi-affine since the sub-{\rm GDD} by deleting  Vertex 2
 with  $q =q^{4} $ is an arithmetic  {\rm  GDD}    by  Type   7.

\subsection* {Quasi-affine over
  {\rm  GDD}  $2$ of Row $15$ }
   {\rm (i) } All quasi-arithmetic {\rm  GDD}s by  adding a vertex  on  Vertex $1$  are listed.
 According  to Lemma \ref {3.1.2} (IV) we have to consider following cases.

 (a)  i.e.  (15.2.2)$, q \in R_3$,  by    {\rm  GDD}  $2$ of Row $15$. It is quasi-affine   by Lemma \ref {3.1.28}.

  {\rm (ii) }  All quasi-arithmetic {\rm  GDD}s by  adding a vertex  on  Vertex $2$  are listed.

 (a)  i.e.  (15.2.1)
 $, q \in R_3$,  by    {\rm  GDD}  $2$ of Row $15$. It is quasi-affine   by Lemma \ref {3.1.28}.

  {\rm (iii) }  All quasi-arithmetic {\rm  GDD}s by  adding a vertex  on  Vertex $2$  are listed.
It  is the same as (i).

  {\rm (iv) }  All quasi-affine circles are listed.

 (nc)
 (a)  i.e.  (15.2.3)  is quasi-affine since the sub-{\rm GDD} by deleting  Vertex 2
  is an arithmetic  {\rm  GDD}   by Lemma \ref {3.1.1} {\rm (X)} or  Type   7.

 (a)
 (a)  i.e.  (15.2.4) is quasi-affine since the sub-{\rm GDD} by deleting  Vertex 2  is  an arithmetic  {\rm  GDD}   by Lemma \ref {3.1.1} (X).
\subsection* {Quasi-affine over
  {\rm  GDD}  $3$ of Row $15$  }
   {\rm (i) } All quasi-arithmetic {\rm  GDD}s by  adding a vertex  on  Vertex $1$  are listed.
 According  to Lemma \ref {3.1.2} (III) we have to consider following cases.

 (a)  i.e.  (15.3.1)$, q \in R_3$,  by    {\rm  GDD}  $3$ of Row $6$.
  The sub-{\rm GDD} by deleting  Vertex  1    is   an arithmetic  {\rm  GDD}   by Lemma \ref {3.1.1} {\rm (IX)}.
 The sub-{\rm GDD} by deleting  Vertex  2   is an arithmetic  {\rm  GDD}   by  Type   7. It  is quasi-affine.

(b)  i.e.  (15.3.3)
 $, q \in R_3$,  by    {\rm  GDD}  $1$ of Row $15$.  The sub-{\rm GDD} by deleting  Vertex  4   is
an arithmetic  {\rm  GDD}   Lemma \ref {3.1.2}{\rm (II)}.
It
is quasi-affine   by Lemma \ref {3.1.3}.

(c)  i.e.  (15.3.2)$, q \in R_3$,  by     {\rm  GDD}  $3$ of Row $15$.
 The sub-{\rm GDD} by deleting  Vertex  2    is  an arithmetic  {\rm  GDD}   by Lemma \ref {3.1.1} {\rm (IV)}.  The sub-{\rm GDD} by deleting  Vertex  1    is  an arithmetic  {\rm  GDD}   by  Type   7. It  is   quasi-affine,
\\

  {\ }\ \ \ \ \ \ \ \ \ \ \ \ $\begin{picture}(100,       15) \put(-68,        -1){ (d)}

\put(60,       1){\makebox(0,       0)[t]{$\bullet$}}
\put(-24,      -12){$-q$}

\put(-7,       -12){$-q^{-1}$}

\put(22,      -12){$q$}

\put(40,       -12){$q^{-1}$}

\put(28,       -1){\line(1,       0){33}}
\put(27,       1){\makebox(0,      0)[t]{$\bullet$}}
\put(58,       -12){$-1$}

\put(-14,       1){\makebox(0,      0)[t]{$\bullet$}}

\put(-18,      -1){\line(1,       0){50}}

\put(59,       0){\line(-1,      1){17}}

\put(28,      -1){\line(1,       1){17}}

\put(43,     18){\makebox(0,      0)[t]{$\bullet$}}

\put(18,     12){$q^{-1}$}

\put(36,     22){$q$}
\put(59,     12){$q^{-1}$}

 \ \ \ \ \ \ \ \ \ \ \ \ \ \ \ \ \ \ \ \ \ \ \ \  \ \ \ \ \ \ \ \
  \ \ \ \ \ \ \ \ \ \ \ \ \ \ \ \ \ \ \ {$, q \in R_3$,  by   {\rm  GDD}  $1$ of Row $16$.}
\put(80,         -1)  {    } \end{picture}$\\
\\
 The sub-{\rm GDD} by deleting  Vertex  4    is not an arithmetic  {\rm  GDD}    by Lemma \ref {3.1.2}{\rm (II)}.\\

  {\ }\ \ \ \ \ \ \ \ \ \ \ \ $\begin{picture}(100,       15) \put(-68,        -1){ (e)}

\put(60,       1){\makebox(0,       0)[t]{$\bullet$}}
\put(58,       -12){$-1$}

\put(40,       -12){$q^{-1}$}
\put(28,       -1){\line(1,       0){33}}
\put(27,       1){\makebox(0,      0)[t]{$\bullet$}}

\put(22,      -12){$q$}
\put(0,       -12){$-1$}

\put(-14,       1){\makebox(0,      0)[t]{$\bullet$}}

\put(-14,      -1){\line(1,       0){50}}

\put(-18,      -12){$-1$}

\put(59,       0){\line(-1,      1){17}}

\put(28,      -1){\line(1,       1){17}}

\put(43,     18){\makebox(0,      0)[t]{$\bullet$}}

\put(18,     12){$q^{-1}$}

\put(36,     22){$q$}
\put(59,     12){$q^{-1}$}

 \ \ \ \ \ \ \ \ \ \ \ \ \ \ \ \ \ \ \ \ \ \ \ \  \ \ \ \ \ \ \ \
  \ \ \ \ \ \ \ \ \ \ \ \ \ \ \ \ \ \ \ {$, q \in R_3$,  by   {\rm  GDD}  $2$ of Row $17$.}
\put(80,         -1)  {    } \end{picture}$\\ \\
 The sub-{\rm GDD} by deleting  Vertex  4    is not an arithmetic  {\rm  GDD}   by Lemma \ref {3.1.2}{\rm (II)}.\\

  {\ }\ \ \ \ \ \ \ \ \ \ \ \ $\begin{picture}(100,       15) \put(-68,        -1){ (f)}

\put(60,       1){\makebox(0,       0)[t]{$\bullet$}}
\put(58,       -12){$-1$}

\put(40,       -12){$q^{-1}$}
\put(28,       -1){\line(1,       0){33}}
\put(27,       1){\makebox(0,      0)[t]{$\bullet$}}

\put(22,      -12){$q$}
\put(0,       -12){$-q$}

\put(-14,       1){\makebox(0,      0)[t]{$\bullet$}}

\put(-14,      -1){\line(1,       0){50}}

\put(-18,      -12){$-1$}

\put(59,       0){\line(-1,      1){17}}

\put(28,      -1){\line(1,       1){17}}

\put(43,     18){\makebox(0,      0)[t]{$\bullet$}}

\put(18,     12){$q^{-1}$}

\put(36,     22){$q$}
\put(59,     12){$q^{-1}$}

 \ \ \ \ \ \ \ \ \ \ \ \ \ \ \ \ \ \ \ \ \ \ \ \  \ \ \ \ \ \ \ \
  \ \ \ \ \ \ \ \ \ \ \ \ \ \ \ \ \ \ \ {$, q \in R_3$,  by   {\rm  GDD}  $7$ of Row $17$.}
\put(80,         -1)  {    } \end{picture}$\\
\\
 The sub-{\rm GDD} by deleting  Vertex  4    is not an arithmetic  {\rm  GDD}   by Lemma \ref {3.1.1} {\rm (VIII)}.

(g)  i.e.  (15.3.9),  $q\in R_3, $ by   {\rm  GDD}  $3$ of Row $15$ in Table A2. It is quasi-affine.

   {\rm (ii) }  All    quasi-arithmetic {\rm  GDD}s by  adding a vertex  on  Vertex $2$  are listed.
 According  to Lemma \ref {3.1.2} (I) we have to consider following cases.

(a)  i.e.  (15.3.4)  by    {\rm  GDD}  $2$ of Row $6$. It is quasi-affine   by Lemma \ref {3.1.34}.

(b)  i.e.  (15.3.5) by  Lemma \ref {3.1.1} {\rm  (IV)} . It is quasi-affine   by Lemma \ref {3.1.34}.
\\ \\ \\ \\  \\

  \ \ \ \ \ \  $\begin{picture}(100,     15)\put(-45,      -1){ }\put(-45,      -1){(i) in  Case {\rm (i) } }

\put(170,    10){\makebox(0,     0)[t]{$\bullet$}}

\put(170,    70){\makebox(0,     0)[t]{$\bullet$}}

\put(230,    10){\makebox(0,     0)[t]{$\bullet$}}
\put(230,    70){\makebox(0,     0)[t]{$\bullet$}}

\put(170,      10){\line(0,     1){60}}

\put(170,      10){\line(1,     1){60}}

\put(230,      10){\line(0,     1){60}}

\put(230,      10){\line(-1,     1){60}}

\put(170,      10){\line(1,      0){60}}
\put(170,      70){\line(1,      0){60}}

\put(150,    10){$q$}
\put(150,    30){$q^{-1}$}

\put(150,    70){$q$}

\put(180,    30){$q^{-1}$}
\put(190,    80){$q^{-1}$}
\put(190,    -10){$-1$}
\put(220,    30){${}$}

\put(250,    10){$-1$}
\put(250,    30){$-q$}

\put(250,    70){$-1$}

\put(80,        -1)  {    } \end{picture}$\\

All quasi-arithmetic {\rm  GDD}s by  deleting   Vertex 2 is not an arithmetic {\rm  GDD}  by Lemma \ref {3.1.1}{\rm (I)}.

  {\rm (iii) }  All    quasi-arithmetic {\rm  GDD}s by  adding a vertex  on  Vertex $3$  are listed. It is the same as  {\rm (i)}.

  {\rm (iv) } All quasi-affine  {\rm  GDD}s which are complete diagrams are listed.

(a)  in  Case {\rm (i) }  i.e.  (15.3.8) with  $q\in R_3.$
If the sub-{\rm GDD}
by deleting  Vertex 2 is  an arithmetic  {\rm  GDD},   then  $ \widetilde{q}_{43} =q^{-1}$
by Lemma \ref {3.1.1}{\rm (I)}. The sub-{\rm GDD} by deleting  Vertex 1  is  an arithmetic  {\rm  GDD}   by
 {\rm  GDD}  $3$ of Row $6$
. It  is   quasi-affine   when  $ \widetilde{q}_{43} =q^{-1}$.
\\ \\ \\ \\ \\

{\ }\!\!\!\!\!\!\!\!\!\!\!\!\!\!\!\!\!\!\!\!\!{\ }\!\!\!\!\!\!\!\!\!\!\!\!\!\!\!\!\!\!\!\!\!{\ }\!\!\!\!\!\!\!\!\!\!\!\!\!\!
  $\begin{picture}(100,       15)\put(-68,       -1){ }\put(68,       -1){  (c) in Case (i)     }

\put(170,     10){\makebox(0,      0)[t]{$\bullet$}}

\put(170,     70){\makebox(0,      0)[t]{$\bullet$}}

\put(230,     10){\makebox(0,      0)[t]{$\bullet$}}
\put(230,     70){\makebox(0,      0)[t]{$\bullet$}}

\put(170,       10){\line(0,      1){60}}

\put(170,       10){\line(1,      1){60}}

\put(230,       10){\line(0,      1){60}}

\put(230,       10){\line(-1,      1){60}}

\put(170,       10){\line(1,       0){60}}
\put(170,       70){\line(1,       0){60}}

\put(150,     10){$q$}
\put(150,     30){$q^{-1}$}

\put(150,     70){$q$}

\put(180,     30){$q^{-1}$}
\put(190,     80){$q^{-1}$}
\put(190,     -10){$q^{-1}$}
\put(220,     30){$$}

\put(250,     10){$q$}
\put(250,     30){$q^{-1}$}

\put(250,     70){$-1$}

\put(270,         -1)  { The sub-{\rm GDD} by deleting  Vertex 2 is not an } \end{picture}$
\\ \\
 arithmetic  {\rm  GDD} by Lemma \ref {3.1.1} {\rm (XV) }.

\subsection* {Quasi-affine over
  {\rm  GDD}  $4$ of Row $15$    }

  {\rm (ii) }  All     quasi-arithmetic {\rm  GDD}s by  adding a vertex  on  Vertex $2$  are listed.

Checking  step by step we have to consider following cases.

(a)  i.e.  (15.4.1),  $q\in R_3$,  by    {\rm  GDD}  $2$ of Row $13$.
 The sub-{\rm GDD} by deleting  Vertex  4   is an arithmetic  {\rm  GDD}   because of  symmetry of this  {\rm  GDD}. It is similar in following of  Case {\rm (ii) }.
It  is quasi-affine   by Lemma \ref {3.1.3}.

(b)  i.e.  (15.4.2),  $q\in R_3$,  by   {\rm  GDD}  $4$ of Row $15$. It is quasi-affine   by Lemma \ref {3.1.3}.

\subsection* {Quasi-affine over
  {\rm  GDD}  $1$ of Row $16$   }
  {\rm (i) }  All    quasi-arithmetic {\rm  GDD}s by  adding a vertex  on  Vertex $1$  are listed.
 According  to Lemma \ref {3.1.2} (I) we have to consider following cases.

(a)  i.e.  (16.1.2)  by   {\rm  GDD}  $2$ of Row $6$.
is quasi-affine   by Lemma \ref {3.1.8}.

(a)  i.e.  (16.1.3)  by   {\rm  GDD}  $2$ of Row $7$.
is quasi-affine   by Lemma \ref {3.1.8} with $q^3\in R_3$.

(c)  i.e.  (16.1.4)   by  {\rm  GDD}  $4$ of Row $7$.
is quasi-affine   by Lemma \ref {3.1.3} with $q^3\in R_3$.

(d)  i.e.  (16.1.5)   by    {\rm  GDD}  $1$ of Row $9$. It  is
an
arithmetic  {\rm  GDD}    when $r= -q$ by Lemma \ref {3.1.8}.
It  is quasi-affine   when $r\not= -q$ by Lemma \ref {3.1.8}.

(e)  i.e.  (16.1.6)$, q \notin R_3$,  by    {\rm  GDD}  $2$ of Row $10$. It
is quasi-affine   by Lemma \ref {3.1.8} with $q\in R_6$.

   {\rm (ii) }  All quasi-arithmetic {\rm  GDD}s by  adding a vertex  on  Vertex $2$  are listed.
 According  to Lemma \ref {3.1.2} (III) we have to consider following cases.\\ \\

  {\ }\ \ \ \ \ \ \ \ \ \ \ \ $\begin{picture}(100,       15) \put(-68,        -1){(a) }

\put(60,       1){\makebox(0,       0)[t]{$\bullet$}}
\put(58,       -12){$-1$}

\put(40,       -12){$q^{-1}$}
\put(28,       -1){\line(1,       0){33}}
\put(27,       1){\makebox(0,      0)[t]{$\bullet$}}

\put(22,      -12){$q$}
\put(0,       -12){$q^{-2}$}

\put(-14,       1){\makebox(0,      0)[t]{$\bullet$}}

\put(-14,      -1){\line(1,       0){50}}

\put(-18,      -12){$q^2$}

\put(27,     38){\makebox(0,      0)[t]{$\bullet$}}

\put(27,       0){\line(0,      1){35}}

\put(30,       30){$-q^{}$}

\put(30,       15){$-q^{-1}$}

\ \ \ \ \ \ \ \ \ \ \ \ \ \ \ \ \ \ \ \ \ \ \ \ \ \ \ \ \ \ \ \ \ \ \ \ \ \   \ \ \ \ \ \ \ \ \ \ \ \ \ \ \ \ \ \ \ \ {$, q \in R_3$,  by   {\rm  GDD}  $1$ of Row $6$.}
\put(80,         -1)  {    } \end{picture}$\\
\\
 The sub-{\rm GDD} by deleting  Vertex  4    is not an arithmetic  {\rm  GDD}    by Lemma \ref {3.1.1} {\rm (V)}.

(b)  i.e.  (16.1.1)$, q \in F^{*}\setminus \{1,  -1\}$,  by    {\rm  GDD}  $3$ of Row $8$.  The sub-{\rm GDD} by deleting
 Vertex  4  is an arithmetic  {\rm  GDD}   by Lemma \ref {3.1.2} {\rm (III)}.
It  is  quasi-affine   by Lemma \ref {3.1.3}.\\ \\

  {\ }\ \ \ \ \ \ \ \ \ \ \ \ $\begin{picture}(100,       15) \put(-68,        -1){ (c)}

\put(60,       1){\makebox(0,       0)[t]{$\bullet$}}
\put(58,       -12){$-1$}

\put(40,       -12){$q^{-1}$}
\put(28,       -1){\line(1,       0){33}}
\put(27,       1){\makebox(0,      0)[t]{$\bullet$}}

\put(22,      -12){$q$}
\put(0,       -12){$q$}

\put(-14,       1){\makebox(0,      0)[t]{$\bullet$}}

\put(-14,      -1){\line(1,       0){50}}

\put(-18,      -12){$-1$}

\put(27,     38){\makebox(0,      0)[t]{$\bullet$}}

\put(27,       0){\line(0,      1){35}}

\put(30,       30){$-q^{}$}

\put(30,       15){$-q^{-1}$}

\ \ \ \ \ \ \ \ \ \ \ \ \ \ \ \ \ \ \ \ \ \ \ \ \ \ \ \ \ \ \ \ \ \ \ \ \ \   \   \ \ \ \ \ \ \ \ \ \ \ \ \ \ \ \ \ \ \ {$, q \in R_3$,  by    {\rm  GDD}  $1$ of Row $15$.}
\put(80,         -1)  {    } \end{picture}$\\
\\
 The sub-{\rm GDD} by deleting  Vertex  4    is not an arithmetic  {\rm  GDD}    by Lemma \ref {3.1.1} {\rm (V)}.\\ \\

    {\ }\ \ \ \ \ \ \ \ \ \ \ \    \ \ \ \ \ \  $\begin{picture}(100,      15) \put(-85,      -1){(d) }

\put(60,      1){\makebox(0,      0)[t]{$\bullet$}}
\put(58,      -12){$-1$}

\put(40,      -12){$q^{-1}$}
\put(28,      -1){\line(1,      0){33}}
\put(27,      1){\makebox(0,     0)[t]{$\bullet$}}

\put(26,     -12){$q$}
\put(0,      -12){$-q^{-1}$}

\put(-14,      1){\makebox(0,     0)[t]{$\bullet$}}

\put(-14,     -1){\line(1,      0){50}}

\put(-18,     -12){$-q$}

\put(27,    38){\makebox(0,     0)[t]{$\bullet$}}

\put(27,      0){\line(0,     1){35}}

\put(30,      30){$-q^{}$}

\put(30,      15){$-q^{-1}$}

\ \ \ \ \ \ \ \ \ \ \ \ \ \ \ \ \ \ \ \ \ \ \ \ \ \ \ \ \ \ \ \ \ \ \ \ \ \   \  \ \ \ \ \ \ \ \ \ \ \ \ \ \ \ \ \ \ \ {$,q \in R_3$.  {\rm GDD} $1$ of Row $16$.}
\put(80,        -1)  {    } \end{picture}$\\
\\
 The sub-{\rm GDD} by deleting  Vertex  4  is not an arithmetic GDD   by Lemma \ref {3.1.1} {\rm (V)}.\\
\\

{\ }\ \ \ \ \ \ \ \ \ \ \ \ $\begin{picture}(100,       15) \put(-68,        -1){ $(e)$}

\put(60,       1){\makebox(0,       0)[t]{$\bullet$}}

\put(28,       -1){\line(1,       0){33}}
\put(27,       1){\makebox(0,      0)[t]{$\bullet$}}
\put(-14,       1){\makebox(0,      0)[t]{$\bullet$}}

\put(-14,      -1){\line(1,       0){50}}

\put(58,       -12){$-1$}

\put(40,       -12){$q^{-1}$}

\put(22,      -12){$q$}
\put(0,       -12){${q^{-1}}$}

\put(-18,      -12){${q}$}

\put(27,     38){\makebox(0,      0)[t]{$\bullet$}}

\put(27,       0){\line(0,      1){35}}

\put(30,       30){$-q^{}$}

\put(30,       15){$-q^{-1}$}

\ \ \ \ \ \ \ \ \ \ \ \ \ \ \ \ \ \ \ \ \ \ \ \ \ \ \ \ \ \ \ \ \ \ \ \ \ \   \   \ \ \ \ \ \ \ \ \ \ \ \ \ \ \ \ \ \ \ { ,  by   {\rm  GDD}  $1$ of Row $4$.}
\put(80,         -1)  {    } \end{picture}$\\
\\
 The sub-{\rm GDD} by deleting  Vertex  4    is not an arithmetic  {\rm  GDD}    by Lemma \ref {3.1.2}{\rm (II)}.\\
\\

{\ }\ \ \ \ \ \ \ \ \ \ \ \ $\begin{picture}(100,       15) \put(-68,        -1){$(f)$ }

\put(60,       1){\makebox(0,       0)[t]{$\bullet$}}

\put(28,       -1){\line(1,       0){33}}
\put(27,       1){\makebox(0,      0)[t]{$\bullet$}}
\put(-14,       1){\makebox(0,      0)[t]{$\bullet$}}

\put(-14,      -1){\line(1,       0){50}}

\put(58,       -12){$-1$}

\put(40,       -12){$q^{-2}$}

\put(22,      -12){$q^2$}
\put(0,       -12){${q^{-2}}$}

\put(-18,      -12){${q}$}

\put(27,     38){\makebox(0,      0)[t]{$\bullet$}}

\put(27,       0){\line(0,      1){35}}

\put(30,       30){$-q^{2}$}

\put(30,       15){$-q^{-2}$}

\ \ \ \ \ \ \ \ \ \ \ \ \ \ \ \ \ \ \ \ \ \ \ \ \ \ \ \ \ \ \ \ \ \ \ \ \ \   \   \ \ \ \ \ \ \ \ \ \ \ \ \ \ \ \ \ \ \ {,  by    {\rm  GDD}  $1$ of Row $5$.}
\put(80,         -1)  {    } \end{picture}$\\
\\
 The sub-{\rm GDD} by deleting  Vertex  4    is not an arithmetic  {\rm  GDD}   by Lemma \ref {3.1.1}{\rm (I)}.

\subsection* {Quasi-affine over
  {\rm  GDD}  $2$ of Row $16$  }
  {\rm (i) }  All quasi-arithmetic {\rm  GDD}s by  adding a vertex  on  Vertex $1$  are listed.
 According  to Lemma \ref {3.1.2} (IV) we have to consider following cases.

(a)  i.e.  (16.2.1)
$, q \in R_3$,  by    {\rm  GDD}  $2$ of Row $15$. It is quasi-affine   by Lemma \ref {3.1.21}.

(b)  i.e.  (16.2.2),
  $q \in R_3$,  by   {\rm  GDD}  $2$ of Row $16$. It is quasi-affine   by Lemma \ref {3.1.21}.

  {\rm (ii) }  All quasi-arithmetic {\rm  GDD}s by  adding a vertex  on  Vertex $2$  are listed.
 According  to Lemma \ref {3.1.2} (IV) we have to consider following cases.

(a)  i.e.  (16.2.3)$, q \in F^{*}\setminus \{1,  -1\}$.   {\rm  GDD}  $2$ of Row $4$,  by   the sub-{\rm GDD} by deleting Vertex  4  is not an arithmetic  {\rm  GDD}   by Lemma \ref {3.1.1} {\rm (IV)}.
It  is  quasi-affine   by Lemma \ref {3.1.16}.\\ \\

  {\ }\ \ \ \ \ \ \ \ \ \ \ \ $\begin{picture}(100,       15) \put(-68,        -1){ $(b)$}

\put(60,       1){\makebox(0,       0)[t]{$\bullet$}}
\put(58,       -12){$-1$}

\put(40,       -12){$q^2$}
\put(28,       -1){\line(1,       0){33}}
\put(27,       1){\makebox(0,      0)[t]{$\bullet$}}

\put(22,      -12){$-1$}
\put(0,       -12){$q^{-2}$}

\put(-14,       1){\makebox(0,      0)[t]{$\bullet$}}

\put(-14,      -1){\line(1,       0){50}}

\put(-18,      -12){$q$}

\put(27,     38){\makebox(0,      0)[t]{$\bullet$}}

\put(27,       0){\line(0,      1){35}}

\put(30,       30){$-q^{2}$}

\put(30,       20){$-q^{-2}$}

\ \ \ \ \ \ \ \ \ \ \ \ \ \ \ \ \ \ \ \ \ \ \ \ \ \ \ \ \ \ \ \   \  \ \ \ \ \ \ \ \ \ \ \ \ \ \ \ \ \ \ \ {$, q \in F^{*}\setminus \{1,  -1\}$,  by  {\rm  GDD}  $2$ of Row $5$.}
\put(80,         -1)  {    } \end{picture}$\\
\\
 The sub-{\rm GDD} by deleting  Vertex  4    is not an arithmetic  {\rm  GDD}   by Lemma \ref {3.1.1}{\rm (I)}.

 (c)  i.e.  (16.2.4)$, q \in F^{*}\setminus \{1,  -1\}$,  by    {\rm  GDD}  $2$ of Row $6$.  The sub-{\rm GDD} by deleting
 Vertex  4  is  an arithmetic  {\rm  GDD}   by Lemma \ref {3.1.1} {\rm (IV)}.
  It  is quasi-affine   by Lemma \ref {3.1.16}.\\ \\

  {\ }\ \ \ \ \ \ \ \ \ \ \ \ $\begin{picture}(100,       15) \put(-68,        -1){(d) }

\put(60,       1){\makebox(0,       0)[t]{$\bullet$}}
\put(58,       -12){$-1$}

\put(40,       -12){$q$}
\put(28,       -1){\line(1,       0){33}}
\put(27,       1){\makebox(0,      0)[t]{$\bullet$}}

\put(22,      -12){$-1$}
\put(0,       -12){$q^{-1}$}

\put(-14,       1){\makebox(0,      0)[t]{$\bullet$}}

\put(-14,      -1){\line(1,       0){50}}

\put(-18,      -12){$-1$}

\put(27,     38){\makebox(0,      0)[t]{$\bullet$}}

\put(27,       0){\line(0,      1){35}}

\put(30,       30){$-q^{}$}

\put(30,       20){$-q^{-1}$}

\ \ \ \ \ \ \ \ \ \ \ \ \ \ \ \ \ \ \ \ \ \ \ \ \ \ \ \ \ \ \   \   \ \ \ \ \ \ \ \ \ \ \ \ \ \ \ \ \ \ \ {$, q \in F^{*}\setminus \{1,  -1\}$,  by   {\rm  GDD}  $2$ of Row $8$.}
\put(80,         -1)  {    } \end{picture}$\\
\\
 The sub-{\rm GDD} by deleting  Vertex  4    is not an arithmetic  {\rm  GDD}   by Lemma \ref {3.1.1}{\rm (II)}.

 (e)  i.e.  (16.2.5)$, q \in R_3$,  by  {\rm  GDD}  $2$ of Row $15$.  The sub-{\rm GDD} by deleting  Vertex  4
   is an arithmetic  {\rm  GDD}   by Lemma \ref {3.1.1}{\rm (II)}.  It  is quasi-affine
by Lemma \ref {3.1.16}.\\ \\

    {\ }\ \ \ \ \ \ \ \ \ \ \ \    \ \ \ \ \ \  $\begin{picture}(100,      15) \put(-85,      -1){ (f)}

\put(60,      1){\makebox(0,      0)[t]{$\bullet$}}
\put(58,      -12){$-1$}

\put(37,      -12){$-q^{-1}$}
\put(28,      -1){\line(1,      0){33}}
\put(27,      1){\makebox(0,     0)[t]{$\bullet$}}

\put(22,     -12){$-1$}
\put(0,      -12){$-q$}

\put(-14,      1){\makebox(0,     0)[t]{$\bullet$}}

\put(-14,     -1){\line(1,      0){50}}

\put(-18,     -12){$q$}

\put(27,    38){\makebox(0,     0)[t]{$\bullet$}}

\put(27,      0){\line(0,     1){35}}

\put(30,      33){$q^{-1}$}

\put(30,      20){$q^{}$}

\ \ \ \ \ \ \ \ \ \ \ \ \ \ \ \ \ \ \ \ \ \ \ \ \ \ \ \ \ \ \ \ \ \ \ \ \ \   \   \ \ \ \ \ \ \ \ \ \ \ \ \ \ \ \ \ \ \ {$,q \in R_3$.  GDD $2$ of Row $14$.}
\put(80,        -1)  {    } \end{picture}$\\ \\
The sub-{\rm GDD} by deleting  Vertex  4    is not an arithmetic  {\rm  GDD}   by Lemma \ref {3.1.1}{\rm (I)}.

 (g)  i.e.  (16.2.6),  $q \in R_3$,  by   {\rm  GDD}  $2$ of Row $16$.

  {\rm (iii) }  All quasi-arithmetic {\rm  GDD}s by  adding a vertex  on  Vertex $3$  are listed.
 According  to Lemma \ref {3.1.2} (III) we have to consider following cases.

 (a)  i.e.  (16.2.7)
 $q=-\xi,  \xi \in R_3 $  $, q \notin R_3$,  by   {\rm  GDD}  $1$ of Row $7$. It is quasi-affine by \cite [Lemma 2.4]{TZ22}.

  {\rm (iv) }  All quasi-affine circles are listed.

 (a)  (nc)  i.e.  (16.2.8)
  is  quasi-affine  since the sub-{\rm GDD} by deleting  Vertex 2

 is  an arithmetic  {\rm  GDD}   by Lemma \ref {3.1.1} (II) or by Lemma \ref {3.1.1} (IX).

  (b)  (b)  i.e.  (16.2.9)  is  quasi-affine  since the sub-{\rm GDD} by deleting  Vertex   is an arithmetic  {\rm  GDD}.

\subsection* {Quasi-affine over
  {\rm  GDD}  $3$ of Row $16$  }
  {\rm (i) }  All quasi-arithmetic {\rm  GDD}s by  adding a vertex  on  Vertex $1$  are listed.
 According  to Lemma \ref {3.1.2} (III) we have to consider following cases.\\

  {\ }\ \ \ \ \ \ \ \ \ \ \ \ $\begin{picture}(100,       15) \put(-68,        -1){ (a)}

\put(60,       1){\makebox(0,       0)[t]{$\bullet$}}
\put(58,       -12){$-1$}

\put(40,       -12){$q^{-1}$}
\put(28,       -1){\line(1,       0){33}}
\put(27,       1){\makebox(0,      0)[t]{$\bullet$}}

\put(22,      -12){$q$}
\put(0,       -12){$q^{-2}$}

\put(-14,       1){\makebox(0,      0)[t]{$\bullet$}}

\put(-14,      -1){\line(1,       0){50}}

\put(-18,      -12){$q^2$}

\put(59,       0){\line(-1,      1){17}}

\put(28,      -1){\line(1,       1){17}}

\put(43,     18){\makebox(0,      0)[t]{$\bullet$}}

\put(18,     12){${-1}$}

\put(36,     22){$-1$}
\put(59,     12){$-q$}

 \ \ \ \ \ \ \ \ \ \ \ \ \ \ \ \ \ \ \ \ \ \ \ \  \ \ \ \ \ \ \ \
  \ \ \ \ \ \ \ \ \ \ \ \ \ \ \ \ \ \ \ {$, q \in R_3$,  by   {\rm  GDD}  $1$ of Row $6$.}
\put(80,         -1)  {    } \end{picture}$\\ \\
 The sub-{\rm GDD} by deleting  Vertex  4    is not an arithmetic  {\rm  GDD}    by Lemma \ref {3.1.1} {\rm (V)}.\\ \\ \\ \\ \\

{\ }\!\!\!\!\!\!\!\!\!\!\!\!\!\!\!\!\!\!\!\!\!{\ }\!\!\!\!\!\!\!\!\!\!\!\!\!\!\!\!\!\!\!\!\!{\ }\!\!\!\!\!\!\!\!\!\!\!\!\!\!\!\!\!\!\!\!\!{\ }\!\!\! \!\!\!\
 $\begin{picture}(100,      15)\put(45,       -1){ }\put(95,       -1){(b) }

\put(170,     10){\makebox(0,      0)[t]{$\bullet$}}

\put(170,     70){\makebox(0,      0)[t]{$\bullet$}}

\put(230,     10){\makebox(0,      0)[t]{$\bullet$}}
\put(230,     70){\makebox(0,      0)[t]{$\bullet$}}

\put(170,       10){\line(0,      1){60}}

\put(170,       10){\line(1,      1){60}}

\put(230,       10){\line(0,      1){60}}


\put(170,       10){\line(1,       0){60}}
\put(170,       70){\line(1,       0){60}}

\put(150,     10){$q$}
\put(150,     30){$q^{-1}$}

\put(150,     70){$-1$}

\put(180,     30){$q^{-1}$}
\put(190,     80){$q^{2}$}
\put(190,     -10){$-1$}
\put(220,     30){$$}

\put(250,     10){$-1$}
\put(250,     30){$-q$}

\put(250,     70){$-1$}

\ \ \ \ \ \ \ \ \ \ \ \ \ \  \ \ \ \ \ \ \ \ \ \ \ \ \ \ \ \ \ \ \ \ \ \ \ \  \ \ \ \ \ \ \ \ \ \ \
 \ \ \ \ \ \ \ \ \ \ \ \ \ \ \ \ \ \ \ \ \ \ \ \ {$, q \in R_3$,  by    {\rm  GDD}  $3$ of Row $6$. The sub-{\rm GDD} by deleting}
\put(80,         -1)  {    } \end{picture}$\\
\\
 Vertex  1   is not  an arithmetic  {\rm  GDD}    by Lemma \ref {3.1.1}{\rm (II)} or by Lemma \ref {3.1.1} {\rm (X)}.

 (c)  i.e.  (16.3.1)$, q \in F^{*}\setminus \{1,  -1\}$,  by    {\rm  GDD}  $3$ of Row $8$.
  The sub-{\rm GDD} by deleting  Vertex 4    is an arithmetic  {\rm  GDD}    by  {\rm  GDD}  2 of Row 17 or by Lemma \ref {3.1.2} {\rm (III)}(l). It  is quasi-affine   by Lemma \ref {3.1.3}.\\

  {\ }\ \ \ \ \ \ \ \ \ \ \ \ $\begin{picture}(100,       15) \put(-68,        -1){ (d)}

\put(60,       1){\makebox(0,       0)[t]{$\bullet$}}
\put(58,       -12){$-1$}

\put(40,       -12){$q^{-1}$}
\put(28,       -1){\line(1,       0){33}}
\put(27,       1){\makebox(0,      0)[t]{$\bullet$}}

\put(22,      -12){$q$}
\put(0,       -12){$q$}

\put(-14,       1){\makebox(0,      0)[t]{$\bullet$}}

\put(-14,      -1){\line(1,       0){50}}

\put(-18,      -12){$-1$}

\put(59,       0){\line(-1,      1){17}}

\put(28,      -1){\line(1,       1){17}}

\put(43,     18){\makebox(0,      0)[t]{$\bullet$}}

\put(18,     12){${-1}$}

\put(36,     22){$-1$}
\put(59,     12){$-q$}

 \ \ \ \ \ \ \ \ \ \ \ \ \ \ \ \ \ \ \ \ \ \ \ \  \ \ \ \ \ \ \ \
  \ \ \ \ \ \ \ \ \ \ \ \ \ \ \ \ \ \ \ {$, q \in R_3$,  by    {\rm  GDD}  $1$ of Row $15$.}
\put(80,         -1)  {    } \end{picture}$\\
\\
 The sub-{\rm GDD} by deleting  Vertex  4    is not an arithmetic  {\rm  GDD}    by Lemma \ref {3.1.1} {\rm (V)}.   \\ \\

  {\ }\\ \\

{\ }\!\!\!\!\!\!\!\!\!\!\!\!\!\!\!\!\!\!\!\!\!{\ }\!\!\!\!\!\!\!\!\!\!\!\!\!\!\!\!\!\!\!\!\!{\ }\!\!\!\!\!\!\!\!\!\!\!\!\!\!\!\!\!\!\!\!\!{\ }\!\!\!\
  $\begin{picture}(100,       15)\put(-68,       -1){ }\put(88,       -1){(e) }

\put(170,     10){\makebox(0,      0)[t]{$\bullet$}}

\put(170,     70){\makebox(0,      0)[t]{$\bullet$}}

\put(230,     10){\makebox(0,      0)[t]{$\bullet$}}
\put(230,     70){\makebox(0,      0)[t]{$\bullet$}}

\put(170,       10){\line(0,      1){60}}

\put(170,       10){\line(1,      1){60}}

\put(230,       10){\line(0,      1){60}}


\put(170,       10){\line(1,       0){60}}
\put(170,       70){\line(1,       0){60}}

\put(150,     10){$q$}
\put(150,     30){$q^{-1}$}

\put(150,     70){$q$}

\put(180,     30){$q^{-1}$}
\put(190,     80){$q^{-1}$}
\put(190,     -10){$-1$}
\put(220,     30){${}$}

\put(250,     10){$-1$}
\put(250,     30){$-q$}

\put(250,     70){$-1$}
\ \ \ \ \ \ \ \ \ \ \ \ \ \
 \ \ \ \ \ \ \ \ \ \ \ \ \ \ \ \ \ \ \ \ \ \ \ \  \ \ \ \ \ \ \ \ \ \ \
 \ \ \ \ \ \ \ \ \ \ \ \ \ \ \ \ \ \ \ \ \ \ \ \ {$, q \in R_3$,  by     {\rm  GDD}  $3$ of Row $15$.}
\put(80,         -1)  {    } \end{picture}$\\ \\
 The sub-{\rm GDD} by deleting  Vertex  2 is not  an arithmetic  {\rm  GDD}      by Lemma \ref {3.1.2} {\rm (II)}.
\\

  {\ }\ \ \ \ \ \ \ \ \ \ \ \ $\begin{picture}(100,       15) \put(-68,        -1){ (f)}

\put(60,       1){\makebox(0,       0)[t]{$\bullet$}}
\put(58,       -12){$-1$}

\put(40,       -12){$q^{-1}$}
\put(28,       -1){\line(1,       0){33}}
\put(27,       1){\makebox(0,      0)[t]{$\bullet$}}

\put(22,      -12){$q$}
\put(0,       -12){$q$}

\put(-14,       1){\makebox(0,      0)[t]{$\bullet$}}

\put(-14,      -1){\line(1,       0){50}}

\put(-18,      -12){$-1$}

\put(59,       0){\line(-1,      1){17}}

\put(28,      -1){\line(1,       1){17}}

\put(43,     18){\makebox(0,      0)[t]{$\bullet$}}

\put(18,     12){${-1}$}

\put(36,     22){$-1$}
\put(59,     12){$-q$}

 \ \ \ \ \ \ \ \ \ \ \ \ \ \ \ \ \ \ \ \ \ \ \ \  \ \ \ \ \ \ \ \
  \ \ \ \ \ \ \ \ \ \ \ \ \ \ \ \ \ \ \ {$, q \in R_3$,  by    {\rm  GDD}  $1$ of Row $15$.}
\put(80,         -1)  {    } \end{picture}$\\
\\
 The sub-{\rm GDD} by deleting  Vertex  4    is not an arithmetic  {\rm  GDD}    by Lemma \ref {3.1.1} {\rm (V)}.   \\ \\ \\ \\ \\ \\ \\

{\ }\!\!\!\!\!\!\!\!\!\!\!\!\!\!\!\!\!\!\!\!\!{\ }\!\!\!\!\!\!\!\!\!\!\!\!\!\!\!\!\!\!\!\!\!{\ }\!\!\!\!\!\!\!\!\!\!\!\!\!\!\!\!\!\!\!\!\!{\ }\!\!\!\
   $\begin{picture}(100,       15)\put(-68,       -1){ }\put(88,       -1){(g) }

\put(170,     10){\makebox(0,      0)[t]{$\bullet$}}

\put(170,     70){\makebox(0,      0)[t]{$\bullet$}}

\put(230,     10){\makebox(0,      0)[t]{$\bullet$}}
\put(230,     70){\makebox(0,      0)[t]{$\bullet$}}

\put(170,       10){\line(0,      1){60}}

\put(170,       10){\line(1,      1){60}}

\put(230,       10){\line(0,      1){60}}


\put(170,       10){\line(1,       0){60}}
\put(170,       70){\line(1,       0){60}}

\put(150,     10){$q$}
\put(150,     30){$q^{-1}$}

\put(150,     70){$q$}

\put(180,     30){$q^{-1}$}
\put(190,     80){$q^{-1}$}
\put(190,     -10){$-1$}
\put(220,     30){${}$}

\put(250,     10){$-1$}
\put(250,     30){$-q$}

\put(250,     70){$-1$}
\ \ \ \ \ \ \ \ \ \ \ \ \ \
 \ \ \ \ \ \ \ \ \ \ \ \ \ \ \ \ \ \ \ \ \ \ \ \  \ \ \ \ \ \ \ \ \ \ \
 \ \ \ \ \ \ \ \ \ \ \ \ \ \ \ \ \ \ \ \ \ \ \ \ {$, q \in R_3$,  by     {\rm  GDD}  $3$ of Row $15$.}
\put(80,         -1)  {    } \end{picture}$\\
\\
 The sub-{\rm GDD} by deleting  Vertex  2    is not an arithmetic  {\rm  GDD}    by Lemma \ref {3.1.2}{\rm (II)}.
\\

  {\ }\ \ \ \ \ \ \ \ \ \ \ \ $\begin{picture}(100,       15) \put(-68,        -1){(h) }

\put(60,       1){\makebox(0,       0)[t]{$\bullet$}}
\put(58,       -12){$-1$}

\put(40,       -12){$q^{-1}$}
\put(28,       -1){\line(1,       0){33}}
\put(27,       1){\makebox(0,      0)[t]{$\bullet$}}

\put(22,      -12){$q$}
\put(-4,       -12){$-q^{-1}$}

\put(-14,       1){\makebox(0,      0)[t]{$\bullet$}}

\put(-14,      -1){\line(1,       0){50}}

\put(-22,      -12){$-q$}

\put(59,       0){\line(-1,      1){17}}

\put(28,      -1){\line(1,       1){17}}

\put(43,     18){\makebox(0,      0)[t]{$\bullet$}}

\put(18,     12){${-1}$}

\put(36,     22){$-1$}
\put(59,     12){$-q$}

 \ \ \ \ \ \ \ \ \ \ \ \ \ \ \ \ \ \ \ \ \ \ \ \  \ \ \ \ \ \ \ \
  \ \ \ \ \ \ \ \ \ \ \ \ \ \ \ \ \ \ \ {$, q \in R_3$,  by    {\rm  GDD}  $1$ of Row $16$.}
\put(80,         -1)  {    } \end{picture}$\\
\\
 The sub-{\rm GDD} by deleting  Vertex  4    is not an arithmetic  {\rm  GDD}    by Lemma \ref {3.1.1} {\rm (V)}.   \\ \\ \\ \\ \\

{\ }\!\!\!\!\!\!\!\!\!\!\!\!\!\!\!\!\!\!\!\!\!{\ }\!\!\!\!\!\!\!\!\!\!\!\!\!\!\!\!\!\!\!\!\!{\ }\!\!\!\!\!\!\!\!\!\!\!\!\!\!\!\!\!\!\!\!\!{\ }\!\!\!\
  $\begin{picture}(100,       15)\put(-68,       -1){ }\put(88,       -1){(i) }

\put(170,     10){\makebox(0,      0)[t]{$\bullet$}}

\put(170,     70){\makebox(0,      0)[t]{$\bullet$}}

\put(230,     10){\makebox(0,      0)[t]{$\bullet$}}
\put(230,     70){\makebox(0,      0)[t]{$\bullet$}}

\put(170,       10){\line(0,      1){60}}

\put(170,       10){\line(1,      1){60}}

\put(230,       10){\line(0,      1){60}}


\put(170,       10){\line(1,       0){60}}
\put(170,       70){\line(1,       0){60}}

\put(150,     10){$q$}
\put(150,     30){${-1}$}

\put(150,     70){$-1$}

\put(180,     30){$q^{-1}$}
\put(190,     80){$-q^{}$}
\put(190,     -10){$-1$}
\put(220,     30){${}$}

\put(250,     10){$-1$}
\put(250,     30){$-q$}

\put(250,     70){$-1$}
\ \ \ \ \ \ \ \ \ \ \ \ \ \
 \ \ \ \ \ \ \ \ \ \ \ \ \ \ \ \ \ \ \ \ \ \ \ \  \ \ \ \ \ \ \ \ \ \
 \ \ \ \ \ \ \ \ \ \ \ \ \ \ \ \ \ \ \ \ \ \ \ \ { $, q \in R_3$,  by    {\rm  GDD}  $3$ of Row $16$.}
\put(80,         -1)  {    } \end{picture}$\\
\\
 The sub-{\rm GDD} by deleting  Vertex  2    is not an arithmetic  {\rm  GDD}     by Lemma \ref {3.1.1} {\rm (V)}.
\\

  {\ }\ \ \ \ \ \ \ \ \ \ \ \ $\begin{picture}(100,       15) \put(-68,        -1){(j) }

\put(60,       1){\makebox(0,       0)[t]{$\bullet$}}

\put(28,       -1){\line(1,       0){33}}
\put(27,       1){\makebox(0,      0)[t]{$\bullet$}}

\put(58,       -12){$-1$}

\put(40,       -12){$q^{-1}$}
\put(22,      -12){$q$}
\put(0,       -12){$-1$}

\put(-18,      -12){$-1$}

\put(-14,       1){\makebox(0,      0)[t]{$\bullet$}}

\put(-14,      -1){\line(1,       0){50}}

\put(59,       0){\line(-1,      1){17}}

\put(28,      -1){\line(1,       1){17}}

\put(43,     18){\makebox(0,      0)[t]{$\bullet$}}

\put(18,     12){${-1}$}

\put(36,     22){$-1$}
\put(59,     12){$-q$}

 \ \ \ \ \ \ \ \ \ \ \ \ \ \ \ \ \ \ \ \ \ \ \ \  \ \ \ \ \ \ \ \
  \ \ \ \ \ \ \ \ \ \ \ \ \ \ \ \ \ \ \ {$, q \in R_3$,  by    {\rm  GDD}  $7$ of Row $17$.}
\put(80,         -1)  {    } \end{picture}$\\
\\
 The sub-{\rm GDD} by deleting  Vertex  4    is not  an arithmetic  {\rm  GDD}    by Lemma \ref {3.1.1} {\rm (V)}. \\

  {\ }\ \ \ \ \ \ \ \ \ \ \ \ $\begin{picture}(100,       15) \put(-68,        -1){(k) }

\put(60,       1){\makebox(0,       0)[t]{$\bullet$}}
\put(58,       -12){$-1$}

\put(40,       -12){$q^{-1}$}
\put(28,       -1){\line(1,       0){33}}
\put(27,       1){\makebox(0,      0)[t]{$\bullet$}}

\put(22,      -12){$q$}
\put(0,       -12){$-q$}

\put(-14,       1){\makebox(0,      0)[t]{$\bullet$}}

\put(-14,      -1){\line(1,       0){50}}

\put(-18,      -12){$-1$}

\put(59,       0){\line(-1,      1){17}}

\put(28,      -1){\line(1,       1){17}}

\put(43,     18){\makebox(0,      0)[t]{$\bullet$}}

\put(18,     12){${-1}$}

\put(36,     22){$-1$}
\put(59,     12){$-q$}

 \ \ \ \ \ \ \ \ \ \ \ \ \ \ \ \ \ \ \ \ \ \ \ \  \ \ \ \ \ \ \ \
  \ \ \ \ \ \ \ \ \ \ \ \ \ \ \ \ \ \ \ {$, q \in R_3$,  by    {\rm  GDD}  $7$ of Row $17$.}
\put(80,         -1)  {    } \end{picture}$\\
\\
 The sub-{\rm GDD} by deleting  Vertex  4    is not  an arithmetic  {\rm  GDD}    by Lemma \ref {3.1.1} {\rm (V)}. \\

{\ }\ \ \ \ \ \ \ \ \ \ \ \ $\begin{picture}(100,       15) \put(-68,        -1){$(l)$ }

\put(60,       1){\makebox(0,       0)[t]{$\bullet$}}

\put(28,       -1){\line(1,       0){33}}
\put(27,       1){\makebox(0,      0)[t]{$\bullet$}}
\put(-14,       1){\makebox(0,      0)[t]{$\bullet$}}

\put(-14,      -1){\line(1,       0){50}}

\put(58,       -12){$-1$}

\put(40,       -12){$q^{-1}$}

\put(22,      -12){$q$}
\put(0,       -12){${q^{-1}}$}

\put(-18,      -12){${q}$}

\put(59,       0){\line(-1,      1){17}}

\put(28,      -1){\line(1,       1){17}}

\put(43,     18){\makebox(0,      0)[t]{$\bullet$}}

\put(18,     12){${-1}$}

\put(36,     22){$-1$}
\put(59,     12){$-q$}

 \ \ \ \ \ \ \ \ \ \ \ \ \ \ \ \ \ \ \ \ \ \ \ \  \ \ \ \ \ \ \ \
 \ \ \ \ \ \ \ \ \ \ \ \ \ \ \ \ \ \ \ { ,  by   {\rm  GDD}  $1$ of Row $4$.}
\put(80,         -1)  {    } \end{picture}$\\
\\
 The sub-{\rm GDD} by deleting  Vertex  4    is not   an arithmetic  {\rm  GDD}    by Lemma \ref {3.1.2}{\rm (II)}.\\ \\

{\ }\ \ \ \ \ \ \ \ \ \ \ \ $\begin{picture}(100,       15) \put(-68,        -1){$(m)$ }

\put(60,       1){\makebox(0,       0)[t]{$\bullet$}}

\put(28,       -1){\line(1,       0){33}}
\put(27,       1){\makebox(0,      0)[t]{$\bullet$}}
\put(-14,       1){\makebox(0,      0)[t]{$\bullet$}}

\put(-14,      -1){\line(1,       0){50}}

\put(58,       -12){$-1$}

\put(40,       -12){$q^{-2}$}

\put(22,      -12){$q^2$}
\put(0,       -12){${q^{-2}}$}

\put(-18,      -12){${q}$}

\put(59,       0){\line(-1,      1){17}}

\put(28,      -1){\line(1,       1){17}}

\put(43,     18){\makebox(0,      0)[t]{$\bullet$}}

\put(18,     12){${-1}$}

\put(36,     22){$-1$}
\put(59,     12){$-q^2$}

 \ \ \ \ \ \ \ \ \ \ \ \ \ \ \ \ \ \ \ \ \ \ \ \  \ \ \ \ \ \ \ \
 \ \ \ \ \ \ \ \ \ \ \ \ \ \ \ \ \ \ \ { , $q\in R_3\cup R_6$,  by   {\rm  GDD}  $1$ of Row $4$.}
\put(80,         -1)  {    } \end{picture}$\\
\\
 The sub-{\rm GDD} by deleting  Vertex  4    is not   an arithmetic  {\rm  GDD}    by Lemma \ref {3.1.1}{\rm (XIII)}.

   {\rm (ii) }  All quasi-arithmetic {\rm  GDD}s by  adding a vertex  on  Vertex $2$  are listed.
 According  to Lemma \ref {3.1.2} (I) we have to consider following cases.

 (a)  i.e.  (16.3.2),  by   {\rm  GDD}  $2$ of Row $4$.
 The sub-{\rm GDD} by deleting  Vertex  4    is  an arithmetic  {\rm  GDD}     by Lemma \ref {3.1.1}{\rm (II)}, \\
  It  is quasi-affine    by Lemma \ref {3.1.3}.
\\

{\ }\ \ \ \ \ \ \ \ \ \ \ \ $\begin{picture}(100,       15) \put(-68,        -1){ (b)}

\put(60,       1){\makebox(0,       0)[t]{$\bullet$}}

\put(28,       -1){\line(1,       0){33}}
\put(27,       1){\makebox(0,      0)[t]{$\bullet$}}
\put(-14,       1){\makebox(0,      0)[t]{$\bullet$}}

\put(-14,      -1){\line(1,       0){50}}

\put(58,       -12){$q^2$}

\put(40,       -12){$q^{-2}$}

\put(-28,      -12){$-q^{-1}$}
\put(0,       -12){$q^{2}$}

\put(18,     12){${-q^2}$}

\put(22,      -12){$-1$}

\put(36,     22){$-1$}
\put(59,     12){$-1$}

 \put(59,       0){\line(-1,      1){17}}

\put(28,      -1){\line(1,       1){17}}

\put(43,     18){\makebox(0,      0)[t]{$\bullet$}}

 \ \ \ \ \ \ \ \ \ \ \ \ \ \ \ \ \ \ \ \ \ \ \ \  \ \ \ \ \ \ \ \
  \ \ \ \ \ \ \ \ \ \ \ \ \ \ \ \ \ \ \ { ,  by   {\rm  GDD}  $3$ of Row $5$.}
\put(80,         -1)  {    } \end{picture}$\\
\\
 The sub-{\rm GDD} by deleting  Vertex  4    is not an arithmetic  {\rm  GDD}    by Lemma \ref {3.1.1}{\rm (I)}.

(c) i.e. (10.3.6), by    {\rm  GDD}  $2$ of Row $6$.
 The sub-{\rm GDD} by deleting  Vertex  4    is not an arithmetic  {\rm  GDD}   $q\in R_3$  by Lemma \ref {3.1.1}{\rm (X)}. The sub-{\rm GDD} by deleting  Vertex  4    is  an arithmetic  {\rm  GDD}   $q\in R_6$  by Lemma \ref {3.1.1}{\rm (X)}. It is  quasi-affine  $q\in R_6$.
\\ \\ \\ \\ \\

{\ }\!\!\!\!\!\!\!\!\!\!\!\!\!\!\!\!\!\!\!\!\!{\ }\!\!\!\!\!\!\!\!\!\!\!\!\!\!\!\!\!\!\!\!\!{\ }\!\!\!\!\!\!\!\!\!\!\!\!\!\!\!\!\!\!\!\!\!{\ }\!\!\!\
  $\begin{picture}(100,       15)\put(88,       -1){ }\put(88,       -1){(d) }

\put(170,     10){\makebox(0,      0)[t]{$\bullet$}}

\put(170,     70){\makebox(0,      0)[t]{$\bullet$}}

\put(230,     10){\makebox(0,      0)[t]{$\bullet$}}
\put(230,     70){\makebox(0,      0)[t]{$\bullet$}}

\put(170,       10){\line(0,      1){60}}

\put(170,       10){\line(1,      1){60}}

\put(230,       10){\line(0,      1){60}}


\put(170,       10){\line(1,       0){60}}
\put(170,       70){\line(1,       0){60}}

\put(150,     10){$-1$}
\put(150,     30){$q^{2}$}

\put(150,     70){$-1$}

\put(180,     30){$q^{-1}$}
\put(190,     80){$q^{-1}$}
\put(190,     -10){$-q$}
\put(220,     30){$$}

\put(250,     10){$-1$}
\put(250,     30){$-1$}

\put(250,     70){$q$}\ \ \ \ \ \ \ \ \ \ \ \ \ \
 \ \ \ \ \ \ \ \ \ \ \ \ \ \ \ \ \ \ \ \ \ \ \ \  \ \ \ \ \ \ \ \ \ \
  \ \ \ \ \ \ \ \ \ \ \ \ \ \ \ \ \ \ \ \ \ \ \ \ { ,  by     {\rm  GDD}  $3$ of Row $6$.}
\put(80,         -1)  {    } \end{picture}$\\
\\
 The sub-{\rm GDD} by deleting  Vertex  1   is not an arithmetic  {\rm  GDD}   by Lemma \ref {3.1.1}{\rm (II)}.
\\

{\ }\ \ \ \ \ \ \ \ \ \ \ \ $\begin{picture}(100,       15) \put(-68,        -1){ (e)}

\put(60,       1){\makebox(0,       0)[t]{$\bullet$}}

\put(28,       -1){\line(1,       0){33}}
\put(27,       1){\makebox(0,      0)[t]{$\bullet$}}
\put(-14,       1){\makebox(0,      0)[t]{$\bullet$}}

\put(-14,      -1){\line(1,       0){50}}

\put(58,       -12){$q^3$}

\put(40,       -12){$q^{-3}$}

\put(22,      -12){$-1$}
\put(0,       -12){${q}$}

\put(-18,      -12){${-1}$}

\put(59,       0){\line(-1,      1){17}}

\put(28,      -1){\line(1,       1){17}}

\put(43,     18){\makebox(0,      0)[t]{$\bullet$}}

\put(18,     12){${-q^3}$}

\put(36,     22){$-1$}
\put(59,     12){$-1$}

 \ \ \ \ \ \ \ \ \ \ \ \ \ \ \ \ \ \ \ \ \ \ \ \  \ \ \ \ \ \ \ \
   \ \ \ \ \ \ \ \ \ \ \ \ \ \ \ \ \ \ \ {, $q \in R_9, $ by    {\rm  GDD}  $2$ of Row $7$.}
\put(80,         -1)  {    } \end{picture}$\\
\\
 The sub-{\rm GDD} by deleting  Vertex  4    is not an arithmetic  {\rm  GDD}    by Lemma \ref {3.1.1}{\rm (X)}.
\\

{\ }\ \ \ \ \ \ \ \ \ \ \ \ $\begin{picture}(100,       15) \put(-68,        -1){ (f)}

\put(60,       1){\makebox(0,       0)[t]{$\bullet$}}

\put(28,       -1){\line(1,       0){33}}
\put(27,       1){\makebox(0,      0)[t]{$\bullet$}}
\put(-14,       1){\makebox(0,      0)[t]{$\bullet$}}

\put(-14,      -1){\line(1,       0){50}}

\put(58,       -12){$q^3$}

\put(40,       -12){$q^{-3}$}

\put(22,      -12){$-1$}
\put(0,       -12){$q^{2}$}

\put(-28,      -12){$-q^{-1}$}

\put(59,       0){\line(-1,      1){17}}

\put(28,      -1){\line(1,       1){17}}

\put(43,     18){\makebox(0,      0)[t]{$\bullet$}}

\put(18,     12){${-q^3}$}

\put(36,     22){$-1$}
\put(59,     12){$-1$}

 \ \ \ \ \ \ \ \ \ \ \ \ \ \ \ \ \ \ \ \ \ \ \ \  \ \ \ \ \ \ \ \
   \ \ \ \ \ \ \ \ \ \ \ \ \ \ \ \ \ \ \ {, $q \in R_9, $ by    {\rm  GDD}  $4$ of Row $7$.}
\put(80,         -1)  {    } \end{picture}$\\
\\
 The sub-{\rm GDD} by deleting  Vertex  4    is not an arithmetic  {\rm  GDD}    by Lemma \ref {3.1.1}{\rm (I)}.

 (g)  i.e.  (16.3.4)   by   {\rm  GDD}  $1$ of Row $8$.  The sub-{\rm GDD} by deleting  Vertex  1
is an arithmetic  {\rm  GDD}   by Lemma \ref {3.1.1} {\rm (IX)}. It  is   quasi-affine   by Lemma \ref {3.1.3}.

 (h)  i.e.  (16.3.3)   by  {\rm  GDD}  $1$ of Row $9$.
 The sub-{\rm GDD} by deleting  Vertex  4    is  an arithmetic  {\rm  GDD}   when $(-q)^3 =r$ or  $(-q)^2 =r$ or
  $-q =r$ by Lemma \ref {3.1.1}{\rm (II)} and  Type   7. It  is   quasi-affine   when $(-q)^3 =r$ or  $q^2 =r$ by Lemma \ref {3.1.3}. It  is an arithmetic  {\rm  GDD}   when $-q =r$ by Lemma \ref {3.1.3}.\\

{\ }\ \ \ \ \ \ \ \ \ \ \ \ $\begin{picture}(100,       15) \put(-68,        -1){ (i)}

\put(60,       1){\makebox(0,       0)[t]{$\bullet$}}

\put(28,       -1){\line(1,       0){33}}
\put(27,       1){\makebox(0,      0)[t]{$\bullet$}}
\put(-14,       1){\makebox(0,      0)[t]{$\bullet$}}

\put(-14,      -1){\line(1,       0){50}}

\put(58,       -12){$q$}

\put(40,       -12){$q^{-1}$}

\put(22,      -12){$-1$}
\put(0,       -12){$q^{-1}$}

\put(-18,      -12){$q{}$}

\put(59,       0){\line(-1,      1){17}}

\put(28,      -1){\line(1,       1){17}}

\put(43,     18){\makebox(0,      0)[t]{$\bullet$}}

\put(18,     12){${-q}$}

\put(36,     22){$-1$}
\put(59,     12){$-1$}

 \ \ \ \ \ \ \ \ \ \ \ \ \ \ \ \ \ \ \ \ \ \ \ \  \ \ \ \ \ \ \ \
   \ \ \ \ \ \ \ \ \ \ \ \ \ \ \ \ \ \ \ {$, q \in R_3$,  by    {\rm  GDD}  $1$ of Row $11$.}
\put(80,         -1)  {    } \end{picture}$\\
\\
  The sub-{\rm GDD} by deleting  Vertex  4    is not an arithmetic  {\rm  GDD}   by Lemma \ref {3.1.1}{\rm (II)}.
\\ \\ \\

 {\ }\\

{\ }\!\!\!\!\!\!\!\!\!\!\!\!\!\!\!\!\!\!\!\!\!{\ }\!\!\!\!\!\!\!\!\!\!\!\!\!\!\!\!\!\!\!\!\!{\ }\!\!\!\!\!\!\!\!\!\!\!\!\!\!\!\!\!\!\!\!\!{\ }\!\!\!\
   $\begin{picture}(100,       15)\put(-68,       -1){ }\put(88,       -1){ (j)}

\put(170,     10){\makebox(0,      0)[t]{$\bullet$}}

\put(170,     70){\makebox(0,      0)[t]{$\bullet$}}

\put(230,     10){\makebox(0,      0)[t]{$\bullet$}}
\put(230,     70){\makebox(0,      0)[t]{$\bullet$}}

\put(170,       10){\line(0,      1){60}}

\put(170,       10){\line(1,      1){60}}

\put(230,       10){\line(0,      1){60}}


\put(170,       10){\line(1,       0){60}}
\put(170,       70){\line(1,       0){60}}

\put(150,     10){$-1$}
\put(150,     30){$q^{-1}$}

\put(150,     70){$q$}

\put(180,     30){$q^{-1}$}
\put(190,     80){$q^{-1}$}
\put(190,     -10){$-q$}
\put(220,     30){$$}

\put(250,     10){$-1$}
\put(250,     30){$-1$}

\put(250,     70){$q$}\ \ \ \ \ \ \ \ \ \ \ \ \ \
 \ \ \ \ \ \ \ \ \ \ \ \ \ \ \ \ \ \ \ \ \ \ \ \  \ \ \ \ \ \ \ \ \ \ \
  \ \ \ \ \ \ \ \ \ \ \ \ \ \ \ \ \ \ \ \ \ \ \ \ { ,  by     {\rm  GDD}  $3$ of Row $15$.}
\put(80,         -1)  {    } \end{picture}$\\
\\
 The sub-{\rm GDD} by deleting  Vertex  1   is not an arithmetic  {\rm  GDD}   by Lemma \ref {3.1.1}{\rm (II)}.\\
\\ \\

  {\ }\\

{\ }\!\!\!\!\!\!\!\!\!\!\!\!\!\!\!\!\!\!\!\!\!{\ }\!\!\!\!\!\!\!\!\!\!\!\!\!\!\!\!\!\!\!\!\!{\ }\!\!\!\!\!\!\!\!\!\!\!\!\!\!\!\!\!\!\!\!\!{\ }\!\!\!\
   $\begin{picture}(100,       15)\put(-68,       -1){ }\put(88,       -1){ (k)}

\put(170,     10){\makebox(0,      0)[t]{$\bullet$}}

\put(170,     70){\makebox(0,      0)[t]{$\bullet$}}

\put(230,     10){\makebox(0,      0)[t]{$\bullet$}}
\put(230,     70){\makebox(0,      0)[t]{$\bullet$}}

\put(170,       10){\line(0,      1){60}}

\put(170,       10){\line(1,      1){60}}

\put(230,       10){\line(0,      1){60}}


\put(170,       10){\line(1,       0){60}}
\put(170,       70){\line(1,       0){60}}

\put(150,     10){$-1$}
\put(150,     30){$-q^{}$}

\put(150,     70){$-1$}

\put(180,     30){$q^{-1}$}
\put(190,     80){$-1$}
\put(190,     -10){$-q$}
\put(220,     30){$$}

\put(250,     10){$-1$}
\put(250,     30){$-1$}

\put(250,     70){$q$}\ \ \ \ \ \ \ \ \ \ \ \ \ \
 \ \ \ \ \ \ \ \ \ \ \ \ \ \ \ \ \ \ \ \ \ \ \ \  \ \ \ \ \ \ \ \ \ \ \
  \ \ \ \ \ \ \ \ \ \ \ \ \ \ \ \ \ \ \ \ \ \ \ \ { ,  by     {\rm  GDD}  $3$ of Row $16$.}
\put(80,         -1)  {    } \end{picture}$\\
\\
 The sub-{\rm GDD} by deleting  Vertex  1   is not an arithmetic  {\rm  GDD}   by Lemma \ref {3.1.1}{\rm (V)}.\\

{\ }\ \ \ \ \ \ \ \ \ \ \ \ $\begin{picture}(100,       15) \put(-68,        -1){(l) }

\put(60,       1){\makebox(0,       0)[t]{$\bullet$}}

\put(28,       -1){\line(1,       0){33}}
\put(27,       1){\makebox(0,      0)[t]{$\bullet$}}
\put(-14,       1){\makebox(0,      0)[t]{$\bullet$}}

\put(-14,      -1){\line(1,       0){50}}

\put(58,       -12){$-q$}

\put(35,       -12){$-q^{-1}$}

\put(22,      -12){$-1$}
\put(0,       -12){${-1}$}

\put(-18,      -12){${q}$}

\put(59,       0){\line(-1,      1){17}}

\put(28,      -1){\line(1,       1){17}}

\put(43,     18){\makebox(0,      0)[t]{$\bullet$}}

\put(18,     12){${q}$}

\put(36,     22){$-1$}
\put(59,     12){$-1$}

 \ \ \ \ \ \ \ \ \ \ \ \ \ \ \ \ \ \ \ \ \ \ \ \  \ \ \ \ \ \ \ \
   \ \ \ \ \ \ \ \ \ \ \ \ \ \ \ \ \ \ \ {$, q \in R_3$,  by   {\rm  GDD}  $4$ of Row $16$.}
\put(80,         -1)  {    } \end{picture}$\\
\\
 The sub-{\rm GDD} by deleting  Vertex 4   is not an arithmetic  {\rm  GDD}   by Lemma \ref {3.1.1}{\rm (I)}.\\

{\ }\ \ \ \ \ \ \ \ \ \ \ \ $\begin{picture}(100,       15) \put(-68,        -1){(m) }

\put(60,       1){\makebox(0,       0)[t]{$\bullet$}}

\put(28,       -1){\line(1,       0){33}}
\put(27,       1){\makebox(0,      0)[t]{$\bullet$}}
\put(-14,       1){\makebox(0,      0)[t]{$\bullet$}}

\put(-14,      -1){\line(1,       0){50}}

\put(58,       -12){$q^2$}

\put(35,       -12){$q^{-2}$}

\put(16,      -12){$-1$}
\put(0,       -12){${q}$}

\put(-18,      -12){${-1}$}

\put(59,       0){\line(-1,      1){17}}

\put(28,      -1){\line(1,       1){17}}

\put(43,     18){\makebox(0,      0)[t]{$\bullet$}}

\put(18,     12){$-q^2$}

\put(36,     22){$-1$}
\put(59,     12){$-1$}

 \ \ \ \ \ \ \ \ \ \ \ \ \ \ \ \ \ \ \ \ \ \ \ \  \ \ \ \ \ \ \ \
   \ \ \ \ \ \ \ \ \ \ \ \ \ \ \ \ \ \ \ {$, q^2 \in R_3$,  by   {\rm  GDD}  $2$ of Row $6$.}
\put(80,         -1)  {    } \end{picture}$\\  \\
  The sub-{\rm GDD} by deleting  Vertex 4   is  an arithmetic  {\rm  GDD} when $ q \in R_6$  by Lemma \ref {3.1.1}{\rm (X)}.

  {\rm (iii) }  All quasi-arithmetic {\rm  GDD}s by  adding a vertex  on  Vertex $3$  are listed.

Checking  step by step we have to consider following cases.

 (a)  i.e.  (16.3.5)
,  by    {\rm  GDD}  $4$ of Row $16$. The sub-{\rm GDD} by deleting  Vertex  4
  is an arithmetic  {\rm  GDD}    by  Type   7.
  It  is quasi-affine    by Lemma \ref {3.1.3}.
\\ \\ \\ \\ \\

{\ }\!\!\!\!\!\!\!\!\!\!\!\!\!\!\!\!\!\!\!\!\!{\ }\!\!\!\!\!\!\!\!\!\!\!\!\!\!\!\!\!\!\!\!\!{\ }\!\!\!\!\!\!\!\!\!\!\!\!\!\!\!\!\!\!\!\!\!{\ }\!\!\!\
   $\begin{picture}(100,       15)\put(-68,       -1){ }\put(88,       -1){ (b)}

\put(170,     10){\makebox(0,      0)[t]{$\bullet$}}

\put(170,     70){\makebox(0,      0)[t]{$\bullet$}}

\put(230,     10){\makebox(0,      0)[t]{$\bullet$}}
\put(230,     70){\makebox(0,      0)[t]{$\bullet$}}

\put(170,       10){\line(0,      1){60}}

\put(170,       10){\line(1,      1){60}}

\put(230,       10){\line(0,      1){60}}


\put(170,       10){\line(1,       0){60}}
\put(170,       70){\line(1,       0){60}}

\put(150,     10){$-1$}
\put(150,     30){$-q$}

\put(150,     70){$-1$}

\put(180,     30){${-1}$}
\put(190,     80){$q^{-1}$}
\put(190,     -10){$-q$}
\put(220,     30){$$}

\put(250,     10){$-1$}
\put(250,     30){$q^{-1}$}

\put(250,     70){$q$}\ \ \ \ \ \ \ \ \ \ \ \ \ \
 \ \ \ \ \ \ \ \ \ \ \ \ \ \ \ \ \ \ \ \ \ \ \ \  \ \ \ \ \ \ \ \ \ \ \
 \ \ \ \ \ \ \ \ \ \ \ \ \ \ \ \ \ \ \ \ \ \ \ \ { ,  by     {\rm  GDD}  $3$ of Row $16$.}
\put(80,         -1)  {    } \end{picture}$\\
\\
 The sub-{\rm GDD} by deleting  Vertex  1   is not an arithmetic  {\rm  GDD}   by Lemma \ref {3.1.1}{\rm (II)}.

  {\rm (iv) } All quasi-affine  {\rm  GDD}s which are complete diagrams are listed.
\\ \\ \\ \\ \\

{\ }\!\!\!\!\!\!\!\!\!\!\!\!\!\!\!\!\!\!\!\!\!{\ }\!\!\!\!\!\!\!\!\!\!\!\!\!\!\!\!\!\!\!\!\!{\ }\!\!\!\!\!\!\!\!\!
  $\begin{picture}(100,       15)\put(-68,       -1){ }\put(68,       -1){ (b) in  Case {\rm (i) }}

\put(170,     10){\makebox(0,      0)[t]{$\bullet$}}

\put(170,     70){\makebox(0,      0)[t]{$\bullet$}}

\put(230,     10){\makebox(0,      0)[t]{$\bullet$}}
\put(230,     70){\makebox(0,      0)[t]{$\bullet$}}

\put(170,       10){\line(0,      1){60}}

\put(170,       10){\line(1,      1){60}}

\put(230,       10){\line(0,      1){60}}

\put(230,       10){\line(-1,      1){60}}

\put(170,       10){\line(1,       0){60}}
\put(170,       70){\line(1,       0){60}}

\put(150,     10){$q$}
\put(150,     30){$q^{-1}$}

\put(150,     70){$-1$}

\put(180,     30){$q^{-1}$}
\put(190,     80){$q^{2}$}
\put(190,     -10){$-1$}
\put(220,     30){$$}

\put(250,     10){$-1$}
\put(250,     30){$-q$}

\put(250,     70){$-1$}

\put(270,         -1)  { with $q\in R_3.$ If the sub-{\rm GDD}   } \end{picture}$\\
\\ by deleting  Vertex 2 is  an arithmetic  {\rm  GDD}   then  $\widetilde{q}_{34} =-q$ by Lemma \ref {3.1.1}{\rm (I)}.
 The sub-{\rm GDD} by deleting  Vertex 1  is not an arithmetic  {\rm  GDD}   since  it is not   {\rm  GDD}  $3$ of Row $9$ and  by  Type   6.
\\ \\ \\ \\ \\ \\ \\ \\

{\ }\!\!\!\!\!\!\!\!\!\!\!\!\!\!\!\!\!\!\!\!\!{\ }\!\!\!\!\!\!\!\!\!\!\!\!\!\!\!\!\!\!\!\!\!{\ }\!\!\!\!\!\!\!\!
   $\begin{picture}(100,       15)\put(-68,       -1){ }\put(68,       -1){(g) in  Case {\rm (i) } }

\put(170,     10){\makebox(0,      0)[t]{$\bullet$}}

\put(170,     70){\makebox(0,      0)[t]{$\bullet$}}

\put(230,     10){\makebox(0,      0)[t]{$\bullet$}}
\put(230,     70){\makebox(0,      0)[t]{$\bullet$}}

\put(170,       10){\line(0,      1){60}}

\put(170,       10){\line(1,      1){60}}

\put(230,       10){\line(0,      1){60}}

\put(230,       10){\line(-1,      1){60}}

\put(170,       10){\line(1,       0){60}}
\put(170,       70){\line(1,       0){60}}

\put(150,     10){$q$}
\put(150,     30){$q^{-1}$}

\put(150,     70){$q$}

\put(180,     30){$q^{-1}$}
\put(190,     80){$q^{-1}$}
\put(190,     -10){$-1$}
\put(220,     30){${}$}

\put(250,     10){$-1$}
\put(250,     30){$-q$}

\put(250,     70){$-1$}

\put(80,         -1)  {    } \end{picture}$\\
\\
The sub-{\rm GDD} by deleting  Vertex 2 is not an arithmetic  {\rm  GDD}   by Lemma \ref {3.1.1}{\rm (I)}.
\\ \\ \\

  {\ }\\

{\ }\!\!\!\!\!\!\!\!\!\!\!\!\!\!\!\!\!\!\!\!\!{\ }\!\!\!\!\!\!\!\!\!\!\!\!\!\!\!\!\!\!\!\!\!{\ }\!\!\!\!\!\!\!\!\!\!\!\!\!
 $\begin{picture}(100,       15)\put(-68,       -1){ }\put(68,       -1){ (i ) in  Case {\rm (i) }}

\put(170,     10){\makebox(0,      0)[t]{$\bullet$}}

\put(170,     70){\makebox(0,      0)[t]{$\bullet$}}

\put(230,     10){\makebox(0,      0)[t]{$\bullet$}}
\put(230,     70){\makebox(0,      0)[t]{$\bullet$}}

\put(170,       10){\line(0,      1){60}}

\put(170,       10){\line(1,      1){60}}

\put(230,       10){\line(0,      1){60}}

\put(230,       10){\line(-1,      1){60}}

\put(170,       10){\line(1,       0){60}}
\put(170,       70){\line(1,       0){60}}

\put(150,     10){$q$}
\put(150,     30){${-1}$}

\put(150,     70){$-1$}

\put(180,     30){$q^{-1}$}
\put(190,     80){${-q}$}
\put(190,     -10){$-1$}
\put(220,     30){$$}

\put(250,     10){$-1$}
\put(250,     30){$-q$}

\put(250,     70){$-1$}

\put(80,         -1)  {    } \end{picture}$\\
\\
The sub-{\rm GDD} by deleting  Vertex 2 is not an arithmetic  {\rm  GDD}   by Lemma \ref {3.1.1}{\rm (I)}.
\\ \\ \\ \\ \\

{\ }\!\!\!\!\!\!\!\!\!\!\!\!\!\!\!\!\!\!\!\!\!{\ }\!\!\!\!\!\!\!\!\!\!\!\!\!\!\!\!\!\!\!\!\!{\ }\!\!\!\!\!\!\!\!\!\!\!\!\!\!\!
  $\begin{picture}(100,       15)\put(-68,       -1){ }\put(68,       -1){(d) in  Case {\rm (ii) } }

\put(170,     10){\makebox(0,      0)[t]{$\bullet$}}

\put(170,     70){\makebox(0,      0)[t]{$\bullet$}}

\put(230,     10){\makebox(0,      0)[t]{$\bullet$}}
\put(230,     70){\makebox(0,      0)[t]{$\bullet$}}

\put(170,       10){\line(0,      1){60}}

\put(170,       10){\line(1,      1){60}}

\put(230,       10){\line(0,      1){60}}

\put(230,       10){\line(-1,      1){60}}

\put(170,       10){\line(1,       0){60}}
\put(170,       70){\line(1,       0){60}}

\put(150,     10){$-1$}
\put(150,     30){$q^{2}$}

\put(150,     70){$-1$}

\put(180,     30){$q^{-1}$}
\put(190,     80){$q^{-1}$}
\put(190,     -10){$-q$}
\put(220,     30){$$}

\put(250,     10){$-1$}
\put(250,     30){$-1$}

\put(250,     70){$q$}

\put(280,         -1)  {  with $q\in R_3.$  If the sub-{\rm GDD} by deleting  Vertex 2 } \end{picture}$\\
\\ is an arithmetic  {\rm  GDD},
then  $\widetilde{q}_{34} = -q$ by Lemma \ref {3.1.1}{\rm (I)}. The sub-{\rm GDD} by deleting  Vertex 1  is not an arithmetic  {\rm  GDD}   checked  step by step.\\ \\ \\ \\ \\ 

 \ \ \ \ \ \  $\begin{picture}(100,     15)\put(-45,      -1){ }\put(-45,      -1){(j) in  Case {\rm (ii) } }

\put(170,    10){\makebox(0,     0)[t]{$\bullet$}}

\put(170,    70){\makebox(0,     0)[t]{$\bullet$}}

\put(230,    10){\makebox(0,     0)[t]{$\bullet$}}
\put(230,    70){\makebox(0,     0)[t]{$\bullet$}}

\put(170,      10){\line(0,     1){60}}

\put(170,      10){\line(1,     1){60}}

\put(230,      10){\line(0,     1){60}}

\put(230,      10){\line(-1,     1){60}}

\put(170,      10){\line(1,      0){60}}
\put(170,      70){\line(1,      0){60}}

\put(150,    10){$-1$}
\put(150,    30){$q^{-1}$}

\put(150,    70){$q$}

\put(180,    30){$q^{-1}$}
\put(190,    80){$q^{-1}$}
\put(190,    -10){$-q$}
\put(220,    30){$$}

\put(250,    10){$-1$}
\put(250,    30){$-1$}

\put(250,    70){$q$}

\put(80,        -1)  {    } \end{picture}$\\
\\
 The sub-{\rm GDD} by deleting Vertex 2 is not an arithmetic {\rm  GDD}  by Lemma \ref {3.1.1}{\rm (I)}.
 \\ \\ \\ \\ \\

  \ \ \ \ \ \  $\begin{picture}(100,     15)\put(-45,      -1){ }\put(-45,      -1){(k) in  Case {\rm (ii) } }

\put(170,    10){\makebox(0,     0)[t]{$\bullet$}}

\put(170,    70){\makebox(0,     0)[t]{$\bullet$}}

\put(230,    10){\makebox(0,     0)[t]{$\bullet$}}
\put(230,    70){\makebox(0,     0)[t]{$\bullet$}}

\put(170,      10){\line(0,     1){60}}

\put(170,      10){\line(1,     1){60}}

\put(230,      10){\line(0,     1){60}}

\put(230,      10){\line(-1,     1){60}}

\put(170,      10){\line(1,      0){60}}
\put(170,      70){\line(1,      0){60}}

\put(150,    10){$-1$}
\put(150,    30){$-q^{}$}

\put(150,    70){$-1$}

\put(180,    30){$q^{-1}$}
\put(190,    80){${-1}$}
\put(190,    -10){$-q$}
\put(220,    30){$$}

\put(250,    10){$-1$}
\put(250,    30){$-1$}

\put(250,    70){$q$}

\put(80,        -1)  {    } \end{picture}$\\
\\
The sub-{\rm GDD} by deleting Vertex 2 is not an arithmetic {\rm  GDD}  by Lemma \ref {3.1.1}{\rm (I)}.
 \\ \\ \\ \\ \\

  \ \ \ \ \ \  $\begin{picture}(100,     15)\put(-45,      -1){ }\put(-45,      -1){ (b) in  Case {\rm (iii) }}

\put(170,    10){\makebox(0,     0)[t]{$\bullet$}}

\put(170,    70){\makebox(0,     0)[t]{$\bullet$}}

\put(230,    10){\makebox(0,     0)[t]{$\bullet$}}
\put(230,    70){\makebox(0,     0)[t]{$\bullet$}}

\put(170,      10){\line(0,     1){60}}

\put(170,      10){\line(1,     1){60}}

\put(230,      10){\line(0,     1){60}}

\put(230,      10){\line(-1,     1){60}}

\put(170,      10){\line(1,      0){60}}
\put(170,      70){\line(1,      0){60}}

\put(150,    10){$-1$}
\put(150,    30){$-q$}

\put(150,    70){$-1$}

\put(180,    30){${-1}$}
\put(190,    80){$q^{-1}$}
\put(190,    -10){$-q$}
\put(220,    30){$$}

\put(250,    10){$-1$}
\put(250,    30){$q^{-1}$}

\put(250,    70){$q$}

\put(80,        -1)  {    } \end{picture}$\\
 \\ with $q\in R_3.$
If The sub-{\rm GDD} by deleting Vertex 1 is   arithmetic  then  $\widetilde{q}_{24} =(-q)^{-2}$ by  Type   6.
If The sub-{\rm GDD} by deleting Vertex 3 is an arithmetic {\rm  GDD},  then $\widetilde{q}_{24} =q^{2}$. Consequently, $ q^{4} =1,$ which implies a contradiction.

  {\rm (iv) } All quasi-affine circles are listed.

\subsection* {Quasi-affine over
  {\rm  GDD}  $4$ of Row $16$   }

  {\rm (ii) } All quasi-arithmetic {\rm  GDD}s by  adding a vertex  on  Vertex $2$  are listed.
 According  to Lemma \ref {3.1.1} (I) we have to consider following cases.

 (a)  i.e.  (16.4.1),   $q=-\xi,  \xi \in R_3 $,  by  {\rm  GDD}  $4$ of Row $16$.  The sub-{\rm GDD} by deleting
 Vertex  4  is an arithmetic  {\rm  GDD}    by Lemma \ref {3.1.1} {\rm (IV)}.
  It  is quasi-affine    by Lemma \ref {3.1.3}.

  {\rm (iii) }  All quasi-arithmetic {\rm  GDD}s by  adding a vertex  on  Vertex $3$  are listed.
 According  to Lemma \ref {3.1.2} (III) we have to consider following cases.

 (a)  i.e.  (16.4.2)
$q=-\xi,  \xi \in R_3 $,  $ q \notin R_3$,  by    {\rm  GDD}  $1$ of Row $7$.
 It is quasi-affine  by  Lemma   \ref {2.63}.

\subsection* {Quasi-affine over
  {\rm  GDD}  $5$ of Row $16$  }

  {\rm (i) }  All quasi-arithmetic {\rm  GDD}s by  adding a vertex  on  Vertex $1$  are listed.

Checking  step by step we have to consider following cases.

 (a)  i.e.  (16.5.1), $q\in R_3, $ by   {\rm  GDD}  $1$ of Row $16$   in Table A2.
is quasi-affine   by Lemma \ref {3.2.33}

  {\rm (ii) }  All quasi-arithmetic {\rm  GDD}s by  adding a vertex  on  Vertex $2$  are listed.

Checking  step by step we have to consider following cases.

 (a)  i.e.  (16.5.2)
 by    {\rm  GDD}  $5$ of Row $16$. The sub-{\rm GDD} by deleting  Vertex  4
is  an arithmetic  {\rm  GDD}     by  {\rm  GDD}  $5$ of Row $16$. It  is quasi-affine  by Lemma \ref {3.1.3}.

  {\rm (iii) }  All quasi-arithmetic {\rm  GDD}s by  adding a vertex  on  Vertex $3$  are listed.
 According  to Lemma \ref {3.1.2} (II) we have to consider following cases.

 (a)  i.e.  (16.5.3)$, q \in  R_6$.   {\rm  GDD}  $1$ of Row $13$.
 It is quasi-affine  by  Lemma   \ref {2.63}.

 (b)  i.e.  (16.5.4),  $q \in R_3.$   {\rm  GDD}  $5$ of Row $16$.
 It is quasi-affine
 by  Lemma   \ref {2.63}.

  {\rm (iv) } All quasi-affine circles are listed.

 (a) (nc)  i.e.  (16.5.5)   is  quasi-affine  since the sub-{\rm GDD} by deleting  Vertex 2
   is  an arithmetic  {\rm  GDD}   by Lemma \ref {3.1.1} {\rm (IV)}.

 (a)  (a)  i.e.  (16.5.6)
 is  quasi-affine  since the sub-{\rm GDD} by deleting  Vertex 2
   is  an arithmetic  {\rm  GDD}   by  Type   3.

\subsection* {Quasi-affine over
  {\rm  GDD}  $1$ of Row $17$  }   {\rm (i) } All quasi-arithmetic {\rm  GDD}s by  adding a vertex  on  Vertex $1$  are listed.
 According  to Lemma \ref {3.1.2} (V) we have to consider following cases.

 (a)  i.e.  (17.1.2)
, $ r\in R_4, $  by   {\rm  GDD}  $2$ of Row $6$;  $ r\in R_3, $  {\rm  GDD}  $1$ of Row $17$;
$ r\in R_6, $  {\rm  GDD}  $2$ of Row $7$. It  is quasi-affine   by Lemma \ref {3.1.4}.

 (b)  i.e.  (17.1.1),   $q,  r \in R_3, $ by   {\rm  GDD}  $9$ of Row $17$  in Table A2.
is quasi-affine   by Lemma \ref {3.1.3}.

 (c)  i.e.  (17.1.3),  $r \in R_6, $ by   {\rm  GDD}  $4$ of Row $7$  in Table A2.
is quasi-affine   by Lemma \ref {3.1.4}.

 (d)  i.e.  (17.1.4),   $r \not=-1, $ by   {\rm  GDD}  $1$ of Row $9$  in Table A2 and  by Lemma \ref {3.1.1} {\rm (IV)}.
is quasi-affine   by Lemma \ref {3.1.4}.

  {\rm (ii) }  All quasi-arithmetic {\rm  GDD}s by  adding a vertex  on  Vertex $2$  are listed.
 According  to Lemma \ref {3.1.2} (V) we have to consider following cases.

 (a)  i.e.  (17.1.5), $ r\in R_4, $  {\rm  GDD}  $2$ of Row $6$;  $ r\in R_3, $ by   {\rm  GDD}  $1$ of Row $17$;
 $ r\in R_6, $ by    {\rm  GDD}  $2$ of Row $7$;
$ r\in R_2, $  {\rm  GDD}  $1$ of Row $1$;
 The sub-{\rm GDD} by deleting  Vertex  4   is an arithmetic  {\rm  GDD}   when
 $r=q$ and  $q\in R_3$,  or $r=q^{-1}$,  or $r=-q$ and $r \in R_6$  by Lemma \ref {3.1.1} {\rm (X)} or $r=-1$.
  It  is quasi-affine   by Lemma \ref {3.1.10}.

 (b)  i.e.  (17.1.6),  $r \in R_6, $ by   {\rm  GDD}  $4$ of Row $7$  in Table A2.
 The sub-{\rm GDD} by deleting  Vertex  4   is  an arithmetic  {\rm  GDD}   when  $ q^{} =r^{-2}$ by  Type   2.
 The sub-{\rm GDD} by deleting  Vertex  4   is not an arithmetic  {\rm  GDD}   when  $ q^{} \not=r^{-2}$ by Lemma \ref {3.1.1}{\rm (I)}.
  It  is quasi-affine   when  $ q^{} =r^{-2}$,    $r \in R_6$ and  $q \in R_3, $ by Lemma \ref {3.1.36}.

 (c)  i.e.  (17.1.7),   $r^2 \not=1, $ by    {\rm  GDD}  $1$ of Row $9$  in Table A2.
 The sub-{\rm GDD} by deleting  Vertex  4   is  an arithmetic  {\rm  GDD}   when $ r= q$ or $ r =q^{2}$ or $ r =-q$  by Lemma \ref {3.1.1} {\rm (IX)}.
  It  is quasi-affine    by Lemma \ref {3.1.36}.\\ \\

     {\ }\ \ \ \ \ \ \ \ \ \ \ \ \ \ $\begin{picture}(100,       15) \put(-68,       -1){(d) }
\put(27,       1){\makebox(0,      0)[t]{$\bullet$}}
\put(60,       1){\makebox(0,       0)[t]{$\bullet$}}
\put(93,      1){\makebox(0,      0)[t]{$\bullet$}}
\put(28,       -1){\line(1,       0){33}}
\put(61,       -1){\line(1,       0){30}}

\put(22,      -12){$r$}
\put(40,       -12){$-1$}
\put(58,       -12){$-1$}
\put(70,       -12){$-1$}

\put(91,       -12){$-1$}

\put(60,     38){\makebox(0,      0)[t]{$\bullet$}}

\put(60,       0){\line(0,      1){35}}

\put(63,       30){$-1$}

\put(63,       18){$q$}

\put(145,         -1)  { , $q \in R_3, $  $r \in R_4, $ by   {\rm  GDD}  $1$ of Row $2$  in Table A2. }
\put(80,         -1)  {    } \end{picture}$\\
\\
 The sub-{\rm GDD} by deleting  Vertex  4   is not  an arithmetic  {\rm  GDD}   by Lemma \ref {3.1.1}{\rm (I)}.

  {\rm (iii) }  All quasi-arithmetic {\rm  GDD}s by  adding a vertex  on  Vertex $3$  are listed.
 According  to Lemma \ref {3.1.2} (IV) we have to consider following cases.

 (a)  i.e.  (17.1.8)
$, q \in R_3$.  by   {\rm  GDD}  $2$ of Row $15$. It is quasi-affine  by  Lemma   \ref {2.63}.

 (b)  i.e.  (17.1.9),  $q \in R_3$. by    {\rm  GDD}  $2$ of Row $16$. It is quasi-affine  by  Lemma   \ref {2.63}.

 (c)  i.e.  (17.1.10),  $q \in R_3$. by    {\rm  GDD}  $1$ of Row $17$. It is quasi-affine  by  Lemma   \ref {3.1.10}.

  {\rm (iv) }  All quasi-affine circles are listed.

 (a) (nc)  i.e.  (17.1.11)  is  quasi-affine  since the sub-{\rm GDD} by deleting  Vertex 2
is  an arithmetic  {\rm  GDD}   when $r= q^{-1}$,  or $r= q^{}$,
or $r= -q^{-1}$ by Lemma \ref {3.1.1} {\rm (X)}.

 (a) (a)  i.e.  (17.1.12)
   is  quasi-affine  since the sub-{\rm GDD} by deleting  Vertex 2
  is  an arithmetic  {\rm  GDD}   when $r= q^{}$,  or $r= q^{-1}$,
or $r= -q^{}$,  or $r= -1^{}$;  $q\in R_3$    by Lemma \ref {3.1.1} {\rm (X)}.

 (a) (c)  i.e.  (17.1.13)
   is  quasi-affine  since the sub-{\rm GDD} by deleting  Vertex 2
  is  an arithmetic  {\rm  GDD}   when $r= q^{}$,  or $r= q^{-1}$,
or $r= -q^{}$,  or $r= -1^{}$;  $q\in R_3$    by Lemma \ref {3.1.1} {\rm (XI)}.

 (b) (nc)  i.e.  (17.1.14)
  is  quasi-affine  since the sub-{\rm GDD} by deleting  Vertex 2
 and $q =r^{-1} $ is an arithmetic  {\rm  GDD}   by     {\rm  GDD}  $7$ of Row $17$.

 (a) (nc)  i.e.  (17.1.15)   is  quasi-affine  since the sub-{\rm GDD} by deleting  Vertex 2
with  $q =r^{-2} $ is an arithmetic  {\rm  GDD}   by Lemma \ref {3.1.1} {\rm (V)}.

 (b) (nc)  i.e.  (17.1.16)
    is  quasi-affine  since the sub-{\rm GDD} by deleting  Vertex 2
with  $q^2= r^{-1}$
is an arithmetic  {\rm  GDD}   by  Lemma \ref {3.1.2} {\rm (III)}.

 (c) (nc)  i.e.  (17.1.17)
  is  quasi-affine  since the sub-{\rm GDD} by deleting  Vertex 2
    with     $q =-r^{-1} $ is  an arithmetic  {\rm  GDD}   by Lemma \ref {3.1.2} (III).

 (c) (nc)  i.e.  (17.1.18),
   is  quasi-affine  since the sub-{\rm GDD} by deleting  Vertex 2   is  an arithmetic  {\rm  GDD}   by Lemma \ref {3.1.2} (III).

 (d) (nc)    i.e.  (17.1.19)
  is  quasi-affine  since the sub-{\rm GDD} by deleting  Vertex 2
   with     $q =r $ is  an arithmetic  {\rm  GDD}.

 (d) (nc)    i.e.  (17.1.20)
 is  quasi-affine  since the sub-{\rm GDD} by deleting  Vertex 2
 with     $q =r^2 $ is  an arithmetic  {\rm  GDD}   by Lemma \ref {3.1.2} {\rm (III)}(h).

 (d) (nc)    i.e.  (17.1.21)   is  quasi-affine  since the sub-{\rm GDD} by deleting  Vertex 2
   with     $r =q^{2} $ is  an arithmetic  {\rm  GDD}   by  Type 1.

 (d) (nc)    i.e.  (17.1.22)  is  quasi-affine  since the sub-{\rm GDD} by deleting  Vertex 2
 with     $r =-q^{} $ is  an arithmetic  {\rm  GDD}   by  Type   7.

\subsection* {Quasi-affine over
  {\rm  GDD}  $2$ of Row $17$  }

  {\rm (i) }  All quasi-arithmetic {\rm  GDD}s by  adding a vertex  on  Vertex $2$  are listed.

Checking  step by step we have to consider following cases.

 (a)    i.e.  (17.2.1), $q\in R_3, $     {\rm  GDD}  $4$ of Row $16$. It is quasi-affine   by Lemma \ref {3.1.8}.

   {\rm (ii) }  All quasi-arithmetic {\rm  GDD}s by  adding a vertex  on  Vertex $2$  are listed.

Checking  step by step we have to consider following cases.

(a)    i.e.  (17.2.2)
,  $q\in R_3$,  by   {\rm  GDD}  $2$ of Row $17$. The sub-{\rm GDD} by deleting
 Vertex  4  is an arithmetic  {\rm  GDD}   by  Type   7.
 It  is quasi-affine   by Lemma \ref {3.1.3}.

  {\rm (iii) }  All quasi-arithmetic {\rm  GDD}s by  adding a vertex  on  Vertex $3$  are listed.
 According  to Lemma \ref {3.1.2} (I) we have to consider following cases.

(a)    i.e.  (17.2.3)$, q \in R_3,  $   by  {\rm  GDD}  $2$ of Row $6$.
 It is quasi-affine  by \cite [Lemma 2.1]{TZ22}.

(b)    i.e.  (17.2.4)$, q \in R_9,  $   by  {\rm  GDD}  $2$ of Row $7$.
 It is quasi-affine  by  by Lemma  \ref {3.1.8}.

(c)    i.e.  (17.2.5),  by Lemma \ref {3.1.1} {\rm (IV)}.
 It is quasi-affine by  Lemma   \ref {2.63}.

  {\rm (iv) }  All quasi-affine circles are listed.

(a) (nc)    i.e.  (17.2.6) is  quasi-affine  since the sub-{\rm GDD} by deleting  Vertex 2
    with  $-q^{-1} =-q^{2} $ is an arithmetic  {\rm  GDD}    by Lemma \ref {3.1.2} (III).

(a) (c)    i.e.  (17.2.7)
  is  quasi-affine  since the sub-{\rm GDD} by deleting  Vertex 2
  with  $r =-q^{} $ is  an arithmetic  {\rm  GDD}   by  Type   7.

\subsection* {Quasi-affine over
  {\rm  GDD}  $3$ of Row $17$  }

  {\rm (i) }  All quasi-arithmetic {\rm  GDD}s by  adding a vertex  on  Vertex $1$  are listed.

Checking  step by step we have to consider following cases.

(a)     i.e.  (17.3.1),   $q\in R_3, $ by   {\rm  GDD}  $5$ of Row $16$ in Table A2. The sub-{\rm GDD}
by deleting  Vertex  4   is an arithmetic  {\rm  GDD}   by  Type   7. It  is   quasi-affine   by Lemma \ref {3.1.3}.

(b)    i.e.  (17.3.2),   $q\in R_3, $ by   {\rm  GDD}  $3$ of Row $17$.  The sub-{\rm GDD} by deleting  Vertex  1   is
 an arithmetic  {\rm  GDD}.    The sub-{\rm GDD} by deleting  Vertex  2   is an arithmetic  {\rm  GDD}.    It  is   quasi-affine.

  {\rm (ii) }  All quasi-arithmetic {\rm  GDD}s by  adding a vertex  on  Vertex $2$  are listed.
Checking  step by step we have to consider following cases.

  {\rm (iii) } All quasi-arithmetic {\rm  GDD}s by  adding a vertex  on  Vertex $3$  are listed.
 According  to Lemma \ref {3.1.2} (I) we have to consider following cases.

(a)    i.e.  (17.3.3), $q\in R_6 $,  by     {\rm  GDD}  $2$ of Row $4$. The sub-{\rm GDD} by deleting  Vertex  4   is
an arithmetic  {\rm  GDD}   by Lemma \ref {3.1.1}{\rm (II)} or  by Lemma \ref {3.1.1} {\rm (IX)}.
  It  is quasi-affine   by Lemma \ref {3.1.3}.
\\

{\ }\ \ \ \ \ \ \ \ \ \ \ \ $\begin{picture}(100,       15) \put(-68,        -1){ (b)}

\put(60,       1){\makebox(0,       0)[t]{$\bullet$}}

\put(28,       -1){\line(1,       0){33}}
\put(27,       1){\makebox(0,      0)[t]{$\bullet$}}
\put(-14,       1){\makebox(0,      0)[t]{$\bullet$}}

\put(-14,      -1){\line(1,       0){50}}

\put(58,       -12){$q^2$}

\put(40,       -12){$q^{-2}$}

\put(22,      -12){$-1$}
\put(0,       -12){$q^{2}$}

\put(-29,      -12){$-q^{-1}$}

 \put(59,       0){\line(-1,      1){17}}

\put(28,      -1){\line(1,       1){17}}

\put(43,     18){\makebox(0,      0)[t]{$\bullet$}}

\put(18,     12){-$q^{-2}$}

\put(36,     22){$-q^2$}
\put(59,     12){$q^{-2}$}

 \ \ \ \ \ \ \ \ \ \ \ \ \ \ \ \ \ \ \ \ \ \ \ \  \ \ \ \ \ \ \ \
 \ \ \ \ \ \ \ \ \ \ \ \ \ \ \ \ \ \ \ { $q^2=-\xi ^2,  \xi ^2\in R_3 $,  by
  {\rm  GDD}  $3$ of Row $5$.}
\put(80,         -1)  {    } \end{picture}$\\ \\
  The sub-{\rm GDD} by deleting  Vertex  4    is not an arithmetic  {\rm  GDD}   by Lemma \ref {3.1.1}{\rm (I)}.
\\ \\

{\ }\ \ \ \ \ \ \ \ \ \ \ \ $\begin{picture}(100,       15) \put(-68,        -1){ (c)}

\put(60,       1){\makebox(0,       0)[t]{$\bullet$}}

\put(28,       -1){\line(1,       0){33}}
\put(27,       1){\makebox(0,      0)[t]{$\bullet$}}
\put(-14,       1){\makebox(0,      0)[t]{$\bullet$}}

\put(-14,      -1){\line(1,       0){50}}

\put(58,       -12){$q^2$}

\put(40,       -12){$q^{-2}$}

\put(22,      -12){$-1$}
\put(0,       -12){${q}$}

\put(-18,      -12){${-1}$}

  \put(59,       0){\line(-1,      1){17}}

\put(28,      -1){\line(1,       1){17}}

\put(43,     18){\makebox(0,      0)[t]{$\bullet$}}

\put(14,     12){$-q^{-2}$}

\put(36,     22){$-q^{2}$}
\put(59,     12){$q^{-2}$}

 \ \ \ \ \ \ \ \ \ \ \ \ \ \ \ \ \ \ \ \ \ \ \ \  \ \ \ \ \ \ \ \
 \ \ \ \ \ \ \ \ \ \ \ \ \ \ \ \ \ \ \ { $q^2=-\xi ^2,  \xi ^2\in R_3 $,  by   {\rm  GDD}  $2$ of Row $6$.}
\put(80,         -1)  {    } \end{picture}$\\
\\
 The sub-{\rm GDD} by deleting  Vertex  4   is not  an arithmetic  {\rm  GDD}    by Lemma \ref {3.1.1}{\rm (II)}.
\\ \\

  {\ }\\  \\

{\ }\!\!\!\!\!\!\!\!\!\!\!\!\!\!\!\!\!\!\!\!\!{\ }\!\!\!\!\!\!\!\!\!\!\!\!\!\!\!\!\!\!\!\!\!{\ }\!\!\!\!\!\!\!\!\!\!\!\!\!\!\!\!\!\!\!\!\!{\ }\!\!\!\!\!\!\!\!\!{\ }\!\!\!\
$\begin{picture}(55,      15)\put(98,       -1){(d) }

\put(170,     10){\makebox(0,      0)[t]{$\bullet$}}

\put(170,     70){\makebox(0,      0)[t]{$\bullet$}}

\put(230,     10){\makebox(0,      0)[t]{$\bullet$}}
\put(230,     70){\makebox(0,      0)[t]{$\bullet$}}

\put(170,       10){\line(0,      1){60}}

\put(170,       10){\line(1,      1){60}}

\put(230,       10){\line(0,      1){60}}


\put(170,       10){\line(1,       0){60}}
\put(170,       70){\line(1,       0){60}}

\put(150,     10){$-1$}
\put(150,     30){$q^{2}$}

\put(150,     70){$-1$}

\put(180,     30){$q^{-1}$}
\put(190,     80){$q^{-1}$}
\put(190,     -10){$-q^{-1}$}
\put(220,     30){$$}

\put(250,     10){$-q$}
\put(250,     30){$q^{-1}$}

\put(250,     70){$q$}

\put(280,         -1)  { $q=-\xi,  \xi \in R_3 $,  by      {\rm  GDD}  $3$ of Row $7$.   } \end{picture}$\\ \\
 The sub-{\rm GDD} by deleting  Vertex  1   is not an  arithmetic  {\rm  GDD}   by Lemma \ref {3.1.1}{\rm (I)}.
\\ \\

{\ }\ \ \ \ \ \ \ \ \ \ \ \ $\begin{picture}(100,       15) \put(-68,        -1){ (e)}

\put(60,       1){\makebox(0,       0)[t]{$\bullet$}}

\put(28,       -1){\line(1,       0){33}}
\put(27,       1){\makebox(0,      0)[t]{$\bullet$}}
\put(-14,       1){\makebox(0,      0)[t]{$\bullet$}}

\put(-14,      -1){\line(1,       0){50}}

\put(58,       -12){$q^3$}

\put(40,       -12){$q^{-3}$}

\put(22,      -12){$-1$}
\put(0,       -12){${q}$}

\put(-18,      -12){${-1}$}

 \put(59,       0){\line(-1,      1){17}}

\put(28,      -1){\line(1,       1){17}}

\put(43,     18){\makebox(0,      0)[t]{$\bullet$}}

\put(18,     12){$-q^{-3}$}

\put(36,     22){$-q^{3}$}
\put(59,     12){$q^{-5}$}

 \ \ \ \ \ \ \ \ \ \ \ \ \ \ \ \ \ \ \ \ \ \ \ \  \ \ \ \ \ \ \ \
  \ \ \ \ \ \ \ \ \ \ \ \ \ \ \ \ \ \ \ {$q^3 \in R_6 $,  by    {\rm  GDD}  $2$ of Row $7$.}
\put(80,         -1)  {    } \end{picture}$\\
\\
 The sub-{\rm GDD} by deleting  Vertex  4   is not an  arithmetic  {\rm  GDD}   by Lemma \ref {3.1.1} {\rm (IX)}.
\\ \\ \\ \\ \\

{\ }\!\!\!\!\!\!\!\!\!\!\!\!\!\!\!\!\!\!\!\!\!{\ }\!\!\!\!\!\!\!\!\!\!\!\!\!\!\!\!\!\!\!\!\!{\ }\!\!\!\!\!\!\!\!\!\!\!\!\!\!\!\!\!\!\!\!\!{\ }\!\!\!\!\! $\begin{picture}(100,       15)\put(-68,       -1){ }\put(88,       -1){ (f)}

\put(170,     10){\makebox(0,      0)[t]{$\bullet$}}

\put(170,     70){\makebox(0,      0)[t]{$\bullet$}}

\put(230,     10){\makebox(0,      0)[t]{$\bullet$}}
\put(230,     70){\makebox(0,      0)[t]{$\bullet$}}

\put(170,       10){\line(0,      1){60}}

\put(170,       10){\line(1,      1){60}}

\put(230,       10){\line(0,      1){60}}


\put(170,       10){\line(1,       0){60}}
\put(170,       70){\line(1,       0){60}}

\put(150,     10){$-1$}
\put(150,     30){$q^{3}$}

\put(150,     70){$-1$}

\put(180,     30){$q^{-1}$}
\put(190,     80){$q^{-2}$}
\put(190,     -10){$-q^{-1}$}
\put(220,     30){$$}

\put(250,     10){$-q$}
\put(250,     30){$q^{-1}$}

\put(250,     70){$q$}

\ \ \ \ \ \ \ \ \ \ \ \ \ \  \ \ \ \ \ \ \ \ \ \ \ \ \ \ \ \ \ \ \ \ \ \ \ \  \ \ \ \ \ \ \ \ \ \ \ \
 \ \ \ \ \ \ \ \ \ \ \ \ \ \ \ \ \ \ \ \ \ \ \ \ { $q=-\xi,  \xi \in R_3 $,   by     {\rm  GDD}  $3$ of Row $7$.}
\put(80,         -1)  {    } \end{picture}$\\
\\
 The sub-{\rm GDD} by deleting  Vertex  1   is not an arithmetic  {\rm  GDD}   by Lemma \ref {3.1.1}{\rm (I)}.

(g)    i.e.  (17.3.4),  $q\in R_6 $,  by    {\rm  GDD}  $1$ of Row $8$. The sub-{\rm GDD} by
deleting  Vertex  4   is an arithmetic  {\rm  GDD}    by Lemma \ref {3.1.1} {\rm (IV)}. It  is   quasi-affine   by Lemma \ref {3.1.3}.

(h)    i.e.  (17.3.5),  $q\in R_6 $,  by    {\rm  GDD}  $1$ of Row $9$. The sub-{\rm GDD} by
 deleting  Vertex  4    is  an arithmetic  {\rm  GDD}    by Lemma \ref {3.1.1} {\rm (IV)}.
  is quasi-affine   by Lemma \ref {3.1.3}.

(i)    i.e.  (17.3.6), $q\in R_6 $,   by   {\rm  GDD}  $2$ of Row $10$. The sub-{\rm GDD} by deleting  Vertex  4   is an arithmetic  {\rm  GDD}   by Lemma \ref {3.1.1} {\rm (IV)}. It  is   quasi-affine   by Lemma \ref {3.1.3}.

(j)    i.e.  (17.3.7), $q=-\xi,  \xi \in R_3 $,  by     {\rm  GDD}  $1$ of Row $10$. The sub-{\rm GDD} by
deleting  Vertex  4   is an arithmetic  {\rm  GDD}   by Lemma \ref {3.1.1} {\rm (IV)}. It  is   quasi-affine   by Lemma \ref {3.1.3}.

(k)    i.e.  (17.3.8), $q \in R_6 $,   by    {\rm  GDD}  $2$ of Row $10$. The sub-{\rm GDD} by
deleting  Vertex  4   is an arithmetic  {\rm  GDD}   by Lemma \ref {3.1.1} {\rm (IV)}. It  is   quasi-affine   by Lemma \ref {3.1.3}.\\

{\ }\ \ \ \ \ \ \ \ \ \ \ \ $\begin{picture}(100,       15) \put(-68,        -1){ (l)}

\put(60,       1){\makebox(0,       0)[t]{$\bullet$}}

\put(28,       -1){\line(1,       0){33}}
\put(27,       1){\makebox(0,      0)[t]{$\bullet$}}
\put(-14,       1){\makebox(0,      0)[t]{$\bullet$}}

\put(-14,      -1){\line(1,       0){50}}

\put(58,       -12){$-q^{-1}$}

\put(40,       -12){$-q$}

\put(22,      -12){$-1$}
\put(0,       -12){$-q^{-1}{}$}

\put(-18,      -12){$q^{-1}$}

 \put(59,       0){\line(-1,      1){17}}

\put(28,      -1){\line(1,       1){17}}

\put(43,     18){\makebox(0,      0)[t]{$\bullet$}}

\put(18,     12){$q^{}$}

\put(36,     22){$q^{-1}$}
\put(59,     12){$-q$}

 \ \ \ \ \ \ \ \ \ \ \ \ \ \ \ \ \ \ \ \ \ \ \ \  \ \ \ \ \ \ \ \
  \ \ \ \ \ \ \ \ \ \ \ \ \ \ \ \ \ \ \ {$, q \in R_3$.  by   {\rm  GDD}  $3$ of Row $14$.}
\put(80,         -1)  {    } \end{picture}$\\
\\
 The sub-{\rm GDD} by deleting  Vertex  4   is not an arithmetic  {\rm  GDD}   by Lemma \ref {3.1.1}{\rm (I)}.

(m)    i.e.  (17.3.9)$, q \in R_3$,  by    {\rm  GDD}  $2$ of Row $16$. The sub-{\rm GDD}
by deleting  Vertex  4   is an arithmetic  {\rm  GDD}    Type   7. It  is   quasi-affine   by Lemma \ref {3.1.3}.\\

{\ }\ \ \ \ \ \ \ \ \ \ \ \ $\begin{picture}(100,       15) \put(-68,        -1){(n) }

\put(60,       1){\makebox(0,       0)[t]{$\bullet$}}

\put(28,       -1){\line(1,       0){33}}
\put(27,       1){\makebox(0,      0)[t]{$\bullet$}}
\put(-14,       1){\makebox(0,      0)[t]{$\bullet$}}

\put(-14,      -1){\line(1,       0){50}}

\put(58,       -12){$-q$}

\put(35,       -12){$-q^{-1}$}

\put(22,      -12){$-1$}
\put(0,       -12){${-1}$}

\put(-18,      -12){${q}$}

 \put(59,       0){\line(-1,      1){17}}

\put(28,      -1){\line(1,       1){17}}

\put(43,     18){\makebox(0,      0)[t]{$\bullet$}}

\put(18,     12){$q^{-1}$}

\put(36,     22){$q$}
\put(59,     12){$-q^{-1}$}

 \ \ \ \ \ \ \ \ \ \ \ \ \ \ \ \ \ \ \ \ \ \ \ \  \ \ \ \ \ \ \ \
  \ \ \ \ \ \ \ \ \ \ \ \ \ \ \ \ \ \ \ {$, q \in R_3$,  by    {\rm  GDD}  $4$ of Row $16$.}
\put(80,         -1)  {    } \end{picture}$\\
\\
 The sub-{\rm GDD} by deleting  Vertex  4   is not  an arithmetic  {\rm  GDD}   by Lemma \ref {3.1.1}{\rm (I)}.

  {\rm (iv) } All quasi-affine  {\rm  GDD}s which are complete diagrams are listed.
\\ \\ \\ \\ \\

{\ }\!\!\!\!\!\!\!\!\!\!\!\!\!\!\!\!\!\!\!\!\!{\ }\!\!\!\!\!\!\!\!\!\!\!\!\!\!\!\!\!\!\!\!\!{\ }\!\!\!\!\!\! $\begin{picture}(100,       15)\put(-68,       -1){ }\put(68,       -1){(b) in  Case {\rm (i) } }

\put(170,     10){\makebox(0,      0)[t]{$\bullet$}}

\put(170,     70){\makebox(0,      0)[t]{$\bullet$}}

\put(230,     10){\makebox(0,      0)[t]{$\bullet$}}
\put(230,     70){\makebox(0,      0)[t]{$\bullet$}}

\put(170,       10){\line(0,      1){60}}

\put(170,       10){\line(1,      1){60}}

\put(230,       10){\line(0,      1){60}}

\put(230,       10){\line(-1,      1){60}}

\put(170,       10){\line(1,       0){60}}
\put(170,       70){\line(1,       0){60}}

\put(150,     10){$-q$}
\put(150,     30){$-q^{-1}$}

\put(150,     70){$-1$}

\put(180,     30){$-q^{-1}$}
\put(190,     80){$q^{-1}$}
\put(190,     -10){$-q^{-1}$}
\put(220,     30){$$}

\put(250,     10){$-q$}
\put(250,     30){$q^{-1}$}

\put(250,     70){$q$}

\put(280,         -1)  { } \end{picture}$\\ \\
 The sub-{\rm GDD} by deleting  Vertex 1  is not  an arithmetic  {\rm  GDD}    by Lemma \ref {3.1.1}{\rm (I)}. \\ \\ \\ \\ \\

 {\ }\!\!\!\!\!\!\!\!\!\!\!\!\!\!\!\!\!\!\!\!\!{\ }\!\!\!\!\!\!\!\!\!\!\!\!\!\!\!\!\!\!\!\!\!{\ }\!\!\!\!\!\!    $\begin{picture}(100,       15)\put(-68,       -1){ }\put(68,       -1){ (f) in  Case {\rm (iii) }}

\put(170,     10){\makebox(0,      0)[t]{$\bullet$}}

\put(170,     70){\makebox(0,      0)[t]{$\bullet$}}

\put(230,     10){\makebox(0,      0)[t]{$\bullet$}}
\put(230,     70){\makebox(0,      0)[t]{$\bullet$}}

\put(170,       10){\line(0,      1){60}}

\put(170,       10){\line(1,      1){60}}

\put(230,       10){\line(0,      1){60}}

\put(230,       10){\line(-1,      1){60}}

\put(170,       10){\line(1,       0){60}}
\put(170,       70){\line(1,       0){60}}

\put(150,     10){$-1$}
\put(150,     30){$q^{3}$}

\put(150,     70){$-1$}

\put(180,     30){$q^{-1}$}
\put(190,     80){$q^{-2}$}
\put(190,     -10){$-q^{-1}$}
\put(220,     30){$$}

\put(250,     10){$-q$}
\put(250,     30){$q^{-1}$}

\put(250,     70){$q$}

\put(280,         -1)  {   with  $q\in R_6.$ } \end{picture}$\\ \\
The sub-{\rm GDD} by deleting  Vertex 1  is not an arithmetic  {\rm  GDD}   by Lemma \ref {3.1.1}{\rm (I)}.

\subsection* {Quasi-affine over
  {\rm  GDD}  $4$ of Row $17$  }
   {\rm (i) } All quasi-arithmetic {\rm  GDD}s by  adding a vertex  on  Vertex $1$  are listed.
 According  to Lemma \ref {3.1.2} (IV) we have to consider following cases.

(a)    i.e.  (17.4.5)$, q \in F^{*}\setminus \{1,  -1\}$,   by    {\rm  GDD}  $2$ of Row $8$.
is quasi-affine   by Lemma \ref {3.1.21}.

(b)    i.e.  (17.4.1)$, q \in R_3$,   by    {\rm  GDD}  $2$ of Row $15$.
    is quasi-affine   by Lemma \ref {3.1.21}.

(c)    i.e.  (17.4.2),  $q \in R_3$,   by    {\rm  GDD}  $2$ of Row $16$.
 is quasi-affine   by Lemma \ref {3.1.21}.

(d)    i.e.  (17.4.3)$, q \in R_3$,   by   {\rm  GDD}  $1$ of Row $17$.
is quasi-affine   by Lemma \ref {3.1.21}.

(e)    i.e.  (17.4.4)$, q \in R_3$,   by    {\rm  GDD}  $4$ of Row $17$.
 is quasi-affine   by Lemma \ref {3.1.21}.

  {\rm (ii) }  All quasi-arithmetic {\rm  GDD}s by  adding a vertex  on  Vertex $2$  are listed.
 According  to Lemma \ref {3.1.2} (IV) we have to consider following cases.\\ \\

  {\ }\ \ \ \ \ \ \ \ \ \ \ \ $\begin{picture}(100,       15) \put(-68,        -1){ $(a)$}

\put(60,       1){\makebox(0,       0)[t]{$\bullet$}}
\put(58,       -12){$-1$}

\put(40,       -12){$q$}
\put(28,       -1){\line(1,       0){33}}
\put(27,       1){\makebox(0,      0)[t]{$\bullet$}}

\put(22,      -12){$-1$}
\put(0,       -12){$q^{-1}$}

\put(-14,       1){\makebox(0,      0)[t]{$\bullet$}}

\put(-14,      -1){\line(1,       0){50}}

\put(-18,      -12){$q$}

\put(27,     38){\makebox(0,      0)[t]{$\bullet$}}

\put(27,       0){\line(0,      1){35}}

\put(30,       30){$-1$}

\put(30,       20){$-q$}

 \ \ \ \ \ \ \ \ \ \ \ \ \ \ \ \ \ \ \ \ \ \ \ \ \ \ \ \ \ \ \ \   \  \ \ \ \ \ \ \ \ \ \ \ \ \ \ \ \ \ \ \ {$, q \in F^{*}\setminus \{1,  -1\}$,   by  {\rm  GDD}  $2$ of Row $4$.}
\put(80,         -1)  {    } \end{picture}$\\
\\
  The sub-{\rm GDD} by deleting  Vertex  4    is not an arithmetic  {\rm  GDD}   by Lemma \ref {3.1.1}{\rm (II)}.\\
\\

  {\ }\ \ \ \ \ \ \ \ \ \ \ \ $\begin{picture}(100,       15) \put(-68,        -1){ $(b)$}

\put(60,       1){\makebox(0,       0)[t]{$\bullet$}}
\put(58,       -12){$-1$}

\put(40,       -12){$q^2$}
\put(28,       -1){\line(1,       0){33}}
\put(27,       1){\makebox(0,      0)[t]{$\bullet$}}

\put(22,      -12){$-1$}
\put(0,       -12){$q^{-2}$}

\put(-14,       1){\makebox(0,      0)[t]{$\bullet$}}

\put(-14,      -1){\line(1,       0){50}}

\put(-18,      -12){$q$}

\put(27,     38){\makebox(0,      0)[t]{$\bullet$}}

\put(27,       0){\line(0,      1){35}}

\put(30,       30){$-1$}

\put(30,       20){$-q^{2}$}

 \ \ \ \ \ \ \ \ \ \ \ \ \ \ \ \ \ \ \ \ \ \ \ \ \ \ \ \ \ \ \ \ \   \   \ \ \ \ \ \ \ \ \ \ \ \ \ \ \ \ \ \ \ {$, q \in F^{*}\setminus \{1,  -1\}$,   by    {\rm  GDD}  $2$ of Row $5$.}
\put(80,         -1)  {    } \end{picture}$\\
\\
  The sub-{\rm GDD} by deleting  Vertex  4    is not an arithmetic  {\rm  GDD}   by Lemma \ref {3.1.1}{\rm (I)}.

(c)    i.e.  (17.4.9)$, q \in F^{*}\setminus \{1,  -1\}$,  by    {\rm  GDD}  $2$ of Row $6$.
  The sub-{\rm GDD} by deleting  Vertex  4    is not  an arithmetic  {\rm  GDD}   by Lemma \ref {3.1.1}{\rm (II)}.  \\ \\

  {\ }\ \ \ \ \ \ \ \ \ \ \ \ $\begin{picture}(100,       15) \put(-68,        -1){ (d)}

\put(60,       1){\makebox(0,       0)[t]{$\bullet$}}
\put(58,       -12){$-1$}

\put(40,       -12){$q$}
\put(28,       -1){\line(1,       0){33}}
\put(27,       1){\makebox(0,      0)[t]{$\bullet$}}

\put(22,      -12){$-1$}
\put(0,       -12){$q^{-1}$}

\put(-14,       1){\makebox(0,      0)[t]{$\bullet$}}

\put(-14,      -1){\line(1,       0){50}}

\put(-18,      -12){$-1$}

\put(27,     38){\makebox(0,      0)[t]{$\bullet$}}

\put(27,       0){\line(0,      1){35}}

\put(30,       30){$-1$}

\put(30,       20){$-q$}

\ \ \ \ \ \ \ \ \ \ \ \ \ \ \ \ \ \ \ \ \ \ \ \ \ \ \ \ \ \ \ \ \ \ \ \ \ \   \   \ \ \ \ \ \ \ \ \ \ \ \ \ \ \ \ \ \ \ {$, q \in F^{*}\setminus \{1,  -1\}$,   by    {\rm  GDD}  $2$ of Row $8$.}
\put(80,         -1)  {    } \end{picture}$\\
\\
 The sub-{\rm GDD} by deleting  Vertex  4   is not an arithmetic  {\rm  GDD}    by Lemma \ref {3.1.1}{\rm (II)}.

 (e)  i.e.  (17.4.6) $, q \in R_3$,   by    {\rm  GDD}  $2$ of Row $15$. The sub-{\rm GDD} by deleting  Vertex  4   is
  an arithmetic  {\rm  GDD}    by Lemma \ref {3.1.1}{\rm (II)} or by Lemma \ref {3.1.1} {\rm (X)}. It  is   quasi-affine   by Lemma \ref {3.1.27}.

   (f)  i.e.  (17.4.7),  $q \in R_3$,   by    {\rm  GDD}  $2$ of Row $16$. The sub-{\rm GDD} by deleting  Vertex  4
  is an arithmetic  {\rm  GDD}    by  Type   7.  It  is quasi-affine   by Lemma \ref {3.1.22}.

 (g)  i.e.  (17.4.8)$, q \in R_3$,   by    {\rm  GDD}  $1$ of Row $17$. The sub-{\rm GDD}
  by deleting  Vertex  4 is an arithmetic  {\rm  GDD}   by Lemma \ref {3.1.2} {\rm (V)}.
  It  is quasi-affine   by Lemma \ref {3.1.27}.\\ \\

  {\ }\ \ \ \ \ \ \ \ \ \ \ \ $\begin{picture}(100,       15) \put(-68,        -1){ (h)}

\put(60,       1){\makebox(0,       0)[t]{$\bullet$}}
\put(58,       -12){$-1$}

\put(40,       -12){$q$}
\put(28,       -1){\line(1,       0){33}}
\put(27,       1){\makebox(0,      0)[t]{$\bullet$}}

\put(22,      -12){$-1$}
\put(0,       -12){$-q$}

\put(-14,       1){\makebox(0,      0)[t]{$\bullet$}}

\put(-14,      -1){\line(1,       0){50}}

\put(-18,      -12){$-1$}

\put(27,     38){\makebox(0,      0)[t]{$\bullet$}}

\put(27,       0){\line(0,      1){35}}

\put(30,       30){$-1$}

\put(30,       20){$-q$}

\ \ \ \ \ \ \ \ \ \ \ \ \ \ \ \ \ \ \ \ \ \ \ \ \ \ \ \   {$, q \in R_3$,   by   {\rm  GDD}  $4$ of Row $17$. The sub-{\rm GDD} by deleting}
\put(80,         -1)  {    } \end{picture}$\\
\\
 Vertex  4 is not an arithmetic  {\rm  GDD}   by Lemma \ref {3.1.1}{\rm (II)} or by Lemma \ref {3.1.1} {\rm (X)}.

  {\rm (iii) }  All quasi-arithmetic {\rm  GDD}s by  adding a vertex  on  Vertex $3$  are listed.
 According  to Lemma \ref {3.1.2} (IV) we have to consider following cases.

 (a)  i.e.  (17.4.11), $q = -\xi,  \xi \in R_3 $,  $ q \notin R_3$,  by  {\rm  GDD}  $2$ of Row $7$.
 It is quasi-affine by  Lemma   \ref {2.63}.

  {\rm (iv) }  All quasi-affine circles are listed.

(a)  (a)    i.e.  (17.4.12)
  is  quasi-affine  since the sub-{\rm GDD} by deleting  Vertex 2
with    $-1 =-q^{3}$ is an arithmetic  {\rm  GDD}   by    {\rm  GDD}  $2$ of Row $6$ or Lemma \ref {3.1.2} (XI).

(a) (nc)    i.e.  (17.4.13)   is  quasi-affine  since the sub-{\rm GDD} by deleting  Vertex 2
  is an arithmetic  {\rm  GDD}   by Lemma \ref {3.1.1} (X).

(b) (a)    i.e.  (17.4.14)
  is   quasi-affine  since the sub-{\rm GDD} by deleting  Vertex 2
  with   $ -q^{3} ={-1}$ is an arithmetic  {\rm  GDD}   by  Type   7.

(c) (nc)    i.e.  (17.4.15)
  is   quasi-affine since the sub-{\rm GDD} by deleting  Vertex 2
 is an arithmetic  {\rm  GDD}.

(c) (a)    i.e.  (17.4.16)    is   quasi-affine  since the sub-{\rm GDD} by deleting  Vertex 2
 with $-q^{3}=-1
$ is
is an arithmetic  {\rm  GDD}.

 (c) (nc)    i.e.  (17.4.17)   is   quasi-affine  since the sub-{\rm GDD} by deleting  Vertex 2
   is  an arithmetic  {\rm  GDD}   by Lemma \ref {3.1.1} (X).

(d) (a)    i.e.  (17.4.18)
  is quasi-affine  since the sub-{\rm GDD} by deleting  Vertex 2
 with $-q^{3}=-1$
is  an arithmetic  {\rm  GDD}   by Lemma \ref {3.1.1} {\rm (X)}.

 (d) (nc)    i.e.  (17.4.19)  is  quasi-affine since the sub-{\rm GDD} by deleting  Vertex 2
  is  an arithmetic  {\rm  GDD}   by  Type   7.

\subsection* {Quasi-affine over
  {\rm  GDD}  $5$ of Row $17$  }
   {\rm (i) } All quasi-arithmetic {\rm  GDD}s by  adding a vertex  on  Vertex $1$  are listed.
According  to   Type   3
 we have to consider following cases.

   {\rm (ii) }  All quasi-arithmetic {\rm  GDD}s by  adding a vertex  on  Vertex $2$  are listed.
  According  to Lemma \ref {3.1.2} (III) we have to consider following cases.

(a)     i.e.  (17.5.1), $q \in R_6$,   by    {\rm  GDD}  $1$ of Row $6$. The sub-{\rm GDD} by deleting
Vertex  4 is  an arithmetic  {\rm  GDD}   by Lemma \ref {3.1.1} {\rm (V)}. It  is   quasi-affine   by Lemma \ref {3.1.3}.\\ \\

  {\ }\ \ \ \ \ \ \ \ \ \ \ \ $\begin{picture}(100,       15) \put(-68,        -1){(b) }

\put(60,       1){\makebox(0,       0)[t]{$\bullet$}}
\put(58,       -12){$-1$}

\put(40,       -12){$q^{-1}$}
\put(28,       -1){\line(1,       0){33}}
\put(27,       1){\makebox(0,      0)[t]{$\bullet$}}

\put(22,      -12){$q$}
\put(0,       -12){$q^{-3}$}

\put(-14,       1){\makebox(0,      0)[t]{$\bullet$}}

\put(-14,      -1){\line(1,       0){50}}

\put(-18,      -12){$q^3$}

\put(27,     38){\makebox(0,      0)[t]{$\bullet$}}

\put(27,       0){\line(0,      1){35}}

\put(30,       30){$-1$}

\put(30,       20){$-q$}

\ \ \ \ \ \ \ \ \ \ \ \ \ \ \ \   \   \ \ \ \ \ \ \ \ \ \   { $q\in R_6$ $, q \in F^{*}\setminus \{1,  -1\}$ $, q \notin R_3$,   by    {\rm  GDD}  $1$ of Row $7$.}
\put(80,         -1)  {    } \end{picture}$\\ \\
 The sub-{\rm GDD} by deleting  Vertex  4   is not an arithmetic  {\rm  GDD}   by Lemma \ref {3.1.1} {\rm (V)}.

 (c) (nc)    i.e.  (17.5.2)$, q \in F^{*}\setminus \{1,  -1\}$. $q\in R_6, $ by    {\rm  GDD}  $3$ of Row $8$.
  The sub-{\rm GDD} by deleting  Vertex  4 is an arithmetic  {\rm  GDD}    by  {\rm  GDD}  5 of Row 17 or by Lemma \ref {3.1.2} {\rm (III)}.   It  is quasi-affine   by Lemma \ref {3.1.3}.\\ \\

  {\ }\ \ \ \ \ \ \ \ \ \ \ \ $\begin{picture}(100,       15) \put(-68,        -1){(d) }

\put(60,       1){\makebox(0,       0)[t]{$\bullet$}}

\put(28,       -1){\line(1,       0){33}}
\put(27,       1){\makebox(0,      0)[t]{$\bullet$}}

\put(58,       -12){$-1$}

\put(34,       -12){$-q$}

\put(6,      -12){$-q^{-1}$}
\put(-10,       -12){$-q$}

\put(-18,      -12){$q$}

\put(-14,       1){\makebox(0,      0)[t]{$\bullet$}}

\put(-14,      -1){\line(1,       0){50}}

\put(27,     38){\makebox(0,      0)[t]{$\bullet$}}

\put(27,       0){\line(0,      1){35}}

\put(30,       30){$-1$}

\put(30,       20){$q^{-1}$}

\ \ \ \ \ \ \ \ \ \ \ \ \ \ \ \ \ \ \ \ \ \ \ \ \ \ \ \ \ \ \ \ \ \ \ \ \ \   \    \ \ \ \ \ \ \ \ \ \ \ \ \ \ \ \ \ \ \ {$, q \in R_3$,   by    {\rm  GDD}  $1$ of Row $14$.}
\put(80,         -1)  {    } \end{picture}$\\
\\
 The sub-{\rm GDD} by deleting  Vertex  4 is not an arithmetic  {\rm  GDD}   by Lemma \ref {3.1.1}{\rm (I)}.

 (e)     i.e.  (17.5.3)$, q \in R_3$,   by    {\rm  GDD}  $5$ of Row $17$.    The sub-{\rm GDD} by deleting  Vertex  4 is an
arithmetic {\rm  GDD} by  {\rm  GDD}  4 of 15 Row or by Lemma \ref {3.1.1} {\rm (V)}.   It  is quasi-affine   by Lemma \ref {3.1.3}.\\ \\

{\ }\ \ \ \ \ \ \ \ \ \ \ \ $\begin{picture}(100,       15) \put(-68,        -1){ $(f)$}

\put(60,       1){\makebox(0,       0)[t]{$\bullet$}}

\put(28,       -1){\line(1,       0){33}}
\put(27,       1){\makebox(0,      0)[t]{$\bullet$}}
\put(-14,       1){\makebox(0,      0)[t]{$\bullet$}}

\put(-14,      -1){\line(1,       0){50}}

\put(58,       -12){$-1$}

\put(40,       -12){$q^{-1}$}

\put(22,      -12){$q$}
\put(0,       -12){${q^{-1}}$}

\put(-18,      -12){${q}$}

\put(27,     38){\makebox(0,      0)[t]{$\bullet$}}

\put(27,       0){\line(0,      1){35}}

\put(30,       30){$-1$}

\put(30,       20){$-q$}

\ \ \ \ \ \ \ \ \ \ \ \ \ \ \ \ \ \ \ \ \ \ \ \ \ \ \ \ \ \ \ \ \ \ \ \ \ \   \    \ \ \ \ \ \ \ \ \ \ \ \ \ \ \ \ \ \ \ {,   $q\in R_6$,  by   {\rm  GDD}  $1$ of Row $4$.}
\put(80,         -1)  {    } \end{picture}$\\
\\
 The sub-{\rm GDD} by deleting  Vertex  4 is not  an arithmetic  {\rm  GDD}   by Lemma \ref {3.1.1} {\rm (VIII)}.\\

{\ }\ \ \ \ \ \ \ \ \ \ \ \ $\begin{picture}(100,       15) \put(-68,        -1){ $(g)$}

\put(60,       1){\makebox(0,       0)[t]{$\bullet$}}

\put(28,       -1){\line(1,       0){33}}
\put(27,       1){\makebox(0,      0)[t]{$\bullet$}}
\put(-14,       1){\makebox(0,      0)[t]{$\bullet$}}

\put(-14,      -1){\line(1,       0){50}}

\put(58,       -12){$-1$}

\put(40,       -12){$q^{-2}$}

\put(22,      -12){$q^2$}
\put(0,       -12){${q^{-2}}$}

\put(-18,      -12){${q}$}

\put(27,     38){\makebox(0,      0)[t]{$\bullet$}}

\put(27,       0){\line(0,      1){35}}

\put(30,       30){$-1$}

\put(30,       20){$-q^2$}

\ \ \ \ \ \ \ \ \ \ \ \ \ \ \ \ \ \ \ \ \ \ \ \ \ \ \ \ \ \ \ \ \ \ \ \ \ \   \    \ \ \ \ \ \ \ \ \ \ \ \ \ \ \ \ \ \ \ {$q^{2}\in R_6, $   by   {\rm  GDD}  $1$ of Row $5$.}
\put(80,         -1)  {    } \end{picture}$\\
\\
 The sub-{\rm GDD} by deleting  Vertex  4 is not an arithmetic  {\rm  GDD}   by Lemma \ref {3.1.1} {\rm (VIII)}.

  {\rm (iii) }  All quasi-arithmetic {\rm  GDD}s by  adding a vertex  on  Vertex $3$  are listed.
 According  to Lemma \ref {3.1.2} (I) we have to consider following cases.

 (a)    i.e.  (17.5.4)
 $q^2 = -\xi,  \xi \in R_3, $ by     {\rm  GDD}  $2$ of Row $6$.
 It is quasi-affine by  Lemma   \ref {2.63}.

 (b)    i.e.  (17.5.5),  $q^3 = -\xi,  \xi \in R_3, $   $q\notin R_3, $  by    {\rm  GDD}  $2$ of Row $7$.
 It is quasi-affine by  Lemma   \ref {2.63}.

 (c)    i.e.  (17.5.6),  $-q\in R_3$,  by Lemma \ref {3.1.1} {\rm (IV)}, $r \not= q^{-1}.$
It is quasi-affine by  Lemma   \ref {2.63}.

(d)    i.e.  (17.5.9)
 $-q^2 \in R_3$ $, q ^2\in R_6$ $, q \notin R_3$,   by    {\rm  GDD}  $2$ of Row $10$.
 It is quasi-affine  by \cite [Lemma 2.1]{TZ22}

(e)    i.e.  (17.5.7)$, q \in R_3$,   by    {\rm  GDD}  $2$ of Row $16$.
 It is quasi-affine
 by  Lemma   \ref {2.63}.

(f)    i.e.  (17.5.8)$, q \in R_3$,   by    {\rm  GDD}  $4$ of Row $16$.
 It is quasi-affine
 by  Lemma   \ref {2.63}.

  {\rm (iv) }  All quasi-affine circles are listed.

(nc)  (a)    i.e.  (17.5.10) is quasi-affine  since the sub-{\rm GDD} by deleting  Vertex 2
   with  $-r^{} =q^{-2} $ is an arithmetic  {\rm  GDD}.

 (nc)
 (f)    i.e.  (17.5.11)
  is quasi-affine  since the sub-{\rm GDD} by deleting  Vertex 2
   with   $q^{-2}=q, $ is an arithmetic  {\rm  GDD}    by Lemma \ref {3.1.2} {\rm (III)}.

(nc)   (a)  is not quasi-affine since\\

  $\begin{picture}(100,       15) \put(-68,        -1){ }

\put(60,       1){\makebox(0,       0)[t]{$\bullet$}}

\put(28,       -1){\line(1,       0){33}}
\put(27,       1){\makebox(0,      0)[t]{$\bullet$}}

\put(-14,       1){\makebox(0,      0)[t]{$\bullet$}}

\put(-14,      -1){\line(1,       0){50}}

\put(-18,      10){$-1$}
\put(0,       5){$q^{-1}$}
\put(22,      10){${-1}$}
\put(40,       5){$\xi$}

\put(58,       10){$-1$}

\put(80,         -1)  { with  $ \xi^{2} =-q^{}$ is not  an arithmetic  {\rm  GDD}    by Lemma \ref {3.1.2} {\rm (XI)} with $q\in R_3.$   } \end{picture}$

 (nc)  (b)  is not quasi-affine since\\

  $\begin{picture}(100,       15) \put(-68,        -1){ }

\put(60,       1){\makebox(0,       0)[t]{$\bullet$}}

\put(28,       -1){\line(1,       0){33}}
\put(27,       1){\makebox(0,      0)[t]{$\bullet$}}

\put(-14,       1){\makebox(0,      0)[t]{$\bullet$}}

\put(-14,      -1){\line(1,       0){50}}

\put(-18,      10){$-1$}
\put(0,       5){$q^{-1}$}
\put(22,      10){${-1}$}
\put(40,       5){$\xi$}

\put(58,       10){$-1$}

\put(80,         -1)  {   with  $ \xi^{3} =-q^{}$ is not an arithmetic  {\rm  GDD}    by Lemma \ref {3.1.2} {\rm (XI)}. } \end{picture}$

 (nc)
 (c)    i.e.  (17.5.12)  is quasi-affine  since the sub-{\rm GDD} by deleting  Vertex 2
 with   $-r=-1, $ is an arithmetic  {\rm  GDD}    by Lemma \ref {3.1.2} {\rm (XI)}.

 (nc)
 (e)    i.e.  (17.5.13)    is quasi-affine  since the sub-{\rm GDD} by deleting  Vertex 2
   is an arithmetic  {\rm  GDD}    by  Type   7.

\subsection* {Quasi-affine over
  {\rm  GDD}  $6$ of Row $17$ }

   {\rm (ii) }  All quasi-arithmetic {\rm  GDD}s by  adding a vertex  on  Vertex $2$  are listed.

 We have to consider following cases step by step.\\

{\ }\!\!\!\!\!\!\!\!\!\!\!\!\!\!\!\!\!\!\!\!\!{\ }\!\!\!\!\!\!\!\!\!\!\!\!\!\!\!\!\!\!\!\!\!{\ }\!\!\!\!\!\!\!\!\!\!\!\!
  $\begin{picture}(100,       15)\put(-68,       -1){ }\put(68,       -1){(a) }
\put(111,      1){\makebox(0,      0)[t]{$\bullet$}}
\put(144,       1){\makebox(0,       0)[t]{$\bullet$}}
\put(170,      -11){\makebox(0,      0)[t]{$\bullet$}}
\put(170,     15){\makebox(0,      0)[t]{$\bullet$}}
\put(113,      -1){\line(1,      0){33}}
\put(142,     -1){\line(2,      1){27}}
\put(170,       -14){\line(-2,      1){27}}

\put(100,       10){$-1$}

\put(115,       5){$-q^{-1}$}

\put(123,       -12){$q^{-1}$}

\put(140,      -20){$q^{-1}$}
\put(145,       15){$-q^{-1}$}

\put(170,       -12){$-1$}
\put(170,       18){$-1$}

\put(195,         -1)  {,  $q\in R_3, $  by   {\rm  GDD}  $6$ of Row $17$. }
\put(80,         -1)  {    } \end{picture}$\\
\\
 The sub-{\rm GDD} by deleting  Vertex  4   is not an arithmetic  {\rm  GDD}   by Lemma \ref {3.1.1} {\rm (V)}.

 (b)    i.e.  (17.6.1),   $q\in R_3, $  by   {\rm  GDD}  $1$ of Row $15$.  The sub-{\rm GDD} by deleting  Vertex  4
is an arithmetic  {\rm  GDD}   by Lemma \ref {3.1.2} {\rm (III)}.
It  is quasi-affine  by Lemma \ref {3.1.3}.

  {\rm (iii) }  All quasi-arithmetic {\rm  GDD}s by  adding a vertex  on  Vertex $3$  are listed.
According  to   Type   4 and  {\rm  GDD}  $8,  9$ of Row $17$
 we have to consider following cases.

 {\rm (iv) }  All quasi-affine circles are listed.

 (a)    i.e.  (17.6.2),   $q\in R_3, $
  is  quasi-affine  since the sub-{\rm GDD} by deleting  Vertex 2
   is  an arithmetic  {\rm  GDD}    $4$ of Row $17$.

\subsection* {Quasi-affine over
  {\rm  GDD}  $7$ of Row $17$  }
   {\rm (i) } All quasi-arithmetic {\rm  GDD}s by  adding a vertex  on  Vertex $1$  are listed.
 According  to Lemma \ref {3.1.2} (I) we have to consider following cases.

 (a)    i.e.  (17.7.1)
  $q^2\in R_3$,   by  {\rm  GDD}  $2$ of Row $6$. It is quasi-affine   by Lemma \ref {3.1.26}.

 (b)    i.e.  (17.7.2), $q^3\in R_3, $  by    {\rm  GDD}  $2$ of Row $7$. It is quasi-affine       by Lemma \ref {3.1.26}.

(c)    i.e.  (17.7.5), $q\in R_9,  $ by   {\rm  GDD}  $4$ of Row $7$. It is quasi-affine       by Lemma \ref {3.1.26}.

 (d)    i.e.  (17.7.3),  by   by Lemma \ref {3.1.1} {\rm (IV)},  $r\not=q^{-1}$ . It is quasi-affine       by Lemma \ref {3.1.26}.

  (e)    i.e.  (17.7.4)$, q \in R_3$,   by    {\rm  GDD}  $1$ of Row $11$. It is quasi-affine       by Lemma \ref {3.1.26}.

  {\rm (ii) }  All quasi-arithmetic {\rm  GDD}s by  adding a vertex  on  Vertex $2$  are listed.
 According  to Lemma \ref {3.1.2} (III) we have to consider following cases.\\
\\

  {\ }\ \ \ \ \ \ \ \ \ \ \ \ $\begin{picture}(100,       15) \put(-68,        -1){ (a)}

\put(60,       1){\makebox(0,       0)[t]{$\bullet$}}
\put(58,       -12){$-1$}

\put(40,       -12){$q^{-1}$}
\put(28,       -1){\line(1,       0){33}}
\put(27,       1){\makebox(0,      0)[t]{$\bullet$}}

\put(22,      -12){$q$}
\put(0,       -12){$q^{-2}$}

\put(-14,       1){\makebox(0,      0)[t]{$\bullet$}}

\put(-14,      -1){\line(1,       0){50}}

\put(-18,      -12){$q^2$}

\put(27,     38){\makebox(0,      0)[t]{$\bullet$}}

\put(27,       0){\line(0,      1){35}}

\put(30,       30){$-1$}

\put(30,       20){$-q$}

\ \ \ \ \ \ \ \ \ \ \ \ \ \ \ \ \ \ \ \ \ \ \ \ \ \ \ \ \ \ \ \ \ \ \ \ \ \   \  \ \ \ \ \ \ \ \ \ \ \ \ \ \ \ \ \ \ \ {$, q \in R_3$,  by    {\rm  GDD}  $1$ of Row $6$.}
\put(80,         -1)  {    } \end{picture}$\\
\\
  The sub-{\rm GDD} by deleting  Vertex  4    is not  an arithmetic  {\rm  GDD} ,    by Lemma \ref {3.1.1} {\rm (V)}.

 (b)    i.e.  (17.7.6)$, q \in F^{*}\setminus \{1,  -1\}$,   by    {\rm  GDD}  $3$ of Row $8$.
  The sub-{\rm GDD} by deleting  Vertex  4    is an arithmetic  {\rm  GDD}   by  {\rm  GDD}  $7$ of Row $17$ or  by Lemma \ref {3.1.1} {\rm (VIII)}. It  is   quasi-affine
 by Lemma \ref {3.1.3}.

 (c)    i.e.  (17.7.7)$, q \in R_3$,   by    {\rm  GDD}  $1$ of Row $15$.
 The sub-{\rm GDD} by deleting  Vertex  4 is an arithmetic  {\rm  GDD}    by  {\rm  GDD}  6 of Row 17 or by Lemma \ref {3.1.1} {\rm (VIII)}.    It  is quasi-affine   by Lemma \ref {3.1.3}.\\ \\

  {\ }\ \ \ \ \ \ \ \ \ \ \ \ $\begin{picture}(100,       15) \put(-68,        -1){ (d)}

\put(60,       1){\makebox(0,       0)[t]{$\bullet$}}
\put(58,       -12){$-1$}

\put(40,       -12){$q^{-1}$}
\put(28,       -1){\line(1,       0){33}}
\put(27,       1){\makebox(0,      0)[t]{$\bullet$}}

\put(22,      -12){$q$}
\put(-4,       -12){$-q^{-1}$}

\put(-14,       1){\makebox(0,      0)[t]{$\bullet$}}

\put(-14,      -1){\line(1,       0){50}}

\put(-22,      -12){$-q$}

 \put(27,     38){\makebox(0,      0)[t]{$\bullet$}}

\put(27,       0){\line(0,      1){35}}

\put(30,       30){$-1$}

\put(30,       20){$-q$}

\ \ \ \ \ \ \ \ \ \ \ \ \ \ \ \ \ \ \ \ \ \ \ \ \ \ \ \ \ \ \   {$, q \in R_3$,   by    {\rm  GDD}  $1$ of Row $16$. The sub-{\rm GDD} by deleting  }
\put(80,         -1)  {    } \end{picture}$\\
\\
Vertex  4  is not an arithmetic  {\rm  GDD}   by Lemma \ref {3.1.1} {\rm (V)} or
 by Lemma \ref {3.1.1} {\rm (VIII)}. \\ \\

{\ }\ \ \ \ \ \ \ \ \ \ \ \ $\begin{picture}(100,       15) \put(-68,        -1){ (e)}

\put(60,       1){\makebox(0,       0)[t]{$\bullet$}}
\put(58,       -12){$-1$}

\put(40,       -12){$q^{-1}$}
\put(28,       -1){\line(1,       0){33}}
\put(27,       1){\makebox(0,      0)[t]{$\bullet$}}

\put(22,      -12){$q$}
\put(0,       -12){$-1$}

\put(-14,       1){\makebox(0,      0)[t]{$\bullet$}}

\put(-14,      -1){\line(1,       0){50}}

\put(-18,      -12){$-1$}

\put(27,     38){\makebox(0,      0)[t]{$\bullet$}}

\put(27,       0){\line(0,      1){35}}

\put(30,       30){$-1$}

\put(30,       20){$-q$}

\ \ \ \ \ \ \ \ \ \ \ \ \ \ \ \ \ \ \ \ \ \ \ \ \ \ \ \ \ \ \ \ \ \ \ \ \ \   \    \ \ \ \ \ \ \ \ \ \ \ \ \ \ \ \ \ \ \ {$, q \in R_3$,   by    {\rm  GDD}  $2$ of Row $17$.}
\put(80,         -1)  {    } \end{picture}$\\
\\
  The sub-{\rm GDD} by deleting  Vertex  4    is not  an arithmetic  {\rm  GDD}   by Lemma \ref {3.1.1} {\rm (V)}.   \\ \\

  {\ }\ \ \ \ \ \ \ \ \ \ \ \ $\begin{picture}(100,       15) \put(-68,        -1){(f) }

\put(60,       1){\makebox(0,       0)[t]{$\bullet$}}
\put(58,       -12){$-1$}

\put(40,       -12){$q^{-1}$}
\put(28,       -1){\line(1,       0){33}}
\put(27,       1){\makebox(0,      0)[t]{$\bullet$}}

\put(22,      -12){$q$}
\put(0,       -12){$-q$}

\put(-14,       1){\makebox(0,      0)[t]{$\bullet$}}

\put(-14,      -1){\line(1,       0){50}}

\put(-18,      -12){$-1$}

 \put(27,     38){\makebox(0,      0)[t]{$\bullet$}}

\put(27,       0){\line(0,      1){35}}

\put(30,       30){$-1$}

\put(30,       20){$-q$}

\ \ \ \ \ \ \ \ \ \ \ \ \ \ \ \ \ \ \ \ \ \ \ \ \ \ \ \ \ \ \ \ \ \ \ \ \ \   \   \ \ \ \ \ \ \ \ \ \ \ \ \ \ \ \ \ \ \ {$, q \in R_3$,   by    {\rm  GDD}  $7$ of Row $17$.}
\put(80,         -1)  {    } \end{picture}$\\
\\
  The sub-{\rm GDD} by deleting  Vertex  4    is not  an arithmetic  {\rm  GDD}   by Lemma \ref {3.1.1} {\rm (V)}.\\ \\

{\ }\ \ \ \ \ \ \ \ \ \ \ \ $\begin{picture}(100,       15) \put(-68,        -1){$(g)$ }

\put(60,       1){\makebox(0,       0)[t]{$\bullet$}}

\put(28,       -1){\line(1,       0){33}}
\put(27,       1){\makebox(0,      0)[t]{$\bullet$}}
\put(-14,       1){\makebox(0,      0)[t]{$\bullet$}}

\put(-14,      -1){\line(1,       0){50}}

\put(58,       -12){$-1$}

\put(40,       -12){$q^{-1}$}

\put(22,      -12){$q$}
\put(0,       -12){${q^{-1}}$}

\put(-18,      -12){${q}$}

 \put(27,     38){\makebox(0,      0)[t]{$\bullet$}}

\put(27,       0){\line(0,      1){35}}

\put(30,       30){$-1$}

\put(30,       20){$-q$}

\ \ \ \ \ \ \ \ \ \ \ \ \ \ \ \ \ \ \ \ \ \ \ \ \ \ \ \ \ \ \ \ \ \ \ \ \ \   \   \ \ \ \ \ \ \ \ \ \ \ \ \ \ \ \ \ \ \ {  by   {\rm  GDD}  $1$ of Row $4$.}
\put(80,         -1)  {    } \end{picture}$\\
\\
  The sub-{\rm GDD} by deleting  Vertex  4    is not an arithmetic  {\rm  GDD}   by Lemma \ref {3.1.1} {\rm (VIII)},  \\ \\

{\ }\ \ \ \ \ \ \ \ \ \ \ \ $\begin{picture}(100,       15) \put(-68,        -1){ $(h)$}

\put(60,       1){\makebox(0,       0)[t]{$\bullet$}}

\put(28,       -1){\line(1,       0){33}}
\put(27,       1){\makebox(0,      0)[t]{$\bullet$}}
\put(-14,       1){\makebox(0,      0)[t]{$\bullet$}}

\put(-14,      -1){\line(1,       0){50}}

\put(58,       -12){$-1$}

\put(40,       -12){$q^{-2}$}

\put(22,      -12){$q^2$}
\put(0,       -12){${q^{-2}}$}

\put(-18,      -12){${q}$}

 \put(27,     38){\makebox(0,      0)[t]{$\bullet$}}

\put(27,       0){\line(0,      1){35}}

\put(30,       30){$-1$}

\put(30,       20){$-q^2$}

\ \ \ \ \ \ \ \ \ \ \ \ \ \ \ \ \ \ \ \ \ \ \ \ \ \ \ \ \ \ \ \ \ \ \ \ \ \   \   \ \ \ \ \ \ \ \ \ \ \ \ \ \ \ \ \ \ \ { by    {\rm  GDD}  $1$ of Row $5$.}
\put(80,         -1)  {    } \end{picture}$\\
\\
  The sub-{\rm GDD} by deleting  Vertex  4    is not an arithmetic  {\rm  GDD}    by Lemma \ref {3.1.1} {\rm (VIII)}.

  {\rm (iv) }  All quasi-affine circles are listed.

(nc)  (nc)  is not quasi-affine since\\

$\begin{picture}(100,       15) \put(-68,        -1){ }

\put(60,       1){\makebox(0,       0)[t]{$\bullet$}}

\put(28,       -1){\line(1,       0){33}}
\put(27,       1){\makebox(0,      0)[t]{$\bullet$}}

\put(-14,       1){\makebox(0,      0)[t]{$\bullet$}}

\put(-14,      -1){\line(1,       0){50}}

\put(-18,      10){$-1$}
\put(0,       5){$q$}
\put(16,      10){$-1$}
\put(34,       5){$-q^{-1}$}

\put(60,       10){$-1$}

\put(80,         -1)  {  is not an arithmetic  {\rm  GDD}   by Lemma \ref {3.1.1} (X) or by Lemma \ref {3.1.1} (II).  } \end{picture}$

  (a) (nc)    i.e.  (17.7.8)    is  quasi-affine  since the sub-{\rm GDD} by deleting  Vertex 2
   is  an arithmetic  {\rm  GDD}   by Lemma \ref {3.1.1}{\rm (II)}(g) or  by Lemma \ref {3.1.1} {\rm (X)}.

 (b)  (nc)  is not quasi-affine since\\

  $\begin{picture}(100,       15) \put(-68,        -1){  }

\put(60,       1){\makebox(0,       0)[t]{$\bullet$}}

\put(28,       -1){\line(1,       0){33}}
\put(27,       1){\makebox(0,      0)[t]{$\bullet$}}

\put(-14,       1){\makebox(0,      0)[t]{$\bullet$}}

\put(-14,      -1){\line(1,       0){50}}

\put(-18,      10){$-1$}
\put(0,       5){$q^{}$}
\put(18,      10){$-1$}
\put(34,       5){$-q^{-3}$}

\put(58,       10){$-1$}

\put(80,         -1)  { is not an arithmetic  {\rm  GDD}   by Lemma \ref {3.1.1} (X).   } \end{picture}$

 (d) (nc)    i.e.  (17.7.9)   is quasi-affine since the sub-{\rm GDD} by deleting  Vertex 2
  with  $r =-1 $   is   an arithmetic  {\rm  GDD}   by Lemma \ref {3.1.1} (X).

  (d) (nc)    i.e.  (17.7.10)  is  quasi-affine  since the sub-{\rm GDD} by deleting  Vertex 2
  with  $r =-q^{} $ is an arithmetic  {\rm  GDD}.
\subsection* {Quasi-affine over
  {\rm  GDD}  $8$ of Row $17$   }
   {\rm (i) } All quasi-arithmetic {\rm  GDD}s by  adding a vertex  on  Vertex $1$  are listed.
 According  to Lemma \ref {3.1.2} (V) we have to consider following cases.

 (a)     i.e.  (17.8.1),     $r \in  R_2 \cup R_4 \cup R_6 \cup R_3 $ by Lemma \ref {3.1.2} (V).  The sub-{\rm GDD} by deleting  Vertex  4
 is  arithmetic  {\rm  GDD}   when    $ r =-q^{-1}$ by  Type   4.
   The sub-{\rm GDD} by deleting  Vertex  4    is not  an arithmetic  {\rm  GDD}   when    $ r \not=-q^{-1}$ by Lemma \ref {3.1.1}{\rm (I)}.
 It  is   quasi-affine   when $ r =-q^{-1}$ by Lemma \ref {3.1.3}.\\

{\ }\!\!\!\!\!\!\!\!\!\!\!\!\!\!\!\!\!\!\!\!\!{\ }\!\!\!\!\!\!\!\!\!\!\!\!\!\!\!
$\begin{picture}(100,       15)\put(-68,       -1){ }\put(38,       -1){(b) }

\put(111,      1){\makebox(0,      0)[t]{$\bullet$}}
\put(144,       1){\makebox(0,       0)[t]{$\bullet$}}
\put(170,      -11){\makebox(0,      0)[t]{$\bullet$}}
\put(170,     15){\makebox(0,      0)[t]{$\bullet$}}
\put(113,      -1){\line(1,      0){33}}
\put(142,     -1){\line(2,      1){27}}
\put(170,       -14){\line(-2,      1){27}}

\put(170,       -14){\line(0,      1){27}}

\put(100,       10){$r^{-1}$}
\put(110,       5){$-r^{-1}$}

\put(127,     - 10){$-1$}

\put(140,      -20){$-1$}
\put(150,       15){$-q$}

\put(175,       10){$q$}

\put(178,       -20){$-1$}

\put(178,       0){$q^{-1}$}

\put(200,         -1)  {,  $r\in R_3, $ by     {\rm  GDD}  $9$ of Row $17$ in Table A2. }
\put(80,         -1)  {    } \end{picture}$\\
\\
  The sub-{\rm GDD} by deleting  Vertex  4    is not an arithmetic  {\rm  GDD}   by Lemma \ref {3.1.1}{\rm (I)}.

 (c)     i.e.  (17.8.2),  $r\in R_3, $ by    {\rm  GDD}  $8$ of Row $17$ in Table A2.
 The sub-{\rm GDD} by deleting  Vertex  1    is  an arithmetic  {\rm  GDD}    when   $ q^{} =-r^{}$ by  Type   4.
 The sub-{\rm GDD} by deleting  Vertex  2    is  an arithmetic  {\rm  GDD}    when   $ q^{} =-r^{}$ by  Type   3.
  It  is quasi-affine   when $ q^{} =-r^{}$.
 \\ \\ \\

  {\ }\\ \\

{\ }\!\!\!\!\!\!\!\!\!\!\!\!\!\!\!\!\!\!\!\!\!{\ }\!\!\!\!\!\!\!\!\!\!\!\!\!\!\!\!\!\!\!\!\!{\ }\!\!\!\!\!\!\!\!\!\!\!
   $\begin{picture}(100,       15)\put(-68,       -1){ }\put(68,       -1){ (d)}

\put(170,     10){\makebox(0,      0)[t]{$\bullet$}}

\put(170,     70){\makebox(0,      0)[t]{$\bullet$}}

\put(230,     10){\makebox(0,      0)[t]{$\bullet$}}
\put(230,     70){\makebox(0,      0)[t]{$\bullet$}}

\put(170,       10){\line(0,      1){60}}

\put(170,       10){\line(1,      1){60}}

\put(230,       10){\line(0,      1){60}}


\put(170,       10){\line(1,       0){60}}
\put(170,       70){\line(1,       0){60}}

\put(150,     30){$r^{-1}$}
\put(150,     10){$-1$}

\put(150,     70){$r$}

\put(190,     80){$-r^{}$}
\put(190,     -10){$q^{-1}$}

\put(180,     30){$-1$}

\put(250,     30){$-q^{}$}

\put(250,     10){$q$}

\put(250,     70){$-1$}


\put(295,         -1)  {,  $r\in R_3, $  by   {\rm  GDD}  $8$ of Row $17$ in Table A2. }
\put(80,         -1)  {    } \end{picture}$\\
\\
  The sub-{\rm GDD} by deleting  Vertex  1    is not an arithmetic  {\rm  GDD}   by Lemma \ref {3.1.1}{\rm (I)}.\\

{\ }\!\!\!\!\!\!\!\!\!\!\!\!\!\!\!\!\!\!\!\!\!{\ }\!\!\!\!\!\!\!\!\!\!\!\!\!\!\!\!\!\!\!\!\!{\ }\!\!\!\!\!\!\!\!\!\!\!
  $\begin{picture}(100,       15)\put(-68,       -1){ }\put(68,       -1){ (e)}

\put(111,      1){\makebox(0,      0)[t]{$\bullet$}}
\put(144,       1){\makebox(0,       0)[t]{$\bullet$}}
\put(170,      -11){\makebox(0,      0)[t]{$\bullet$}}
\put(170,     15){\makebox(0,      0)[t]{$\bullet$}}
\put(113,      -1){\line(1,      0){33}}
\put(142,     -1){\line(2,      1){27}}
\put(170,       -14){\line(-2,      1){27}}

\put(170,       -14){\line(0,      1){27}}

\put(100,       10){$-r^{-1}$}
\put(125,       5){$r^{2}$}

\put(130,      -10){$-1$}

\put(145,      -20){$-1$}
\put(150,       15){$-q^{}$}

\put(175,       10){$q$}

\put(178,       -20){$-1$}

\put(178,       0){$q^{-1}$}

\put(200,         -1)  {,  $r\in R_6, $  by   {\rm  GDD}  $4$ of Row $7$ in Table A2. }
\put(80,         -1)  {    } \end{picture}$\\
\\
  The sub-{\rm GDD} by deleting  Vertex  4    is not   an arithmetic  {\rm  GDD}   by Lemma \ref {3.1.1}{\rm (I)}.

 (f)     i.e.  (17.8.3),  $r^2 \not=1, $ by    {\rm  GDD}  $1$ of Row $9$ in Table A2. The sub-{\rm GDD} by deleting
 Vertex  4  is not  an arithmetic  {\rm  GDD}    when $r \not=-q$  by Lemma \ref {3.1.1} {\rm (IV)}.
  The sub-{\rm GDD} by deleting  Vertex  4    is   an arithmetic  {\rm  GDD}    when $r =-q$  by  Type   4.
  It  is  quasi-affine
 by Lemma \ref {3.1.3}.\\

{\ }\!\!\!\!\!\!\!\!\!\!\!\!\!\!\!\!\!\!\!\!\!{\ }\!\!\!\!\!\!\!\!\!\!\!\!\!\!\!\!\!\!\!\!\!{\ }\!\!\!\!\!\!\!\!\!\!\!
   $\begin{picture}(100,       15)\put(-68,       -1){ }\put(68,       -1){ (g)}

\put(111,      1){\makebox(0,      0)[t]{$\bullet$}}
\put(144,       1){\makebox(0,       0)[t]{$\bullet$}}
\put(170,      -11){\makebox(0,      0)[t]{$\bullet$}}
\put(170,     15){\makebox(0,      0)[t]{$\bullet$}}
\put(113,      -1){\line(1,      0){33}}
\put(142,     -1){\line(2,      1){27}}
\put(170,       -14){\line(-2,      1){27}}

\put(170,       -14){\line(0,      1){27}}

\put(100,       10){$r^{}$}
\put(120,       5){$-1$}

\put(130,      -10){$-1$}

\put(145,      -20){$-1$}
\put(150,       15){$-q^{}$}

\put(175,       10){$q$}

\put(178,       -20){$-1$}

\put(178,       0){$q^{-1}$}

\put(200,         -1)  {,  $r\in R_4, $ by   {\rm  GDD}  $1$ of Row $2$ in Table A2. }
\put(80,         -1)  {    } \end{picture}$\\
\\
  The sub-{\rm GDD} by deleting  Vertex  4    is not an arithmetic  {\rm  GDD}   by Lemma \ref {3.1.1}{\rm (I)}.

   {\rm (ii) }  All quasi-arithmetic {\rm  GDD}s by  adding a vertex  on  Vertex $2$  are listed.
 According  to Lemma \ref {3.1.2} (V) we have to consider following cases.

 (a)     i.e.  (17.8.4),     $r\in R_3, $ by   {\rm  GDD}  $1$ of Row $17$.  The sub-{\rm GDD} by deleting  Vertex  4    is an arithmetic  {\rm  GDD}
when   $ q^{} =r^{}$ or  $ q^{} =r^{2}$ by Lemma \ref {3.1.1} {\rm (IX)}. It  is   quasi-affine   by Lemma \ref {3.1.3}.
\\

{\ }\!\!\!\!\!\!\!\!\!\!\!\!\!\!\!\!\!\!\!\!\!{\ }\!\!\!\!\!\!\!\!\!\!\!\!\!\!\!\!\!\!\!\!\!{\ }\!\!\!\!\!\!\!\!\!\!\!
   $\begin{picture}(100,       15)\put(-68,       -1){ }\put(68,       -1){ (b)}

\put(111,      1){\makebox(0,      0)[t]{$\bullet$}}
\put(144,       1){\makebox(0,       0)[t]{$\bullet$}}
\put(170,      -11){\makebox(0,      0)[t]{$\bullet$}}
\put(170,     15){\makebox(0,      0)[t]{$\bullet$}}
\put(113,      -1){\line(1,      0){33}}
\put(142,     -1){\line(2,      1){27}}
\put(170,       -14){\line(-2,      1){27}}

\put(170,       -14){\line(0,      1){27}}

\put(100,       10){$r^{-1}$}
\put(110,       5){$-r^{-1}$}

\put(127,      -15){$-1$}

\put(148,      -20){$-1$}
\put(155,       15){$q^{-1}$}

\put(175,       10){$q$}

\put(178,       -20){$-1$}

\put(178,       0){$-q^{}$}

\put(195,         -1)  {,  $r\in R_3, $   by   {\rm  GDD}  $9$ of Row $17$ in Table A2. }
\put(80,         -1)  {    } \end{picture}$\\
\\
  The sub-{\rm GDD} by deleting  Vertex  4    is not an arithmetic  {\rm  GDD}   by Lemma \ref {3.1.1}{\rm (I)}.\\ \\ \\ \\ \\

{\ }\!\!\!\!\!\!\!\!\!\!\!\!\!\!\!\!\!\!\!\!\!{\ }\!\!\!\!\!\!\!\!\!\!\!\!\!\!\!\!\!\!\!\!\!{\ }\!\!\!\!\!\!\!\!\!\!\!
  $\begin{picture}(100,       15)\put(-68,       -1){ }\put(68,       -1){ (c) }

\put(170,     10){\makebox(0,      0)[t]{$\bullet$}}

\put(170,     70){\makebox(0,      0)[t]{$\bullet$}}

\put(230,     10){\makebox(0,      0)[t]{$\bullet$}}
\put(230,     70){\makebox(0,      0)[t]{$\bullet$}}

\put(170,       10){\line(0,      1){60}}

\put(170,       10){\line(1,      1){60}}

\put(230,       10){\line(0,      1){60}}


\put(170,       10){\line(1,       0){60}}
\put(170,       70){\line(1,       0){60}}

\put(150,     30){$-r^{}$}
\put(150,     10){$-1$}

\put(150,     70){$r$}

\put(250,     30){$-q^{}$}
\put(250,     10){$q$}

\put(250,     70){$-1$}

\put(190,     80){$r^{-1}$}
\put(190,     -10){$q^{-1}$}

\put(180,     30){$-1$}

\put(280,         -1)  {,  $r\in R_3, $  by   {\rm  GDD}  $8$ of Row $17$ in Table A2. }
\put(80,         -1)  {    } \end{picture}$\\
\\
 The sub-{\rm GDD} by deleting  Vertex  1    is not  an arithmetic  {\rm  GDD}    by Lemma \ref {3.1.1}{\rm (I)}.
 \\ \\ \\  \\ \\

{\ }\!\!\!\!\!\!\!\!\!\!\!\!\!\!\!\!\!\!\!\!\!{\ }\!\!\!\!\!\!\!\!\!\!\!\!\!\!\!\!\!\!\!\!\!{\ }\!\!\!\!\!\!\!\!\!\!\!
  $\begin{picture}(100,       15)\put(-68,       -1){ }\put(68,       -1){ (d)}

\put(170,     10){\makebox(0,      0)[t]{$\bullet$}}

\put(170,     70){\makebox(0,      0)[t]{$\bullet$}}

\put(230,     10){\makebox(0,      0)[t]{$\bullet$}}
\put(230,     70){\makebox(0,      0)[t]{$\bullet$}}

\put(170,       10){\line(0,      1){60}}

\put(170,       10){\line(1,      1){60}}

\put(230,       10){\line(0,      1){60}}


\put(170,       10){\line(1,       0){60}}
\put(170,       70){\line(1,       0){60}}

\put(150,     30){$r^{-1}$}
\put(150,     10){$-1$}

\put(150,     70){$r$}

\put(250,     30){$-q^{}$}
\put(250,     10){$q$}

\put(250,     70){$-1$}

\put(190,     80){$-r^{}$}
\put(190,     -10){$q^{-1}$}

\put(180,     30){$-1$}

\put(295,         -1)  {,  $r\in R_3, $  by   {\rm  GDD}  $8$ of Row $17$ in Table A2. }
\put(80,         -1)  {    } \end{picture}$\\
\\
  The sub-{\rm GDD} by deleting  Vertex  2    is not an arithmetic  {\rm  GDD}   by Lemma \ref {3.1.1}{\rm (I)}.

 (e)  i.e.  (17.8.5),  $r\in R_6, $ by    {\rm  GDD}  $4$ of Row $7$ in Table B.
  The sub-{\rm GDD} by deleting  Vertex  4    is not an arithmetic  {\rm  GDD}   when $ q^{} \not=r^{2}$ by Lemma \ref {3.1.1}{\rm (I)}.  The sub-{\rm GDD} by deleting  Vertex  4    is  an arithmetic  {\rm  GDD}   when $ q^{} =r^{2}$ by  Type   3. It  is   quasi-affine.
 by Lemma \ref {3.1.3}.

(f)  i.e.  (17.8.6),  by  Lemma \ref {3.1.2}(IV). The sub-{\rm GDD} by
 deleting  Vertex  4    is  an arithmetic  {\rm  GDD}   by Lemma \ref {3.1.1} {\rm (IV)}. It  is   quasi-affine
 by Lemma \ref {3.1.3}.\\

{\ }\!\!\!\!\!\!\!\!\!\!\!\!\!\!\!\!\!\!\!\!\!{\ }\!\!\!\!\!\!\!\!\!\!\!\!\!\!\!\!\!\!\!\!\!{\ }\!\!\!\!\!\!\!\!\!\!\!
$\begin{picture}(100,       15)\put(-68,       -1){ }\put(68,       -1){ (g)}

\put(111,      1){\makebox(0,      0)[t]{$\bullet$}}
\put(144,       1){\makebox(0,       0)[t]{$\bullet$}}
\put(170,      -11){\makebox(0,      0)[t]{$\bullet$}}
\put(170,     15){\makebox(0,      0)[t]{$\bullet$}}
\put(113,      -1){\line(1,      0){33}}
\put(142,     -1){\line(2,      1){27}}
\put(170,       -14){\line(-2,      1){27}}

\put(170,       -14){\line(0,      1){27}}

\put(100,       10){$r^{}$}
\put(125,       5){$-1$}

\put(130,      -15){$-1$}

\put(145,      -20){$-1$}
\put(155,       15){$q^{-1}$}

\put(175,       10){$q$}

\put(178,       -20){$-1$}

\put(178,       0){$-q^{}$}

\put(195,         -1)  {,  $r\in R_4, $ by    {\rm  GDD}  $1$ of Row $2$ in Table B. }
\put(80,         -1)  {    } \end{picture}$\\
\\
  The sub-{\rm GDD} by deleting  Vertex  4    is not an arithmetic  {\rm  GDD}   by Lemma \ref {3.1.1}{\rm (I)}.

  {\rm (iii) } All quasi-arithmetic {\rm  GDD}s by  adding a vertex  on  Vertex $3$  are listed.

Checking  step by step we have to consider following cases.
\\

{\ }\!\!\!\!\!\!\!\!\!\!\!\!\!\!\!\!\!\!\!\!\!{\ }\!\!\!\!\!\!\!\!\!\!\!\!\!\!\!\!\!\!\!\!\!{\ }\!\!\!\!\!\!\!\!\!\!\!
  $\begin{picture}(100,       15)\put(-68,       -1){ }\put(68,       -1){ (a)}

\put(111,      1){\makebox(0,      0)[t]{$\bullet$}}
\put(144,       1){\makebox(0,       0)[t]{$\bullet$}}
\put(170,      -11){\makebox(0,      0)[t]{$\bullet$}}
\put(170,     15){\makebox(0,      0)[t]{$\bullet$}}
\put(113,      -1){\line(1,      0){33}}
\put(142,     -1){\line(2,      1){27}}
\put(170,       -14){\line(-2,      1){27}}

\put(170,       -14){\line(0,      1){27}}

\put(100,       10){$-q$}
\put(115,       5){$-q^{-1}$}

\put(127,       -12){$q^{-1}$}

\put(140,      -20){$-q^{-1}$}
\put(145,       15){$q^{}$}

\put(175,       10){$-1$}

\put(178,       -12){$-1$}

\put(178,       0){$-1$}

\put(195,         -1)  {,  $q\in R_3, $  by   {\rm  GDD}  $2$ of Row $13$. }
\put(80,         -1)  {    } \end{picture}$\\
\\
 The sub-{\rm GDD} by deleting  Vertex  4   is not an  arithmetic  {\rm  GDD}   by Lemma \ref {3.1.2} {\rm (III)}.

(b)  i.e.  (17.8.7),  $q\in R_3, $ by   {\rm  GDD}  $4$ of Row $15$.
 The sub-{\rm GDD} by deleting  Vertex  4   is an arithmetic  {\rm  GDD}   by  {\rm  GDD}  $7$ of Row $17$ or by Lemma \ref {3.1.2} {\rm (III)}.
  It  is quasi-affine   by Lemma \ref {3.1.3}.

 (c)     i.e.  (17.8.8),  $q\in R_3$ by    {\rm  GDD}  $6$ of Row $17$.  The sub-{\rm GDD} by deleting  Vertex  4
 is an arithmetic  {\rm  GDD}    by Lemma \ref {3.1.2} (III). It  is   quasi-affine   by Lemma \ref {3.1.3}.

 (d)     i.e.  (17.8.9),  $q\in R_3, $  by   {\rm  GDD}  $7$ of Row $17$.
 The sub-{\rm GDD} by deleting  Vertex  4   is an arithmetic  {\rm  GDD}.
  It  is quasi-affine   by Lemma \ref {3.1.3}.

  {\rm (iv) } All quasi-affine  {\rm  GDD}s which are complete diagrams are listed.
 \\
\\ \\ \\ \\

{\ }\!\!\!\!\!\!\!\!\!\!\!\!\!\!\!\!\!\!\!\!\!{\ }\!\!\!\!\!\!\!\!\!\!\!\!\!\!\!\!\!\!\!\!\!\!\!\!\!\! $\begin{picture}(100,       15)\put(55,       -1){ (d) in  Case {\rm (i) }}\put(-68,       -1){ }

\put(170,     10){\makebox(0,      0)[t]{$\bullet$}}

\put(170,     70){\makebox(0,      0)[t]{$\bullet$}}

\put(230,     10){\makebox(0,      0)[t]{$\bullet$}}
\put(230,     70){\makebox(0,      0)[t]{$\bullet$}}

\put(170,       10){\line(0,      1){60}}

\put(170,       10){\line(1,      1){60}}

\put(230,       10){\line(0,      1){60}}

\put(230,       10){\line(-1,      1){60}}

\put(170,       10){\line(1,       0){60}}
\put(170,       70){\line(1,       0){60}}

\put(150,     10){$-1$}
\put(150,     30){$r^{-1}$}

\put(150,     70){$r$}

\put(180,     30){$-1$}
\put(190,     80){$-r$}
\put(190,     -10){$q^{-1}$}
\put(220,     30){$$}

\put(250,     10){$q$}
\put(250,     30){$-q^{}$}

\put(250,     70){$-1$}

\put(280,         -1)  {  } \end{picture}$\\
\\ The sub-{\rm GDD} by deleting  Vertex 1  is not an arithmetic  {\rm  GDD}   by Lemma \ref {3.1.1}{\rm (I)}.
\\ \\ \\ \\ \\

 {\ }\!\!\!\!\!\!\!\!\!\!\!\!\!\!\!\!\!\!\!\!\!{\ }\!\!\!\!\!\!\!\!\!\!\!\!\!\!\!\!\!\!\!\!\!\!\!\!   $\begin{picture}(100,       15)\put(-68,       -1){ }\put(48,       -1){ (c) in  Case {\rm (i) }}

\put(170,     10){\makebox(0,      0)[t]{$\bullet$}}

\put(170,     70){\makebox(0,      0)[t]{$\bullet$}}

\put(230,     10){\makebox(0,      0)[t]{$\bullet$}}
\put(230,     70){\makebox(0,      0)[t]{$\bullet$}}

\put(170,       10){\line(0,      1){60}}

\put(170,       10){\line(1,      1){60}}

\put(230,       10){\line(0,      1){60}}

\put(230,       10){\line(-1,      1){60}}

\put(170,       10){\line(1,       0){60}}
\put(170,       70){\line(1,       0){60}}

\put(150,     10){$-1$}
\put(150,     30){$-r$}

\put(150,     70){$r$}

\put(180,     30){$-1$}
\put(190,     80){$r^{-1}$}
\put(190,     -10){$q^{-1}$}
\put(220,     30){$$}

\put(250,     10){$q$}
\put(250,     30){$-q^{}$}

\put(250,     70){$-1$}

\put(280,         -1)  { } \end{picture}$\\
The sub-{\rm GDD} by deleting  Vertex 1  is not an arithmetic  {\rm  GDD}   by Lemma \ref {3.1.1}{\rm (I)}.

(c) in  Case {\rm (ii) } and (d) in  Case {\rm (ii) } are the same as (c) in  Case {\rm (i) }and (d) in Case {\rm (i) }.

\subsection* {Quasi-affine over
  {\rm  GDD}  $9$ of Row $17$  }
  {\rm (i)} All quasi-arithmetic {\rm  GDD}s by  adding a vertex  on  Vertex $1$  are listed.
 According  to Lemma \ref {3.1.2} (V) we have to consider following cases.

 (a)     i.e.  (17.9.1), $q\in R_3, $ $r \in  R_3 \cup R_4 \cup R_6 $,     by    {\rm  GDD}  $1$ of Row $17$  in Table A2.
It is quasi-affine   by Lemma \ref {3.1.3}.

 (b)     i.e.  (17.9.2),  $r\in R_3, $  by   {\rm  GDD}  $9$ of Row $17$. It is quasi-affine   by Lemma \ref {3.1.3}.

 (c)     i.e.  (17.9.3),  $r \in R_6, $ by   {\rm  GDD}  $4$ of Row $7$  in Table A2.
It is quasi-affine   by Lemma \ref {3.1.3}.

 (d)     i.e.  (17.9.4),   $r \in   R_3 \cup R_4 \cup R_6 $ by Lemma \ref {3.1.2} (IV).
It is quasi-affine   by Lemma \ref {3.1.3}.

 (e)     i.e.  (17.9.5), $r \in R_4, $ by   {\rm  GDD}  $1$ of Row $2$  in Table A2.
It is quasi-affine   by Lemma \ref {3.1.3}.

  {\rm (ii) }  All quasi-arithmetic {\rm  GDD}s by  adding a vertex  on  Vertex $2$  are listed.
According  to   Type   4
 we have to consider following cases.

 (a)     i.e.  (17.9.6),   Type   4.  The sub-{\rm GDD} by deleting  Vertex  4   is an arithmetic  {\rm  GDD} }
by Lemma \ref {3.1.2} {\rm (V)}. It  is   quasi-affine    by Lemma \ref {3.1.3}.

 (b)     i.e.  (17.9.7)
   by   Type   4. The sub-{\rm GDD} by deleting  Vertex  4   is
an arithmetic  {\rm  GDD}   by Lemma \ref {3.1.2} {\rm (IV)}. It  is   quasi-affine    by Lemma \ref {3.1.3}.

 (c)     i.e.  (17.9.8),  $q \in R_3, $ by   {\rm  GDD}  $9$ of Row $17$  in Table A2. The sub-{\rm GDD} by
 deleting  Vertex  4   is an arithmetic  {\rm  GDD}. It  is   quasi-affine    by Lemma \ref {3.1.22} or by Lemma \ref {3.1.3}.

  {\rm (iii) }  All quasi-arithmetic {\rm  GDD}s by  adding a vertex  on  Vertex $3$  are listed.
According  to  Type   4
 we have to consider following cases.

 (a)     i.e.  (17.9.9)   by   {\rm  GDD}  $6$ of Row $17$. It is quasi-affine
 by  Lemma   \ref {2.63}.

 (b)     i.e.  (17.9.10)  by   {\rm  GDD}  $7$ of Row $17$. It is quasi-affine
 by  Lemma   \ref {2.63}.

  {\rm (iv) } All quasi-affine circles are listed.

 (a)
 (a)     i.e.  (17.9.11)      is quasi-affine  since the sub-{\rm GDD} by deleting  Vertex 2 is not an arithmetic  {\rm  GDD}   by Lemma \ref {3.1.1} (I) when  $ q^{} \not=r^{}$.  It is  an arithmetic  {\rm  GDD}   by  Type   2 when  $ q^{} =r^{}$.

 (a)
 (b)       i.e.  (17.9.12)  is  quasi-affine  since the sub-{\rm GDD} by deleting  Vertex 2
 is  an arithmetic  {\rm  GDD}   when $r=q^{-1}$ or $ q^{-1}= r^{2}$ by Lemma \ref {3.1.1} (IX).

\subsection* {Quasi-affine over
 {\rm  GDD}  $1$ of Row $18$   }
  {\rm (i) }  All quasi-arithmetic {\rm  GDD}s by  adding a vertex  on  Vertex $1$  are listed.
 According  to Lemma \ref {3.1.2} (II) we have to consider following cases.

 (a)     i.e.  (18.1.1)$, q \in R_9$,   by    {\rm  GDD}  $1$ of Row $18$. It is quasi-affine   by Lemma \ref {3.1.3}.

  {\rm (ii) }  All quasi-arithmetic {\rm  GDD}s by  adding a vertex  on  Vertex $2$  are listed.
 According  to Lemma \ref {3.1.2} (II) we have to consider following cases.

(a)     i.e.  (18.1.2)$, q \in R_9$,   by    {\rm  GDD}  $1$ of Row $18$. It is quasi-affine   by Lemma \ref {3.1.3}.
\\ \\

  {\ }\ \ \ \ \ \ \ \ \ \ \ \ $
$\\
\\
  The sub-{\rm GDD} by deleting  Vertex  4    is not an arithmetic  {\rm  GDD}    by Lemma \ref {3.1.1}{\rm (I)}.\\
\\

  {\ }\ \ \ \ \ \ \ \ \ \ \ \ $%
$\\
\\
  The sub-{\rm GDD} by deleting  Vertex 4    is not an arithmetic  {\rm  GDD}    by Lemma \ref {3.1.1}{\rm (I)}.\\ \\

  {\ }\ \ \ \ \ \ \ \ \ \ \ \ $%
$\\
\\
  The sub-{\rm GDD} by deleting  Vertex  4    is not an arithmetic  {\rm  GDD}   by Lemma \ref {3.1.1}{\rm (I)}.
 \\ \\

  {\ }\ \ \ \ \ \ \ \ \ \ \ \ $%
$\\
\\
  The sub-{\rm GDD} by deleting  Vertex  4    is not an arithmetic  {\rm  GDD}    by Lemma \ref {3.1.1}{\rm (I)}.

\subsection* {Quasi-affine over
  {\rm  GDD}  $2$ of Row $18$  }

  {\rm (ii) }  All quasi-arithmetic {\rm  GDD}s by  adding a vertex  on  Vertex $2$  are listed.
Checking  step by step we have to consider following cases.\\

{\ }\!\!\!\!\!\!\!\!\!\!\!\!\!\!\!\!\!\!\!\!\!{\ }\!\!\!\!\!\!\!\!\!\!\!\!\!\!\!\!\!\!\!\!\!\!\!\!\!\!\!\!\!\!\!\!
 $%
$\\ \\
 The sub-{\rm GDD} by deleting  Vertex  4   is not an arithmetic  {\rm  GDD}    by Lemma \ref {3.1.1}{\rm (I)}.
 \hfill $\Box$

\section {Appendix}


\begin {Lemma} \label {a1} All simple circles and semi-simple circles are listed.\\ \\

 {\ }\ \ \ \ \ \ \ \ \ \ \ \ $%
$\\

\end {Lemma}

\begin {Lemma}\label {2.85}   Assume rank $n = 4.$  Continuation of all  {\rm  GDD}s of Table A2 are listed.

 {\rm  {\rm  GDD} } $1$ of Row $7$ when  $q \in  R_{4}$ is continual on  1  via T5 by  {\rm  {\rm  GDD} } $1$ of Row $22$ in  Table B.

 {\rm  {\rm  GDD} } $2$ of Row $7$ when  $q \in  R_{4}$ is continual on  1  via T6 by
 {\rm  {\rm  GDD} } $2$ of Row $22$.

  {\rm  {\rm  GDD} } $4$ of Row $7$ when  $q \in  R_{4}$ is continual on  3  via T6 by
 {\rm  {\rm  GDD} } $6$ of Row $22$.

 {\rm  {\rm  GDD} } $1$ of Row $9$ when  $q \in  R_{4}$ is continual on  3  via T6 by
 {\rm  {\rm  GDD} } $7$ of Row $22$ with   $r =-1.$

 {\rm  {\rm  GDD} } $1$ of Row $9$ when  $q \notin  R_{5}   \cup R_3\cup R_2$  is continual on  1  via T5 by
  {\rm  {\rm  GDD} } $5$ of Row $9$ with   $r =q^{-3}.$

  {\rm  {\rm  GDD} } $1$ of Row $9$ when  $q \notin  R_{4}$ is continual on  1  via T5 by
   {\rm  {\rm  GDD} } $1$ of Row $14$.

{\rm  {\rm  GDD} } $3$ of Row $9$ when  $q \notin  R_{4}$ is continual on  1  via T5 by
  {\rm  {\rm  GDD} } $5$ of Row $14$.

 {\rm  {\rm  GDD} } $1$ of Row $13$ when  $q \in  R_{3}$  is continual on  1  via T6 by
  {\rm  {\rm  GDD} } $1$ of Row $21$.

 {\rm  {\rm  GDD} } $2$ of Row $13$ when  $q \in  R_{3}$ is continual on  1  via T6 by
{\rm  {\rm  GDD} } $2$ of Row $21$.

 {\rm  {\rm  GDD} } $1$ of Row $15$ when  $q \in  R_{3}$   is continual on 1  via T5 by    {\rm  {\rm  GDD} } $1$ of Row $20$.

 {\rm  {\rm  GDD} } $1$ of Row $15$ when  $q \in  R_{3}$ is continual on  1  via T6 by
 {\rm  {\rm  GDD} } $3$ of Row $21$.

 {\rm  {\rm  GDD} } $1$ of Row $15$ when  $q \in  R_{3}$  is continual on  3  via T6 by   {\rm  {\rm  GDD} } $5$ of Row $20$.

 {\rm  {\rm  GDD} } $1$ of Row $15$ when  $q \in  R_{3}$ is continual on  3  via T5 by
  {\rm  {\rm  GDD} } $5$ of Row $21$.

{\rm  {\rm  GDD} } $2$ of Row $15$ when  $q \in  R_{3}$ is continual on  1  via T5 by
  {\rm  {\rm  GDD} } $6$ of Row $21$.

 {\rm  {\rm  GDD} } $2$ of Row $15$ when  $q \in  R_{3}$ is continual on  1  via T6 by
 {\rm  {\rm  GDD} } $3$ of Row $20$.

  {\rm  {\rm  GDD} } $3$ of Row $15$ when  $q \in  R_{3}$ is continual on  1  via T6 by
  {\rm  {\rm  GDD} } $7$ of Row $20$.

 {\rm  {\rm  GDD} } $4$ of Row $15$ when  $q \in  R_{3}$     is continual on  1  via T5 by      {\rm  {\rm  GDD} } $2$ of Row $20$.

 {\rm  {\rm  GDD} } $4$ of Row $15$ when  $q \in  R_{3}$ is continual on  1  via T6 by
 {\rm  {\rm  GDD} } $4$ of Row $21$.

 {\rm  {\rm  GDD} } $3$ of Row $15$ when  $q \in  R_{3}$ is continual on  1  via T5 by
 {\rm  {\rm  GDD} } $7$ of Row $21$.

 {\rm  {\rm  GDD} } $5$ of Row $16$ when  $q \in  R_{3}$   is continual on  3  via T5 by    {\rm  {\rm  GDD} } $1$ of Row $17$.

 {\rm  {\rm  GDD} } $4$ of Row $16$ when  $q \in  R_{3}$ is continual on  3  via T5 by
 {\rm  {\rm  GDD} } $2$ of Row $17$.



 {\rm  {\rm  GDD} } $1$ of Row $2$ when   $q^2 \not=1$ is continual on 3 via T6 by {\rm  {\rm  GDD} } $1$ of Row $9$.

 {\rm  {\rm  GDD} } $1$ of Row $4$ when   $q\in R_3$ is continual on 2 via T6 by {\rm  {\rm  GDD} } $6$ of Row $18$.

 {\rm  {\rm  GDD} } $1$ of Row $5$ when  $q\in R_3$ is continual on 3 via T6 by {\rm  {\rm  GDD} } $2$ of Row $6$.

 {\rm  {\rm  GDD} } $2$ of Row $5$ when  $q\in R_3$ is continual on 3 via T6 by {\rm  {\rm  GDD} } $6$ of Row $20$.

 {\rm  {\rm  GDD} } $1$ of Row $6$ when  $q\in R_3$ is continual on 3 via T6 by {\rm  {\rm  GDD} } $6$ of Row $20$.

 {\rm  {\rm  GDD} } $1$ of Row $6$ when  $q\in R_3$ is continual on 3 via T5 by {\rm  {\rm  GDD} } $1$ of Row $18$.

 {\rm  {\rm  GDD} } $1$ of Row $6$ when  $q\notin R_3 \cup R_2$ is continual on 3 via T5 by {\rm  {\rm  GDD} } $1$ of Row $9$.

 {\rm  {\rm  GDD} } $2$ of Row $6$ when   $q^2 \not=1$ is continual on 3 via T5 by {\rm  {\rm  GDD} } $2$ of Row $9$.

 {\rm  {\rm  GDD} } $3$ of Row $6$ when   $q\in R_3$  is continual on 2 via T6 by {\rm  {\rm  GDD} } $8$ of Row $20$.

 {\rm  {\rm  GDD} } $3$ of Row $8$ when   $q\in R_3$  is continual on 2 via T5 by {\rm  {\rm  GDD} } $6$ of Row $18$.

 {\rm  {\rm  GDD} } $1$ of Row $10$ when   $q\in R_4$  is continual on 1 via T6 by {\rm  {\rm  GDD} } $3$ of Row $22$.

 {\rm  {\rm  GDD} } $1$ of Row $10$ when   $q\in R_6$  is continual on 1 via T5 by {\rm  {\rm  GDD} } $6$ of Row $17$.

 {\rm  {\rm  GDD} } $1$ of Row $10$ when   $q\in R_5$  is continual on 1 via T5 by {\rm  {\rm  GDD} } $5$ of Row $9$.

 {\rm  {\rm  GDD} } $1$ of Row $10$ when   $q\in R_4$  is continual on 1 via T5 by {\rm  {\rm  GDD} } $1$ of Row $14$.

 {\rm  {\rm  GDD} } $1$ of Row $11$ when   $q\in R_3$  is continual on 1 via T6 by {\rm  {\rm  GDD} } $3$ of Row $12$.

 {\rm  {\rm  GDD} } $1$ of Row $11$ when   $q\in R_3$  is continual on 1 via T5 by {\rm  {\rm  GDD} } $1$ of Row $13$.

\end {Lemma}

We write the order of  {\rm  GDD}s as follows: \\ \\ \\

{\ }\!\!\!\!\!\!\!\!\!\!\!\!\!\!\!\!\!\!\!\!\!{\ }\!\!\!\!\!\!\!\!\!\!\!\!\!\!\!\!\!\!\!\!\!\!\!\!\!\!\!\!\!\!\!\!\!\!\!\!\!\!\!\!\!\!\!
   $\begin{picture}(100,       15)\put(-68,       -1){ }\put(68,       -1){ }

\put(111,      1){\makebox(0,      0)[t]{$\bullet$}}
\put(144,       1){\makebox(0,       0)[t]{$\bullet$}}
\put(170,      -11){\makebox(0,      0)[t]{$\bullet$}}
\put(170,     15){\makebox(0,      0)[t]{$\bullet$}}
\put(113,      -1){\line(1,      0){33}}
\put(142,     -1){\line(2,      1){27}}
\put(170,       -14){\line(-2,      1){27}}

\put(170,       -14){\line(0,      1){27}}

\put(100,       10){$q_{11}$}
\put(120,       5){$$}

\put(130,     10){$q_{22}$}

\put(145,      -20){$$}
\put(150,       15){$$}

\put(175,       10){$q_{33}$}

\put(178,       -20){$q_{44}$}

\put(178,       0){$$}

\put(200,         -1)  {   }
\put(80,         -1)  {    } \end{picture}$  {\ }\ \ \ \ \ \ \ \ \ \ \ \   {\ }\ \ \ \ \ \ \ \ \ \ \ \
  {\ }\ \ \ \ \ \ \ \ \ \ \ \ $\begin{picture}(100,       15) \put(-68,        -1){  }

\put(27,     38){\makebox(0,      0)[t]{$\bullet$}}

\put(60,       1){\makebox(0,       0)[t]{$\bullet$}}
\put(58,       -12){$q_{44}$}

\put(40,       -12){$$}
\put(28,       -1){\line(1,       0){33}}
\put(27,       1){\makebox(0,      0)[t]{$\bullet$}}

\put(-14,       1){\makebox(0,      0)[t]{$\bullet$}}

\put(22,      -12){$q_{22}$}
\put(0,       -12){$$}

\put(-14,      -1){\line(1,       0){50}}

\put(-18,      -12){$q_{11}$}

\put(27,       0){\line(0,      1){35}}

\put(30,       30){$q_{33}$}

\put(30,       20){$$}

\ \ \ \ \ \ \ \ \ \ \ \ \ \ \ \ \ \ \ \ \ \ \ \ \ \ \ \ \ \ \ \ \ \ \ \ \ \   \   \ \ \ \ \ \ \ \ \ \ \ \ \ \ \ \ \ \ \ {}
\put(80,         -1)  {    } \end{picture}$
{\ }\!\!\!\!\!\!\!\!\!\!\!\!\!\!\!\!\!\!\!\!\!{\ }\!\!\!\!\!\!\!\!\!\!\!\!\!\!\!\!\!\!\!\!\!{\ }\!\!\!\!\!\!\!\!\!\!\!
   $\begin{picture}(100,       15)\put(-68,       -1){ }\put(68,       -1){ }

\put(170,     10){\makebox(0,      0)[t]{$\bullet$}}

\put(170,     70){\makebox(0,      0)[t]{$\bullet$}}

\put(230,     10){\makebox(0,      0)[t]{$\bullet$}}
\put(230,     70){\makebox(0,      0)[t]{$\bullet$}}

\put(170,       10){\line(0,      1){60}}

\put(170,       10){\line(1,      1){60}}

\put(230,       10){\line(0,      1){60}}


\put(170,       10){\line(1,       0){60}}
\put(170,       70){\line(1,       0){60}}

\put(150,     30){$$}
\put(150,     10){$q_{11}$}

\put(150,     70){$q_{44}$}

\put(190,     80){$$}
\put(190,     -10){$$}

\put(180,     30){$$}

\put(250,     30){$$}

\put(250,     10){$q_{33}$}

\put(250,     70){$q_{22}$}


\put(295,         -1)  { }
\put(80,         -1)  {    } \end{picture}$\\

\end {document}